\documentclass[
	12pt,				
	openright,			
    oneside,
	a4paper,			
	english,			
	french,				
	spanish,			
	brazil				
	]{abntex2}

\usepackage{amssymb}
\usepackage{amsmath}
\usepackage{mathrsfs}
\usepackage{epstopdf}
\usepackage{enumerate}
\usepackage{amsfonts}
\newtheorem{prop}{Proposition}[section]
\newtheorem{remark}{Remark}[section]
\newtheorem{lemma}{Lemma}[section]
\newtheorem{theorem}{Theorem}[section]
\newtheorem{definition}{Definition}[section]
\newtheorem{corollary}{Corollary}[section]

\newtheorem{example}{Example}[section]

\newcommand{\qed}{\begin{flushright}
$\Box$
\end{flushright}}
\newcommand{\cosg}{\mathrm{cosg}}
\newcommand{\sing}{\mathrm{sing}}
\newcommand{\vecpe}{\times_{e}}

\usepackage{lmodern}			
\usepackage[T1]{fontenc}		
\usepackage[utf8]{inputenc}		
\usepackage{lastpage}			
\usepackage{indentfirst}		
\usepackage{color}				
\usepackage{graphicx}			
\usepackage{microtype} 			
		
\usepackage{lipsum}				

\usepackage[english,hyperpageref]{backref}	 
\usepackage[alf,abnt-etal-list=0,abnt-etal-cite=3]{abntex2cite}	

\renewcommand*{\backrefalt}[4]{
	\ifcase #1 %
		Nenhuma citação no texto.%
	\or
		Cited on page #2.%
	\else
		Cited #1 times on pages #2.%
	\fi}

\titulo{Differential Geometry of Rotation Minimizing Frames, Spherical Curves, and Quantum Mechanics of a Constrained Particle}
\autor{Luiz Carlos Barbosa da Silva}
\local{Recife}
\data{2017}
\preambulo{Trabalho apresentado ao Programa de Pós-graduação em Matemática do Departamento de Matemática da Universidade Federal de Pernambuco como requisito parcial para a
obtenção do grau de Doutor em Matemática.
 \newline
 \newline
 \newline
 \newline
 \newline
 \newline
Supervisor: Prof. Dr. Fernando Antônio Nóbrega Santos}

\definecolor{blue}{RGB}{41,5,195}

\makeatletter
\hypersetup{
		pdftitle={\@title}, 
		pdfauthor={\@author},
    	pdfsubject={\imprimirpreambulo},
	    pdfcreator={LaTeX with abnTeX2},
		pdfkeywords={abnt}{latex}{abntex}{abntex2}{trabalho acadêmico}, 
		colorlinks=true,       		
    	linkcolor=black,          	
    	citecolor=blue,        		
    	filecolor=magenta,      		
		urlcolor=blue,
		bookmarksdepth=4
}
\makeatother

\makeindex

\begin{document}

\frenchspacing 
\selectlanguage{english}


  \begin{capa}%

\vspace{-0.5cm}

\begin{figure}[h!]%
        \centering%
        \includegraphics[width=0.25\linewidth]{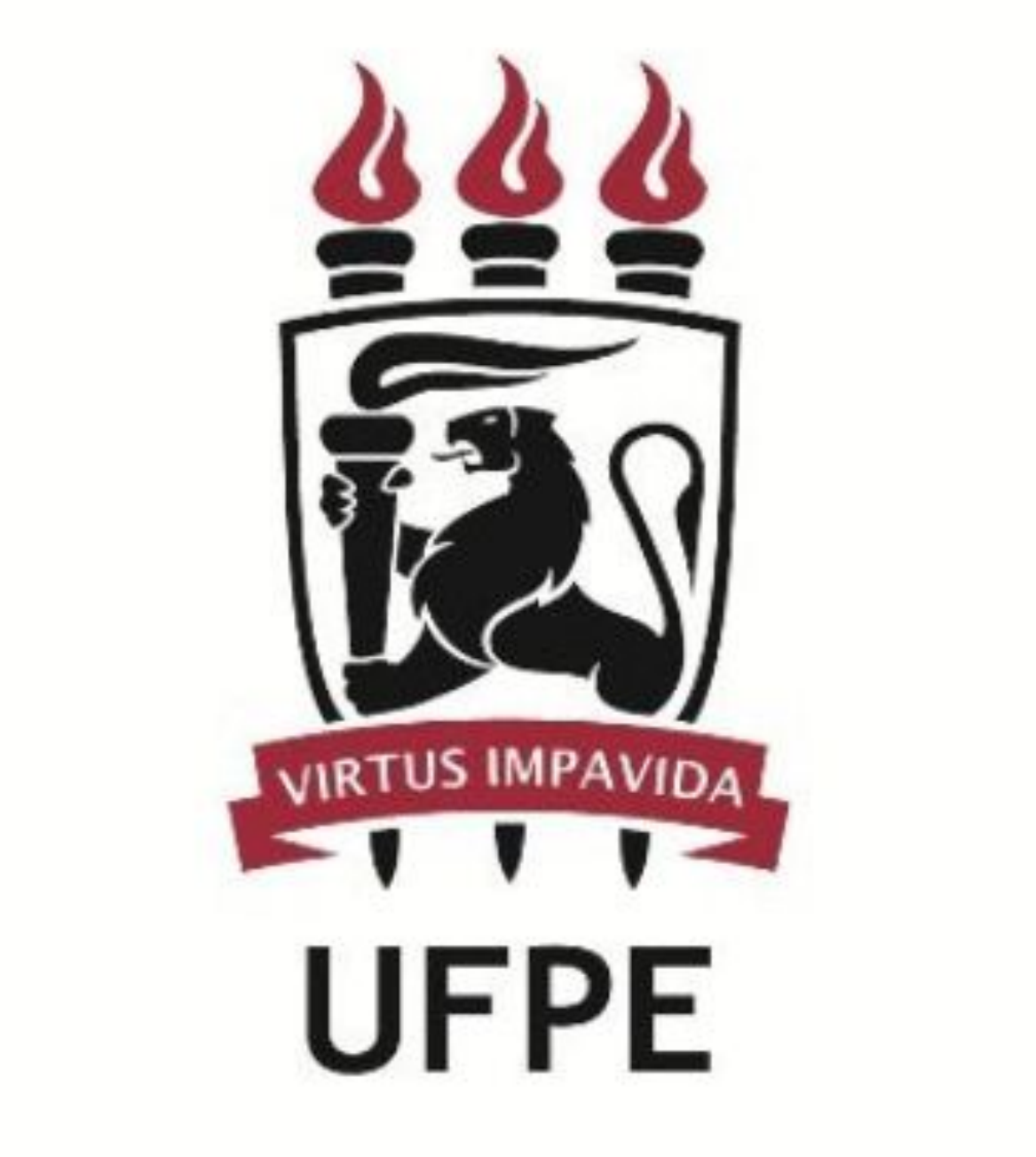}%
      \end{figure}%
    \center
\ABNTEXchapterfont\large{Universidade Federal de Pernambuco\\
    Centro de Ciências Exatas e da Natureza\\
    Departamento de Matemática\\
    Programa de Doutorado em Matemática}

\vfill

\ABNTEXchapterfont\large\imprimirautor

\vfill
    \ABNTEXchapterfont\bfseries\LARGE\imprimirtitulo
    \vfill

	
%
    \large\imprimirlocal 
    
    \large\imprimirdata

    \vspace*{1cm}
  \end{capa}


\imprimirfolhaderosto*


%
 \begin{fichacatalografica}
 \begin{figure}[t] 
     \begin{flushleft} 
     \includegraphics[width=\linewidth]{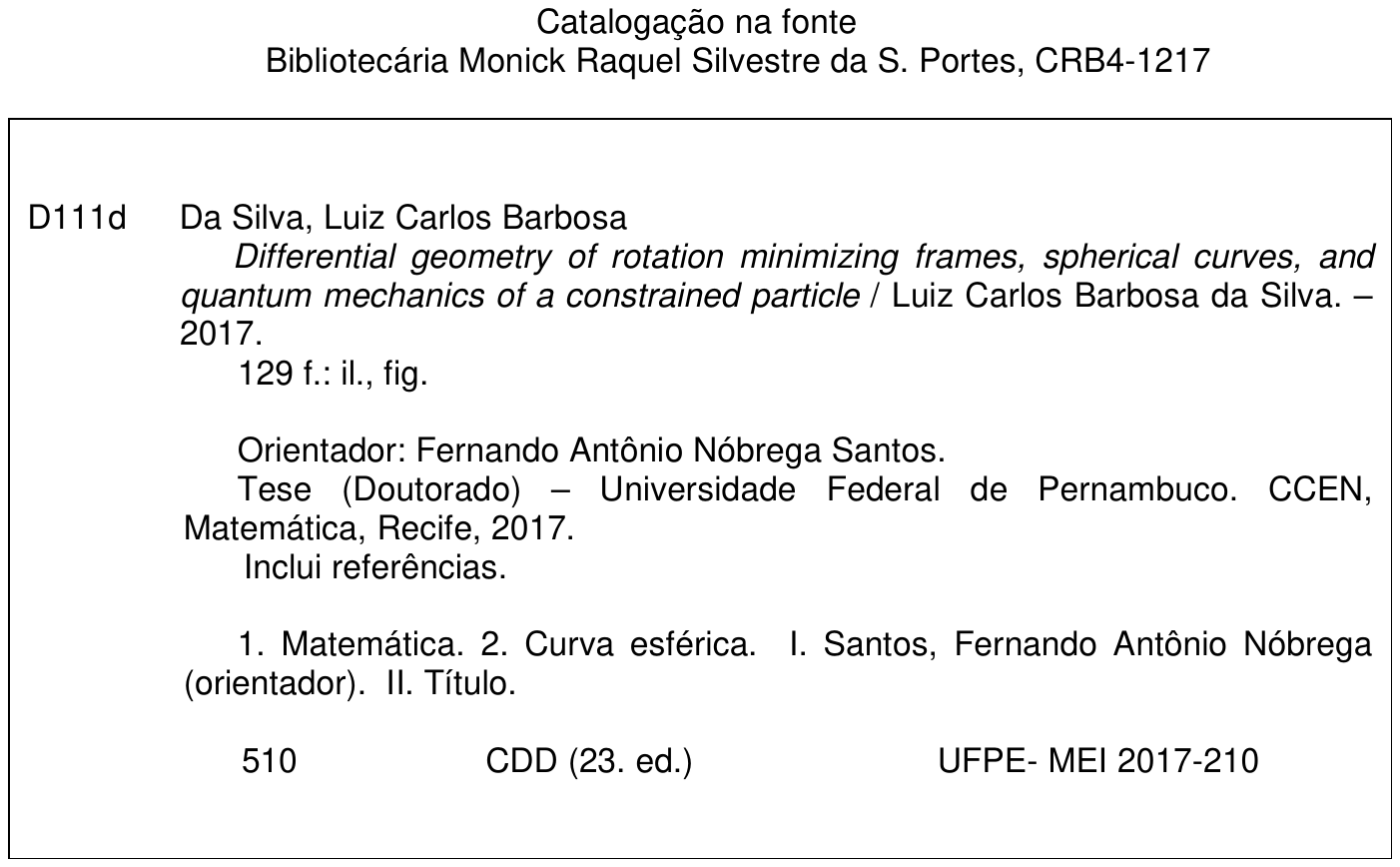} 
     \end{flushleft}
\end{figure}
 \end{fichacatalografica}

\begin{folhadeaprovacao}

%
%
%
%

\begin{figure}[t] 
     \begin{flushleft} 
     \includegraphics[width=\linewidth]{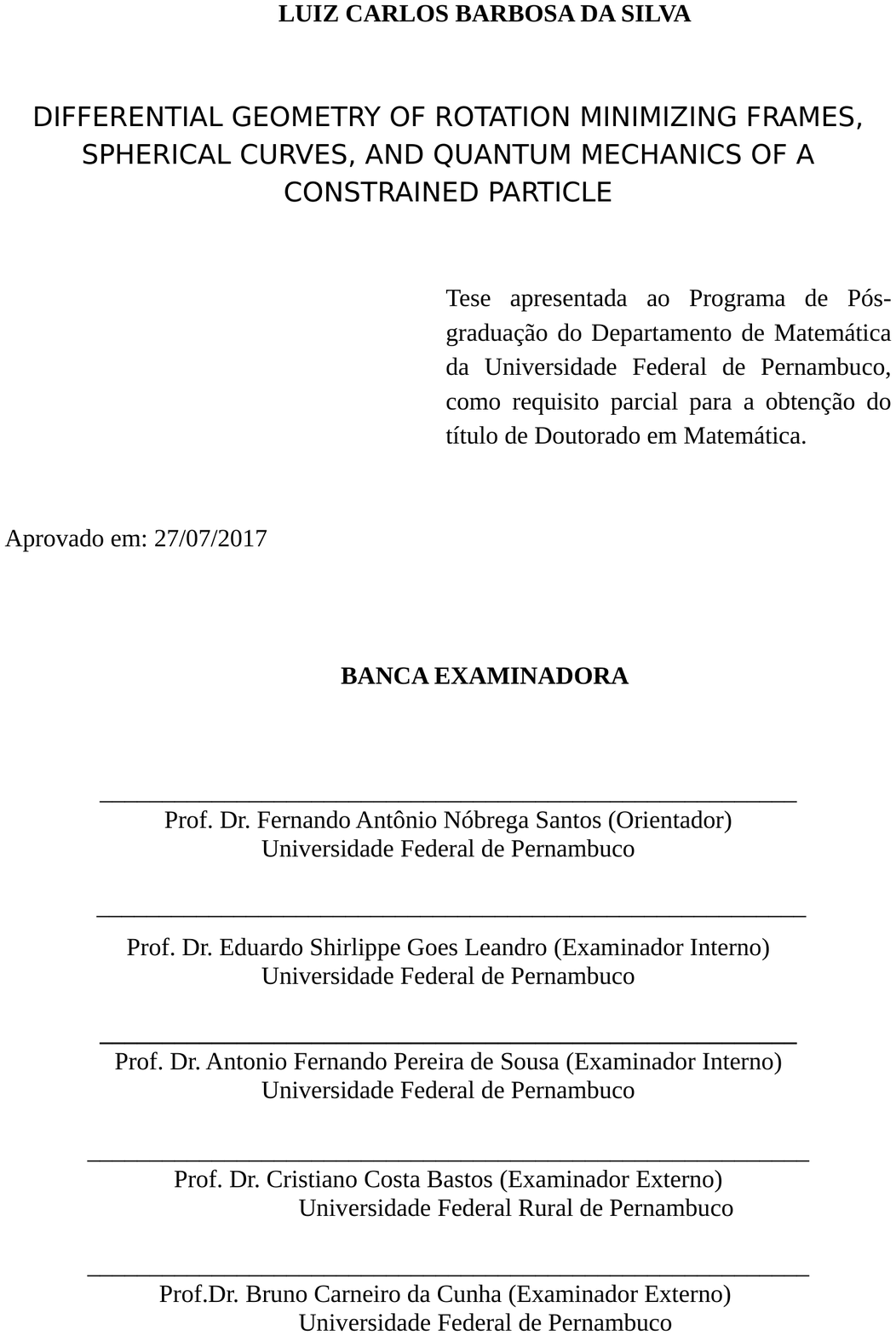} 
     \end{flushleft}
\end{figure}

\end{folhadeaprovacao}


\begin{dedicatoria}
   \vspace*{\fill}
   \centering
   \noindent
   \textit{ To my mother Noemia\\
   for all the love and support.} \vspace*{\fill}
\end{dedicatoria}

\begin{agradecimentos}

It was March of 2008, the beginning of my journey as an undergrad student at the Universidade Federal de Pernambuco (UFPE). My dreams at that point? Well, I just wanted to learn as much as possible and then work as a scientist... I'm kidding, of course I also dreamed about winning a Nobel prize, mom would be really proud of me... but I promptly realized that it was a better idea to invest my time and efforts on solving the homework assignments! What a wonderful period... I remember that the beginning of each new semester at those initials years was accompanied by the same recurrent questions: am I smart enough to learn all these subjects? will I be happy doing this? is this my dream after all???... maybe I should throw it all away and go to Bollywood ... With time questions changed a little: does that professor hate me? will I have a job???... Well, ten years have passed since then and I do not ask myself such questions anymore, or at least not all of them. I am a person who finds pleasure in doing what I do and yes, I had a dream and I made it real... I am a scientist!
\begin{figure*}[h]
\centering
  {\includegraphics[width=0.8\linewidth]{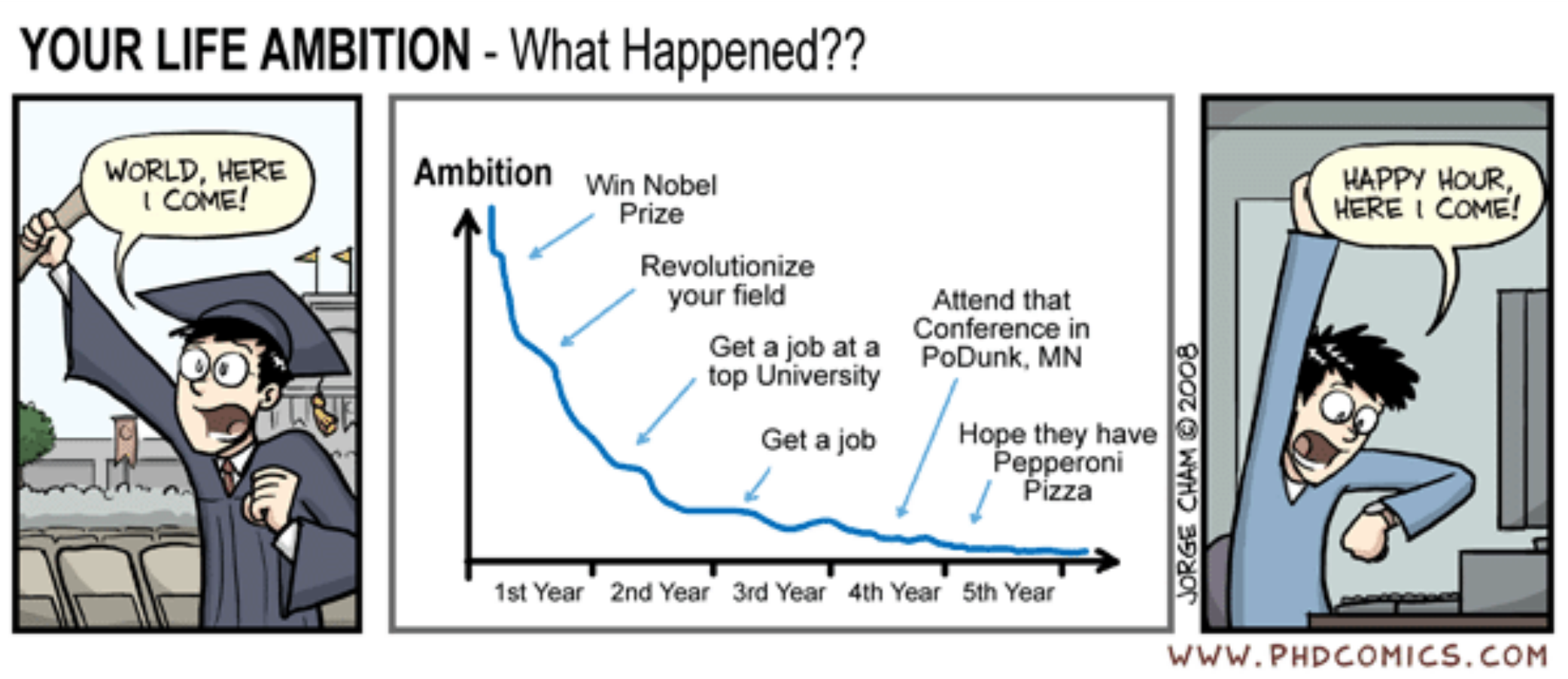}}
          \caption{From PHD Comics \url{http://phdcomics.com/}.}
\end{figure*}

Now it is the time for new goals, new dreams, and new disappointments also. Maybe, I should try working and living in that not so beautiful city where one of my heroes worked too... or maybe teach that infamous course the way I think is the best one. Anyway, the point here is that I didn't walk alone. On the road I met friends, good professors, not so good professors... advisors, collaborators, crazy and also reasonable referees, other researches... Well, science is also a human endeavor after all. In short, this is the moment to say thanks!  

I know it is a cliché to thank all the professors that one encounters along the way, but I will do it anyway. I would like to thank all the individuals that helped me along these long, wonderful, illuminating, and sometimes very difficult 10 years at UFPE. It was a great pleasure to meet and learn from you all. In particular, I would like to mention professors Antônio de Sousa, César Castilho, Eduardo Leandro, Francisco Brito, Hamid Hassanzadeh, Henrique Araújo, Maurício Coutinho-Filho, Ramón Mendoza, Sérgio Santa Cruz, Solange Rutz, and naturally my advisor Fernando Santos. I also thank the staff from the Department of Mathematics. 

I would like to take the opportunity to thank all the good discussions with professor Carmine Ortix from Utrecht University, the Netherlands, and professor Fernando Etayo from Universidad de Cantabria, Spain. I would also like to thank all my collaborators: Bertrand Berche and Sebastien Fumeron from Université de Lorraine, France; Fábio G. Ribeiro from Instituto Federal da Paraíba, Brazil; Cristiano C. Bastos, Fernando Moraes, and José Deibsom da Silva from Universidade Federal Rural de Pernambuco, Brazil; and finally my advisor Fernando A. N. Santos and former advisor Maurício D. Coutinho-Filho from UFPE. They all contributed a lot for my formation. In particular, professor Coutinho-Filho and his pursuit of a research done with all the rigor and care that it deserves influenced me a lot. I do my best to be as good and as professional as he is. Besides, I also want to express my gratitude to my advisor, professor Santos. We work together since I was an undergrad and then I must say that he decisively helped me in achieving my goals. In addition, his enthusiasm, support, friendship, and efforts to teach as much as possible about how it is like to be a scientist were decisive to me. I don't have words to say how much I am grateful and I am not able to imagine a better advisor. 

An attentive reader will probably note that our Bibliography contains texts written in 5 different languages. This reflects the importance of language skills in the daily life of a scientist. So, I would like to take this opportunity and thank all the teachers that I had and have: in particular, Mariano Hebenbrock (German), Alexandra Kasky (German), Laura Trutt (French), and Alisa Zingerman (Hebrew). In this respect ,one should not forget the important role played by some institutions. So, I would like to thank the Núcleo de Línguas e Culturas from UFPE, the Centro Cultural Brasil Alemanha, and also Italki, HelloTalk, Memrise, and Duolingo.

Probably, the ultimate cliché is that of acknowledging the family. But sometimes life is not that easy and the support from our relatives can not be taken for granted. So, I must say that I am a very luck guy for the opportunity of having a very supportive and special mother. She didn't have the opportunity to proceed further with her education, but did everything she could in order to offer her children a different fate. Then, thank you mom for everything, your support was crucial to give me the necessary strength to go on in pursuit of my dreams. Last but not least, I would like to thank the rest of my family. In particular, my sister Luciana, my nephew Kléber, my grandparents Josefa, José and Severina. I also thank my girlfriend Bruna, the last two years were like a dream. In addition, some people say that friends are the family that we can choose. I would like to thank all my friends. Unfortunately, I don't have a good memory to remember and mention here all their names. So, let me just thank some friends from the high school Leonardo, Leandro, and Thony; the friends who graduated with me in 2011 (BSc. in Math.) Edgar, Gilson, and João; Gabriel for being my tourist guide in Rome and for the useful discussion about geometry and mathematics in general; David and Frank for all the beers in Nancy and useful discussions about physics; and Renato, Alan, Jaime, Sally, and Rúbia for all the useful discussion about mathematics. 

Finally, let me thank the staff from Overleaf where I wrote this thesis and all my other manuscripts, Researchgate which keeps me updated about what my peers are doing, arXiv for its essential preprint service, and also Google Scholar for the excellent search engine and its track of citations. They definitely help in making science more democratic (in this respect, I would also like to thank other initiatives that fight for making research available for all). Lastly, let me also thank the financial support from the Brazilian agency Conselho Nacional de Desenvolvimento Científico e Tecnológico (CNPq) for the master and doctoral scholarships and also the Instituto Nacional de Matemática Pura e Aplicada (IMPA) for the financial support to participate in the XIX Brazilian School on Differential Geometry (2016) and XXXI Brazilian Mathematical Colloquium (2017).
\end{agradecimentos}

\begin{epigrafe}
    \vspace*{\fill}
	\begin{flushright}
		\textit{``La filosofia naturale è scritta in questo grandissimo libro che continuamente ci sta aperto innanzi agli occhi, io dico l’universo, ma non si può intendere se prima non s’impara a intender la lingua e conoscer i caratteri nei quali è scritto. Egli è scritto in lingua matematica, e i caratteri son triangoli, cerchi ed altre figure geometriche, senza i quali mezzi è impossibile a intenderne umanamente parola; senza questi è un aggirarsi vanamente per un oscuro labirinto''\\
		- Il Saggiatore, Galileo Galilei (1564-1642)}
	\end{flushright}
\end{epigrafe}


\begin{resumo}[Abstract]
 \begin{otherlanguage*}{english}
$\,\,\,\,\,\,\,\,\,\,\,\,\,\,\,$ This thesis is devoted to the differential geometry of curves and surfaces along with applications in quantum mechanics. In the 1st part we initially introduce the well known Frenet frame and then discuss on plane curves with a power-law curvature. Later, we show that the curvature function is a lower bound for the scalar angular velocity of any moving frame, from which one defines Rotation Minimizing (RM) frames as those frames that achieve this minimum. Remarkably, RM frames are ideal to study spherical curves and allow us to characterize them through a linear equation, in contrast with a differential equation from a Frenet approach. We also apply these ideas to curves that lie on level surfaces, $\Sigma=F^{-1}(c)$, by reinterpreting the problem in the context of a metric induced by $\mbox{Hess}\,F$, which may fail to be positive or non-degenerate and naturally leads us to a Lorentz-Minkowski $\mathbb{E}_1^3$ or isotropic $\mathbb{I}^3$ space. We then develop a systematic approach to construct RM frames and characterize spherical curves in $\mathbb{E}_1^3$ and $\mathbb{I}^3$ and furnish a criterion for a curve to lie on a level set surface. Finally, we extend these investigations to characterize curves that lie on the (hyper)surface of geodesic spheres in a Riemannian manifold. Using that for geodesic spherical curves the (radial) geodesics connecting the curve to a fixed point induce a normal vector field, we are able to characterize geodesic spherical curves in hyperbolic and spherical geometries through a linear equation. In the 2nd part we apply some of the previous ideas in the quantum dynamics of a constrained particle, where differential geometry is a relevant and timing tool due to the possibility of synthesizing nanostructures with non-trivial shapes. After describing the confining potential formalism, from which emerges a geometry-induced potential (GIP), we devote our attention to tubular surfaces as a mean to model curved nanotubes. The use of RM frames offers a simpler description for the constrained dynamics and allows us to show that the torsion of the centerline of a curved tube gives rise to a geometric phase. Later, we study the problem of prescribed GIP for curves and surfaces in Euclidean space: for curves it is solved by integrating Frenet equations, while for surfaces it involves a non-linear 2nd order PDE. Here we explore the GIP for surfaces invariant by a 1-parameter group of isometries, which turns the PDE into an ODE and leads to cylindrical, revolution, and helicoidal surfaces. The latter class is an important candidate to establish a link with chirality. Here we devote a special attention to helicoidal minimal surfaces and prove the existence of geometry-induced bound and localized states and the possibility of controlling the change in the probability density when the surface is subjected to an extra charge.   

   \noindent 
   \textbf{Key-words}: Frenet frame. Rotation minimizing frame. Spherical curve. 
   Constrained dynamics. Curved nanotube. 
   Prescribed curvature. 
 \end{otherlanguage*}
\end{resumo}

\setlength{\absparsep}{18pt} 
\begin{resumo}[Resumo]
$\,\,\,\,\,\,\,\,\,\,\,\,\,\,\,$ Esta tese \'e dedicada à geometria diferencial de curvas e superfícies e aplicações na mecânica quântica. Na 1ª parte introduzimos o conhecido triedro de Frenet e então estudamos curvas planas com curvatura dada por potências. Adiante, mostramos que a função curvatura é uma cota inferior para a velocidade de rotação de um referencial móvel qualquer, de onde se define Referenciais que Minimizam Rotação (RMR) como aqueles que atingem essa cota. Notavelmente, RMR são ideais no estudo de curvas esféricas e nos permitem caracterizá-las através de uma equação linear, em contraste com uma EDO em abordagens \`a Frenet. Também aplicamos essas ideias na caracterização de curvas em surperfícies de nível, $\Sigma=F^{-1}(c)$, reinterpretando o problema no contexto de uma métrica induzida por $\mbox{Hess}\,F$, que pode não ser positiva ou não-degenerada e então nos levar a um espaço de Lorentz-Minkowski $\mathbb{E}_1^3$ ou isotrópico $\mathbb{I}^3$. De forma unificada, construímos RMR e caracterizamos curvas esféricas em $\mathbb{E}_1^3$ e $\mathbb{I}^3$ e então fornecemos um critério para que curvas estejam em superfícies de nível. Finalmente, estedemos essas investigações a fim de caracterizar curvas na hiperfície de esferas geodésicas de uma variedade Riemanniana. Usando que para tais curvas as geodésicas (radiais) que ligam a curva a um determinado ponto fixo induz um campo de vetores normais, somos capazes de caracterizar curvas em esferas geodésicas em geometrias hiperbólica e esférica através de uma equação linear. Na 2ª parte aplicamos à dinâmica quântica de uma partícula confinada alguns dos conceitos já discutidos, onde a geometria diferencial é uma ferramenta relevante e atual devido à possibilidade de se sintetizar estruturas com formas não-triviais. Após descrever o formalismo do potencial confinante, de onde emerge um potencial induzido por geometria (PIG), nos dedicamos às superfícies tubulares a fim de modelar nanotubos curvos. O uso de RMR fornece uma descrição simples e nos permite mostrar que a torção do eixo do nanotubo dá origem a uma fase geométrica. Adiante, estudamos o problema de PIG prescrito: para curvas ele é solucionado integrando-se as equações de Frenet, enquanto para superfícies ele se escreve em termos de uma EDP não-linear de 2ª ordem. Aqui exploramos o PIG em superfícies invariantes por um grupo a 1 parâmetro de isometrias, que transforma a EDP do problema em uma EDO e nos leva ao estudo de superfícies cilíndricas, de revolução e helicoidais. Estas últimas são candidatas naturais para se estabelecer um link com quiralidade. Aqui dedicamos uma atenção especial às superfícies helicoidais mínimas e mostramos a existência de estados ligados e localizados induzidos por geometria e também a possibilidade de se controlar a distribuição de probabilidade ao submeter a superfície a uma carga extra.
 
   \noindent 
   \textbf{Palavras-chave}: Triedro de Frenet. Triedro que minimiza rotação. Curva esférica. 
   Dinâmica confinada. Nanotubo curvo.  
   Curvatura prescrita. 
\end{resumo}


\pdfbookmark[0]{\listfigurename}{lof}
\listoffigures*
\cleardoublepage



\pdfbookmark[0]{\contentsname}{toc}
\tableofcontents*
\cleardoublepage

\textual


\chapter[INTRODUCTION]{INTRODUCTION}

Before presenting the content of this thesis, let us first discuss on some philosophical viewpoints. In the course of his(er) studies in mathematics, the student is systematically exposed to more and more abstract ideas. With time, (s)he acquires some \textit{mathematical maturity} and has the possibility of doing useful and interesting things after mastering these abstract concepts \cite{Steen1983}. Our goal is not to criticize this learning process, probably this is a necessary evil, but we may wonder if this could explain why some (professional) mathematicians do think that research is an one-way road toward generality and abstractness: if a theorem does not offer very general informations about very general objects assuming the less possible assumptions, then it is not worth of our attention. Naturally, we are exaggerating, but the example above has its merits: it is natural to expect that the initials years have an influence on the professionals' viewpoints\footnote{Maybe, Stewart's book \cite{Stewart} may offer a better view into how mathematics and mathematicians function. See also \cite{DambeckSO2012,DambeckSO2016}.}. In this introduction we are neither interested in building a philosophy of mathematical learning/research nor in criticizing the abstraction in mathematics \cite{FerrariPTB2003}. Abstraction and the ability to handle it are essential to mathematics! Our point here is that research does not always function like this, sometimes it is important to revisit and work on classical themes and problems, e.g., on Fran\c cois Tr\`eves' words \cite{Treves} in the preface of one of his books\footnote{Tr\`eves was talking about the gap between what is taught to students and what they need in order to understand recent developments in research. Anyway, this passage can also serve our purposes.} ``[...] progress comes not only from pushing further and further into new territory but also from frequent return to the familiar grounds, from seeking an ever-deeper understanding of their nature, and finding there new inspiration and guidance''.

Such a return to ``well known'' subjects may have many origins, such as the development/improvement of techniques that allow for (i) a solution of unsolved problems or (ii) a new solution of already solved problems in classical subjects, and also (iii) applications, which may foster new advances and problems in research. Restricting ourselves to the differential geometry realm, we may mention the renewed interested on the theory of curves and surfaces in both three dimensional Euclidean and Non-Euclidean geometries\footnote{A look at the abstracts of the plenary talks in the last Brazilian Schools on Differential Geometry, for example, can testify this.}. In the applied arena, we can mention the advances in nanotechnology, that make the synthesis of curved structures a reality \cite{CastroNetoRepProgPhys,TerronesNewJPhys} and, therefore, requires the use of differential geometry of curves and surfaces tools for an appropriate modeling. In addition, we may mention applications in architecture \cite{KasapNNJ2016,LawrenceNNJ2011,Pottmann2007}, in computer graphics \cite{Farin2001,Pottmann2001}, and also in the emergent field of experimental mathematics \cite{Borwein2004}, in which case experiments using numerical examples or graphical images can help mathematicians in solving problems. 

This thesis is highly influenced by this viewpoint. Indeed, as the reader will be able to testify, here we are mainly interested in the differential geometry of curves and surfaces along with some applications in quantum mechanics.

\section{About the content of this thesis}

This work will be divided in two parts: part I - ``On the Differential Geometry of Rotation Minimizing Frames and Spherical Curves'' (from chapter 2 to 6); and part II - ``Applications in Physics: On the Quantum Mechanics of a Constrained Particle'' (from chapter 7 to 10). In essence, chapters 2 and 3 correspond to the content in \cite{daSilvaArXiv2017}; chapters 3, 4, and 5 correspond to \cite{daSilvaArXiv}; part of chapter 5 corresponds to \cite{daSilvaArXivIso2017};  chapter 6 corresponds to \cite{DaSilvaDeibsomDaSilvaArXiv2017}; chapters 7, 9, and 10 correspond to \cite{daSilvaAP2017}; and, finally, chapter 8 corresponds to a manuscript in progress with Fernando Santos (advisor), Fernando Moraes (UFRPE\footnote{Universidade Federal Rural de Pernambuco, Recife - Brazil.}), Bertrand Berche, and Sebastien Fumeron (both from Universit\'e de Lorraine, France). In the following we give a general picture about the content of this thesis.

\subsection{Differential geometry of curves}

The geometry of spheres is certainly one of the most important topic of investigation in differential geometry; the search for necessary and/or sufficient conditions for a submanifold be a sphere being one of its major pursuit. In this respect, a related and interesting problem then is that of characterizing curves that lie on the surface of a sphere. The focus of part I is on the study of the geometry of curves by means of moving frames defined along them\footnote{It is worth mentioning that moving frames here are not meant in the more general sense of Cartan \cite{GriffithsDMJ1973}. Indeed, we shall adopt a more elementary approach.} and applications to spherical curves. In chapter \ref{chap::FrenetCurves} we first review the well known Frenet frames and Frenet equations. Motivated by the search for curves whose geometry-induced potential is Hydrogen-like (see chapter \ref{chap::Vgip}) we investigate plane curves with power-law curvature function. We also describe the curvature and torsion of space curves in terms of spherical analogs, i.e., in terms of their osculating spheres: the use of osculating spheres will prove to be very useful when dealing with curves in isotropic space (sections \ref{sec::GeomIso} and \ref{Sec::DiffGeomIsoSpaces}). In chapter \ref{chap::RMframes} we introduce a general description of adapted frames along curves and show that the curvature function is a lower bound for the scalar angular velocity of any moving frame, from which we define Rotation Minimizing (RM) frames as those frames that achieve this minimum\footnote{In practice we identify this property by observing that $\{\mathbf{t},\mathbf{n}_1,\mathbf{n}_2\}$ is an RM frame if and only if $\mathbf{n}_i'$ is parallel to the tangent $\mathbf{t}$.}. Remarkably, these frames fit like a glove in the study of spherical curves\footnote{RM frames have other remarkable geometric properties, e.g., they are parallel transported along a curve with respect to the normal connection \cite{Etayo2016}; and a ruled surface along a curve is developable (zero Gaussian curvature) if and only if the rulings point in the direction of an RM vector field \cite{EtayoTJM2017,TuncerGMN2015}.}. Indeed, while the characterization of spherical curves in terms of a Frenet frame is made through a differential equation involving curvature and torsion, when we use an RM frame the characterization can be made through a linear equation and, in addition, such a characterization remains the same in higher dimensions, something that is not true in a Frenet-like approach.

Motived by the quest of a characterization of curves that lie on a given surface in Euclidean space, we study in chapter \ref{chap::curvesE13} the problem of defining RM frames and characterizing spherical curves in a Lorentz-Minkowski space, i.e., $\mathbb{R}^3$ equipped with an index 1 metric such as $\langle x,y\rangle_1=x_1y_1+x_2y_2-x_3y_3$. By equipping the neighborhood of a level set surface $\Sigma=F^{-1}(c)$ with a Hessian metric $h(\cdot,\cdot)_p=\langle\mbox{Hess}_pF\,\cdot,\cdot\rangle$ one is naturally led to the study of the differential geometry of curves in non-Riemannian spaces. Indeed, in general, a Hessian  $\mbox{Hess}\,F=\partial^2F/\partial x^i\partial x^j$ may fail to be positive
or non-degenerate, and then leads us to the study of the differential geometry of curves in the just mentioned Lorentz-Minkowski space (chapter \ref{chap::curvesE13}) and in isotropic  (section \ref{sec::GeomIso}) space\footnote{It is worth mentioning that such geometries are particular instances of the so called Cayley-Klein geometries \cite{Sulanke2006}, as mentioned in section \ref{sec::GeomIso}.}, i.e., $\mathbb{R}^3$ equipped with a degenerate metric $\langle x,y\rangle_i=x_1y_1+\delta\,x_2y_2$, where $\delta\in\{-1,0,+1\}$. Although simple, this idea proves to be very useful \cite{daSilvaArXiv}, as we will made clear in chapter \ref{Chap_CurvInSurf}. 

Finally, in chapter \ref{chap::RMFinRiemGeom} we extend these investigations for curves on geodesic spheres in $\mathbb{S}^{m+1}(r)$ and $\mathbb{H}^{m+1}(r)$, the $(m+1)$-dimensional sphere and hyperbolic space of radius $r$, respectively. An important observation is that in Euclidean space spherical curves are normal curves, and vice-versa: since $\langle\alpha-p,\alpha-p\rangle=\,\mbox{constant}\,\Leftrightarrow\langle\mathbf{t},\alpha-p\rangle=0$ it follows that, up to a translation, the position vector of a spherical curve lies on the normal plane. Such an equivalence makes sense due to the double nature of $\mathbb{R}^{m+1}$ as both a manifold and as a tangent space\footnote{This problem has to do with the more general quest of studying curves that lie on a given (moving) plane generated by two chosen vectors of a moving trihedron, e.g., one would define osculating, normal or rectifying curves as those curves whose position vector, up to a translation, lies on their osculating, normal or rectifying planes, respectively \cite{ChenMonthly2003,ChenAJMS2017}. It is known that (i) osculating curves are precisely the plane curves (if we substitute the principal normal by an RM vector field, we still have a characterization for plane curves \cite{daSilvaArXiv2017}: see chapter \ref{chap::RMframes}), (ii) normal curves are precisely the spherical curves, and (iii) rectifying curves are precisely geodesics on a cone \cite{ChenManuscript,ChenAJMS2017}.}. In order to extend these notions to a Riemannian setting one should replace the line segment $\alpha(s)-p$ by a geodesic connecting $p$ to a point $\alpha(s)$, as pointed out by Lucas and Ortega-Yag\"ues in the study of rectifying curves \cite{LucasJMAA2015,LucasMJM2016}(\footnote{They proved that rectifying curves in the 3d sphere and hyperbolic space are geodesics on a conical surface (in analogy with the Euclidean case).}). We show, as a consequence of the Gauss lemma for the exponential map in a Riemannian manifold $M^{m+1}$, that on a sufficiently small neighborhood of $p\in M^{m+1}$ a curve $\alpha:I\to M^{m+1}$ is normal (with center $p$) if and only if it lies on a geodesic sphere (with center $p$) in $M^{m+1}$. Using this equivalence in $\mathbb{S}^{m+1}(r)$ and $\mathbb{H}^{m+1}(r)$ we are able to characterize those curves that lie on the (hyper)surface of a geodesic sphere through a linear equation involving the coefficients (curvatures) that dictate an RM frame motion. For completeness, we also discuss in this work the characterization of geodesic spherical curves in terms of a Frenet frame (theorem \ref{thr::FrenetChar3DSphCurv}) and show that the characterization of (geodesic) spherical curves is the same as in Euclidean space. Lastly, the relation between totally geodesic submanifolds, which play the role of planes in Riemannian geometry, and curves with a normal development curve  $(\kappa_1,...\,,\kappa_m)$ lying on a line passing through the origin is more delicate, since in general a manifold has no totally geodesic submanifold up to the trivial ones \cite{MurphyArXiv2017,NikolauevskyIJM2015,Tsukada1996}. Nonetheless, in this work we are able to show that if a Riemannian manifold contains totally geodesic submanifolds, then any curve on a totally geodesic submanifold is associated with a normal development that lies on a line passing through the origin (theorem \ref{thr::NormDevelopRiemPlaneCurv}). We show in addition that a curve in $\mathbb{S}^{m+1}(r)$ and $\mathbb{H}^{m+1}(r)$ lies on a totally geodesic submanifold if and only if its normal development is a line passing through the origin (theorem \ref{thr::NormDevelopS3andH3PlaneCurv}).

\subsection{Quantum dynamics of a constrained particle}

In part II we apply some of the theoretical framework developed in the first part in the quantum dynamics of a constrained particle. This is an important topic in contemporary research due to the many advances in the experimental techniques in nanotechnology, which demand a better understanding of how the geometry of nanostructures may influence its chemical and physical properties \cite{LahiffABC2010,NovoselovNature2005}. In this respect, the differential geometry of curves and surfaces may be a valuable tool in describing the dynamics of a particle constrained to a curve/surface \cite{DaCostaPRA1981,JensenKoppeAnnPhys}, which possibly models a given nanostructure, and also in understanding how the nanostructure geometry and their properties interact. This offers the possibility of engineering structures with certain physical properties prescribed \textit{a priori} through a geometric approach \cite{delCampoSciRep2014,daSilvaAP2017,FernandosEfranceses}, which is crucial in order to create new technologies. Chapter \ref{chap::Vgip} is devoted to discussing the fundamentals of the dynamics of a particle confined on a curve/surface. The history begins with De Witt's attempt to approach the problem through a quantization procedure in the intrinsic coordinates of the constraint region. The resulting equations suffer however
from an ordering ambiguity \cite{DeWittRMP1957}. The point here is that in order to do a more realist and ``correct'' modeling for the confinement it is necessary to take into account that the curved region is embedded somewhere. Indeed, a confining potential approach taking into account how the constraint region is embedded in ambient space does not suffer
from such a problem: the confining potential (or extrinsic) approach gives a unique effective Hamiltonian to the constrained dynamics \cite{DaCostaPRA1981,JensenKoppeAnnPhys}. It is shown that a geometry-induced potential (GIP) acts upon the dynamics and that it depends on both intrinsic and extrinsic geometric quantities, e.g., for the constrained dynamics on a surface the GIP is \cite{DaCostaPRA1981,JensenKoppeAnnPhys}
$$V_{gip}=-\frac{\hbar^2}{2m^*}(H^2-K),$$
where $m^*$ is the mass of the constrained particle, $\hbar=h/2\pi$ ($h$ being the Planck constant), and $H$ and $K$ are the mean and Gaussian curvatures of the constraint surface, respectively. 

In chapter \ref{chap::TubularSurf}, we then devote our attention to tubular surfaces. This can serve as a mean to model nanotubes, which play an important role in modern nanotechnology \cite{TakeuchiNANOSYST2014}. A tubular surface can be geometrically constructed by moving a circle of fixed radius (located on the normal plane) along a given curve. The use of RM frames offers a simpler description for the equations of motion and, in addition, it allows us to show that the curve torsion gives rise to a geometric phase \cite{BerryPRSL1984,BerryPT1990,ColinDeVerdiere2006}, which is an important ingredient in phenomena such as the Aharanov-Bohm effect \cite{AharanovBohm,FechnerPUZ1998,TwistedRingProgTheorPhys}. 

In chapter \ref{chap::PrescVgip} we address the problem of prescribed GIP for curves and
surfaces in Euclidean space $\mathbb{R}^3$, i.e., how to find a curved region with
a potential given \textit{a priori}. For curves this is easily solved
by integrating Frenet equations, which is a system of linear 1st order ODE's, while the problem for surfaces involves a non-linear 2nd order PDE. A comprehensive study of these PDE's is not a trivial task and, in addition, it can encode in its generality useless examples. In this respect, the study of particular classes can turn to be more useful and insightful than a general analysis. In fact, in most physical systems of interest it is always supposed some kind of symmetry. Here, we explore the GIP for surfaces invariant by a 1-parameter group of isometries of $\mathbb{R}^3$, which leads to cylindrical, revolution, and
helicoidal surfaces \cite{DoCarmoTohoku1982,MedeirosRMU1991}. It also has the advantage of turning the PDE's into ODE's and then making a general study a realistic goal \cite{daSilvaAP2017}. Finally, in chapter \ref{chap_ConstDynHelSurf} we discuss in more detail the important class of helicoidal surfaces. This represents the most general kind of invariant surfaces and is particularly important
since helicoidal surfaces are natural candidates to establish a link between chirality and a GIP. For the family of helicoidal minimal surfaces, we prove the existence of geometry-induced bound and localized states and the possibility of controlling the change in the distribution of the probability density when the surface is subjected to an extra charge \cite{daSilvaAP2017}. 

In short, in this 2nd part of the thesis we believe we contributed to a better understanding of the geometrical aspects of the quantum constrained dynamics in showing how to control the geometry-induced potential, which is a fundamental step toward future potential applications of this formalism.

\section{List of publications}

Finally, let us mention that a large part of this thesis is based on the content of a few manuscripts. Under the supervision of Fernando A. N. Santos we did:
\begin{itemize}
\item \textbf{DA SILVA, L. C. B.}; BASTOS, C. C.; RIBEIRO, F. G. Quantum mechanics of a constrained particle and the problem of prescribed geometry-induced potential. Annals of Physics, v. 379, p. 13,
2017. \url{http://dx.doi.org/10.1016/j.aop.2017.02.012}
\item \textbf{DA SILVA, L. C. B.}; SANTOS, F. A. N.; MORAES, F.; FUMERON, S.; BERCHE, B. Quantum mechanics of particles confined on surfaces: symbolic approach and application to curved tubes. (Working paper).
\end{itemize}

In addition, there was the opportunity of conducting some independent research:
\begin{itemize}
\item \textbf{DA SILVA, L. C. B.} Moving frames and the characterization of curves that lie on a surface. Journal of Geometry, 2017. \url{http://dx.org/10.1007/s00022-017-0398-7}
\item \textbf{DA SILVA, L. C. B.} Characterization of spherical and plane curves using rotation
minimizing frames. e-print, 2017. \url{https://arxiv.org/abs/1706.01577v3}.
\item \textbf{DA SILVA, L. C. B.} Rotation minimizing frames and spherical curves in simply isotropic and semi-isotropic 3-spaces. e-print, 2017. \url{https://arxiv.org/abs/1707.06321}.
\item \textbf{DA SILVA, L. C. B.}; DEIBSOM DA SILVA, J. Characterization of curves that lie on a geodesic sphere in the $(m+1)$-dimensional sphere and hyperbolic space. e-print, 2017. \url{https://arxiv.org/abs/1707.07335}.
\end{itemize}

Finally, in parallel to this thesis we have completed under the supervision of F. A. N. Santos and Maur\'icio D. Coutinho-Filho (co-advisor for the master degree at UFPE):
\begin{itemize}
\item BASTOS, C. C.; \textbf{DA SILVA, L. C. B.}; SANTOS, F. A. N. Semi-empirical and \textit{ab initio} calculations for twisted M\"obius strips molecular models. In: Anais do XIII Encontro da SBPMat, João Pessoa, v. 1, p. 99, 2014.
\item SANTOS, F. A. N.; \textbf{DA SILVA, L. C. B.}; COUTINHO-FILHO, M. D. Topological approach to microcanonical
thermodynamics and phase transition
of interacting classical spins. Journal of Statistical Mechanics: Theory and Experiments, v. 2017, p. 013202, 2017. \url{http://dx.doi.org/10.1088/1742-5468/2017/1/013202}.
\end{itemize}

\part{ON THE DIFFERENTIAL GEOMETRY OF ROTATION MINIMIZING FRAMES AND SPHERICAL CURVES}
\chapter{DIFFERENTIAL GEOMETRY OF CURVES IN EUCLIDEAN SPACE}
\label{chap::FrenetCurves}

Let $\mathbb{E}^3$ denote the three dimensional \emph{Euclidean space}, i.e., $\mathbb{R}^3$ equipped with the standard metric $\langle x,y\rangle=\sum_{i=1}^3x_iy_i$. In addition, we equip $\mathbb{E}^3$ with a norm $\Vert \mathbf{x}\Vert=\sqrt{\langle\mathbf{x},\mathbf{x}\rangle}$.

A regular curve (of class $C^k$) is a function $\alpha:I\rightarrow\mathbb{E}^{3}$ satisfying $\alpha'\not=0$ and that has a continuous derivative of order $k$. Here, $I$ is an interval of $\mathbb{R}$, possibly infinite, or $\mathbb{S}^1$ in the case of a closed curve. We say that $\alpha(t)$ is \emph{parametrized by arc-length} if $\langle\alpha(u),\alpha(u)\rangle=1$. If this is not the case, we should reparameterize it according to
\begin{equation}
s(t) = \int_{t_0}^t\,\sqrt{\langle\alpha(u),\alpha(u)\rangle}\,\mathrm{d}u\,.
\end{equation}
We say that $s$ is an \emph{arc-length parameter} and we have $\Vert\alpha'(s)\Vert=1$.

\section{Frenet frame and Frenet equations for a plane curve}

It can be easily proved that a curve $\alpha$ is a straight line if and only if its acceleration vector, $\mathbf{a}=\alpha''$, vanishes identically. When parametrized by an arc-length parameter $s$, the acceleration vector still gives information about the bending of $\alpha$. Indeed, along $\alpha$ we may introduce the \emph{unit tangent}
\begin{equation}
\mathbf{t}(s) = \frac{\mathrm{d}\alpha(s)}{\mathrm{d}s}\,.
\end{equation}
Since $\langle\mathbf{t},\mathbf{t}\rangle=1$, we have that $\mathbf{t}$ is orthogonal to $\mathbf{t}'$. If $\Vert\mathbf{t}'\Vert\not=0$, we can define a normal vector field along $\alpha$ as
\begin{equation}
\mathbf{n}(s)=\frac{\mathbf{t}'(s)}{\Vert\mathbf{t}'(s)\Vert}\,.
\end{equation}
We call $\mathbf{n}$ the (\emph{principal}) \emph{normal} and $\kappa(s)=\Vert\mathbf{t}'(s)\Vert=\Vert\alpha''(s)\Vert$ the \emph{curvature function} of the curve $\alpha$.

The unit tangent and principal normal are linearly independent. So, if $\alpha$ is a plane curve, i.e., $\alpha\subset\mathbb{E}^2$, we can write $\mathbf{n}'=a\,\mathbf{t}+b\,\mathbf{n}$. As $\Vert\mathbf{n}\Vert=1$, it follows that $b=0$. On the other hand, since $\langle\mathbf{t},\mathbf{n}\rangle=0$, one has $a=\langle\mathbf{n}',\mathbf{t}\rangle=-\langle\mathbf{n},\mathbf{t}'\rangle=-\kappa$. In short, the equation of motion for the plane \emph{Frenet (moving) frame} $\{\mathbf{t},\mathbf{n}\}$ is
\begin{equation}
\frac{\mathrm{d}}{\mathrm{d}s}\left(
\begin{array}{c}
\mathbf{t}\\
\mathbf{n}\\
\end{array}
\right)=\left(
\begin{array}{cc}
0 & \kappa\\
-\kappa & 0\\
\end{array}
\right)\left(
\begin{array}{c}
\mathbf{t}\\
\mathbf{n}\\
\end{array}
\right).\label{eq::FrenetEqsPlaneCurv}
\end{equation}

It is not difficult to prove that the curvature $\kappa$ is invariant by rotations and translations in $\mathbb{E}^2$, i.e., it is invariant by rigid motions. Straight lines have $\kappa=0$, while circles are the only curves with constant curvature, $\kappa=1/R$ (sometimes it is convenient to see a line as a circle with infinity radius). In addition, given a function $\tilde{\kappa}(\tilde{s})>0$ there exists a unique curve $\alpha$, up to rigid motions, with curvature $\kappa=\tilde{\kappa}$ and arc-length parameter $s=\tilde{s}$. The parametrization of the solution curve for the Frenet equations is \cite{Struik}
\begin{equation}
\left\{
\begin{array}{c}
x(s)=z_{1}\,C(s)-z_{2}\,S(s)+x_0\\[6pt]
y(s)=z_{1}\,S(s)+z_{2}\,C(s)+y_0\\
\end{array}
\right.\,,
\end{equation}
where $x_0,\,y_0,$ and $z_{i}$ are constants to be specified by the initial conditions and
\begin{equation}
\left\{
\begin{array}{c}
S(s)=+\displaystyle\int_{s_0}^s\cos\Big(\int_{s_0}^v\kappa(u)\,\mathrm{d}u\Big)\mathrm{d}v\\[8pt]
C(s)=-\displaystyle\int_{s_0}^s\sin\Big(\int_{s_0}^v\kappa(u)\,\mathrm{d}u\Big)\mathrm{d}v\\
\end{array}
\right.\,.\label{eq::SenoCosFrenetFrame}
\end{equation}

The idea of approximating a curve by a simpler one is very fruitful. This gives rise to the concept of order of contact.

\begin{definition}
 We say that two regular curves $\alpha$ and $\beta$ in $\mathbb{E}^n$ have a \emph{contact of order $k$} at $\alpha(s_0)=\beta(s_0^*)$ if all the higher order derivatives, up to order $k$, also coincide: 
 \begin{equation}
 \forall\,i\in\{1,...\,,k\},\,\frac{\mathrm{d}^i\alpha(s_0)}{\mathrm{d}s^i}=\frac{\mathrm{d}^i\beta(s_0^*)}{\mathrm{d}(s^*)^i}.
 \end{equation}
\end{definition}

For example, the tangent line has a contact of order 1 with its reference curve. The circle that has a contact of order 2 at $\alpha(s_0)$ is called the \emph{osculating circle} and its radius is $\rho_c(s_0)=1/\kappa(s_0)$ \cite{Kreyszig1991,Struik}. The center of the osculating circle $P_c(s)$ at $\alpha(s)$ is
\begin{equation}
P_c(s) = \alpha(s) +\frac{1}{\kappa(s)}\mathbf{n}(s)\,.\label{eq::OscCircle}
\end{equation}
The \emph{radius of curvature} is $\rho_c=1/\kappa$. When $\kappa(s_0)=0$, it means that the tangent line has a contact of order 2 at $\alpha(s_0)$ and in this case we may say that $\rho_c=+\infty$. 

\subsection{Plane curves with power-law curvature function}
\label{subsec::PowerLawCurvature}

Let us consider plane curves with a power law curvature function, i.e., $\tau\equiv0$ and $\kappa(s)=c_0/s^p$, where $c_0>0$ and $p\in\mathbb{R}$ are constants\footnote{Theses curve can be used to find a plane ``Hydrogen curve'', i.e., a curve whose geometry-induced potential is $V_{gip}=e^2(4\pi\varepsilon_0\vert s\vert)^{-1}$, where $p=1/2$ and $c_0=\sqrt{8me^2/4\pi\varepsilon_0\hbar^2}$. See chapter \ref{chap::PrescVgip} for details.}. As can be easily verified, the solutions of Eq. (\ref{eq::FrenetEqsPlaneCurv}) for the power-law case is given by 
\begin{equation}
\left\{\begin{array}{c}
t_i(s) = a_i\, C_p(s) + b_i\, S_p(s)\\[5pt]
n_i(s) = (s^p/c_0)t'_i(s) = -a_i\,S_p(s)+b_i\,C_p(s)\\
\end{array}
\right.,\,i\in\{1,2\},
\end{equation}
where $\mathbf{t}=(t_1,t_2)$, $\mathbf{n}=(n_1,n_2)$ are the tangent and principal normal, $a_i,b_i$ are constants for all $i$, and we have defined
\begin{equation}
C_p(s) = \left\{
\begin{array}{ccc}
\displaystyle\cos\left(\frac{c_0s^{1-p}}{1-p}\right) & , & p \not= 1\\[8pt]
\cos\left(c_0\ln s\right)  & , & p = 1\\
\end{array}
\right.;
\,\,
S_p(s) = \left\{
\begin{array}{ccc}
\displaystyle\sin\left(\frac{c_0s^{1-p}}{1-p}\right) & , & p \not= 1\\[8pt]
\sin\left(c_0\ln s\right)  & , & p = 1\\
\end{array}
\right..
\end{equation}
By imposing initial conditions $\mathbf{t}_0=\mathbf{t}(s_0)=(t_{1,0},t_{2,0})$ and $\mathbf{n}_0=\mathbf{n}(s_0)=(n_{1,0},n_{2,0})$, we find
\begin{equation}
\left(
\begin{array}{c}
t_i \\
n_i \\
\end{array}
\right)=\mathcal{R}_p(s)^{T}\,\mathcal{R}_p(s_0)\left(
\begin{array}{c}
t_{i,0}\\
n_{i,0}\\
\end{array}
\right);\,\mathcal{R}_p(s) = \left(
\begin{array}{cc}
C_p(s) & -S_p(s)\\
S_p(s) & \,C_p(s)\\
\end{array}
\right),
\end{equation}
where we have used a ``rotation'' matrix $\mathcal{R}_p(s)$.

\begin{figure*}[t]
\centering
\includegraphics[width=0.5\textwidth]{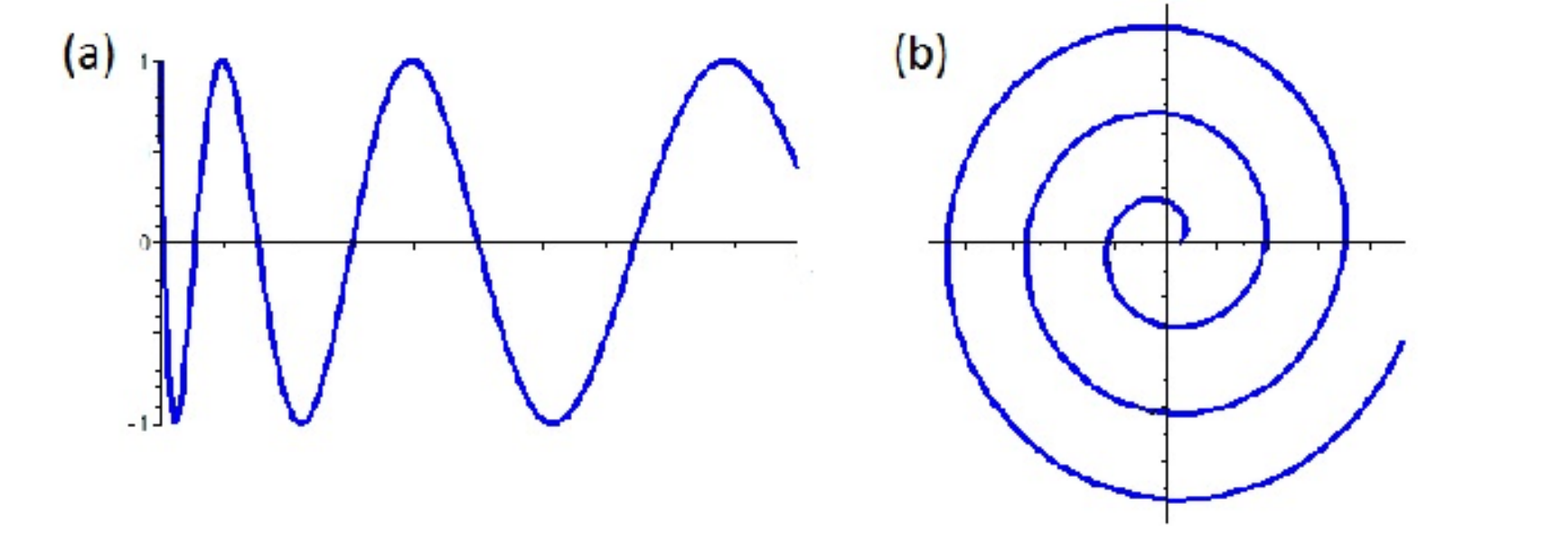}
\caption{\textbf{(a)} Plot of $C_{1/2}(s)$; \textbf{(b)} Plot of a plane Hydrogen curve.}
\label{fig:Hcurve}
\end{figure*}

These results show us that the Frenet frame $\{\mathbf{t}(s),\mathbf{n}(s)\}$ of such a plane curve rotates with $s$. Although the functions $C_p(s)$ and $S_p(s)$ are not periodic, the point described by $(C_p(s),S_p(s))$ moves along a circle which rotates with a non-constant angular velocity, then one would see $C_p(s)$ and $S_p(s)$ as an almost periodic functions whose period depends on $s$, see Fig. \ref{fig:Hcurve}.(a). So, any curve with a power-law curvature $\kappa(s)=c_0/s^{p}$ will displays a kind of almost periodic character.

To obtain the curve parametrization, we must integrate the functions $C_p(s)$ and $S_p(s)$, since $\alpha(s)=\int_{s_0}^s\,\mathbf{t}(u)\,\mathrm{d}u$. To the best of our knowledge, for a general value of $p$, this integration can not be expressed in terms of elementary functions, so from now on, we restrict the discussion to the case of our main interest, namely, a Hydrogen potential $\Leftrightarrow p=1/2$:
\begin{eqnarray}
\int C_{1/2}(u)\,du & = & +\sqrt{\frac{u}{c_0^2}}\,S_{1/2}(u)+\frac{C_{1/2}(u)}{2c_0^2}+c_1\,;\\
\int S_{1/2}(u)\,du & = & -\sqrt{\frac{u}{c_0^2}}\,C_{1/2}(u)+\frac{S_{1/2}(u)}{2c_0^2}+c_2\,,
\end{eqnarray}
where $c_1,c_2$ are arbitrary constants.

Since the integration of Frenet equations gives a unique curve up to rigid motions we are free to choose the initial conditions. Assuming for simplicity $\mathbf{t}(s_0)=(1,0)$ and $\mathbf{n}(s_0)=(0,1)$, we can then write
\begin{equation}
\mathbf{t}(s)=\mathcal{R}_{1/2}(s_0)
\left(
\begin{array}{c}
C_{1/2}(s)\\
S_{1/2}(s)\\
\end{array}
\right)
\,\Rightarrow\,
\alpha(s)=\mathcal{R}_{1/2}(s_0)\left(
\begin{array}{cc}
\frac{1}{2c_0^2} & \frac{\sqrt{s}}{c_0}\\[8pt]
-\frac{\sqrt{s}}{c_0} & \frac{1}{2c_0^2}\\
\end{array}
\right)\left(
\begin{array}{c}
C_{1/2}(s)\\
S_{1/2}(s)\\
\end{array}
\right)+c_3\,,
\end{equation}
where $c_3$ is a constant. By the uniqueness up to rigid motions, we can choose $c_3=0$ (by translating the curve) and ignore the factor $\mathcal{R}_{1/2}(s_0)$ (by rotating the curve). In essence, such a curve $\alpha(s)$ is obtained by the superposition of a circle of radius $1/2c_0^2$ and a spiral curve whose distance to the origin increases as $\sqrt{s}/c_0$, see Fig. \ref{fig:Hcurve}.(b).

\section{Frenet frame and Frenet equations for a space curve}

The definition of the unit tangent $\mathbf{t}$ and the principal normal $\mathbf{n}$ is the same for a space curve. The \emph{binormal} vector is defined as
\begin{equation}
\mathbf{b}(s) = \mathbf{t}(s)\times \mathbf{n}(s),
\end{equation}
where $\times$ is the usual cross product in $\mathbb{E}^3$. The trihedron $\{\mathbf{t},\mathbf{n},\mathbf{b}\}$ is the \emph{Frenet frame} of $\alpha$ . As the Frenet frame is an orthonormal base for $\mathbb{E}^3$, its equation of motion is characterized by a screw-symmetric matrix. By definition, we have $\mathbf{t}'=\kappa\mathbf{n}$. On the other hand, let us write $\mathbf{b}'=a\mathbf{t}+b\mathbf{n}$ (since $\Vert\mathbf{b}\Vert=1$, the vector $\mathbf{b}'$ has no component in $\mathbf{b}$). Using the definition of the binormal vector, we have
\begin{equation}
\mathbf{b}' = \mathbf{t}'\times\mathbf{n}+\mathbf{t}\times\mathbf{n}'=\mathbf{t}\times\mathbf{n}'\Rightarrow\langle\mathbf{b}',\mathbf{t}\rangle=\langle\mathbf{t}\times\mathbf{n}',\mathbf{t}\rangle=0.
\end{equation}
So, the derivative of the binormal must be parallel to the principal normal: $\mathbf{b}'=-\tau\mathbf{n}$. We call $\tau$ the \emph{torsion} of $\alpha$.

The Frenet equations, i.e., the equation of motion of the Frenet frame, are written as
\begin{equation}
\frac{\mathrm{d}}{\mathrm{d}s}\left(
\begin{array}{c}
\mathbf{t}\\
\mathbf{n}\\
\mathbf{b}\\
\end{array}
\right)=\left(
\begin{array}{ccc}
0 & \kappa & 0\\
-\kappa & 0 & \tau\\
0 & -\tau & 0\\ 
\end{array}
\right)\left(
\begin{array}{c}
\mathbf{t}\\
\mathbf{n}\\
\mathbf{b}\\
\end{array}
\right).\label{eq::FrenetFrameEqs}
\end{equation}

As in the case of the curvature, the torsion $\tau$ is invariant by rotations and translations, i.e., by rigid motions. In addition, it measures how much the curve deviates from being planar, i.e., $\alpha$ is a plane curve if and only if $\tau=0$ \cite{Kreyszig1991,Struik}. We say that a point $\alpha(s)$ is \emph{twisted} if at $s$ we have a non-zero curvature and a non-zero torsion (intuitively, the curve is truly three dimensional at $s$).

The Frenet equation uniquely determines a curve up to rigid motions, i.e., given two functions $\tilde{\kappa}(\tilde{s})>0$ and $\tilde{\tau}(\tilde{s})$, integration of the Frenet equations determines a unique curve $\alpha$, up to rigid motions, with arc-length parameter $s=\tilde{s}$, curvature function $\kappa=\tilde{\kappa}$, and torsion $\tau=\tilde{\tau}$: the curve parametrization is given by $\alpha(s)=\alpha_0+\int_{s_0}^s\,\mathbf{t}(u)\mathrm{d}u$.

\begin{definition}
 We say that a regular curve $\alpha$ and a surface $\Sigma$ in $\mathbb{E}^3$ have a \emph{contact of order $k$} at $\alpha(s_0)\in\Sigma$ if there exists a curve $\beta\subset\Sigma$ such that they have a contact of order $k$ with  $\alpha$. For a level set surface $\Sigma=G^{-1}(c)$, where $\mathrm{grad}_p\,G\not=0$ for $p\in \Sigma$, this condition is equivalent to say that
 \begin{equation}
 \forall\,i\in\{1,...\,,k\},\,\frac{\mathrm{d}^i\gamma(s^*_0)}{\mathrm{d}(s^*)^i}=0\,,\mathrm{ where }\,\gamma(s^*)=(G\circ\alpha)(s^*)\,,\,\gamma(s_0^*)=\alpha(s_0).
 \end{equation}
\end{definition}

The plane $\mbox{span}\{\mathbf{t}(s),\mathbf{n}(s)\}$ is known as the \emph{osculating plane}. It has a contact of order 2 with $\alpha$ at $\alpha(s)$. The osculating plane contains the velocity and acceleration vectors $\alpha'$ and $\alpha''$, regardless of the curve parameter. Indeed, let $t$ be a generic regular parameter for $\alpha$. Then $\mathrm{d}s/\mathrm{d}t=v(t)=\Vert\alpha'(t)\Vert$ and  $\mathbf{t}(t)=\mathbf{v}(t)/v(t)$. It follows that
\begin{equation}
\frac{\mathrm{d}\mathbf{t}}{\mathrm{d}s}=\frac{\mathrm{d}t}{\mathrm{d}s}\frac{\mathrm{d}\mathbf{t}}{\mathrm{d}t}=\frac{\mathbf{a}}{v^2}-\dot{v}\frac{\mathbf{v}}{v^3}\,,
\end{equation}
where $\mathbf{v}(t)=\dot{\alpha}(t)$ and $\mathbf{a}(t)=\ddot{\alpha}(t)$ are the velocity and acceleration vectors with respect to $t$. We can alternatively write
\begin{equation}
\mathbf{a} = \dot{v}\,\mathbf{t}+v^2\kappa\,\mathbf{n}\,.
\end{equation}
Kinetically speaking, the above relation means that every motion is locally a circular motion with respect to the osculating plane: $\dot{v}$ is the instantaneous scalar velocity and $v^2\kappa=v^2/\rho$ the instantaneous centripetal acceleration.

Going back to geometry, we can use that the binormal vector is orthogonal to the osculating plane in order to write, for a generic regular parameter, the following expression
\begin{equation}
\mathbf{b}=\frac{\dot{\alpha}\times\ddot{\alpha}}{\Vert\dot{\alpha}\times\ddot{\alpha}\Vert}.
\end{equation}
Then, the Frenet frame according to a generic regular parameter is 
\begin{equation}
\mathbf{t} = \frac{\dot{\alpha}}{\Vert\dot{\alpha}\Vert},\,\mathbf{b}=\frac{\dot{\alpha}\times\ddot{\alpha}}{\Vert\dot{\alpha}\times\ddot{\alpha}\Vert},\mbox{ and }\,\mathbf{n}=\mathbf{b}\times\mathbf{t}.
\end{equation}
It is possible to show that the curvature function and torsion are
\begin{equation}
\kappa=\frac{\Vert\dot{\alpha}\times\ddot{\alpha}\Vert}{\Vert\dot{\alpha}\Vert^3}\,\,\mathrm{ and }\,\,\tau=\frac{\langle\dot{\alpha}\times\ddot{\alpha},\dddot{\alpha}\rangle}{\Vert\dot{\alpha}\times\ddot{\alpha}\Vert^2}\,.\label{eq::CurvAndTorsion}
\end{equation}

\section{Spherical curves and osculating spheres}

An osculating sphere at $\alpha(s)$ is a sphere that has an order 3 contact with its reference curve. At a zero torsion point the osculating plane has a contact of order 3 and so we may say that the osculating sphere at this point has an infinity radius. At a twisted point, the center $P_S(s)$ and radius $R_S(s)$ of the osculating sphere are respectively given by \cite{Kreyszig1991}
\begin{equation}
P_S = \alpha+\rho\,\mathbf{n}+\frac{1}{\tau}\frac{\mathrm{d}\rho}{\mathrm{d}s}\,\mathbf{b}\,\mbox{ and }\,R_S = \sqrt{\rho^2+\frac{1}{\tau^2}\left[\rho'\right]^2}\,\,,
\end{equation}
where $\rho=1/\kappa$ is the radius of curvature.

It is possible to characterize spherical curves through their Frenet frame. Indeed, since the osculating spheres must be always the same, we have \cite{Kreyszig1991,Kuhnel2010}

\begin{theorem}
Let $\alpha:I\to \mathbb{E}^3$ be a $C^4$ regular curve with a non-zero torsion. It lies on a sphere of radius $r$ if and only if
\begin{equation}
\varsigma(s)=\tau(s)\rho(s)+\frac{\mathrm{d}}{\mathrm{d}s}\left(\frac{\rho'(s)}{\tau(s)}\right)=0.\label{eq::SphAnnihilator}
\end{equation}
\end{theorem}

In the following we express the curvature function and torsion of a spherical curve in term of its spherical curvature \cite{SabanRendLincei1958}. This will be essential in the study of a generic curve near an osculating sphere.
 
\begin{theorem}
Let $\alpha:I\to \mathbb{S}^2(p,r)$ be a spherical $C^3$ curve parametrized by an arc-length $s$. Then the curvature and torsion  are respectively given by
\begin{equation}
\kappa = \frac{1}{r}\sqrt{1+J^2}\,\mbox{ and }\,\tau = \frac{J'}{1+J^2}\,,
\end{equation}
where $J=\langle\alpha-p,\alpha'\times\alpha''\rangle$ is the \emph{spherical curvature} \cite{SabanRendLincei1958}.
\label{thr::KappaTauJforSpheres}
\end{theorem}
\textit{Proof. } Assume $p$ to be the origin (the general case is reduced to this one by studying $\tilde{\alpha}=\alpha-p$). The vectors $\alpha/r,\,\alpha'$, and $(\alpha/r)\times\alpha'$ form an orthonormal frame along the curve\footnote{Such frames are also known in the mathematical literature as \emph{Saban frames} \cite{IzumiyaRM2017}. They were introduced by Giacomo Saban in the characterization of spheres through the vanishing of $\oint\rho^n\tau\,\mathrm{d}s$ for every closed curve on the surface \cite{SabanRendLincei1958}, which generalizes the case $n=0$, i.e., $\oint\tau\,\mathrm{d}s$, in Ref. \cite{ScherrerVNGZ1940}.}. Write
\begin{equation}
\alpha'' = \frac{1}{r}\langle\alpha'',\alpha\rangle\frac{\alpha}{r}+\langle\alpha'',\alpha'\rangle\,\alpha'+\frac{1}{r}\langle\alpha'',\alpha\times\alpha'\rangle\frac{\alpha}{r}\times\alpha'\,.
\end{equation}
Since $\alpha$ is parametrized by arc-length, we have $\langle\alpha'',\alpha'\rangle=0$. In addition, from $\langle\alpha,\alpha\rangle=r^2$, it follows that $\langle\alpha'',\alpha\rangle=-\langle\alpha',\alpha'\rangle=-1$. In conclusion, the acceleration vector gives
\begin{equation}
\alpha'' = -\frac{1}{r}\,\frac{\alpha}{r}+\frac{J}{r}\,\frac{\alpha}{r}\times\alpha'\,\Rightarrow\,\kappa=\Vert\alpha''\Vert=\frac{1}{r}\sqrt{1+J^2}\,.
\end{equation}

Now, writing the normal and binormal vectors as $\mathbf{n}=\alpha''/\kappa$ and $\mathbf{b}=\alpha'\times\mathbf{n}=\alpha'\times\alpha''/\kappa$, and using the Frenet equation $\tau=-\langle\mathbf{b}',\mathbf{n}\rangle$, we have
\begin{eqnarray}
\tau & = & -\left\langle\frac{\textrm{d}}{\textrm{d}s}\left(\frac{\alpha'\times\alpha''}{\kappa}\right),\frac{\alpha''}{\kappa}\right\rangle\nonumber\\
& = & -\frac{1}{\kappa^2}\langle\alpha'\times\alpha''',\alpha''\rangle\nonumber\\
& = & -\frac{1}{r^2\kappa^2}\langle\alpha'\times\alpha''',-\alpha+J\,\alpha\times\alpha'\rangle.
\end{eqnarray}
Finally, using the expression for $\kappa$ above,  $J'=\langle\alpha,\alpha'\times\alpha'''\rangle$, and the vector identity 
$
\langle \mathbf{A}\times\mathbf{B},\mathbf{C}\times\mathbf{D}\rangle = \langle \mathbf{A},\mathbf{C}\rangle\langle\mathbf{B},\mathbf{D}\rangle-\langle \mathbf{A},\mathbf{D}\rangle\langle\mathbf{B},\mathbf{C}\rangle
$, we get the desired result for the torsion: $\tau = J'/(1+J^2)$.
\qed

\begin{remark}
It follows that a $C^3$ curve on a sphere has no inflection points. So, its torsion is always well defined. If we drop the $C^3$ condition and allows the torsion to possibly have discontinuities, it follows from the Darboux theorem \cite{OlsenMonthly2004} that on a point $\tau(s_0)$ it is not possible to have distinct values for the lateral limits $\tau_0^-=\lim_{s\to s_0^{-}}\tau(s)$ and $\tau_0^+=\lim_{s\to s_0^{+}}\tau(s)$, since $\tau$ is the derivative of another function: i.e., for a spherical curve $\tau_0^+$ and $\tau_0^-$ either exist and coincide or one of them does not exist\footnote{For a generic analytic curve, the lateral limits $\tau_0^-$ and $\tau_0^+$ do exist and coincide \cite{HordMonthly1972}. But if we drop the analyticity assumption, one of the lateral limit may diverge even for a $C^{\infty}$ curve \cite{HordMonthly1972}. Indeed, analytic curves are well behaved with respect to inflection points: given two analytic functions $K\geq0$ and $\tau$, there exists an analytic curve, up to a rigid motion in $\mathbb{E}^3$, with curvature $\kappa=\sqrt{K}$ and torsion $\tau$ \cite{SasaiTMJ1984}.}.
\end{remark}

Using the concept of osculating spheres, we would intuitively say that every curve is locally spherical. Now we investigate how to extend Theorem \ref{thr::KappaTauJforSpheres} for a generic curve. Let $\alpha$ be a regular twisted curve and $\Sigma_s = \mathbb{S}^2(P_S(s),R_S(s))$ be its osculating sphere at $\alpha(s)$. Near a fixed point $\alpha(s_0)$ we can obtain a spherical curve $\beta: (s_0-\epsilon,s_0+\epsilon)\to\Sigma_{s_0}$ by projecting $\alpha$ on $\Sigma_{s_0}$ according to
\begin{equation}
\beta(t)=r_0\frac{\alpha(t)-a_0}{\Vert\alpha(t)-a_0\Vert}\,,\label{eq::OscSphericalProjection}
\end{equation}
where $a_0=P_S(s_0)$ and $r_0=R_S(s_0)$. So we have \cite{daSilvaArXiv2017}
\begin{theorem}
The torsion $\tau_{\alpha}$ and the curvature $\kappa_{\alpha}$ of a $C^3$ regular twisted curve $\alpha$ and the torsion $\tau_{\beta}$ and the curvature $\kappa_{\beta}$ of its (osculating) spherical projection $\beta$ coincide at $s_0$: 
\begin{equation}
\kappa_{\alpha}(s_0) = \kappa_{\beta}(s_0)\,\mbox{ and }\,\tau_{\alpha}(s_0) = \tau_{\beta}(s_0).
\end{equation}
In addition, we can write the curvature function as
\begin{equation}
\kappa_{\alpha}(s_0) = \frac{1}{R_S(s_0)}\sqrt{1+J^2(s_0)}\,,
\end{equation}
and the torsion as
\begin{equation}
\tau_{\alpha}(s_0) = \frac{\langle\alpha(s_0)-P_S(s_0),\alpha'(s_0
)\times\alpha'''(s_0)\rangle}{1+J^2(s_0)}=\frac{J'(s_0)}{1+J^2(s_0)}+\frac{\kappa(s_0)\varsigma(s_0)}{1+J^2(s_0)}\,,
\end{equation}
where $J(s)=\langle\alpha(s)-P_S(s),\alpha'(s)\times\alpha''(s)\rangle$ and $\varsigma$ is defined in Eq. (\ref{eq::SphAnnihilator}).
\label{thr::tauUsingOscSphere}
\end{theorem}
\textit{Proof. } In order to compute $\tau_{\beta}$ and $\kappa_{\beta}$ it is enough to find $\beta'$, $\beta''$, and $\beta'''$. Calculating the derivatives of $\beta$ and 
taking into account the relations 
\begin{equation}
\left\{
\begin{array}{ccc}
\langle\,\alpha(s_0)-a_0,\mathbf{t}(s_0)\,\rangle & = &0\\ 1+\kappa_{\alpha}(s_0)\,\langle\,\alpha(s_0)-a_0,\mathbf{n}(s_0)\,\rangle & = & 0\\ \kappa'_{\alpha}(s_0)-\kappa_{\alpha}^2(s_0)\tau_{\alpha}(s_0)\,\langle\,\alpha(s_0)-a_0,\mathbf{b}(s_0)\,\rangle & = & 0\\
\end{array}
\right.
\end{equation}
satisfied by an osculating sphere \cite{Kreyszig1991}, we obtain after some lengthy but straightforward calculations the following relations at $s=s_0$
\begin{equation}
\left\{
\begin{array}{ccc}
\beta'(s_0) &=& \mathbf{t}(s_0)\\
\beta''(s_0) &=& \kappa_{\alpha}(s_0)\mathbf{n}(s_0)\\
\beta'''(s_0)&=& -\kappa_{\alpha}^2(s_0)\mathbf{t}(s_0)+\kappa_{\alpha}'(s_0)\mathbf{n}(s_0)+\tau_{\alpha}(s_0)\kappa_{\alpha}(s_0)\mathbf{b}(s_0)\\
\end{array}
\right.\,,
\end{equation}
where $\{\mathbf{t},\mathbf{n},\mathbf{b}\}$ are the Frenet frame of $\alpha$. It follows that $\alpha(s_0)=\beta(s_0)$, $\alpha'(s_0)=\beta'(s_0)$, $\alpha''(s_0)=\beta''(s_0)$, and $\alpha'''(s_0)=\beta'''(s_0)$ (this is not a surprise, since an osculating sphere has a  contact of order 3 with its reference curve).

Substituting the expressions above for $\beta',\beta''$, and $\beta''$, in the equation for the torsion and curvature and using the equalities $\alpha^{(i)}(s_0)=\beta^{(i)}(s_0)$ ($i=0,1,2,3$) gives the desired result: $\tau_{\beta}=\tau_{\alpha}$, $\kappa_{\beta}=\kappa_{\alpha}$. 

Now let us express $\kappa$ and $\tau$ in terms of $J$. First, observe that $\beta$ is not necessarily parametrized by arc-length, so we must adapt the expressions in Theorem \ref{thr::KappaTauJforSpheres}. For a curve with a generic regular parameter $t$, we can write
\begin{equation}
\frac{\mathrm{d}\beta}{\mathrm{d}s}=\frac{1}{v}\frac{\mathrm{d}\beta}{\mathrm{d}t}\,\mbox{ and }\,\frac{\mathrm{d}^2\beta}{\mathrm{d}s^2}=-\frac{1}{v^3}\frac{\mathrm{d}v}{\mathrm{d}t}\frac{\mathrm{d}\beta}{\mathrm{d}t}+\frac{1}{v^2}\frac{\mathrm{d}^2\beta}{\mathrm{d}t^2}\,,
\end{equation}
where $v(t)=\Vert\beta'(t)\Vert$. Thus, it follows that
\begin{equation}
j(t)=j(s_{\beta}(t))=\frac{1}{v^3(t)}\langle\,\beta(t)-P_S(s_0),\,\beta'(t)\times\beta''(t)\,\rangle,
\end{equation}
where $s_{\beta}$ is the arc-length parameter of $\beta$ and $j(s_{\beta})=\langle\,\beta(s_{\beta})-P_S(s_0),\,\beta'(s_{\beta})\times\beta''(s_{\beta})\,\rangle$. 

Finally, applying the expression from Theorem \ref{thr::KappaTauJforSpheres} to the spherical curve $\beta$ and using that at $t=s_0$ one has $j(t=s_0)=J(s_0)$ and $v(t=s_0)=1$, we find the expressions for $\kappa_{\alpha}$ and $\tau_{\alpha}$: in the second equality for $\tau_{\alpha}$, we should use the fact that $\textrm{d}P_S/\textrm{d}s=\varsigma\,\mathbf{b}$ in order to conclude that $J'=\langle\alpha-P_S,\alpha'\times\alpha'''\rangle-\kappa\varsigma$.
\qed

\chapter{ROTATION MINIMIZING FRAMES}
\label{chap::RMframes}

In the previous chapter we studied the geometry of curves in $\mathbb{E}^3$ through the use of Frenet frames.  This is the usual way of studying the geometry of curves, but by no means it represents the only choice. Indeed, since the principal normal always points to the center of curvature, it may result in unnecessary rotation and then making the use of a Frenet frame unsuitable in some contexts. In this respect, the consideration of \emph{rotation minimizing frames} may represent a more suitable choice. Due to their minimal twist, rotation minimizing frames are of fundamental importance in many branches, such as in camera \cite{FaroukiCAVW2009,JaklicCAGD2013} and rigid body motions \cite{FaroukiMC2012,FaroukiCAGD2014}, robotics \cite{WebsterIJRR2010}, fluid flow \cite{GermanoJFM1982,HuttlIJHFF2000,VashisthIECR2008}, quantum mechanics \cite{DaCostaPRA1981,HaagAHP2015,TwistedRingProgTheorPhys}, integrable systems \cite{SandersMMJ2003}, visualization \cite{BanksIEEETVCG,HansonTechrep1995} and deformation of tubes \cite{Li2009,LiCGF2010}, in mathematical biology in the study of DNA \cite{ChirikjianBST2013,ClauvelinJCTC2012} and protein folding \cite{HuPRE2011}, sweep surface modeling \cite{Bloomenthal1991,PottmannIJSM1998,SiltanenCGF1992,WangCAD1997}, and in differential geometry as well \cite{BishopMonthly,daSilvaArXiv,daSilvaArXiv2017,EtayoTJM2017,OzdemirMJMS2008}, just to name a few. In this thesis we will be interested in these frames and their applications in geometry and quantum mechanics\footnote{In chapter \ref{chap::curvesE13} and in section \ref{sec::GeomIso} we discuss the extension of such a concept in Lorentz-Minkowski space, non-degenerate and index one metric, and isotropic space, degenerate metric, respectively. In chapter \ref{chap::RMFinRiemGeom} we extend it to a Riemannian ambient space.}.

\section{Adapted moving frames and velocity of rotation}

Besides the Frenet frame, we may consider any other orthonormal trihedron $\{\mathbf{e}_0,\mathbf{e}_1,\mathbf{e}_2\}$ along $\alpha$. We say that a frame is \emph{adapted} to $\alpha$ if $\mathbf{e}_0=\mathbf{t}$ and $\langle\mathbf{e}_i,\mathbf{e}_j\rangle=\delta_{ij}$. We suppose in the remaining of this work that all the frames are adapted. The orthonormality condition implies that the frame motion is characterized by a skew-symmetric matrix, i.e., 
\begin{equation}
\frac{\mathrm{d}}{\mathrm{d}s}\left(
\begin{array}{c}
\mathbf{e}_0\\
\mathbf{e}_1\\
\mathbf{e}_2\\
\end{array}
\right)=\left(
\begin{array}{ccc}
 0        & \chi_{2} & -\chi_{1}\\
-\chi_{2} & 0        & \omega   \\
 \chi_{1} & -\omega  & 0        \\ 
\end{array}
\right)\left(
\begin{array}{c}
\mathbf{e}_0\\
\mathbf{e}_1\\
\mathbf{e}_2\\
\end{array}
\right).\label{eq::FrameEqs}
\end{equation}
By defining the \emph{Darboux vector} $\mathbf{w}=\omega\mathbf{e}_0+\chi_1\mathbf{e}_1+\chi_2\mathbf{e}_2$, the equation of motion above can be written in a compact form as
\begin{equation}
\frac{\mathrm{d}\mathbf{e}_i}{\mathrm{d}s}=\mathbf{w}\times\mathbf{e}_i\,,\,i\in\{0,1,2\}\,.
\end{equation}
The above equation means that the Darboux vector is the angular velocity of a rigid body with axes $\mathbf{e}_0$, $\mathbf{e}_1$, and $\mathbf{e}_2$ \cite{Thornton2008}. Then, the scalar angular velocity of the trihedron $\{\mathbf{e}_i\}_i$ is $w=\sqrt{\omega^2+\chi_1^2+\chi_2^2}$. 

We may ask how does the trihedron $\{\mathbf{e}_i\}_i$ relates to the Frenet one. Since the curvature function $\kappa$ only depends on the derivatives $\mathbf{e}_0=\alpha'$ and $\mathbf{e}_0'=\alpha''$, Eq. (\ref{eq::CurvAndTorsion}), we have the relation
\begin{equation}
\kappa=\Vert\mathbf{e}_0'\Vert=\Vert\chi_2\mathbf{e}_1-\chi_1\mathbf{e}_2\Vert=\sqrt{\chi_1^2+\chi_2^2}\,.
\end{equation}
It follows from the expression above that
\begin{prop}
The scalar angular velocity $w$ of \emph{any} adapted moving trihedron satisfies
\begin{equation}
w=\sqrt{\omega^2+\chi_1^2+\chi_2^2}\geq\,\kappa\,.
\end{equation}
\label{prop::relationAngVelAdaptedFrame}
\end{prop}
This proposition tells us that $\kappa$ is a minimum for the rotation of any frame. On the other hand, for the torsion we also need $\mathbf{e}''_0$, Eq. (\ref{eq::CurvAndTorsion}), which can be written as
\begin{equation}
\mathbf{e}_0'' = -\kappa^2\mathbf{e}_0+(\chi_2'+\omega\chi_1)\mathbf{e}_1-(\chi_1'-\omega\chi_2)\mathbf{e}_2  \,.
\end{equation}
Then
\begin{equation}
\tau=\frac{\omega\kappa^2+\chi_1\chi_2'-\chi_1'\chi_2}{\kappa^2}\Rightarrow\omega=\tau-\frac{\chi_1\chi_2'-\chi_1'\chi_2}{\kappa^2}.
\end{equation}
Now, defining $\kappa\mathrm{e}^{-\mathrm{i}\theta}=\chi_2+\mathrm{i}\,\chi_1$, we deduce that $\theta'=(\chi_1\chi_2'-\chi_1'\chi_2)\kappa^{-2}$ and, therefore,
\begin{prop} \cite{TwistedRingProgTheorPhys}
A generic adapted moving frame $\{\mathbf{e}_i\}$ can be obtained from the Frenet one by a rotation:
\begin{equation}
\left\{
\begin{array}{ccc}
\mathbf{e}_1(s) &=& \cos\theta(s)\,\mathbf{n}(s)-\sin\theta(s)\,\mathbf{b}(s)\\
\mathbf{e}_2(s) &=& \sin\theta(s)\,\mathbf{n}(s)+\cos\theta(s)\,\mathbf{b}(s)\\
\end{array}
\right.\,.
\end{equation}
In addition, one has the following relations
\begin{equation}
\left\{
\begin{array}{ccc}
\chi_2(s) & = &  \kappa(s)\cos\theta(s)\\
\chi_1(s) & = & -\kappa(s)\sin\theta(s)\\
\omega(s) & = &  \tau(s)-\theta'(s)\\
\end{array}
\right.\,.
\end{equation}
\label{prop::RelationBetwAdaptedFrameAndFrenet}
\end{prop}

Combining Propositions \ref{prop::relationAngVelAdaptedFrame} and \ref{prop::RelationBetwAdaptedFrameAndFrenet}, it follows that a moving frame with $\theta=\int\tau$ minimizes rotation, i.e., one must discount the unnecessary rotation associated with $\tau$ in order to minimize rotation.
\begin{definition}
A moving frame $\{\mathbf{t},\mathbf{n}_1,\mathbf{n}_2\}$ is said to be a \emph{rotation minimizing} (RM) \emph{frame} if $\theta'(s)=\tau(s)$. The equation of motion for an RM must be
\begin{equation}
\frac{\mathrm{d}}{\mathrm{d}s}\left(
\begin{array}{c}
\mathbf{t} \\
\mathbf{n}_1\\
\mathbf{n}_2\\
\end{array}
\right)=\left(
\begin{array}{ccc}
 0        & \kappa_{1} & \kappa_{2}\\
-\kappa_{1} & 0        & 0   \\
 -\kappa_{2} & 0  & 0        \\ 
\end{array}
\right)\left(
\begin{array}{c}
\mathbf{t}\\
\mathbf{n}_1\\
\mathbf{n}_2\\
\end{array}
\right),\label{eq::BishopEqs}
\end{equation}
i.e., $\chi_2=\kappa_1$, $\chi_1=-\kappa_2$. and $\omega=0$. A normal vector field $\mathbf{u}$ satisfying $\mathbf{u}'=\lambda\mathbf{t}$ is said to be a \emph{rotation minimizing vector field}. 
\end{definition}

\begin{remark}
Due to their remarkable properties, RM frames have been independently discovered several times\footnote{Here we do not attempt to furnish a complete list of ``discoveries''. Our list is probably incomplete.}, such as in the study of PDE's on tubular neighborhoods  \cite{DaCostaPRA1981,GermanoJFM1982,TangIEEE1970} and in computer graphics \cite{KlokCAGD1986}. However, Bishop seems to be the first to exploit their geometric implications \cite{BishopMonthly} (albeit he named them \emph{relatively parallel frames}). In addition, it can be proved that an RM vector field is parallel transported along $\alpha(s)$ with respect to the normal connection of the curve \cite{Etayo2016}: if $\alpha(s_i)=\alpha(s_f)$, $\mathbf{n}_1(s_i)$ will differ from $\mathbf{n}_1(s_f)$ by an angular amount of $\Delta \theta = \int_{s_i}^{s_f} \tau(x)\,\mathrm{d}x$. Due to this last properties, some authors call RM frames \emph{parallel frames}, see e.g. \cite{HansonTechrep1995,OzdemirMJMS2008}.
\label{remark::TerminologyRM}
\end{remark}

It is worth mentioning that RM frames can be globally defined even if the curve has points with zero curvature \cite{BishopMonthly}. In addition, RM frames are not uniquely defined, since any rotation of $\mathbf{n}_i$ on the normal plane still gives an RM field, i.e., the angle $\theta$ is well defined up to an additive constant. But most importantly, the prescription of curvatures $\kappa_1,\kappa_2$ still determines a curve up to rigid motions \cite{BishopMonthly}. Two RM frames $\{\mathbf{t},\mathbf{n}_1,\mathbf{n}_2\}$ and $\{\mathbf{t},\tilde{\mathbf{n}}_1,\tilde{\mathbf{n}}_2\}$  are related by
\begin{equation}
\left\{
\begin{array}{ccc}
\tilde{\mathbf{n}}_1(s) &=& \cos\theta_0\,\mathbf{n}_1(s)-\sin\theta_0\,\mathbf{n}_2(s)\\
\tilde{\mathbf{n}}_2(s) &=& \sin\theta_0\,\mathbf{n}_1(s)+\cos\theta_0\,\mathbf{n}_2(s)\\
\tilde{\kappa}_1(s) &=& \cos\theta_0\,\kappa_1(s)-\sin\theta_0\,\kappa_2(s)\\
\tilde{\kappa}_2(s) &=& \sin\theta_0\,\kappa_1(s)+\cos\theta_0\,\kappa_2(s)
\end{array}
\right.\,.
\end{equation}
Consequently, it follows from the two last equations above that the \emph{normal development curve} $(\kappa_1,\kappa_2)$ preserves its shape when we choose a new RM frame.

\section{Characterization of spherical curves}

Interestingly, RM frames allows for a simple characterization of spherical curves. Indeed\footnote{An attempt to extend these ideas in order to characterize curves that lie on a surface was devised in \cite{daSilvaArXiv}. See Chapter \ref{Chap_CurvInSurf}.}
\begin{theorem} \cite{BishopMonthly}
A regular $C^2$ curve $\alpha:I\to \mathbb{E}^3$ lies on a sphere of radius $r$ if and only if its normal development, i.e., the curve $(\kappa_1(s),\kappa_2(s))$, lies on a line not passing through the origin. In addition, the distance of this line to the origin is $r^{-1}$. \label{thr::charSphCurves}
\end{theorem}
\textit{Proof. } Let $\alpha$ be a spherical curve and $\{\mathbf{t},\mathbf{n}_1,\mathbf{n}_2\}$ an RM frame along it. Since $\langle\alpha-P,\alpha-P\rangle=r^2$, it follows that $\langle\mathbf{t},\alpha-P\rangle=0$. Then, we can write
\begin{equation}
\alpha-P = a_1\mathbf{n}_1+a_2\mathbf{n}_2.
\end{equation}
We have the relation $a_i=\langle\alpha-P,\mathbf{n}_i\rangle$, whose derivative is $a_i'=\langle\mathbf{t},\mathbf{n}_i\rangle+\langle\alpha-P,-\kappa_i\mathbf{t}\rangle=0$. Thus, $a_1$ and $a_2$ are constants. Now, deriving $\langle\alpha-P,\mathbf{t}\rangle=0$ gives $\langle\mathbf{t},\mathbf{t}\rangle+\langle\kappa_1\mathbf{n}_1+\kappa_2\mathbf{n}_2,a_1\mathbf{n}_1+a_2\mathbf{n}_2\rangle=0$ and, therefore, $a_1\kappa_1+a_2\kappa_2+1=0$. Finally, we have $r^2=\langle\alpha-P,\alpha-P\rangle=a_1+a_2=\mathrm{dist}\left(\{\sum a_i\kappa_i+1\},(0,0)\right)^{-2}$.

Conversely, let $\alpha$ be a curve such that its RM curvatures satisfy $\sum a_i\kappa_i+1=0$. Defining $P(s)=\alpha(s)-a_1\mathbf{n}_1-a_2\mathbf{n}_2$, one has $P'=\mathbf{t}+a_1\kappa_1\mathbf{t}+a_2\kappa_2\mathbf{t}=0$ and, therefore, $P$ is a fixed point. We can easily deduce that $\langle\alpha-P,\alpha-P\rangle=a_1^2+a_2^2=$ constant. Then, $\alpha$ is a spherical curve.
\qed

Here we also furnish a proof for the above result by using osculating spheres. But first, let us describe its parametrization by using an RM frame. Indeed, we can write
\begin{equation}
P_S(s_0)=\alpha(s_0)+\beta_0\mathbf{t}(s_0)+\beta_1 \mathbf{n}_1(s_0)+\beta_2\mathbf{n}_2(s_0).
\end{equation}
Now, defining a function $g(s)=\langle P_S-\alpha(s),P_S-\alpha(s)\rangle-r^2$, we have
\begin{eqnarray}
g ' & = & -2\langle P_S-\alpha,\mathbf{t}\rangle = -2\beta_0\,,\\
g'' & = & 2\langle\mathbf{t},\mathbf{t}\rangle-2\langle P_S-\alpha,\kappa_1\mathbf{n}_1+\kappa_2\mathbf{n}_2\rangle = -2(-1+\kappa_1\beta_1+\kappa_2\beta_2)\,,\\
g''' & = & -2\langle P_S-\alpha,\sum_i(\kappa_i'\mathbf{n}_i-\kappa_i^2\mathbf{t})\,\rangle = -2\sum_i(\kappa_i'\beta_i-\kappa_i^2\beta_0)\,.
\end{eqnarray}
Imposing the order 3 contact condition leads to $g'(s_0)=g''(s_0)=g'''(s_0)=0$ and gives
\begin{equation}
\beta_0 = 0,\,\kappa_1(s_0)\beta
_1+\kappa_2(s_0)\beta_2-1=0,\,\mbox{ and }\kappa_1'(s_0)\beta_1+\kappa_2'(s_0)\beta_2=0.\label{eq::coefOfOscSphUsingRMF}
\end{equation}
Thus, the coefficients $\beta_0$, $\beta_1$, and $\beta_2$ as functions of $s_0$ are
\begin{equation}
\beta_0=0,\,\beta_1=\frac{\kappa_2'}{\kappa_1\kappa_2'-\kappa_1'\kappa_2}=\frac{\kappa_2'}{\tau\kappa^2},\,\mbox{ and }\beta_2=-\frac{\kappa_1'}{\kappa_1\kappa_2'-\kappa_1'\kappa_2}=-\frac{\kappa_1'}{\tau\kappa^2}\,,
\end{equation}
where in the equalities above we used the relation between $(\kappa_1,\kappa_2)$ and $(\kappa,\tau)$.
\newline
\newline
\textit{Proof of Theorem \ref{thr::charSphCurves} for $C^4$ curves\footnote{We need a $C^4$ condition in order to compute $\kappa_i''$: $C^1$ is enough to have $\mathbf{t}$; $C^2$ to have $\mathbf{t}'$ and then $\kappa_i$; and $C^3$ ($C^4$) to have $\kappa_i'$ ($\kappa_i''$).}.} Taking the derivative of the osculating center gives
\begin{equation}
P_S' = \frac{\mathrm{d}}{\mathrm{d}s}\left(\alpha+\frac{\kappa_2'}{\tau\kappa^2}\mathbf{n}_1-\frac{\kappa_1'}{\tau\kappa^2}\mathbf{n}_2\right)\nonumber\\
=\left(\frac{\mathrm{d}}{\mathrm{d}s}\frac{\kappa_2'}{\tau\kappa^2}\right)\mathbf{n}_1-\left(\frac{\mathrm{d}}{\mathrm{d}s}\frac{\kappa_1'}{\tau\kappa^2}\right)\mathbf{n}_2\,.
\end{equation}
From the linear independence of $\{\mathbf{n}_1,\mathbf{n}_2\}$ we conclude that $P'_S=0$, i.e., $\alpha$ is spherical, if and only if $\beta_1$ and $\beta_2$ are constants. From Eq. (\ref{eq::coefOfOscSphUsingRMF}), this is equivalent to say that the normal development lies on a line not passing through the origin.
\qed

\begin{remark}
The approach above has some weaknesses when compared with that of Bishop \cite{BishopMonthly}. Indeed, the use of osculating spheres demands that the curve must be $C^4$ and also that $\tau\not=0$, while in Bishop's approach one needs just a $C^2$ condition and no restriction on the torsion. However, the use of osculating spheres will prove to be very useful when dealing with curves in isotropic space, see section \ref{sec::GeomIso}. 
\end{remark}

By using the chain rule and that $(\arctan x)'=(1+x^2)^{-1}$, we can use the results in Theorem \ref{thr::KappaTauJforSpheres} in order to find $\theta$, the angle between the principal normal and an RM vector, for a spherical curve.

\begin{prop}
Let $\alpha:I\to \mathbb{S}^2(p,r)$ be a spherical $C^3$ curve parametrized by arc-length $s$, then the angle $\theta$ between a rotation minimizing vector and the principal normal satisfies
\begin{equation}
\theta(s_2) - \theta(s_1) = \arctan\,J(s_2)-\arctan\,J(s_1)\,,
\end{equation}
where $J=\langle\alpha-p,\alpha'\times\alpha''\rangle$ is the \emph{spherical curvature} \cite{SabanRendLincei1958}.
\end{prop}

Another important observation is that for a spherical curve $\alpha\subset \mathbb{S}^2(p,r)$ the normals to the sphere along $\alpha$ minimize rotation, i.e., the normalized position vector $\mathbf{N}=(\alpha-p)/r$ is an RM vector field. Indeed, $\frac{\mathrm{d}}{\mathrm{d}s}(\alpha-p)/r=-(-1/r)\,\alpha'$ (this is an important step in the implementation of the \textit{double reflection method} for computing approximations of RM frames \cite{WangACMTOG}). The curvature $\kappa_1$ associated with $\mathbf{n}_1=(\alpha-p)/r$ is then $\kappa_1=-1/r$. For the other RM vector field, i.e., $\mathbf{n}_2=\alpha'\times(\alpha-p)/r$, one has $\kappa_2=-J/r$. 

The observation above furnishes an alternative proof for the fact that the total torsion $\oint\tau(x)\mathrm{d}x$ vanishes for all closed curves on a sphere \cite{SabanRendLincei1958,ScherrerVNGZ1940}. Indeed, for a closed curve one has
\begin{equation}
\oint\tau(x)\mathrm{d}x=\theta(s_f)-\theta(s_i)\,.
\end{equation}
As the position vector $\mathbf{N}=(\alpha-P)/r$ is the same at the initial and final points $\alpha(s=s_i)$ and $\alpha(s=s_f)$, respectively, the result follows.

In addition, using the concept of osculating spheres, we would intuitively say that every curve is locally spherical. In this case, it is tempting to ask if the normals to the osculating spheres minimize rotation. Unfortunately, this strategy does not work unless the curve is spherical:
\begin{prop}
If $\alpha:I\to\mathbb{E}^3$ is a regular twisted curve of class $C^4$, then
\begin{equation}
\frac{\mathrm{d}}{\mathrm{d}s}\left(\frac{\alpha(s)-P_S(s)}{R_S(s)}\right)=\frac{1}{R_S}\left[\mathbf{t}+\frac{\rho\rho'}{\tau R_S^2}\varsigma\,\mathbf{n}+\left(\frac{\rho'\,^2}{\tau^2R_S^2}-1\right)\varsigma\,\mathbf{b}\right]\,,\label{eq::DerOfNormalToOscSphere}
\end{equation}
where $\rho=\kappa^{-1}$ and 
\begin{equation}
\varsigma(s)=\tau(s)\rho(s) +\frac{\mathrm{d}}{\mathrm{d}s}\left(\frac{\rho'(s)}{\tau(s)}\right).
\end{equation}
In addition, the normal vector field to the curve given by the normals to the osculating sphere along $\alpha(s)$ minimizes rotation if and only if $\alpha$ is spherical, i.e., when $\varsigma\equiv0$.
\end{prop}
\textit{Proof. } Using that $R_S'=\rho'\varsigma/\tau R_S$ and $P_S'=\varsigma\mathbf{b}$, direct computation of the derivative of $\mathbf{N}=(\alpha-P_S)/R_S$ leads to Eq. (\ref{eq::DerOfNormalToOscSphere}). Finally, by a known result of geometry, the condition to be spherical leads to $\varsigma\equiv0$ \cite{Kreyszig1991,Kuhnel2010}, which by direct examination of Eq. (\ref{eq::DerOfNormalToOscSphere}) is a necessary and sufficient condition to have $\mathbf{N}$ and $\mathbf{t}$ parallel.
\qed

\subsection{Higher dimensional spherical curves}
\label{subsec::DdimSphCurvFrenet}

One can extrapolate the definition of RM frames from three dimensions to $\mathbb{E}^{m+1}$. Then, we say that $\{\mathbf{t},\mathbf{n}_1,...\,,\mathbf{n}_m\}$ is an RM frame along $\alpha:I\to\mathbb{E}^{m+1}$ if 
\begin{equation}
\frac{\mathrm{d}}{\mathrm{d}s}\left(
\begin{array}{c}
\mathbf{t}\\
\mathbf{n}_1\\
\vdots\\
\mathbf{n}_m\\
\end{array}
\right)=\left(
\begin{array}{cccc}
0 & \kappa_{1} & \cdots & \kappa_{m}\\
-\kappa_{1} & 0 & \cdots & 0\\
\vdots & \vdots & \ddots & \vdots \\
-\kappa_{m} & 0 & \cdots & 0\\ 
\end{array}
\right)\left(
\begin{array}{c}
\mathbf{t}\\
\mathbf{n}_1\\
\vdots \\
\mathbf{n}_m\\
\end{array}
\right).\label{eq::nDBishopEqs}
\end{equation}

It is not difficult to see that we can easily generalize theorem \ref{thr::charSphCurves} to higher dimensions. The proof is left as an exercise to the reader:
\begin{theorem}
A regular $C^2$ curve $\alpha:I\to \mathbb{E}^{m+1}$ lies on a sphere of radius $r$ if and only if its normal development, i.e., the curve $(\kappa_1(s),...\,,\kappa_m(s))$, lies on a line not passing through the origin. In addition, the distance of this line to the origin is $r^{-1}$. \label{thr::charSphCurvesNd}
\end{theorem}

On the other hand, this is not the case for a Frenet frame approach. The approach to higher dimensions is not straightforward and, in addition, the higher the dimension, the higher the differentiability the curve must have. A $C^{m+2}$ regular curve $\alpha:I\to\mathbb{R}^{m+1}$ is said to be a \emph{Frenet curve} if $\{\alpha',\alpha'',...\,,\alpha^{(m)}\}$ is a linearly independent set along all the points of $\alpha$. The Frenet frame $\{\mathbf{e}_0,\mathbf{e}_1,...\,,\mathbf{e}_m\}$ is such that
\begin{enumerate}
\item $\{\mathbf{e}_0,\mathbf{e}_1,...\,,\mathbf{e}_m\}$ is an adapted orthonormal moving frame;
\item for each $k=0,...\,,m$ one has $\mbox{span}\{\mathbf{e}_0,...\,,\mathbf{e}_k\}=\mbox{span}\{\alpha',\alpha'',...\,,\alpha^{(k+1)}\}$; and
\item $\langle \alpha^{(k+1)},\mathbf{e}_{k+1}\rangle>0$ for  $k=0,...\,,m$.
\end{enumerate}
The Frenet equations are
\begin{equation}
\left\{
\begin{array}{ccc}
\mathbf{e}_0' & = & \tau_0\mathbf{e}_1\\
\mathbf{e}_i' & = &-\tau_{i-1}\mathbf{e}_{i-1}+\tau_{i}\mathbf{e}_{i+1}\\
\end{array}
\right.,\,\,1\leq i \leq m,
\end{equation}
where $\tau_0 =\kappa$ and $\tau_{m}=0$. As in $\mathbb{E}^3$, it is possible to find expressions for the Frenet apparatus in $\mathbb{E}^{m+1}$ \cite{GluckMonthly1966,GutkinJGP2011}.

In order to study spherical curves in $\mathbb{E}^{m+1}$, the crucial observation is that they are normal curves
\begin{definition}
A regular curve $\alpha:I\to\mathbb{E}^{m+1}$ is said to be a \emph{normal curve} if its position vector lies, up to a translation, on the normal plane to the curve, i.e., 
\begin{equation}
\alpha(s)-p\in\mathrm{span}\{\mathbf{t}(s)\}^{\perp}\,,
\end{equation}
where $p$ is a fixed point.
\end{definition}
\begin{prop}
A regular curve in $\mathbb{E}^{m+1}$ is a normal curve if and only if it is a spherical curve\footnote{This equivalence in Euclidean space is quite silly, but it is very useful in the study of spherical curves in non-Euclidean settings, such as in affine geometry \cite{KreyszigPAMS1975} and also in the $(m+1)$-dimensional sphere $\mathbb{S}^{m+1}(r)$ and hyperbolic space $\mathbb{H}^{m+1}(r)$, as we will made clear in chapter \ref{chap::RMFinRiemGeom}.}.
\end{prop}
\textit{Proof.} In $\mathbb{E}^{m+1}$ we have the relation $\langle\alpha-p,\alpha-p\rangle=R^2\Leftrightarrow\langle\alpha-p,\mathbf{t}\rangle=0$. So, normal and spherical curves are equivalent concepts. 
\qed

Based on the remarks above, we have for a spherical curve $\alpha$ that
\begin{equation}
\alpha(s)-p = c_1(s)\mathbf{e}_1(s)+\cdots+c_m(s)\mathbf{e}_m(s)\,.
\end{equation}
Then,
\begin{eqnarray}
\mathbf{t} & = & \sum_{i=1}^mc_i'\mathbf{e}_i+\sum_{i=1}^mc_i(-\tau_{i-1}\mathbf{e}_{i-1}+\tau_{i}\mathbf{e}_{i+1})\nonumber\\
& = & -\kappa c_1\mathbf{t}+(c_1'-\tau_1 c_2)\mathbf{e}_1+\sum_{i=2}^{m}(c_i'-\tau_ic_{i+1}+\tau_{i-1}c_{i-1})\mathbf{e}_i\,.
\end{eqnarray}
Comparison of coefficients leads to
\begin{equation}
\left\{
\begin{array}{ccc}
c_1 +\kappa^{-1} & = & 0\\
c_1'-\tau_1 c_2  & = & 0\\
c_i'-\tau_i c_{i+1}+\tau_{i-1}c_{i-1} &=& 0\\
c'_m+\tau_{m-1}c_{m-1} &= & 0\\
\end{array}
\right.,\,2\leq i\leq m-1\,.
\end{equation}

The general expression to characterize spherical curves in $\mathbb{E}^{m+1}$ is quite cumbersome and we will not attempt to write it here. For spherical curves in $\mathbb{E}^4$ and $\mathbb{E}^5$ with non-zero curvature and torsions, we have
\begin{equation}
\left\{
\begin{array}{c}
\displaystyle\frac{\mathrm{d}}{\mathrm{d}s}\left\{\frac{1}{\tau_2}\frac{\mathrm{d}}{\mathrm{d}s}\left[\frac{1}{\tau_1}\frac{\mathrm{d}}{\mathrm{d}s}\left(\frac{1}{\kappa}\right)\right]+\frac{\tau_1}{\kappa}\right\}+\frac{\tau_2}{\tau_1}\frac{\mathrm{d}}{\mathrm{d}s}\left(\frac{1}{\kappa}\right)=0\\[10pt]
\displaystyle\frac{\mathrm{d}}{\mathrm{d}s}\left\{\frac{1}{\tau_3}\frac{\mathrm{d}}{\mathrm{d}s}\left[\frac{1}{\tau_2}\frac{\mathrm{d}}{\mathrm{d}s}\left(\frac{1}{\tau_1}\frac{\mathrm{d}}{\mathrm{d}s}\frac{1}{\kappa}\right)\right]+\frac{\tau_2}{\tau_1\tau_3}\frac{\mathrm{d}}{\mathrm{d}s}\frac{1}{\kappa}+\frac{1}{\tau_3}\frac{\mathrm{d}}{\mathrm{d}s}\frac{\tau_1}{\kappa\tau_2}\right\}+\frac{\tau_3}{\tau_2}\frac{\mathrm{d}}{\mathrm{d}s}\left[\frac{1}{\tau_1}\frac{\mathrm{d}}{\mathrm{d}s}\frac{1}{\kappa}\right]+\frac{\tau_1\tau_3}{\kappa\tau_2}=0
\end{array}
\right.\,.
\end{equation}
Needless to say, the approach via RM frames is much simpler and, in addition, it only demands a $C^2$ condition and no restrictions on the torsions and curvature.

\section{Characterization of plane curves}

In the previous section we used the concept of normal curves in order to study spherical curves in $\mathbb{E}^{m+1}$. Now we also address the problem of characterizing those curves in $\mathbb{E}^{3}$ whose position vector lies, up to a translation, on a (moving) plane spanned by the unit tangent and a rotation minimizing vector and prove that they are precisely the plane curves. This problem has to do with the more general quest of studying curves that lie on a given (moving) plane generated by two chosen vectors of a moving trihedron, e.g., one would define \emph{osculating}, \emph{normal} or \emph{rectifying curves} as those curves whose position vector, up to a translation, lies on their osculating, normal or rectifying planes, respectively \cite{ChenMonthly2003,ChenAJMS2017}. It is known that (i) osculating curves are precisely the plane curves\footnote{In fact, every curve is locally contained, up to second order, in its osculating plane. Thus, a plane curve must satisfy this globally.}, (ii) normal curves are precisely the spherical curves, and (iii) rectifying curves are precisely geodesics on a cone \cite{ChenManuscript,ChenAJMS2017} (\footnote{This is also valid for curves in the 3d sphere and hyperbolic space \cite{LucasJMAA2015,LucasMJM2016}.}): these investigations can be generalized to higher dimensions
\cite{CambieTJM2016,ChenAJMS2017}, 
to moving frames adapted to surfaces \cite{CamciASUOC2011} or to other ambient spaces \cite{BozkurtLSJ2013,IlarslanNSJM2007,IlarslanNSJM2003,LucasJMAA2015,LucasMJM2016} as well.

Plane curves are characterized by a vanishing torsion. This means that the curve must be contained on the plane orthogonal to the binormal vector. Next we show that this remains valid if we change the binormal vector by the second vector field of an RM frame.

\begin{theorem}
Up to a translation, the position vector of a $C^2$ regular curve $\alpha:I\to \mathbb{E}^3$ lies on a plane spanned by the unit tangent and a rotation minimizing vector field if and only if $\alpha$ is a plane curve.
\end{theorem}
\textit{Proof. } Since a plane curve $\alpha$ lies on its osculating plane, then it lies on an RM moving plane: the principal normal vector of a plane curve is an RM vector. Conversely, let $\alpha$ lies on an RM moving plane $\mbox{span}\{\mathbf{t},\mathbf{n}_1\}$, i.e.,
\begin{equation}
\alpha(s)-p = A(s)\mathbf{t}(s)+B(s)\mathbf{n}_1(s),
\end{equation}
where $p$ is constant. Taking the derivative gives
\begin{equation}
\mathbf{t}  =  (A'-\kappa_1B)\mathbf{t}+(B'+\kappa_1A)\mathbf{n}_1+\kappa_2A\,\mathbf{n}_2,
\end{equation}
and then
\begin{equation}
 \left\{
    \begin{array}{ccc}
      A'-\kappa_1B & = & 1\\
      B'+\kappa_1A & = & 0\\
      \kappa_2A    & = & 0\\
      \end{array}
  \right.\,.\label{eq::coefRMplaneCurve}
\end{equation}
If $\kappa_2(s)=0$ for all $s$, then $\kappa=0$ or $\tau=0$, since $\kappa_2=\kappa\sin\theta$. In any case, the curve is planar. On the other hand, if $A(s)=0$ for all $s$, it follows from the second equation in (\ref{eq::coefRMplaneCurve}) that $B$ is a constant. In this case, $\alpha(s)-p = B\,\mathbf{n}_1(s)$ and the curve must be spherical: it lies on a sphere of radius $\vert B\vert$ and center $p$. In addition, from the first equation in (\ref{eq::coefRMplaneCurve}), it follows that $\kappa_1=-B^{-1}$ is a constant. Now, using that for a spherical curve the normal development $(\kappa_1,\kappa_2)$ lies on a straight line not passing through the origin, theorem \ref{thr::charSphCurves}, we find that $\kappa_2$ is also a constant. Finally, since $\kappa_1$ and $\kappa_2$ are constants, we conclude that $\kappa$ is a constant and $\tau=0$, i.e., $\alpha$ is a circle (of radius $\vert B\vert$ and center $p$).
\qed
\chapter{DIFFERENTIAL GEOMETRY OF CURVES IN LORENTZ-MINKOWSKI SPACE}
\label{chap::curvesE13}

In previous chapters we were interested in the study of spherical curves in Euclidean space $(\mathbb{E}^3,\langle\cdot,\cdot\rangle)$, however we should not restrict ourselves to this context. More precisely, we can consider the more general setting of a Lorentz-Minkowski space, $(\mathbb{E}_1^3,\langle\cdot,\cdot\rangle_1)$, where one has to deal with three types of spheres: pseudo-spheres $\mathbb{S}_1^2(P;r)=F_P^{-1}(r^2)$; pseudo-hyperbolic spaces $\mathbb{H}_0^2(P;r)=F_P^{-1}(-r^2)$; and light-cones $\mathcal{C}^2(P)=F_P^{-1}(0)$, where $F_P(x)=\langle x-P,x-P\rangle_1$ and $\langle\cdot,\cdot\rangle_1$ has index 1, e.g., $\langle\mathbf{u},\mathbf{v}\rangle_1=u_1v_1+u_2v_2-u_3v_3$. In this respect, it is possible to find characterizations of some classes of spherical curves scattered among a few papers. In fact, we can find descriptions of pseudo-spherical \cite{BektasBMMS1998,
IlarslamJII-PP2003,PekmenMM1999,Petrovic-TorgasevMM2000,Petrovic-TorgasevMV2001} and pseudo-hyperbolic curves \cite{IlarslamJII-PP2003,Petrovic-TorgasevKJM2000}
via Frenet frames, and also curves on light-cones \cite{ErdoganJST2009,LiuJG2016,LiuRM2011} 
by exploiting conformal invariants and the concept of cone curvature \cite{LiuBAG2004}. It is also possible to find constructions of RM-like frames in $\mathbb{E}_1^3$ for spacelike curves \cite{BukcuCFSUA2008,BukcuSJAM2010,LowJGSP2012,OzdemirMJMS2008} with a non-lightlike normal, and timelike curves \cite{KaracanSDUJS2008,LowJGSP2012,OzdemirMJMS2008}, along with several characterizations of spherical curves through a linear equation via RM frames \cite{BukcuCFSUA2008,BukcuSJAM2010,KaracanSDUJS2008,OzdemirMJMS2008}. All the above mentioned studies in $\mathbb{E}_1^3$ have in common that much attention is paid on the possible combinations of causal characters of the tangent and normal vectors, which makes necessary the consideration of several instances of the investigation of RM frames and spherical curves. Moreover, none of them take into account the possibility of a lightlike tangent or a lightlike normal. Naturally, this reflects in the incompleteness of the available characterizations of spherical curves in $\mathbb{E}_1^3$.

In this chapter, we present a systematic approach to moving frames on curves in $\mathbb{E}_1^3$. The turning point is that one should exploit the causal character of the tangent vector and the induced causal character on the normal plane only. In this way, we are able to furnish a systematic approach to the construction of RM frames in $\mathbb{E}_1^3$ and a complete characterization of its spherical curves.

\begin{remark}
An interesting problem which we will not consider in this thesis is to consider the possibility of a curve changing its causal character. Since the property of being space- or timelike is open, i.e., if it is valid at a point it must be valid on a neighborhood of that point, the real problem is to understand what happens near lightlike points. In addition, we will neither consider curves in higher dimensional semi-Riemannian spaces nor applications in relativity, see e.g. \cite{FormigaAJP2006}.
\end{remark}

\section{Preliminaries}

Let us denote by $\mathbb{E}_{1}^3$ the vector space $\mathbb{R}^3$ equipped with a pseudo-metric $\langle\cdot,\cdot\rangle_1$ of index $1$. In fact, the concepts below, and the construction of rotation minimizing frames as well, are still valid in the context of a 3-dimensional semi-Riemannian manifold, but to help intuition, the reader may keep in mind the particular setting of $\mathbb{R}^3$ equipped with the standard Minkowski metric, i.e. $\langle x,y\rangle_1=x_1y_1+x_2y_2-x_3y_3$. Naturally, in a more general context, the derivative of a vector field along a curve should be understood as a covariant derivative. Before discussing the construction of moving frame along curves in $\mathbb{E}_1^3$, let us first introduce some terminology and geometric properties associated with $\mathbb{E}_1^3$: for more details, we refer to \cite{LopesIEJG2014,ONeill}.

One property that makes the geometry in Lorentz-Minkowski spaces $\mathbb{E}_1^3$ more difficult and richer than the geometry in $\mathbb{E}^3$ is that curves and vector subspaces may assume different causal characters:
\begin{definition}
A vector $v\in \mathbb{E}_1^3$ assumes one of the following \textit{causal characters}:
\begin{enumerate}
\item $v$ is \emph{spacelike}, if $\langle v,v\rangle_1>0$ or $v=0$;
\item $v$ is \emph{timelike}, if $\langle v,v\rangle_1<0$;
\item $v$ is \emph{lightlike}, if $\langle v,v\rangle_1=0$ and $v\not=0$.
\end{enumerate}
\end{definition}
The inner product $\langle\cdot,\cdot\rangle_1$ induces a pseudo-norm defined by $\Vert x\Vert=\sqrt{\vert\langle x,x\rangle_1\vert}$. Given a vector subspace $U\subseteq\mathbb{R}^3$, we define the orthogonal complement $U^{\perp}$ in the usual way: $U^{\perp}=\{v\in \mathbb{E}_1^3:\forall\,u\in U,\,\langle v,u\rangle_1=0\}$. Moreover, we can consider the restriction of $\langle\cdot,\cdot\rangle_1$ to $U$, $\langle\cdot,\cdot\rangle_1|_{U}$.
\begin{definition}
Let $U$ be a vector subspace, then
\begin{enumerate}
\item $U$ is \emph{spacelike} if $\langle\cdot,\cdot\rangle_1|_{U}$ is positive definite;
\item $U$ is \emph{timelike} if $\langle\cdot,\cdot\rangle_1|_{U}$ has index 1;
\item $U$ is \emph{lightlike} if $\langle\cdot,\cdot\rangle_1|_{U}$ is degenerate.
\end{enumerate}
\end{definition}
We have the following useful properties related to the causal characters of vector subspaces:
\begin{prop}
Let $U\subseteq \mathbb{E}_1^3$ be a vector subspace. Then,
\begin{enumerate}
\item $\dim U^{\perp} = 3-\dim U$ and $(U^{\perp})^{\perp}=U$;
\item $U$ is lightlike if and only if $U^{\perp}$ is lightlike;
\item $U$ is spacelike (timelike) if and only if $U^{\perp}$ is timelike (spacelike).
\item $U$ is lightlike if and only if $U$ contains a lightlike vector but not a timelike one. Moreover, $U$ admits an orthogonal basis formed  by a lightlike and a spacelike vectors.
\end{enumerate}
\end{prop}

Given two vectors $u,v\in \mathbb{E}_1^3$, the Lorentzian vector product of $u,v$ is the only vector $u\times v$ that satisfies
\begin{equation}
\forall\,w\in \mathbb{E}_1^3,\,\langle u\times v,w\rangle_1=\det(u,v,w),\label{eq::LorentzVectorProd}
\end{equation}
where the columns of $(u,v,w)$ are formed by the entries of $u,v$, and $w$. Mnemonically, we can calculate $u\times v$ for the standard Minkowski metric $\langle x,y\rangle_1=x_1y_1+x_2y_2-x_3y_3$ as
\begin{equation}
u\times v = \det\left(
\begin{array}{ccr}
\mathbf{i} & \mathbf{j} & -\mathbf{k}\\
u_1 & u_2 & u_3\\
v_1 & v_2 & v_3\\
\end{array}
\right)\,.
\end{equation}

From these definitions, we say that a curve $\alpha:I\to \mathbb{E}_1^3$ is \emph{spacelike}, \emph{timelike}, or \emph{lightlike}, if its velocity vector $\alpha'$ is spacelike, timelike, or lightlike, respectively. Analogously, we say that a surface is \emph{spacelike}, \emph{timelike}, or \emph{lightlike}, if its tangent planes are spacelike, timelike, or lightlike, respectively.

If a curve is lightlike we can not define an arc-length parameter (in $\mathbb{E}^3$ this is always possible). In this case, one must introduce the notion of a \textit{pseudo arc-length parameter}, i.e., a parameter $s$ such that $\langle\alpha''(s),\alpha''(s)\rangle_1=1$. More precisely, if $\alpha$ is a lightlike curve and $\langle\alpha'',\alpha''\rangle_1\not=0$ (otherwise $\alpha''$ and $\alpha'$ will be linearly dependent and the curve is a straight line), we define the \textit{pseudo arc-length parameter} as
\begin{equation}
s =\displaystyle\int_a^t\Vert\alpha''(u)\Vert\,\textrm{d}u\,.\label{eq:PseudoArclength}
\end{equation}
On the other hand, if $\alpha$ is not a lightlike curve, then the \textit{arc-length parameter} is defined as usual
\begin{equation}
s = \int_a^t\Vert\alpha'(u)\Vert\,\textrm{d}u\,.\label{eq:arclength}
\end{equation}
In the following we will assume every curve is regular, i.e., $\alpha'\not=0$, and parametrized by an arc-length or pseudo arc-length parameter\footnote{In the Physics literature, the arc-length parameter is sometimes referred as the \emph{proper time} \cite{ONeill}. On the other hand, the pseudo arc-length should not be confused with the affine parameter for lightlike geodesics.}.

\section{Frenet frames in Lorentz-Minkowski space}

The study of the local properties of a curve $\alpha\subset \mathbb{E}_1^3$ in a Frenet frame fashion may become quite cumbersome due to the various possibility for the causal characters of the tangent and its derivative: in essence, there is a construction for each combination of the causal characters of $\mathbf{t}$ and $\mathbf{t}'$. Indeed, let $\mathbf{t}(s)=\alpha'(s)$ be the (unit) tangent. If $\mathbf{t}'$ is not a lightlike vector, let $\mathbf{n}=\mathbf{t}'/\Vert\mathbf{t}'\Vert$ be the normal vector. We shall denote by $\epsilon=\langle\mathbf{t},\mathbf{t}\rangle_1$ and $\eta=\langle\mathbf{n},\mathbf{n}\rangle_1$ the parameters that enclose the causal character of the tangent and normal vectors. If $\mathbf{t}$ and $\mathbf{n}$ are not lightlike, then the  Frenet equation reads
\begin{equation}
\frac{\textrm{d}}{\textrm{d}s}\left(
\begin{array}{c}
\mathbf{t}\\
\mathbf{n}\\
\mathbf{b}\\
\end{array}
\right)=\left(
\begin{array}{ccc}
0 & \eta\,\kappa & 0\\
-\epsilon\,\kappa & 0 & -\epsilon\eta\,\tau\\
0 & -\eta\,\tau & 0\\ 
\end{array}
\right)\left(
\begin{array}{c}
\mathbf{t}\\
\mathbf{n}\\
\mathbf{b}\\
\end{array}
\right)=\left(
\begin{array}{ccc}
0 & \kappa & 0\\
-\kappa & 0 & \tau\\
0 & -\tau & 0\\ 
\end{array}
\right)E_{\mathbf{t},\mathbf{n},\mathbf{b}}\left(
\begin{array}{c}
\mathbf{t}\\
\mathbf{n}\\
\mathbf{b}\\
\end{array}
\right),\label{eq::FrenetEqsInE13}
\end{equation}
where $\mathbf{b}=\mathbf{t}\times\mathbf{n}$, and $\kappa = \langle\mathbf{t}',\mathbf{n}\rangle_1$ and $\tau=\langle\mathbf{n}',\mathbf{b}\rangle_1$ are the curvature function and torsion of $\alpha$, respectively\footnote{Our definition for $\kappa$ is slightly different from that of L\'opez \cite{LopesIEJG2014}, but coincides with that of K\"uhnel \cite{Kuhnel2010}. Despite the fact that our definition is formally identical to the Euclidean version, our $\kappa$ has a signal that encloses the causal character of the curve in a natural manner.}. Here $E_{\mathbf{t},\mathbf{n},\mathbf{b}}=\mbox{diag}(\epsilon,\eta,-\epsilon\eta)=[\langle\mathbf{e}_i,\mathbf{e}_j\rangle_1]_{ij}$ denotes the matrix associated with the frame $\{\mathbf{e}_0=\mathbf{t},\mathbf{e}_1=\mathbf{n},\mathbf{e}_2=\mathbf{b}\}$.

If $\mathbf{t}$ is spacelike and $\mathbf{t}'$ is lightlike, we define $\mathbf{n}=\mathbf{t}'$, while $\mathbf{b}$ is the unique lightlike vector orthonormal to $\mathbf{t}$ that satisfies $\langle\mathbf{n},\mathbf{b}\rangle_1=-1$. The Frenet equations are 
\begin{equation}
\frac{\textrm{d}}{\textrm{d}s}\left(
\begin{array}{c}
\mathbf{t}\\
\mathbf{n}\\
\mathbf{b}\\
\end{array}
\right)=\left(
\begin{array}{ccc}
0 & \,1 & 0\\
0 & \,\tau & 0\\
1 & \,0 & -\tau\\ 
\end{array}
\right)\left(
\begin{array}{c}
\mathbf{t}\\
\mathbf{n}\\
\mathbf{b}\\
\end{array}
\right)=E_{\mathbf{t},\mathbf{n},\mathbf{b}}\left(
\begin{array}{ccc}
0 & 1 & 0\\
-1 & 0 & \tau\\
0 & -\tau & 0\\ 
\end{array}
\right)\left(
\begin{array}{c}
\mathbf{t}\\
\mathbf{n}\\
\mathbf{b}\\
\end{array}
\right),\label{eq::FrenetEqsInE13tprimelightlike}
\end{equation}
where $\tau=-\langle\mathbf{n}',\mathbf{b}\rangle_1$ is the pseudo-torsion. Here $E_{\mathbf{t},\mathbf{n},\mathbf{b}}=[\langle\mathbf{e}_i,\mathbf{e}_j\rangle_1]_{ij}$ denotes the matrix associated with the null frame $\{\mathbf{e}_0=\mathbf{t},\mathbf{e}_1=\mathbf{n},\mathbf{e}_2=\mathbf{b}\}$.

Finally, if $\mathbf{t}$ is lightlike, we define $\mathbf{n}=\mathbf{t}'$ (we assume this normal vector to be spacelike, otherwise $\alpha$ is a straight line), while $\mathbf{b}$ is the unique lightlike vector that satisfies $\langle\mathbf{n},\mathbf{b}\rangle_1=0$ and $\langle\mathbf{t},\mathbf{b}\rangle_1=-1$. The Frenet equations are then
\begin{equation}
\frac{\mathrm{d}}{\mathrm{d}s}\left(
\begin{array}{c}
\mathbf{t}\\
\mathbf{n}\\
\mathbf{b}\\
\end{array}
\right)=\left(
\begin{array}{ccc}
0 & 1 & 0\\
-\tau & 0 & 1\\
0 & -\tau & 0\\ 
\end{array}
\right)\left(
\begin{array}{c}
\mathbf{t}\\
\mathbf{n}\\
\mathbf{b}\\
\end{array}
\right)=\left(
\begin{array}{ccc}
0 & 1 & 0\\
-1 & 0 & \tau\\
0 & -\tau & 0\\ 
\end{array}
\right)E_{\mathbf{t},\mathbf{n},\mathbf{b}}\left(
\begin{array}{c}
\mathbf{t}\\
\mathbf{n}\\
\mathbf{b}\\
\end{array}
\right),\label{eq::FrenetEqsInE13tangentlightlike}
\end{equation}
where $\tau=\langle\mathbf{n}',\mathbf{b}\rangle_1$ is the pseudo-torsion. Here $E_{\mathbf{t},\mathbf{n},\mathbf{b}}=[\langle\mathbf{e}_i,\mathbf{e}_j\rangle_1]_{ij}$ denotes the matrix associated with the null frame $\{\mathbf{e}_0=\mathbf{t},\mathbf{e}_1=\mathbf{n},\mathbf{e}_2=\mathbf{b}\}$.
\begin{remark}
In $\mathbb{E}^3$ the coefficient matrix of a Frenet frame is always skew-symmetric. On the other hand, this does not happen in $\mathbb{E}_1^3$ \cite{LowJGSP2012}. However, the above expressions show that the coefficient matrix can be obtained from a skew-symmetric matrix through a right-multiplication, or a left one if $\mathbf{t}'$ is lightlike, by the matrix $E_{\mathbf{t},\mathbf{n},\mathbf{b}}=[\langle\mathbf{e}_i,\mathbf{e}_j\rangle_1]_{ij}$ associated with the respective Frenet frame $\{\mathbf{e}_0=\mathbf{t},\mathbf{e}_1=\mathbf{n},\mathbf{e}_2=\mathbf{b}\}$ in $\mathbb{E}_1^3$. This skew-symmetric matrix is precisely the coefficient matrix that we would obtain for a Frenet frame in $\mathbb{E}^3$. Let us mention that when $\mathbf{t}'$ is lightlike it does not mean that the curvature function is $\kappa=1$; a curvature is not well defined for such curves \cite{LopesIEJG2014}.
\end{remark}
\begin{remark}
In the following, when discussing RM frames in $\mathbb{E}_1^3$ along non-lightlike curves and null frames along lightlike curves, we will see that the coefficient matrix can be obtained from a skew-symmetric matrix (precisely the matrix that we would obtain for a Bishop frame in $\mathbb{E}^3$) through a right-multiplication by the matrix associated with a convenient basis.
\end{remark}

\section{Rotation minimizing frames along spacelike or lightlike curves}

A quite complete and systematic approach to the problem of the existence of RM-like frames along curves in $\mathbb{E}_1^3$ was presented by \"Ozdemir and Ergin \cite{OzdemirMJMS2008}, where they build RM-like frames on timelike and spacelike curves with a non-lightlike normal. However, as in the Frenet frame case, they also paid much attention to the causal character of $\mathbf{t}'$. Here, we show that one must exploit the structure of the normal plane inherited from the causal character of $\mathbf{t}$ in order to build a unified treatment of the problem. More precisely, instead of considering the problem for each combination of the causal character of $\mathbf{t}$ and $\mathbf{t}'$, one must pay attention to the symmetry associated with the problem, which is reflected in an ambiguity in the definition of an RM frame. The study of moving frames along curves in $\mathbb{E}_1^3$ is then divided in three cases only: (i) timelike curves; (ii) spacelike curves; and (iii) lightlike curves. As a direct consequence, the characterization of spherical curves can be split along three theorems only.

\begin{definition}
A vector field $\mathbf{e}(s)$ along a regular curve $\alpha:I\to \mathbb{E}_1^3$ is a \emph{rotation minimizing} field if the derivative of its normal component is a multiple of the unit tangent vector $\mathbf{t}=\alpha'$ and its tangent component is a constant multiple of $\mathbf{t}$.
\end{definition}

Let $\alpha:I\to \mathbb{E}_1^3$ be a timelike curve. Since $\mathbf{t}$ is a timelike vector, the normal plane $N_{\alpha(s)}=\mbox{span}\{\mathbf{t}(s)\}^{\perp}$ is spacelike. To prove the existence of rotation minimizing moving frames, let $\mathbf{x}_1$ and $\mathbf{x}_2=\mathbf{t}\times\mathbf{x}_1$ be an orthonormal basis of $N_{\alpha}$. The frame $\{\mathbf{t},\mathbf{x}_1,\mathbf{x}_2\}$ satisfies the following equation
\begin{equation}
\frac{\textrm{d}}{\textrm{d}s}\left(
\begin{array}{c}
\mathbf{t}\\
\mathbf{x}_1\\
\mathbf{x}_2\\
\end{array}
\right)=\left(
\begin{array}{ccc}
0 & p_{01} & p_{02}\\
p_{01} & 0 & p_{12}\\
p_{02} & -p_{12} & 0\\ 
\end{array}
\right)\left(
\begin{array}{c}
\mathbf{t}\\
\mathbf{x}_1\\
\mathbf{x}_2\\
\end{array}
\right)=\left(
\begin{array}{ccc}
0 & p_{01} & p_{02}\\
-p_{01} & 0 & p_{12}\\
-p_{02} & -p_{12} & 0\\ 
\end{array}
\right)E_{\mathbf{t},\mathbf{x}_1,\mathbf{x}_2}\left(
\begin{array}{c}
\mathbf{t}\\
\mathbf{x}_1\\
\mathbf{x}_2\\
\end{array}
\right),
\end{equation}
for some functions $p_{ij}$, where $E_{\mathbf{t},\mathbf{x}_1,\mathbf{x}_2}=[\langle\mathbf{e}_i,\mathbf{e}_j\rangle]_{ij}$ denotes the matrix associated with the time-oriented frame $\{\mathbf{e}_0=\mathbf{t},\mathbf{e}_k=\mathbf{x}_k\}$. Let $\theta$ be a smooth function such that $\mathbf{x}=L\cos\theta\,\mathbf{x}_1+L\sin\theta\,\mathbf{x}_2$, where $L$ is a constant. Then,
\begin{equation}
\mathbf{x}'=L(p_{01}\cos\theta+p_{02}\sin\theta)\mathbf{t}+L(\theta'+p_{12})(-\sin\theta\mathbf{x}_1+\cos\theta\mathbf{x}_2).
\end{equation}
Thus, it follows that $\mathbf{x}$ is rotation minimizing if and only if $\theta'+p_{12}=0$. By the existence of a solution $\theta(s)$ for any initial condition, this shows that rotation minimizing vector fields do exist along timelike curves. Observe that RM frames are not unique. Indeed, any rotation of the normal vectors still gives two RM vector fields, i.e., there is an ambiguity associated with the group $SO(2)$.

On the other hand, if $\alpha:I\to \mathbb{E}_1^3$ is a spacelike curve, $\mathbf{t}$ is a spacelike vector and then the normal plane $N_{\alpha(s)}=\mbox{span}\{\mathbf{t}(s)\}^{\perp}$ is timelike. In a Frenet frame fashion, the study is divided into three cases, depending on the causal character of $\mathbf{t}'\in N_{\alpha}$, i.e., if $\mathbf{t}'$ is a space-, time-, or lightlike vector. But, if we only take into account the structure of $N_{\alpha}$, this is no longer necessary.

To prove the existence of rotation minimizing moving frames along spacelike curves, let $\mathbf{y}_1\in N_{\alpha}$ be a timelike vector and let $\mathbf{y}_2=\mathbf{t}\times\mathbf{y}_1$ be spacelike. Then, the frame $\{\mathbf{t},\mathbf{y}_1,\mathbf{y}_2\}$ is an orthonormal time-oriented basis of $\mathbb{E}_1^3$ along $\alpha$. The frame $\{\mathbf{t},\mathbf{y}_1,\mathbf{y}_2\}$ satisfies the following equation of motion
\begin{equation}
\frac{\textrm{d}}{\textrm{d}s}\left(
\begin{array}{c}
\mathbf{t}\\
\mathbf{y}_1\\
\mathbf{y}_2\\
\end{array}
\right)=\left(
\begin{array}{ccc}
0 & -p_{01} & p_{02}\\
-p_{01} & 0 & p_{12}\\
-p_{02} & p_{12} & 0\\ 
\end{array}
\right)\left(
\begin{array}{c}
\mathbf{t}\\
\mathbf{y}_1\\
\mathbf{y}_2\\
\end{array}
\right)=\left(
\begin{array}{ccc}
0 & p_{01} & p_{02}\\
-p_{01} & 0 & p_{12}\\
-p_{02} & -p_{12} & 0\\ 
\end{array}
\right)E_{\mathbf{t},\mathbf{y}_1,\mathbf{y}_2}\left(
\begin{array}{c}
\mathbf{t}\\
\mathbf{y}_1\\
\mathbf{y}_2\\
\end{array}
\right),
\end{equation}
for some functions $p_{ij}$, where $E_{\mathbf{t},\mathbf{y}_1,\mathbf{y}_2}=[(\mathbf{e}_i,\mathbf{e}_j)]_{ij}$ denotes the matrix associated with the time-oriented frame $\{\mathbf{e}_0=\mathbf{t},\mathbf{e}_k=\mathbf{y}_k\}$. Let $\theta$ be a smooth function such that $\mathbf{y}=L\cosh\theta\,\mathbf{y}_1+L\sinh\theta\,\mathbf{y}_2$, where it is used hyperbolic trigonometric functions because the normal plane is timelike. Then, we have
\begin{equation}
\mathbf{y}'=L(-p_{01}\cosh\theta-p_{02}\sinh\theta)\mathbf{t}+L(\theta'+p_{12})(\sinh\theta\mathbf{y}_1+\cosh\theta\mathbf{y}_2).
\end{equation}
Thus, it follows that $\mathbf{y}$ is rotation minimizing if and only if $\theta'+p_{12}=0$. By the existence of a solution $\theta(s)$ for any initial condition, this shows that rotation minimizing vector fields do exist along spacelike curves. As in the previous case, observe that RM frames are not unique. Indeed, any (hyperbolic) rotation of the normal vectors still gives two rotation minimizing vector fields, i.e., there is an ambiguity associated with the group $SO_1(2)$, which is a
component of the symmetry group of a Lorentzian plane $\mathbb{E}^2_1$ \cite{LopesIEJG2014,ONeill}. 

When $\mathbf{n}$ has a distinct causal character from that of $\mathbf{n}_1$, then we can not obtain $\mathbf{n},\mathbf{b}$ from a $SO_1(2)$-rotation of $\mathbf{n}_1,\mathbf{n}_2$, i.e., there exists no $ M\in SO_1(2)$ such that $M(\mathbf{n})=\mathbf{n}_1$ and $M(\mathbf{b})=\mathbf{n}_2$. In this case, we must first exchange $\mathbf{n}_1$ and $\mathbf{n}_2$ and then rotate them \cite{OzdemirMJMS2008}. However, we can still read the information about the causal character of $\mathbf{n}$, including the lightlike case, from the ``circles'' of the normal plane, i.e., the orbits of $O_1(2)$, see figure \ref{fig::diagramNormalDevelpCurves} and Proposition \ref{prop::geomNormalDevelopm} below.

Now we put together the above mentioned existence results of rotation minimizing frames on non-lightlike curves. Let $\{\mathbf{n}_1,\mathbf{n}_2\}$ be a basis for $N_{\alpha}$ formed by rotation minimizing vectors such that
\begin{equation}
\mathbf{n}'_i(s) = -\epsilon \kappa_i\,\mathbf{t}(s),
\end{equation}
where $\epsilon = \langle\mathbf{t},\mathbf{t}\rangle_1=\pm1$ and we have defined the RM curvatures
\begin{equation}
\kappa_i = \langle\mathbf{t}',\mathbf{n}_i\rangle_1,\,i=1,2\,.
\end{equation}
Then, defining $\epsilon_1=\langle\mathbf{n}_1,\mathbf{n}_1\rangle_1=\pm1$, we can write the following equation of motion
\begin{equation}
\frac{\textrm{d}}{\textrm{d}s}\left(
\begin{array}{c}
\mathbf{t}\\
\mathbf{n}_1\\
\mathbf{n}_2\\
\end{array}
\right)=\left(
\begin{array}{ccc}
0 & \epsilon_1\kappa_{1} & \kappa_{2}\\
-\epsilon\kappa_{1} & 0 & 0\\
-\epsilon\kappa_{2} & 0 & 0\\ 
\end{array}
\right)\left(
\begin{array}{c}
\mathbf{t}\\
\mathbf{n}_1\\
\mathbf{n}_2\\
\end{array}
\right)=\left(
\begin{array}{ccc}
0 & \kappa_{1} & \kappa_{2}\\
-\kappa_{1} & 0 & 0\\
-\kappa_{2} & 0 & 0\\ 
\end{array}
\right)E_{\mathbf{t},\mathbf{n}_1,\mathbf{n}_2}\left(
\begin{array}{c}
\mathbf{t}\\
\mathbf{n}_1\\
\mathbf{n}_2\\
\end{array}
\right),\label{eq::GenBishopEqsNonlightCurves}
\end{equation}
 where $E_{\mathbf{t},\mathbf{n}_1,\mathbf{n}_2}=[\langle\mathbf{e}_i,\mathbf{e}_j\rangle_1]_{ij}$ denotes the matrix associated with the time-oriented frame $\{\mathbf{e}_0=\mathbf{t},\mathbf{e}_k=\mathbf{n}_k\}$. The numbers $\epsilon$ and $\epsilon_1$ determine the causal character of $\mathbf{t}$ and $\mathbf{n}_1$, respectively, and since $\mathbf{n}_2=\mathbf{t}\times\mathbf{n}_1$, we have $\epsilon_2=\langle\mathbf{n}_2,\mathbf{n}_2\rangle_1=-\epsilon\epsilon_1=+1$. So, in this case $E_{\mathbf{t},\mathbf{n}_1,\mathbf{n}_2}=\mbox{diag}(\epsilon,\epsilon_1,-\epsilon\epsilon_1)$.

\subsection{Geometry of the normal development for spacelike and timelike curves}

\begin{figure*}[tbp]
\centering
    {\includegraphics[width=0.33\linewidth]{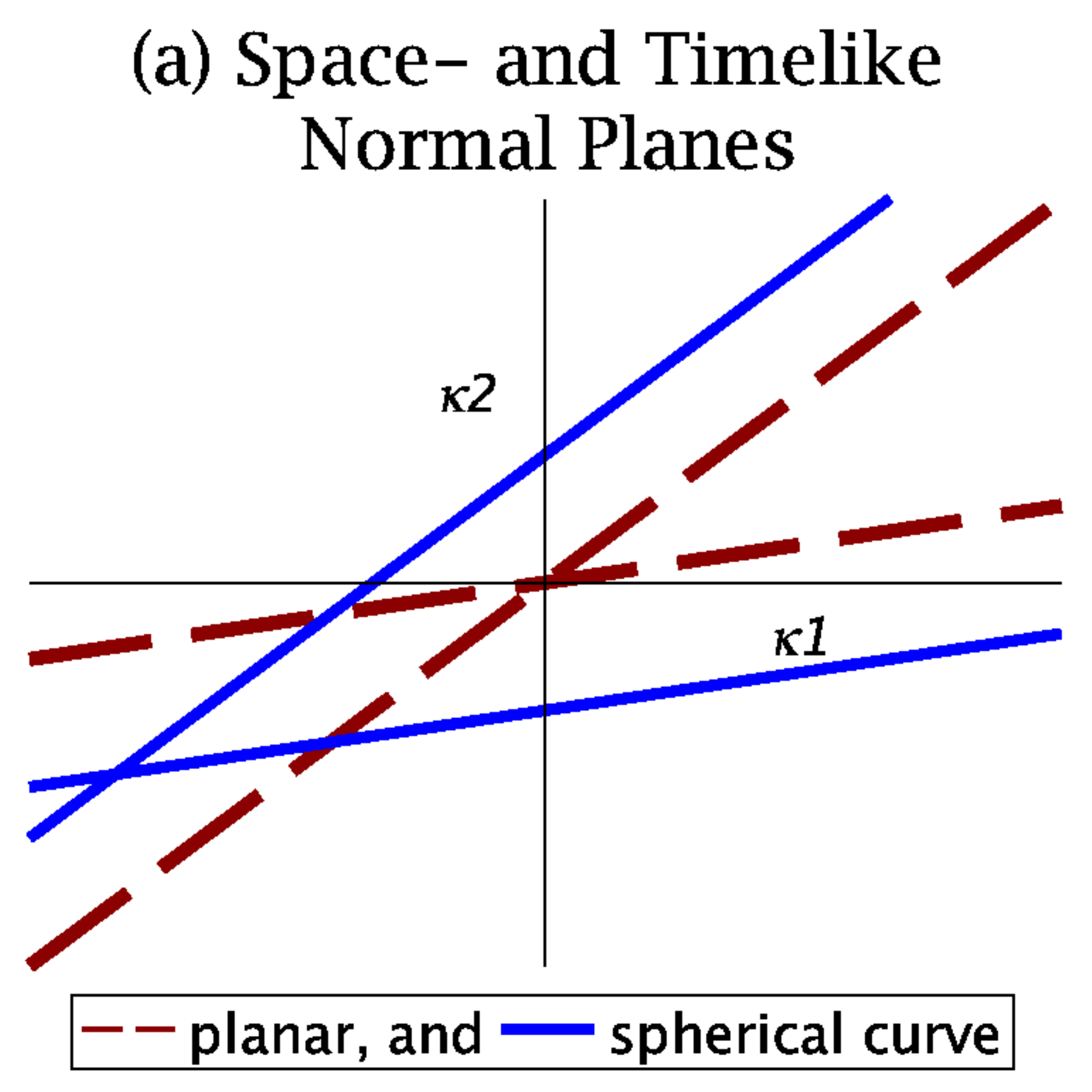}}
    {\includegraphics[width=0.32\linewidth]{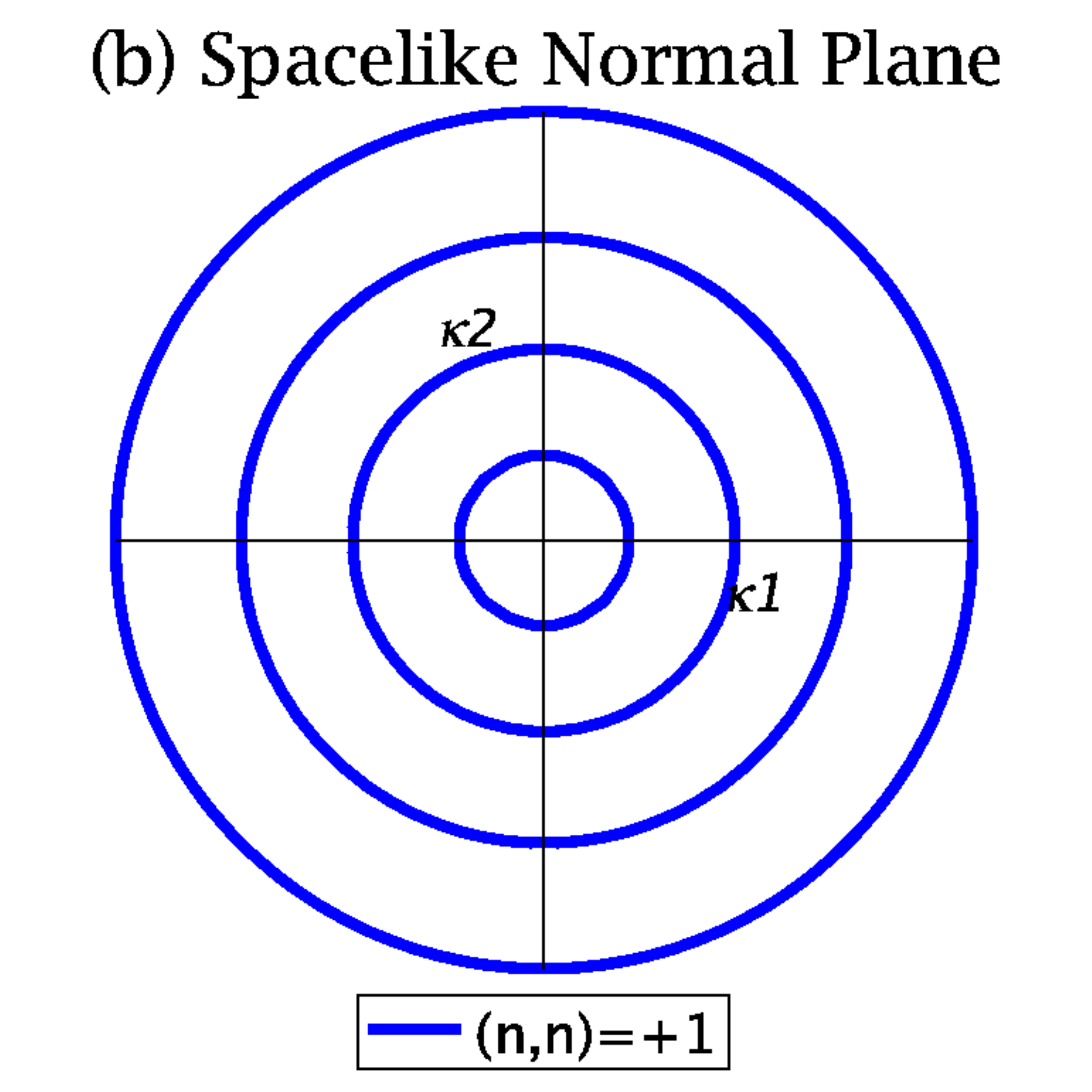}}
    {\includegraphics[width=0.32\linewidth]{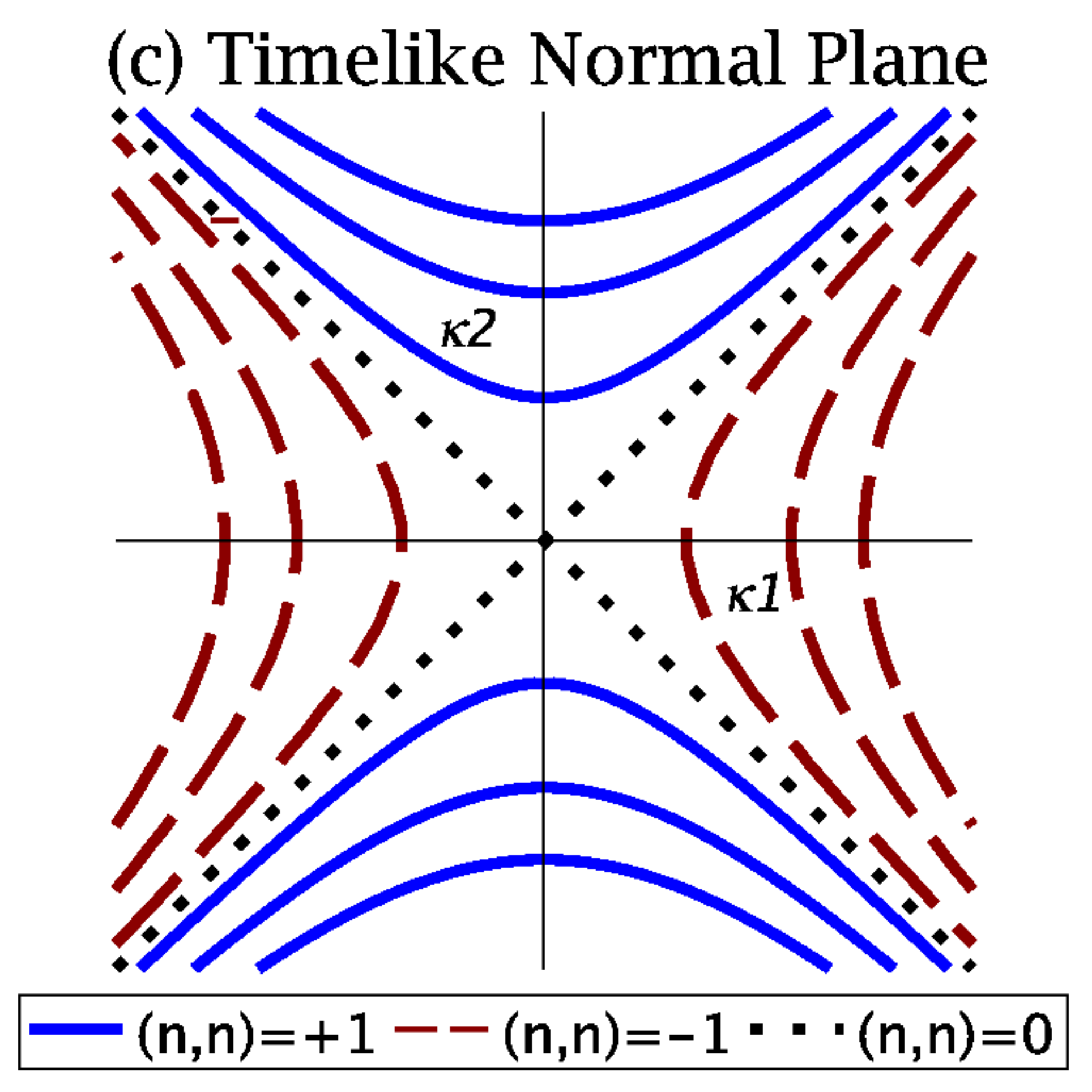}}
          \caption{The geometry of the normal development $(\kappa_1,\kappa_2)$: \textbf{(a)} On a space- or timelike normal plane, lines through the origin (dashed red line) represent plane curves (Proposition \ref{prop::CharacPlaneCurves}), and lines not passing through the origin (solid blue line) represent spherical curves (section \ref{sec::SphCurvE13}); \textbf{(b)} On a spacelike normal plane, circles represent $\kappa$-constant curves; and \textbf{(c)} On a timelike normal plane, hyperbolas represent $\kappa$-constant curves with spacelike normal vector (solid blue line) or timelike normal vector (dashed red line), and the degenerate hyperbola $\kappa_1=\pm \kappa_2$ represents curves with a lightlike normal vector (dotted black line).}
          \label{fig::diagramNormalDevelpCurves}
\end{figure*}

The \emph{normal development} of $\alpha(s)$ is the plane curve $(\kappa_1(s),\kappa_2(s))$. The normal development plane is a plane where we put the values of $\kappa_1$ and $\kappa_2$ in the first and second axis, respectively. We may equip the normal development plane with a structure of $\mathbb{E}^2$ (or $\mathtt{E}_1^2$) if the curve is timelike (or spacelike, respectively).

After proving the existence of RM moving frames for non-lightlike curves the natural question is how to relate the curvatures $\kappa_1,\kappa_2$ to the geometry of the curve which defines them. From the Frenet equations we have
\begin{equation}
\eta\mathbf{n}=\frac{\mathbf{t}'}{\kappa}=\frac{\epsilon_1\kappa_1\mathbf{n}_1+\kappa_2\mathbf{n}_2}{\kappa}\Rightarrow \eta = \epsilon_1\frac{\kappa_1^2}{\kappa^2}+\frac{\kappa_2^2}{\kappa^2}\,,
\end{equation}
where $\eta=\langle\mathbf{n},\mathbf{n}\rangle_1\in\{-1,0,+1\}$. Then we have the following relations (see figure \ref{fig::diagramNormalDevelpCurves}):
\begin{prop}
For a fixed value of the parameter $s$, the point $(\kappa_1(s),\kappa_2(s))$ lies on a conic. More precisely,
\begin{enumerate}
\item If $\mathbf{t}(s)$ is timelike (so $\mathbf{n}(s)$ must be spacelike), then $(\kappa_1(s),\kappa_2(s))$ lies on a circle of radius $\kappa(s)$: $\kappa^2=X^2+Y^2$;
\item If $\mathbf{t}(s)$ is spacelike and $\mathbf{n}(s)$ is timelike (spacelike), then $(\kappa_1(s),\kappa_2(s))$ lies on a hyperbola with foci on the $x$ axis ($y$ axis): $\kappa^2=\pm X^2\mp Y^2$;
\item If $\mathbf{t}(s)$ is spacelike and $\mathbf{n}(s)$ is lightlike, then $(\kappa_1(s),\kappa_2(s))$ lies on the line $X=\pm Y$, which form the asymptotes lines of the hyperbolas from item (2).
\end{enumerate}
\label{prop::geomNormalDevelopm}
\end{prop}

\begin{remark}
Observe that $\kappa$-constant curves correspond precisely to the orbits of the symmetry group $O(2)$ or $O_1(2)$ of an Euclidean or a Lorentzian plane, respectively, on the normal plane.
\end{remark}

Going further with the geometry of the normal development, we investigate lines passing through the origin. This characterizes plane curves (we will see in the following that straight lines not passing through the origin correspond to spherical curves, as happens in $\mathbb{E}^3$).

\begin{prop}
Let $\alpha:I\to \mathbb{E}_1^3$ be a $C^2$ regular curve which is not spherical. Then, the curve is planar if and only if its normal development $(\kappa_1(s),\kappa_2(s))$ lies on a straight line passing through the origin.
\label{prop::CharacPlaneCurves}
\end{prop}
\textit{Proof.} Suppose that $a\kappa_1+b\kappa_2=0$, with $a,\,b$ constants. Defining $\mathbf{x}(s)=a\mathbf{n}_1(s)+b\mathbf{n}_2(s)\in N_{\alpha(s)}=\mbox{span}\{\mathbf{t}(s)\}^{\perp}$, it follows that $\mathbf{x}$ is constant, $\mathbf{x}'=-a\epsilon\kappa_1\mathbf{t}-b\epsilon\kappa_2\mathbf{t}=0$, and also that
\begin{equation}
\langle\alpha,\mathbf{x}\rangle_1' = \langle\alpha,-(a\epsilon\kappa_1+b\epsilon\kappa_2)\mathbf{t}\rangle_1=-\epsilon(a\kappa_1+b\kappa_2)\,\langle\alpha,\mathbf{t}\rangle_1=0\,.
\end{equation}
Thus, $\langle\alpha(s),\mathbf{x}\rangle_1$ is constant and then $\langle\alpha(s)-\alpha(s_0),\mathbf{x}\rangle_1=0$. So, $\alpha$ is a plane curve.

Conversely, suppose $\langle\alpha(s)-\alpha(s_0),\mathbf{x}\rangle_1=0$. Since the tangent $\mathbf{t}$ also belongs to this plane, we have $\mathbf{x}=a\mathbf{n}_1+b\mathbf{n}_2$. Moreover, $a=\langle\mathbf{x},\mathbf{n}_1\rangle_1$ and then $a'=\langle\mathbf{x},-\epsilon\mathbf{t}\rangle_1=0$. Analogously, $b'=0$. In short, $a,b$ are constants. In addition
\begin{equation}
0=\langle\alpha-\alpha_0,\mathbf{x}\rangle_1' = \langle\alpha-\alpha_0,-(a\epsilon\kappa_1+b\epsilon\kappa_2)\mathbf{t}\rangle_1=-\epsilon(a\kappa_1+b\kappa_2)\,\langle\alpha-\alpha_0,\mathbf{t}\rangle_1\,.
\end{equation}
Thus, $a\kappa_1+b\kappa_2=0$ (if it were $\langle\alpha-\alpha_0,\mathbf{t}\rangle_1=0$, the curve would be spherical, since $\langle\alpha-\alpha_0,\alpha-\alpha_0\rangle_1'=2\langle\alpha-\alpha_0,\mathbf{t}\rangle_1=0\Leftrightarrow\langle\alpha-\alpha_0,\alpha-\alpha_0\rangle_1\,=\mbox{ constant}$).
\qed

\begin{remark}
If the pseudo-torsion of a spacelike curve with a lightlike normal vanishes, then the curve is planar, but the converse is not true: indeed, L\'opez \cite{LopesIEJG2014} gives an example of a curve which is planar and has a non-zero pseudo-torsion. It follows from the above propositions that all spacelike curves with a lightlike normal are planar, no matter the value of the pseudo-torsion.
\end{remark}


\section{Moving frames along lightlike curves}
\label{Sec::MovFrmNullCrv}

It is not possible to define RM frames along lightlike curves, we can not even define an orthonormal frame. In this case we must work with the concept of a null frame (see e.g. \cite{InoguchiIEJG2008} for a survey on the geometry of lightlike curves and null frames along them). As in the previous case, we will introduce along the curve $\alpha$ a (null) frame by exploiting the structure of the normal plane only.

Let $\alpha:I\to \mathbb{E}_1^3$ be a lightlike curve. In this case, since $\alpha'$ is a lightlike vector, the normal plane $N_{\alpha(s)}=\mbox{span}\{\alpha'(s)\}^{\perp}$ is lightlike and $\alpha'\in N_{\alpha}$. Then, we have $N_{\alpha(s)}=\mbox{span}\{\alpha'(s),\mathbf{z}_1(s)\}$, where $\mathbf{z}_1$ is a unit spacelike vector. Denote by $\mathbf{t}=\alpha'$ the tangent vector. If $\mathbf{t}'$ is spacelike, then we can assume $\alpha$ to parametrized by a pseudo arc-length. Let $\mathbf{z}_2$ be the lightlike vector orthogonal to $\mathbf{z}_1$ and satisfying $\langle\mathbf{t},\mathbf{z}_2\rangle_1=-1$. In this case, the equation of motion is
\begin{equation}
\frac{\textrm{d}}{\textrm{d}s}\left(
\begin{array}{c}
\mathbf{t}\\
\mathbf{z}_1\\
\mathbf{z}_2\\
\end{array}
\right)=\left(
\begin{array}{ccc}
\kappa_3 & \kappa_1 & 0\\
-\kappa_2 & 0 & \kappa_1\\
0 & -\kappa_{2} & -\kappa_3\\ 
\end{array}
\right)\left(
\begin{array}{c}
\mathbf{t}\\
\mathbf{z}_1\\
\mathbf{z}_2\\
\end{array}
\right)=\left(
\begin{array}{ccc}
0 & \kappa_1 & -\kappa_3\\
-\kappa_1 & 0 & \kappa_2\\
\kappa_3 & -\kappa_{2} & 0\\ 
\end{array}
\right)E_{\mathbf{t},\mathbf{z}_1,\mathbf{z}_2}\left(
\begin{array}{c}
\mathbf{t}\\
\mathbf{z}_1\\
\mathbf{z}_2\\
\end{array}
\right),\label{eq::FrameEqsLightCurves}
\end{equation}
where $\kappa_1=\langle\mathbf{t}',\mathbf{z}_1\rangle_1$, $\kappa_2=\langle\mathbf{z}_1',\mathbf{z}_2\rangle_1$, and $\kappa_3=\langle\mathbf{z}_2',\mathbf{t}\rangle_1$. Here $E_{\mathbf{t},\mathbf{n},\mathbf{b}}=[\langle\mathbf{e}_i,\mathbf{e}_j\rangle_1]_{ij}$ denotes the matrix associated with the null frame $\{\mathbf{e}_0=\mathbf{t},\mathbf{e}_1=\mathbf{n},\mathbf{e}_2=\mathbf{b}\}$. The coefficient $\kappa_1$ plays a significant role on the theory of moving frames along lightlike curves.
\begin{remark}
If $\mathbf{t}'$ is spacelike and if we take $\mathbf{z}_1=\mathbf{t}'=\mathbf{n}$, then $\mathbf{z}_2=\mathbf{b}$ and $\kappa_1=1$, $\kappa_2=\tau$ and $\kappa_3=0$. However, the Frenet frame is not defined when $\mathbf{t}'$ is lightlike. Here, the presence of $\kappa_3$ allows for a description of lightlike curves regardless of the causal character of $\mathbf{t}'$. 
\end{remark}

\begin{prop}
A lightlike curve $\alpha:I\to \mathbb{E}_1^3$ is a straight line if and only if $\kappa_1=0$. Moreover, if $\alpha$ is not a straight line and is parametrized by the pseudo arc-length, then $\kappa_1^2=1$.
\label{prop::CharacLightlikeLine}
\end{prop}
\textit{Proof.} If $\kappa_1=0$, then $\mathbf{t}'=\kappa_3\mathbf{t}$. Integration of this equation gives $\alpha=\alpha_0+(\int \textrm{e}^{\int\kappa_3})\mathbf{t}_0$, where $\alpha_0$ and $\mathbf{t}_0$ are constants. Then, $\alpha$ is a straight line. Conversely, let $\alpha=\alpha_0+f\,\mathbf{t}_0$, with $f$ a smooth function. Taking derivatives, it is easy to verify that $\kappa_1=0$.

Now, suppose $\kappa_1\not=0$, so if $\alpha$ is parametrized by pseudo arc-length we have
\begin{equation}
1=\langle\mathbf{t}',\mathbf{t}'\rangle_1=\kappa_1^2,
\end{equation}
as expected.
\qed

\section{Characterization of spherical curves in Lorentz-Minkowski space}
\label{sec::SphCurvE13}
\begin{figure*}[tbp]
\centering
    {\includegraphics[width=\linewidth]{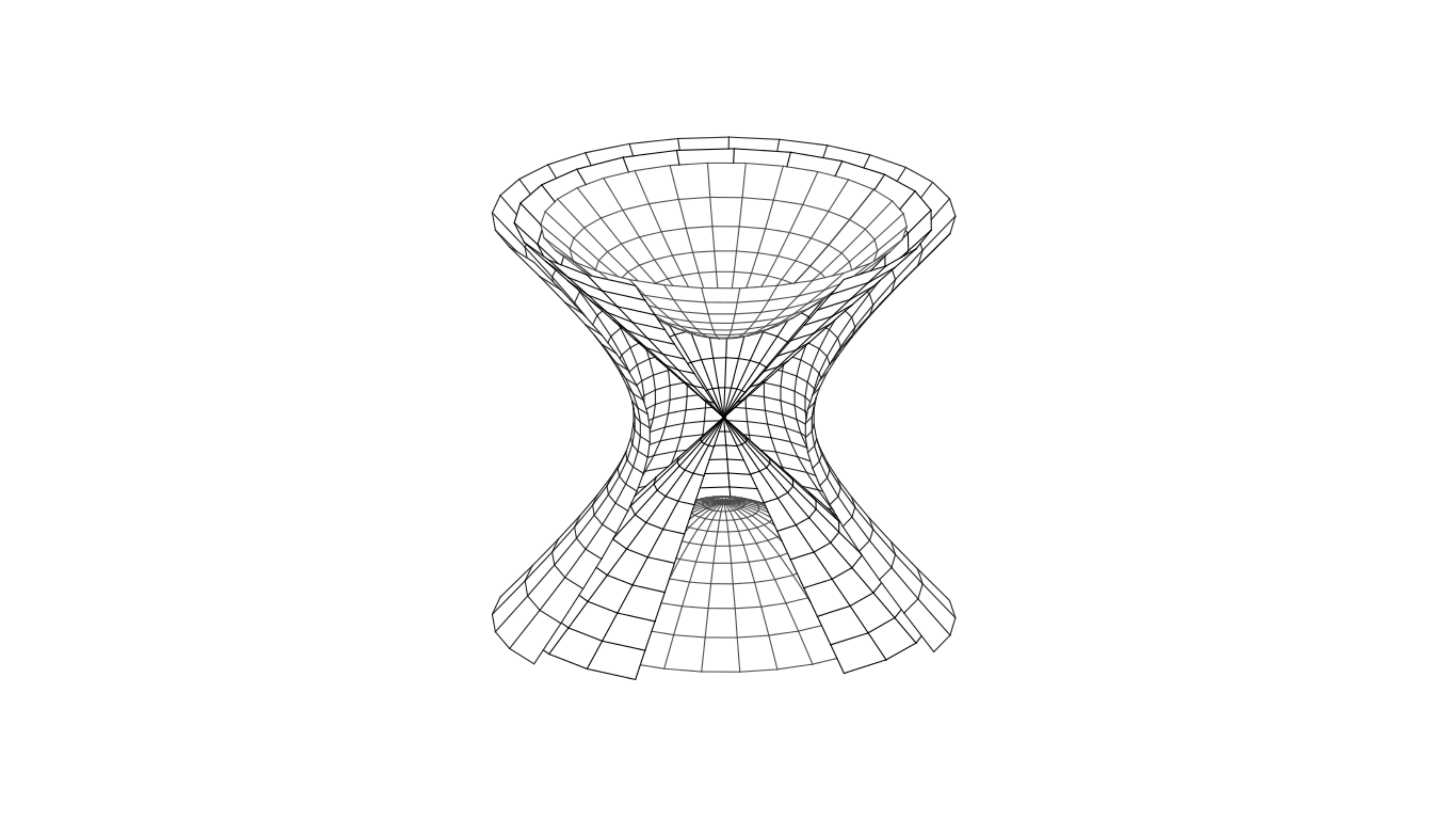}}
          \caption{The three types of spheres in $\mathbb{E}^3_1$. Pseudo-spheres $\mathbb{S}^2_1$ are represented by a one-sheeted hyperboloid, pseudo-hyperbolic spaces $\mathbb{H}_0^2$ by a two-sheeted hyperboloid, and light-cones $\mathcal{C}^2$ by a cone. From \cite{Kuhnel2010}.}
          \label{fig::pseudospheres}
\end{figure*}

In $\mathbb{E}^3$ the function $F(x)=\langle x-P,x-P\rangle$ is non-negative. A sphere of radius $r$ and center $P$ in $\mathbb{E}^3$, $\mathbb{S}^2(P;r)$, is then defined as the level sets of $F$, i.e. $\langle x-P,x-P\rangle=r^2$ (if $r=0$ the sphere degenerates to a single point). On the other hand, in $\mathbb{E}_1^3$ the function $F_1(x)=\langle x-P,x-P\rangle_1$ may assume any value on the real numbers. Then, in $\mathbb{E}_1^3$ we can still define spheres as the level sets of $F_1$, but one must consider three types of spheres, depending on the sign of $F_1$. We shall adopt the following standard notations:
\begin{equation}
\mathbb{S}_1^2(P;r) =\{x\in \mathbb{E}_1^3\,:\, \langle x-P,x-P\rangle_1=r^2\},\\
\end{equation}
\begin{equation}
\mathcal{C}^2(P) = \{x\in \mathbb{E}_1^3\,:\,\langle x-P,x-P\rangle_1=0\},
\end{equation}
and
\begin{equation}
\mathbb{H}_0^2(P;r) =\{x\in \mathbb{E}_1^3\,:\, \langle x-P,x-P\rangle_1=-r^2\},
\end{equation}
where $r\in (0,\infty)$. These spheres are known as \emph{pseudo-sphere}, \emph{light-cone}, and \emph{pseudo-hyperbolic space}\footnote{The sphere $\mathbb{H}^2_0$ is precisely the hyperboloid model for the hyperbolic geometry \cite{ReynoldsMonthly1993}. In this case we just write $\mathbb{H}^2$ and omit the term ``pseudo'' from the name.}, respectively. As surfaces in $\mathbb{E}_1^3$ pseudo-spheres and pseudo-hyperbolic spaces have constant Gaussian curvature $1/r^2$ and $-1/r^2$ \cite{LopesIEJG2014}, respectively\footnote{If we see them as surfaces in $\mathbb{E}^3$, their Gaussian curvatures are not constant and, additionally, for $\mathbb{S}_1^2(P;r)$ it is negative, while for $\mathbb{H}_0^2(P;r)$ it is positive.}. 

It is well known that the Minkowski metric restricted to $\mathbb{H}_0^2(P;r)$ is a positive definite metric. Then, it follows that $\mathbb{H}_0^2(P;r)$ is a spacelike surface and, consequently, there is no lightlike or timelike curves in $\mathbb{H}_0^2(P;r)$. On the other hand, light-cones are lightlike surfaces \cite{LopesIEJG2014} and, consequently, there is no lightlike curves on them. The pseudo-sphere is the only one that has the three types of curves \cite{InoguchiIEJG2008,LopesIEJG2014}:
\begin{lemma}
There exist no time- and lightlike curves in $\mathbb{H}_0^2(P;r)$ and no timelike curves in $\mathcal{C}^2(P)$.
\label{lemma::CasualCharSphCurv}
\end{lemma}

Now we generalize Bishop's characterization of spherical curves in $\mathbb{E}^3$ \cite{BishopMonthly} to the context of spheres in $\mathbb{E}_1^3$.

\begin{theorem}
A $C^2$ regular spacelike or timelike curve $\alpha:I\to \mathbb{E}_1^3$ lies on a sphere of nonzero radius, i.e., $\alpha\subseteq \mathbb{H}_0^2(P;r)$ or $\mathbb{S}_1^2(P;r)$, if and only if its normal development, i.e., the curve $(\kappa_1(s),\kappa_2(s))$, lies on a line not passing through the origin. Moreover, the distance of this line from the origin, $d$, and the radius of the sphere are reciprocals: $d=1/r$.
\label{theo::characSpaceAndLightCurves}
\end{theorem}
\begin{remark}
When a curve is spacelike the normal plane is timelike and then the distance in the normal development plane should be understood as the distance induced by the restriction of $\langle\cdot,\cdot\rangle_1$ on the normal plane. So, circles in this plane are hyperbolas.
\end{remark}
\textit{Proof of Theorem} \ref{theo::characSpaceAndLightCurves}. Denote by $\mathcal{Q}$ a sphere $\mathbb{H}_0^2(P;r)$ or $\mathbb{S}_1^2(P;r)$. If $\alpha$ lies in $\mathcal{Q}$, then taking the derivative of $\langle\alpha-P,\alpha-P\rangle_1=\pm\, r^2$ gives 
\begin{equation}
\langle\alpha-P,\mathbf{t}\rangle_1=0.\label{eq::aux1} 
\end{equation}
This implies that $\alpha-P=a_1\mathbf{n}_1+a_2\mathbf{n}_2$. Now, let us investigate the coefficients $a_i$. Since $a_i=\epsilon_i\langle\alpha-P,\mathbf{n}_i\rangle_1$, where $\epsilon_i=\langle\mathbf{n}_i,\mathbf{n}_i\rangle_1$, we have
\begin{eqnarray}
a_i' & = & \epsilon_i\langle\mathbf{t},\mathbf{n}_i\rangle_1+\epsilon_i\langle\alpha-P,\mathbf{n}_i'\rangle_1=0\,.
\end{eqnarray}
Therefore, the coefficients $a_1$ and $a_2$ are constants. Finally, taking the derivative of Eq. (\ref{eq::aux1}), we find
\begin{equation}
0=\langle\mathbf{t},\mathbf{t}\rangle_1+\langle\alpha-P,\epsilon_1\kappa_1\mathbf{n}_1+\kappa_2\mathbf{n}_2\rangle_1=\epsilon+a_1\kappa_1+a_2\kappa_2.
\end{equation}
Thus, the normal  development $(\kappa_1,\kappa_2)$ lies on a straight line $1+\epsilon a_1X+\epsilon a_2Y=0$ not passing through the origin. If $Q=\mathbb{S}_1^2(P;r)$, then $r^2=\langle\alpha-P,\alpha-P\rangle_1=\epsilon_1 a_1^2+a_2^2=1/d^2$, where $d$ is the distance of the line from the origin. On the other hand, if $Q=\mathbb{H}_0^2(P;r)$, then the curve is necessarily spacelike and $\epsilon_1=-1$, since $\mathbf{n}_1$ is timelike (as mentioned before, $\mathbb{H}_0^2(P;r)$ is a spacelike surface). So, we have $r^2=-\langle\alpha-P,\alpha-P\rangle_1= a_1^2-a_2^2=\pm1/d^2$ (the orientation of the hyperbolas will depend on the causal character of the normal vector $\mathbf{n}$ according to Proposition \ref{prop::geomNormalDevelopm}: see figure 1).

Conversely, assume that $0=1+\epsilon a_1\kappa_1+\epsilon a_2\kappa_2$ for some constants $a_1$ and $a_2$.  Define the function $P(s)=\alpha(s)-a_1\mathbf{n}_1(s)-a_2\mathbf{n}_2(s)$. Then $P'=\mathbf{t}+(a_1\epsilon\kappa_1+a_2\epsilon\kappa_2)\mathbf{t}=0$ and, therefore, $P$ is a fixed point. It follows that $\alpha$ lies on a sphere of nonzero radius and center $P$: $\langle\alpha-P,\alpha-P\rangle_1=\epsilon_1a_1^2+a_2^2$.
\qed

For spacelike curves on light-cones (as mentioned before there is no timelike curve on light-cones: Lemma \ref{lemma::CasualCharSphCurv}) we have an analogous characterization:
\begin{theorem}
A $C^2$ regular spacelike curve $\alpha:I\to \mathbb{E}_1^3$ lies on a light-cone $\mathcal{C}^2(P)$, i.e., lies on a sphere of zero radius, if and only if its normal development, i.e., the curve $(\kappa_1(s),\kappa_2(s))$, lies on a line $\{a_1X\pm a_1Y+1=0\}$ not passing through the origin.
\end{theorem}
\textit{Proof.} Let $\alpha$ be a curve in $\mathcal{C}^2(P)$ with $\langle\mathbf{t},\mathbf{t}\rangle_1=1$ and $\langle\mathbf{n}_1,\mathbf{n}_1\rangle_1=-1$, i.e., $\epsilon=1$ and $\epsilon_1=-1$. Now taking the derivative of $\langle\alpha-P,\alpha-P\rangle_1=0$ gives 
\begin{equation}
\langle\alpha-P,\mathbf{t}\rangle_1=0.\label{eq::aux2} 
\end{equation}
This implies that $\alpha-P=a_1\mathbf{n}_1+a_2\mathbf{n}_2$. Since $a_i=\epsilon_i\langle\alpha-P,\mathbf{n}_i\rangle_1$, where $\epsilon_i=\langle\mathbf{n}_i,\mathbf{n}_i\rangle_1$, we have
\begin{eqnarray}
a_i' & = & \epsilon_i\langle\mathbf{t},\mathbf{n}_i\rangle_1+\langle\alpha-P,\mathbf{n}_i'\rangle_1=0\,.
\end{eqnarray}
Therefore, the coefficients $a_1$ and $a_2$ are constants. Finally, taking the derivative of Eq. (\ref{eq::aux2}), we find
\begin{equation}
0=\langle\mathbf{t},\mathbf{t}\rangle_1+\langle\alpha-P,-\kappa_1\mathbf{n}_1+\kappa_2\mathbf{n}_2\rangle_1=1+a_1\kappa_1+a_2\kappa_2.
\end{equation}
Thus, the normal  development $(\kappa_1(s),\kappa_2(s))$ lies on a straight line $1+a_1X+a_2Y=0$ not passing through the origin. Moreover, $0=\langle\alpha-P,\alpha-P\rangle_1=- a_1^2+a_2^2$, which implies $a_2=\pm a_1$.

Conversely, assume that $0=1+ a_1\kappa_1\pm a_1\kappa_2$ for some constant $a_1$. Define the function $P(s)=\alpha(s)-a_1\mathbf{n}_1(s)\mp a_1\mathbf{n}_2(s)$, which satisfies $P'=\mathbf{t}+(a_1\kappa_1\pm a_1\kappa_2)\mathbf{t}=0$. In other words, $P$ is a fixed point and it follows that $\alpha$ lies on a light-cone $\mathcal{C}^2(P)$ of center $P$.
\qed

For lightlike curves we are not able to use an RM frame approach.  However, by using null frames (section \ref{Sec::MovFrmNullCrv}), we can still state a criterion for a lightlike curve be contained on pseudo-spheres or light-cones (trying to follow steps as in the previous cases does not work, due to the lack of good orthogonality properties). In fact, the following results are generalizations of those of Inoguchi and Lee \cite{InoguchiIEJG2008} for pseudo-spherical lightlike curves.
\begin{theorem}
If a $C^2$ regular lightlike curve $\alpha:I\to \mathbb{E}_1^3$ lies on a pseudo-sphere or a light-cone, then $\kappa_1=0$ or, equivalently, $\alpha$ is a straight line.
\label{theo::SphrLightcurves}
\end{theorem}
\textit{Proof.} Let $\mathcal{Q}$ be a sphere of non-negative radius denoted by $\mathcal{Q}=\{x\,:\,\langle x-P,x-P\rangle_1=\rho\}$ where $\rho=r^2$ ($r>0$) or $0$, i.e., $\mathcal{Q}$ is a pseudo-sphere $\mathbb{S}_1^2(P;r)$ or a light-cone $\mathcal{C}^2(P)$. If $\alpha\subseteq\mathcal{Q}$, taking the derivative of $\langle x-P,x-P\rangle_1=\rho$ gives
\begin{equation}
\langle\mathbf{t},x-P\rangle_1=0.\label{eq::auxtx-pzero}
\end{equation}
Deriving the above equation gives
\begin{equation}
\kappa_1\langle\mathbf{z}_1,x-P\rangle_1=0.
\end{equation}
If $\kappa_1$ were not zero, then we would find $\langle\mathbf{z}_1,x-P\rangle_1=0$, which by taking a derivative again gives $\langle\mathbf{z}_2,x-P\rangle_1=0$. From these two last equations, and from Eq. (\ref{eq::auxtx-pzero}), we would conclude that $x-P=0$, which is not possible. In short, the curve must satisfy $\kappa_1=0$. Finally, by Proposition \ref{prop::CharacLightlikeLine} it follows that $\alpha$ must be a straight line.
\qed

\begin{remark}
Surfaces in a semi-Riemannian manifold $M_1^3$ have an interesting property: a lightlike curve is always a pregeodesic, i.e., there exists a parametrization that makes the curve a parametrized geodesic \cite{ONeill}. Theorem \ref{theo::SphrLightcurves} shows that all lightlike geodesics of $\mathbb{S}_1^2$ and $\mathcal{C}^2$ are straight lines.
\end{remark}

The converse of the above theorem is not true. In fact, taking $\langle\cdot,\cdot\rangle_1$ as the standard Minkowski metric, the straight line $\alpha(\sigma)=(0,0,\sigma)$ does not lie on any pseudo-sphere or light-cone. However, we have the following partial converse: 
\begin{prop}
Let $\alpha_0\in \mathcal{Q}(P;\rho)=\{x:\langle x-P,x-P\rangle_1=\rho\}$ be a point on a pseudo-sphere or light-cone, i.e., $\rho=r^2$ ($r>0$) or $=0$. If $\mathbf{u}\in T_{\alpha_0}\mathcal{Q}(P;\rho)$ is a lightlike vector, then for any smooth function $f(\sigma)$ the curve $\alpha(\sigma)=\alpha_0+f(\sigma)\,\mathbf{u}$ is a lightlike straight line that lies on $\mathcal{Q}(P;\rho)$.
\end{prop}
\textit{Proof.} Using that $\mathbf{u}\in T_{\alpha_0}\mathcal{Q}(P;\rho)$ implies $\langle\alpha_0-P,\mathbf{u}\rangle_1=0$, we find
\begin{eqnarray}
\langle\alpha-P,\alpha-P\rangle_1 & = & \langle\,(\alpha_0-P)+f\,\mathbf{u},(\alpha_0-P)+f\,\mathbf{u}\rangle_1\nonumber\\
& = & \langle\alpha_0-P,\alpha_0-P\rangle_1=\rho\,.
\end{eqnarray}
So, the desired result follows.
\qed
\chapter{CHARACTERIZATION OF CURVES THAT LIE ON A LEVEL SURFACE IN EUCLIDEAN SPACE}
\label{Chap_CurvInSurf}

Here we apply the ideas presented in previous chapters in the investigation of spherical curves in order to characterize those spatial curves that belong to surfaces implicitly defined by a smooth function, $\Sigma=F^{-1}(c)$, by reinterpreting the problem in the new geometric setting of an inner product induced by the Hessian, $\mbox{Hess}\,F=\partial^2F/\partial x^i\partial x^j$. Since a Hessian may fail to be positive or non-degenerate, one is naturally led to the study of the differential geometry of curves in Lorentz-Minkowski and isotropic spaces. Here we present a necessary and sufficient criterion for a curve to lie on a level surface of a smooth function. More precisely, we present a functional relationship involving the coefficients of an RM frame with respect to the Hessian metric along a curve on $\Sigma=F^{-1}(c)$, which reduces to a linear relation when $\mbox{Hess}\,F$ is constant. In this last case, we are able to characterize spatial curves that belong to a given Euclidean quadric $\mathcal{Q}=\{x:\langle B(x-P),(x-P)\rangle=\rho\}$, $\rho\in\mathbb{R}$ constant, by using $h(\cdot,\cdot)=\langle B\cdot,\cdot\rangle$. We also furnish an interpretation for the casual character that a curve may assume when we pass from $\mathbb{E}^3$ to $\mathbb{E}_1^3$, which also allows us to understand why certain types of curves do not exist on a given quadric or on a given Lorentzian sphere, if we reinterpret the problem from $\mathbb{E}_1^3$ in $\mathbb{E}^3$. Finally, in the end of this chapter, we devote our attention to the geometry of isotropic spaces, which naturally appear in the case of a degenerate Hessian and allows for a characterization of curves on degenerate quadrics.

To the best of our knowledge, this is the first time that this characterization problem is considered in a general context.

\section{Characterization of curves on Euclidean quadrics}

Quadrics are the simplest examples of level set surfaces and understanding how the characterization works in this particular instance will prove very useful. Indeed, it will become clear in the following that the proper geometric setting to attack the characterization problem on a surface $\Sigma=F^{-1}(c)$ is that of a metric induced by the Hessian of $F$.

Points on a quadratic surface $\mathcal{Q}\subset\mathbb{R}^3$ can be characterized by a symmetric matrix $B\in\mbox{M}_{3\times3}(\mathbb{R})$ as
\begin{equation}
x\in\mathcal{Q}\Leftrightarrow\Big\langle \,B(x-P),x-P\,\Big\rangle=r^2,\label{eq:QuadraticSurface}
\end{equation}
where $P$ is a fixed point (the center of $\mathcal{Q}$), $r>0$ is a constant, and $\langle\cdot,\cdot\rangle$ is the canonical inner product on $\mathbb{R}^3$. Naturally, if the symmetric matrix $B$ has a non-zero determinant, then this non-degenerate quadric induces a metric or a pseudo-metric on $\mathbb{R}^3$ by defining
\begin{equation}
\langle\cdot,\cdot\rangle_i=\pm\,\langle B\,\cdot,\cdot\rangle\,,
\end{equation}
where $i$ stands for the index of the matrix $B$. More precisely, if the matrix $B$ has index 0, then $\mathcal{Q}$ is an ellipsoid and it can be seen as a sphere on the 3-dimensional Riemannian manifold $M^3=(\mathbb{R}^3,\langle\cdot,\cdot\rangle_0=\langle B\,\cdot,\cdot\rangle)$. The characterization of those spatial curves that belong to an ellipsoid can be made through a direct adaption of Bishop's characterization of spherical curves in $\mathbb{E}^3$ (see chapter \ref{chap::RMframes}). Indeed, one just uses the metric $\langle B\,\cdot,\cdot\rangle$ instead of $\langle\cdot,\cdot\rangle$ and then follows the steps on the construction of an RM frame in $\mathbb{E}^3$. On the other hand, if the matrix $B$ has index 1, then $\mathcal{Q}$ is a one-sheeted hyperboloid and can be seen as a pseudo-sphere on a Lorentz-Minkowski space $\mathbb{E}^3_1=(\mathbb{R}^3,\langle\cdot,\cdot\rangle_1=\langle B\,\cdot,\cdot\rangle)$. If $B$ has index 2, $\mathcal{Q}$ is then a two-sheeted hyperboloid and can be seen as a pseudo-hyperbolic plane on a Lorentz-Minkowski space $\mathbb{E}^3_1=(\mathbb{R}^3,\langle\cdot,\cdot\rangle_2=-\langle B\,\cdot,\cdot\rangle)$. This way, the results on the previous chapters can be applied in order to characterize those spatial curves that belong to a (one or two-sheeted) hyperboloid.

Since the characterization of curves on a quadric is made be reinterpreting the problem on a new geometric setting, a natural question then arises: \textit{How do we interpret the casual character that a spatial curve assumes when we pass from $\mathbb{E}^3$ to $\mathbb{E}_1^3$?}

This question can be answered if we take into account the following expression for the normal curvature on a level surface $\Sigma=F^{-1}(c)$ \cite{DombrowskiMN1968}
\begin{equation}
\kappa_n(p,\mathbf{v}) = \frac{\langle \mbox{Hess}_pF\,\mathbf{v},\mathbf{v}\rangle}{\Vert\nabla_p F\Vert},\label{eqNormalCurvLevelSets}
\end{equation}
where $\mathbf{v}\in T_p\Sigma$, and $\mbox{Hess}\,F$ and  $\nabla F$ are the Hessian and the gradient vector of $F$, respectively (for more details involving the expressions for the curvatures of level set surfaces see \cite{GoldmanCAGD2005}). Then, we have the following interpretation:
\begin{prop}
If $\alpha:I\to\mathbb{R}^3$ is a curve on a non-degenerate quadric $\mathcal{Q}$, then asymptotic directions (in $\mathcal{Q}\subseteq \mathbb{E}^3$) correspond to lightlike directions (in $\mathcal{Q}\subseteq \mathbb{E}_1^3$).
\label{prop::InterpretCasualChar}
\end{prop}
\textit{Proof.} Quadrics are level sets of $F(x)=\langle B\,(x-P),x-P\rangle$ and $\mbox{Hess}\,F=B$. Now, since the quadric is non-degenerate, we have that $\mathcal{Q}$ is the inverse image of a regular value of $F$. Thus, we can apply Eq. (\ref{eqNormalCurvLevelSets}).
\qed

Based on these constructions we can better interpret why pseudo-spheres $\mathbb{S}_1^2$ have both space- and timelike tangent vectors, while pseudo-hyperbolic planes $\mathbb{H}_0^2$ only have spacelike ones. Indeed, Eq. (\ref{eqNormalCurvLevelSets}) shows that the sign of the Gaussian curvature in $\mathbb{E}^3$, $K_{\mathbb{E}^3}$, has an impact on the casual character of the tangent plane:  points with $K_{\mathbb{E}^3}>0$ have spacelike tangent planes, while points with $K_{\mathbb{E}^3}<0$ have timelike tangent planes. 

Finally, observe that quadrics are level sets of  $F(x)=\langle B(x-P),x-P\rangle$, which has a constant Hessian: $\mbox{Hess}\,F=B$. This motivates us to consider this procedure for any level surface.

\section{Curves on level surfaces of a smooth function}

Let $\Sigma$ be a surface implicitly defined by a smooth function $F:U\subseteq\mathbb{R}^3\to\mathbb{R}$. Then, the Hessian of $F$ induces on $\mathbb{R}^3$ a (pseudo-) metric 
\begin{equation}
h(\cdot,\cdot)_p = \pm\langle\mbox{Hess}_p\,F\,\cdot\,,\cdot\rangle=\pm\left\langle\frac{\partial^2F(p)}{\partial x^i\partial x^j}\,\cdot\,,\cdot\right\rangle\,.\label{eq::HessMetric}
\end{equation}
By using Eq. (\ref{eqNormalCurvLevelSets}), Proposition \ref{prop::InterpretCasualChar} is still valid for $\Sigma$ in the context of a Hessian pseudo-metric. Moreover, if $\det(\mbox{Hess}_pF)\not=0$, then $\mbox{Hess}\,F$ is non-degenerate on a neighborhood of $p$. Likewise, since the eigenvalues vary continuously \cite{SerreMatrixBook} and the index can be seen as the number of negative eigenvalues, the Hessian $\mbox{Hess}\,F$ has a constant index on an open neighborhood. Then, $h(\cdot,\cdot)$ in Eq. (\ref{eq::HessMetric}) is well defined on a neighborhood of a non-degenerate point $p$ (for an index 2 or 3 we take $h(\cdot,\cdot)_p=-\langle\mathrm{Hess}_pF\,\cdot,\cdot\rangle$). 

Now we ask ourselves if the techniques developed in the previous sections can be applied to characterize curves that lie on a level surface. Unhappily, we are not able to establish a characterization via a linear equation as previously done. Nonetheless, we can still exhibit a functional relationship between the curvatures $\kappa_1$ and $\kappa_2$ of an RM frame of the corresponding curves with respect to the Hessian metric. Before that, let us try to understand the technical difficulties involved in the study of level surfaces:

\begin{example}[index 1 Hessian]
Suppose that $\mathrm{index}(\mathrm{Hess}\,F)=1$ on a certain neighborhood of a non-degenerate point $p$. Let $\alpha:I\to \mathbb{R}^3$ be a curve on a regular level surface $\Sigma = F^{-1}(c)$ whose velocity vector $\alpha'\in T_{\alpha(s)}\Sigma$ is not an asymptotic direction for all $s\in I$, i.e. $\kappa_n(\alpha(s),\mathbf{\alpha}'(s))\not=0$. This means that the curve is timelike or spacelike. Denote by $\{\mathbf{t},\mathbf{n}_1,\mathbf{n}_2\}$ an RM frame along $\alpha$, with respect to Eq. (\ref{eq::HessMetric}), and denote by $D$ the covariant derivative and by a prime $'$ the usual one. 

From $F(\alpha(s))=c$ it follows that
\begin{equation}
h(\mathrm{grad}_{\alpha(s)}F,\mathbf{t}) = 0\Rightarrow\mathrm{grad}_{\alpha}F=a_1\mathbf{n}_1+a_2\mathbf{n}_2\,,\label{eq:gradF_hessMet_InNormalPlane}
\end{equation}
where $\mathrm{grad}_{\alpha}F$ denotes the gradient vector with respect to $h(\cdot,\cdot)$. The coefficients $a_1$ and $a_2$ satisfy $a_i=\epsilon_i\,h(\mathrm{grad}_{\alpha}F,\mathbf{n}_i)$ and, therefore,
\begin{eqnarray}
\epsilon_ia_i' & = & h(D\,\mathrm{grad}_{\alpha}F, \mathbf{n}_i)+h(\mathrm{grad}_{\alpha}F,D\,\mathbf{n}_i)\nonumber\\
& = & H^F(\mathbf{t}, \mathbf{n}_i)-\epsilon\kappa_i\,h(\mathrm{grad}_{\alpha}F,\mathbf{t})\nonumber\\
& = & H^F(\mathbf{t}, \mathbf{n}_i),
\end{eqnarray}
where $H^F$ denotes the Hessian with respect to $h(\cdot,\cdot)$, whose coefficients can be expressed as \cite{ONeill}
\begin{equation}
H^F_{ij} = \left(\frac{\partial^2F}{\partial x^i\partial x^j}-\sum_k\Gamma_{ij}^k\frac{\partial F}{\partial x^k}\right)\,.
\end{equation}
From this expression we see that $a_i'$ does not need to be zero and then we can not apply the same steps as in the previous sections. Indeed, the orthogonality of the RM frame $\{\mathbf{t},\mathbf{n}_1,\mathbf{n}_2\}$ with respect to $\mathrm{Hess}\,F=\partial^2F/\partial x^i\partial x^j$ and $H^F$ does not coincide, unless $\mathrm{Hess}\,F$ is constant (here we use that $h_{ij}=\partial^2F/\partial x^i\partial^j$ and $h_{ij,k}=\partial_kh_{ij}=\partial^3F/\partial x^i\partial^j\partial x^k$). \label{exe::index1Hessian}
\qed
\end{example}

\begin{theorem}
Let $\mathcal{U}_p\subseteq\mathbb{R}^3$ be a neighborhood of a non-degenerate point $p\in\Sigma=F^{-1}(c)$ where the index is constant. Let $H^F$ denotes the Hessian with respect to the Hessian metric $h(\cdot,\cdot)_q=\langle\mathrm{Hess}_qF\,\cdot,\cdot\rangle$. If $\alpha:I\to\mathcal{U}_p\cap\Sigma$ is a $C^2$ regular curve, with no asymptotic direction if $\mathrm{index}(\mathrm{Hess}\,F
)\not\in\{0,3\}$, i.e., $\kappa_n(\alpha,\alpha')\not=0$, then its normal development $(\kappa_1(s),\kappa_2(s))$ satisfies
\begin{equation}
a_2(s)\kappa_2(s)+a_1(s)\kappa_1(s)+a_0(s)=0, \label{eq::characLevelSurfacesCurves}
\end{equation}
where $a_0=H^F(\mathbf{t},\mathbf{t})$, $a_i=h(\mathrm{grad}_{\alpha}F,\mathbf{n}_i)$, and $a_i'(s)=H^F(\mathbf{t},\mathbf{n}_i)$: or $\epsilon_i H^F(\mathbf{t},\mathbf{n}_i)$, $\epsilon_i=h(\mathbf{n}_i,\mathbf{n}_i)=\pm1$, if $\mathrm{index}(\mathrm{Hess}\,F
)\not\in\{0,3\}$. Here, the RM frame is defined with respect to the Hessian metric. 

Conversely, if Eq. (\ref{eq::characLevelSurfacesCurves}) is valid and $h(\mathrm{grad}_{\alpha(s_0)}F,\mathbf{t}(s_0))=0$ at some point $\alpha(s_0)$, then $\alpha$ lies in a level surface of $F$.
\label{theo::CurvesInLevelSets}
\end{theorem}
\begin{remark}
If $\Sigma=F^{-1}(c)$, where $c$ is a regular value of $F$, then $\Sigma$ is an orientable surface. The reciprocal of this result is also valid, i.e., every orientable surface is the inverse image of a regular value of some smooth function \cite{Guillemin}. Then, the above theorem can be applied to any orientable surface (we still have to exclude those points where the Hessian has a zero determinant).
\end{remark}
\textit{Proof of theorem }\ref{theo::CurvesInLevelSets}.
If the index is $0$, then the Hessian metric defines a Riemannian metric: if $\mbox{index}(\mbox{Hess}\,F)=-3$, then its negative defines a metric. On the other hand, the construction of an RM frame for a pseudo-metric with index 2 in dimension 3 is completely analogous to the case of index 1. Moreover, when the index of $\mbox{Hess}\,F$ is 1 (or 2), the assumption that $\alpha'$ is not an asymptotic direction means that $\alpha$ must be a space- or a timelike curve.  

In the following, let us assume that $\mathrm{index}(\mathrm{Hess}\,F)=1$, the other cases being analogous. In this case, Eq. (\ref{eq::HessMetric}) defines a pseudo-metric in $\mathcal{U}_p\subseteq\mathbb{R}^3$.

Since $F(\alpha(s))=c$, we have
\begin{equation}
h(\mbox{grad}_{\alpha(s)}F,\mathbf{t}) = 0\Rightarrow\mbox{grad}_{\alpha}F=a_1\mathbf{n}_1+a_2\mathbf{n}_2\,,\label{eqGradFHessMetric}
\end{equation}
where $\mbox{grad}_{\alpha}F$ denotes the gradient vector with respect to $h(\cdot,\cdot)$. The coefficients $a_1$ and $a_2$ satisfy $a_i=\epsilon_i\,h(\mbox{grad}_{\alpha}F,\mathbf{n}_i)$ and, therefore,
\begin{eqnarray}
a_i' & = & \epsilon_i\,h(D\,\mbox{grad}_{\alpha}F, \mathbf{n}_i)+\epsilon_i\,h(\mbox{grad}_{\alpha}F,D\,\mathbf{n}_i)\nonumber\\
& = & \epsilon_iH^F(\mathbf{t}, \mathbf{n}_i)-\epsilon_i\epsilon\kappa_i\,h(\mbox{grad}_{\alpha}F,\mathbf{t})\nonumber\\
& = & \epsilon_iH^F(\mathbf{t}, \mathbf{n}_i),
\end{eqnarray}
where $H^F$ denotes the Hessian with respect to $h(\cdot,\cdot)_p$ \cite{ONeill}. Taking the derivative of Eq. (\ref{eqGradFHessMetric}) gives
\begin{eqnarray}
0 & = & h(D\,\mbox{grad}_{\alpha}F,\mathbf{t})+h(\mbox{grad}_{\alpha}F,D\,\mathbf{t})\nonumber\\
& = & H^F(\mathbf{t},\mathbf{t})+h(a_1\mathbf{n}_1+a_2\mathbf{n}_2,\epsilon_1\kappa_1\mathbf{n}_1+\kappa_2\mathbf{n}_2)\nonumber\\
& = & H^F(\mathbf{t},\mathbf{t})+a_1\kappa_1+a_2\kappa_2\,.
\end{eqnarray}
Then, Eq. (\ref{eq::characLevelSurfacesCurves}) is satisfied. 

Conversely, suppose that Eq. (\ref{eq::characLevelSurfacesCurves}) is satisfied. Let us define the function $f(s)=F(\alpha(s))$. We must show that $f$ is constant, i.e., $f'(s)=0$. Taking the derivative of $f$ twice gives
\begin{equation}
f' = h(\mbox{grad}_{\alpha}F,\mathbf{t}),
\end{equation}
and
\begin{eqnarray}
f'' & = & h(D\,\mbox{grad}_{\alpha}F,\mathbf{t})+h(\mbox{grad}_{\alpha}F,D\,\mathbf{t})\nonumber\\
& = & H^F(\mathbf{t},\mathbf{t})+\epsilon_1\kappa_1\,h(\mbox{grad}_{\alpha}F,\mathbf{n}_1)+\kappa_2\,h(\mbox{grad}_{\alpha}F,\mathbf{n}_2)\nonumber\\
& = & 0.
\end{eqnarray}
Then, $f'(s)=h(\mbox{grad}_{\alpha(s)}F(s),\mathbf{t}(s))$ is constant. By assumption, we have $f'(s_0)=0$, then $f(s)=F(\alpha(s))$ is constant on an open neighborhood of $s_0$, i.e., $\alpha$ lies on a level surface of $F$.
\qed

\begin{remark}
The Christoffel symbols $\Gamma_{ij}^k$ of a Hessian metric $g_{ij}=\partial^2F/\partial x^i\partial x^j$ vanish if and only if $\mbox{Hess}\,F$ is constant; which is valid for a quadratic surface, this case being treated in the previous section. 
\end{remark}

If $\mbox{Hess}\,F$ degenerates, i.e., $\det(\mbox{Hess}_pF)=0$ at some points, then the Hessian matrix does not define a metric. Nonetheless, it is still possible to characterize curves on a level surface by using the standard metric of $\mathbb{R}^3$. In fact, it can be used even if $\mbox{Hess}\,F$ is non-degenerate, but in this case we do not have non-degenerate quadrics as a particular instance. The obtained criterion is completely analogous to the previous one in Theorem \ref{theo::CurvesInLevelSets} (\footnote{An alternative approach will be investigated in the next sections through the use of isotropic spaces.}). Indeed, we have

\begin{theorem}
If $\alpha:I\to \mathbb{E}^3\cap \Sigma$ is a $C^2$ regular curve, where $\Sigma=F^{-1}(c)$, then its normal development $(\kappa_1(s),\kappa_2(s))$ satisfies
\begin{equation}
b_2(s)\kappa_2(s)+b_1(s)\kappa_1(s)+b_0(s)=0, \label{eq::characLevelSurfacesCurves2}
\end{equation}
where $b_0=\langle (\mathrm{Hess}\,F)\,\mathbf{t},\mathbf{t}\rangle$, $b_i=\langle\nabla_{\alpha}F,\mathbf{n}_i\rangle$, and $b_i'(s)=\langle (\mathrm{Hess}\,F)\,\mathbf{t},\mathbf{n}_i\rangle$. Here, the rotation minimizing frame is defined with respect to the usual metric in $\mathbb{E}^3$. Conversely, if Eq. (\ref{eq::characLevelSurfacesCurves2}) is valid and $\langle\nabla_{\alpha(s_0)}F,\mathbf{t}(s_0)\rangle=0$ at some point $\alpha(s_0)$, then $\alpha$ lies in a level surface of $F$.
\end{theorem}
\textit{Proof.} Let $\{\mathbf{t},\mathbf{n}_1,\mathbf{n}_2\}$ be an RM frame along $\alpha:I\to \mathbb{E}^3$. If $F(\alpha(s))=c$, then we have
\begin{equation}
\langle\nabla_{\alpha(s)}F,\mathbf{t}\rangle = 0\Rightarrow\nabla_{\alpha}F=b_1\mathbf{n}_1+b_2\mathbf{n}_2\,,\label{eqGradFUsualMetric}
\end{equation}
where $\nabla_{\alpha}F$ denotes the gradient vector with respect to usual metric in $\mathbb{E}^3$. The coefficients $b_1$ and $b_2$ satisfy $b_i=\langle\nabla_{\alpha}F,\mathbf{n}_i\rangle$ and, therefore,
\begin{eqnarray}
b_i' = \langle (\mbox{Hess}\,F)\,\mathbf{t}, \mathbf{n}_i\rangle-\kappa_i\langle\nabla_{\alpha}F,\mathbf{t}\rangle=\langle (\mbox{Hess}\,F)\,\mathbf{t}, \mathbf{n}_i\rangle.
\end{eqnarray}
Taking the derivative of Eq. (\ref{eqGradFUsualMetric}) gives
\begin{eqnarray}
0 & = & \langle (\mbox{Hess}\,F)\,\mathbf{t},\mathbf{t}\rangle+\langle b_1\mathbf{n}_1+b_2\mathbf{n}_2,\kappa_1\mathbf{n}_1+\kappa_2\mathbf{n}_2\rangle\nonumber\\
& = & \langle (\mbox{Hess}\,F)\,\mathbf{t},\mathbf{t}\rangle+b_1\kappa_1+b_2\kappa_2\,.
\end{eqnarray}
Thus, Eq. (\ref{eq::characLevelSurfacesCurves2}) is valid. 

Conversely, suppose that Eq. (\ref{eq::characLevelSurfacesCurves}) is satisfied. Let us define the function $f(s)=F(\alpha(s))$. Taking the derivative of $f$ twice gives
\begin{equation}
f' = \langle\nabla_{\alpha}F,\mathbf{t}\rangle=0,\,
f''
 =  \langle (\mbox{Hess}\,F)\,\mathbf{t},\mathbf{t}\rangle+\kappa_1\langle\nabla_{\alpha}F,\mathbf{n}_1\rangle+\kappa_2\langle\nabla_{\alpha}F,\mathbf{n}_2\rangle= 0.
\end{equation}
Then, $f'(s)=\langle\nabla_{\alpha(s)}F(s),\mathbf{t}(s)\rangle$ is constant. By assumption, we have $f'(s_0)=0$, then $f(s)=F(\alpha(s))$ is constant on an open neighborhood of $s_0$, i.e., $\alpha$ lies on a level surface of $F$.
\qed

\section{Curves on degenerate quadrics: geometry of isotropic spaces}
\label{sec::GeomIso}

In the first section of this chapter we illustrated how the framework developed above applies for non-degenerate quadrics. In such cases the Hessian is constant and this makes the analysis easier, as discussed in example \ref{exe::index1Hessian}. Indeed, such Hessians give rise to a metric, or pseudo-metric, in $\mathbb{R}^3$ that allows for a characterization of curves on (non-degenerate) quadrics through a linear equation (intuitively, these curves are seen as spherical curves in the new geometry setting). However, there are other cases where the Hessian is constant: for cylindrical quadrics, i.e., translation surfaces with a cross section given by a quadric, we basically have
\begin{equation}
\mbox{Hess}_p\,F = \left(
\begin{array}{ccc}
1 & 0 & 0 \\
0 & \delta & 0 \\
0 & 0 &0\\
\end{array}
\right),
\end{equation}
where $\delta\in\{-1,0,1\}$. Here the Hessian is degenerate, but we can still apply the framework developed above to a (semi)metric induced by such a degenerate symmetric matrix. In fact, this leads us to the geometry of isotropic spaces \cite{Giering1982,SipusGM1998,VoglerGMB1989}: simply isotropic space $\mathbb{I}^3$ if $\delta=1$ \cite{Sachs1990,StrubeckerSOA1941}; semi-isotropic space $\mathbb{SI}^3$ if $\delta=-1$ \cite{AydinArXiv2016,daSilvaArXivIso2017}; and doubly isotropic space $\mathbb{I}_2^3$ if $\delta=0$ \cite{BraunerCrelle1967}.  

The just mentioned three dimensional isotropic geometries are examples of 3d Cayley-Klein (CK) geometries \cite{Giering1982,Sulanke2006}. The basic idea behind a CK geometry is the study of those properties in projective space $\mathbb{P}^3$ that preserves a certain configuration, the so called \emph{absolute figure}, i.e., in the spirit of the Klein ``Erlangen Program'' \cite{birkhoff1988felix,klein1893vergleichende}, it is the study of those properties invariant by the action of the subgroup of projective collineations that fix the absolute figure. There are $27$ types of 3d CK geometries \cite{Giering1982}. For example, the Euclidean (Minkowski) space $\mathbb{E}^3$ ($\mathbb{E}_1^3$) is modeled through an absolute figure given by a plane at infinity, identified in homogeneous coordinates with $x_0=0$, and a non-degenerate quadric of index zero (index one), identified with $x_0^2+\cdots+x_3^2=0$ ($x_0^2+x_1^2+x_2^2-x_3^2=0$, respectively) \cite{Giering1982,StruveJG2010}. In our cases of interest, i.e., isotropic space geometries, the absolute figure is given by a plane at infinity, identified with $x_0=0$, and a degenerate quadric of index one or zero, identified with $x_0^2+x_1^2+\delta\, x_2^2=0$. 

Thus, in order to study curves on a degenerate quadric, we may equip $\mathbb{R}^3$ with the following degenerate metric
\begin{equation}
\langle\mathbf{u},\mathbf{v}\rangle_i = u_1v_1 + \delta\,u_2v_2\,.
\end{equation}

Let us now outline how to characterize curves on a (degenerate) quadric $\langle\mathbf{x}-P,\mathbf{x}-P\rangle_i=r^2$ for $\delta=\pm1$: we shall study in more detail the differential geometry of curves on the simply isotropic space $\mathbb{I}^3$ in the next subsection\footnote{Sometimes the simply isotropic space is denoted as $I_3^1$ or $I_3^{(1)}$, see e.g., \cite{Sachs1990,SipusGM1998}.}, which is the geometric setting associated withe a circular cylinder.

The definition of an arc-length parameter, unit tangent, principal normal, and curvature are done as usual (however, we shall avoid curves with isotropic velocity vectors, i.e., $\langle\alpha',\alpha'\rangle_i=0$ but $\alpha'\not=0$):
\begin{equation}
s(t)=\int_{t_0}^t\sqrt{\vert\langle\alpha'(u),\alpha'(u)\rangle_i\vert}\,\mathrm{d}u\,,\,\,\mathbf{t}=\alpha'(s)\,,\,\,\mathbf{n}=\eta\frac{\mathbf{t}'}{\Vert\mathbf{t}'\Vert_i}\,,\,\mbox{ and }\,\kappa(s)=\eta\Vert\alpha''(s)\Vert_i\,,
\end{equation}
where $\eta=\langle\mathbf{n},\mathbf{n}\rangle_i$ when $\delta=-1$, while for $\delta=+1$ we always have $\eta=+1$.
Since the metric $\langle\cdot,\cdot\rangle_i$ is degenerate, we have that the curvature is just the curvature of the projection of $\alpha$ onto the $z=0$ plane (this projection is called the \emph{top view}), which is just the curvature function for a curve in $\mathbb{E}^2$ if $\delta=1$ or $\mathbb{E}_1^2$ if $\delta=-1$.  

It happens that we can not introduce a vector product in isotropic spaces with the same invariance properties as in $\mathbb{E}^3$. Then, in order to complete a trihedron along a curve, we just define the binormal vector $\mathbb{b}$ as the isotropic vector $(0,0,1)$: a vector in the $z$-direction is an \emph{isotropic vector}. Suppose we can introduce RM vector fields (since $\mathbf{b}'=0$, the binormal is RM and then we just need to find a new vector field in substitution to principal normal in order to build an RM frame $\{\mathbf{t},\mathbf{n}_1,\mathbf{n}_2=\mathbf{b}\}$). If $\alpha$ is a curve such that $\langle\alpha(s)-P,\alpha(s)-P\rangle_i=\pm r^2$ and $\{\mathbf{t},\mathbf{n}_1,\mathbf{n}_2=\mathbf{b}\}$ an RM frame along it, then taking the derivative gives
\begin{equation}
\langle\mathbf{t},\alpha-P\rangle_i =0 \Rightarrow \alpha-P = a_1\mathbf{n}_1+a_2\mathbf{b}\,.
\end{equation}
The first coefficient is $a_1=\eta\langle\alpha-P,\mathbf{n}_1\rangle_i$, where $\eta=\langle\mathbf{n}_1,\mathbf{n}_1\rangle_i$. Its derivative is $a_1'=\eta\langle\mathbf{t},\mathbf{n}_1\rangle_i-\eta\kappa_1\langle\alpha-P,\mathbf{t}\rangle_i=0$ and then $a_1$ is a constant\footnote{We can not apply this same strategy to $a_2$, since $\mathbf{b}$ is isotropic, i.e., $\langle\mathbf{b},\mathbf{b}\rangle_i=0$. Nonetheless, we do not need to know $a_2$ in order to find basic informations about $\alpha$.}. Now, taking the derivative of $\langle\mathbf{t},\alpha-P\rangle_i=0$ furnishes
\begin{eqnarray}
0 & = & \langle\mathbf{t},\mathbf{t}\rangle_i+\langle \kappa_1\mathbf{n}_1+\kappa_2\mathbf{b},\alpha-P\rangle_i\nonumber\\
& = & \epsilon +\langle \kappa_1\mathbf{n}_1+\kappa_2\mathbf{b},a_1\mathbf{n}_1+a_2\mathbf{b}\rangle_i\nonumber\\
& = & \epsilon + \eta a_1\kappa_1,  
\end{eqnarray}
where $\epsilon=\langle\mathbf{t},\mathbf{t}\rangle_i=-\eta$. The above result shows that for a cylindrical curve $\alpha$ the first curvature $\kappa_1$ must be a constant. It can be shown that $\kappa_1=\kappa$ (see next subsection). 

On the other hand, if $\kappa=\kappa_1$ is constant, then defining $P=\alpha+\kappa^{-1}\mathbf{n}_1$, we have $P'=\mathbf{t}-\kappa^{-1}\kappa\mathbf{t}=0$. Then, $P$ is a fixed point and $\langle\alpha-P,\alpha-P\rangle_i=\mbox{constant}$. In short, we have
\begin{theorem}
A regular curve $\alpha:I\to\mathbb{R}^3$ lies on a degenerate quadric $\mathcal{Q}=\{(x,y,z):\mathcal{P}(x,y,z)=x^2\pm y^2=r^2\}$ if and only if the (isotropic) curvature function of $\alpha$ with respect to the (isotropic) metric induced by $\mathcal{P}$ in $\mathbb{R}^3$ is a constant.
\end{theorem}

\section{Differential geometry of curves in the simply isotropic space}
\label{Sec::DiffGeomIsoSpaces}

We now discuss in more detail the differential geometry of curves in the simply isotropic space $\mathbb{I}^3$, which is the proper geometric setting to describe curves in a circular cylinder. Motived by the success of RM frames in the study of spherical curves in both $\mathbb{E}^3$ and $\mathbb{E}_1^3$  \cite{BishopMonthly,daSilvaArXiv,Etayo2016,OzdemirMJMS2008}, as described in the previous chapters of this thesis, we develop here the fundamentals of RM frames in the simply isotropic space. We also show, by using the Galilean trigonometric functions \cite{Yaglom1979}, how to relate RM and Frenet frames in $\mathbb{I}^3$. In addition, despite that isotropic spheres can not be always defined as the set of equidistant points from a given center\footnote{As we will see in the following, there are two types of spheres and only one of them corresponds to an equidistant definition.}, we are able to prove that spherical curves can be characterized through a linear equation by conveniently using osculating spheres, in analogy with what happens in $\mathbb{E}^3$ and $\mathbb{E}_1^3$. 

Besides its mathematical interest \cite{AydinJG2015,KaracanTJM2017,SipusPMH2014,YoonJG2017}, see also \cite{Sachs1990} and references therein, isotropic geometry also finds applications in image processing and shape interrogation \cite{koenderink2002image,PottmannCAGD1994}, elasticity \cite{pottmann2009laguerre}, and in economics \cite{AydinTJM2016,chenKJM2014}, just to name a few. Let us now introduce some basic terminology (we refer the reader to Sachs' monograph \cite{Sachs1990} for more details).

Isotropic geometry is the study of those properties in $\mathbb{R}^3$ invariant by the action of the 6-parameter group $\mathcal{B}_6$
\begin{equation}
\left\{
\begin{array}{ccc}
\bar{x} & = & a + x\,\cos\phi-y\,\sin\phi\\ 
\bar{y} & = & b + x\,\sin\phi+y\,\cos\phi\\
\bar{z} & = & c + c_1x+c_2y+z\\
\end{array}
\right.\,.\label{eq::IsoGroupB6}
\end{equation}
So, $\mathcal{B}_6$ forms the set of rigid motions of $\mathbb{I}^3$. In addition, observe that on the $z=0$ plane this geometry looks exactly like the plane Euclidean geometry. The projection of a vector $\mathbf{u}=(u_1,u_2,u_3)\in \mathbb{I}^3$ on the $xy$ plane is called the \emph{top view} of $\mathbf{u}$ and we shall denote it by $\tilde{\mathbf{u}}=(u_1,u_2,0)$. The top view concept plays a fundamental role in the simply isotropic space $\mathbb{I}^3$. In fact, the $z$-direction is preserved under the action of $\mathcal{B}_6$ (\footnote{Maybe, it would be interesting to mention that, from a Physics perspective, such a space is not isotropic. Indeed, the $z$-direction is a distinguished direction and gives rise to an anisotropy (in the physics jargon). Thus, anisotropic geometry would be a better name. Anyway, this is a well established nomenclature and we will not attempt to change it.}). A line with this direction is called an \emph{isotropic line} and a plane that contains an isotropic line is said to be an \emph{isotropic plane}.

The isotropic inner product  between two vectors $\mathbf{u}=(u_1,u_2,u_3)$ and $\mathbf{v}=(v_1,v_2,v_3)$ is defined as
\begin{equation}
\langle\mathbf{u},\mathbf{v}\rangle_{z} = u_1v_1+u_2v_2\,,
\end{equation}
from which we define an isotropic distance as usual\footnote{The index $z$ is here just to emphasize that $z$ is the isotropic (degenerate) direction. Note, in addition, that the isotropic inner product induces in fact a semi-distance in $\mathbb{R}^3$, since points in the isotropic line have zero distance.}: 
\begin{equation}
\mathrm{d}_z(A,B)=\sqrt{\langle B-A,B-A\rangle_z}\,.
\end{equation}
 Note that the inner product and distance above are just the plane Euclidean counterparts of the top views. Since the metric is degenerate, the distance from $(u_1,u_2,u_3)$ to $(u_1,u_2,v_3)$ is zero. In such cases, one may define a codistance by $\mathrm{cd}_z(A,B)=\vert b_3-a_3\vert$ (the codistance is preserved by $\mathcal{B}_6$ and then is an isotropic invariant: it can be used to define angles involving isotropic lines and planes \cite{PottmannCAGD1994,Sachs1990}).

Now we introduce some terminology related to curves. A regular curve $\alpha:I\to \mathbb{I}^3$, i.e., $\alpha'\not=0$, is parametrized by an arc-length $s$ if $\Vert\alpha'(s)\Vert_z\stackrel{\mathrm{def}}{=}\Vert\tilde{\alpha}'(s)\Vert=1$. In the following we assume that all the curves are parametrized by an arc-length $s$ (in particular, this excludes the possibility of an isotropic velocity vector). In addition, a point $\alpha(s_0)$ where $\{\alpha'(s_0),\alpha''(s_0)\}$ is linearly dependent is an \emph{inflection point} and a regular unit speed curve $\alpha(s)=(x(s),y(s),z(s))$ with no inflection point is called an \emph{admissible curve} if $x'y''-x''y'\not=0$.
\begin{remark}
The admissible condition implies that the osculating planes, i.e., the planes that have a contact of order 2 with the reference curve\footnote{For a level set surface $\Sigma=G^{-1}(c)$, a contact of order $k$ with $\alpha$ at $\alpha(s_0)$ is equivalent to say that $\beta^{(i)}(s_0)=0$ ($1\leq i\leq k$), where $\beta=G\circ\alpha$ and $c=\beta(s_0)=\alpha(s_0)$ \cite{Kreyszig1991}.}, can not be isotropic. Moreover, the only curves with $x'y''-x''y'=0$ are precisely the isotropic lines \cite{Sachs1990}.
\end{remark}

\subsection{Isotropic Frenet frame}

The (isotropic) unit tangent, principal normal, and curvature function are defined as usual
\begin{equation}
\mathbf{t}(s)=\alpha'(s),\,\mathbf{n}(s)=\frac{\mathbf{t}'(s)}{\kappa(s)},\,\,\mathrm{ and }\,\,\kappa(s)=\Vert\mathbf{t}'(s)\Vert_z=\Vert\tilde{\mathbf{t}}'(s)\Vert,
\end{equation}
respectively. As usually happens in isotropic geometry, the curvature $\kappa$ is just the curvature function of its top view $\tilde{\alpha}$ and then $\kappa=x'y''-x''y'$. To complete the moving trihedron, we define the binormal vector as the (co)unit vector $\mathbf{b}=(0,0,1)$ in the isotropic direction. The three vectors $\mathbf{t},\mathbf{n},\mathbf{b}$ are linearly independent:
\begin{equation}
\det(\mathbf{t},\mathbf{n},\mathbf{b})=\frac{1}{\kappa}(x'y''-x''y')=1\,.
\end{equation}

The Frenet equations corresponding to the isotropic Frenet frame $\{\mathbf{t},\mathbf{n},\mathbf{b}\}$ can be written as 
\begin{equation}
\frac{\mathrm{d}}{\mathrm{d}s}\left(
\begin{array}{c}
\mathbf{t}\\
\mathbf{n}\\
\mathbf{b}\\
\end{array}
\right)=\left(
\begin{array}{ccr}
0\,\, & \kappa &\,\, 0\\
-\kappa\,\, & 0 & \,\,\tau\\
0\,\, & 0 &\,\, 0 \\
\end{array}
\right)\left(
\begin{array}{c}
\mathbf{t}\\
\mathbf{n}\\
\mathbf{b}\\
\end{array}
\right),
\end{equation}
where $\tau$ is the (isotropic) torsion:
\begin{equation}
\tau=\frac{\det(\alpha',\alpha'',\alpha''')}{\det(\tilde{\alpha}',\tilde{\alpha}'')}\,;\,\,\kappa=\frac{\det(\tilde{\alpha}',\tilde{\alpha}'')}{\sqrt{\langle\alpha',\alpha'\rangle_z}\,^3}\,.
\end{equation}
The above expressions for the torsion and curvature are also valid for a generic regular parameter for $\alpha$ and, in addition, they are invariant by rigid motions in $\mathbb{I}^3$. Contrary to the Euclidean space $\mathbb{E}^3$, we can not define the torsion through the derivative of the binormal vector. However, we should remember that the idea behind such a definition is that in $\mathbb{E}^3$ one can measure the variation of the osculating plane by measuring $\mathbf{b}'$. It can be shown that the isotropic torsion is directly associated with the velocity of variation of the osculating plane, see \cite{Sachs1990}, pp. 112-113. On the other hand, contrary to the isotropic curvature, the torsion is not defined as the torsion of the top view (this would result in $\tau=0$). The isotropic torsion is an intermediate concept depending on its top view behavior and on how much the curve leaves a plane. In fact, an admissible curve lies on a non-isotropic plane if and only if its torsion vanishes identically. 

\subsection{Isotropic osculating spheres}

Due to the degeneracy of the isotropic metric, some geometric concepts can not be uniquely defined using 
$\langle\cdot,\cdot\rangle_z$. This is the case for spheres.

\begin{definition}
We define \emph{isotropic spheres} as connected and irreducible surface of degree 2 given by the 4-parameter family\footnote{Rigorously speaking, isotropic spheres are connected and irreducible surfaces of degree 2 in $\mathbb{P}^3$ that contains the absolute figure (in fact, this definition applies to any CK geometry). One then shows that in $\mathbb{I}^3$ this condition is satisfied by the 4-parameter family in Eq. (\ref{eq::IsoSpheres}) \cite{Sachs1990}.}
\begin{equation}
(x^2+y^2)+2c_1x+2c_2y+2c_3z+c_4=0\,.\label{eq::IsoSpheres}
\end{equation}
In addition, up to a rigid motion (in $\mathbb{I}^3$), we can express a sphere in one of the two normal forms below
\begin{enumerate}
\item (sphere of parabolic type) 
\begin{equation}
z = \frac{1}{2p}(x^2+y^2)\,\mbox{ with }\,p\not=0;
\end{equation}
\item (sphere of cylindrical type) 
\begin{equation}
x^2+y^2=r^2\,\mbox{ with }\,r>0.
\end{equation}
\end{enumerate}
\end{definition}

It can be shown that the quantities $R=1/2p$ and $r$ are isotropic invariants. Moreover, spheres of cylindrical type are precisely the set of points equidistant from a given center
\begin{equation}
\langle\mathbf{x}-P,\mathbf{x}-P\rangle_z = r^2.
\end{equation}
Observe however, that the center $P$ of a cylindrical sphere is not well defined. More precisely, any other point $Q$ with the same top view as $P$, i.e., $\tilde{Q}=\tilde{P}$, would do the same job! We can remedy this by assuming the center located on the $z=0$ plane.

An osculating sphere of an admissible curve $\alpha$ at a point $\alpha(s_0)$ is the (isotropic) sphere that has a contact of order 3 with $\alpha$. The position vector $\mathbf{x}$ of an osculating sphere can be conveniently written as (Eq. (7.18) of  \cite{Sachs1990})
\begin{equation}
\lambda\langle\mathbf{x}-\mathbf{x}_0,\mathbf{x}-\mathbf{x}_0\rangle_z+\langle\mathbf{u},\mathbf{x}-\mathbf{x}_0\rangle = 0\,,\label{eq::IsoOscSphere}
\end{equation}
where $\mathbf{x}_0=\alpha(s_0)$, $\langle\cdot,\cdot\rangle$ is the usual inner product in Euclidean space $\mathbb{E}^3$, and $\lambda\in\mathbb{R}$ and $\mathbf{u}\in\mathbb{R}^3$ are constants.

\subsection{Rotation minimizing vector fields}

Let $\alpha:I\mapsto\mathbb{I}^3$ be an admissible curve parametrized by an arc-length $s$. A normal vector field $\mathbf{v}$ is an RM vector field if $\mathbf{v}'=\mu\, \mathbf{t}$, for some function $\mu$. We easily see that the binormal $\mathbf{b}$ is an RM field, $\mathbf{b}'=\mathbf{0}$. On the other hand, in general the principal normal is not RM, since $\mathbf{n}'=-\kappa\mathbf{t}+\tau\mathbf{b}$.

If $\mathbf{v}$ is a normal vector, then we may write
\begin{equation}
\mathbf{v} = \mu\mathbf{n}+\nu\mathbf{b}\,,
\end{equation}
where we suppose $\mu\not=0$ (otherwise $\mathbf{v}$ is just a multiple of $\mathbf{b}$). Now, imposing $\langle\mathbf{v},\mathbf{v}\rangle_z=1$ implies that
\begin{equation}
1=\langle\mathbf{v},\mathbf{v}\rangle_z=\mu^2\langle\mathbf{n},\mathbf{n}\rangle_z\Rightarrow \mu = \pm 1\,.
\end{equation}
The derivative of $\mathbf{v}$ is 
\begin{eqnarray}
\mathbf{v}' 
& = & \mu\mathbf{n}' +\nu'\, \mathbf{b}\nonumber\\
& = & -\mu\kappa\,\mathbf{t}+\mu\tau\,\mathbf{b}+\nu'\,\mathbf{b}\nonumber\\
& = & -\mu\kappa\,\mathbf{t}+(\mu\tau+\nu')\,\mathbf{b}\,.
\end{eqnarray}
Thus, if we assume $\mathbf{v}$ to be an RM vector field, it follows that
\begin{equation}
\mathbf{v}'\parallel \mathbf{t} \Rightarrow \nu = -\mu \int \tau+\mbox{constant}\,.
\end{equation}
Finally, if impose that $\{\mathbf{t},\mathbf{v},\mathbf{b}\}$ has the same orientation as  $\{\mathbf{t},\mathbf{n},\mathbf{b}\}$, we conclude that 
\begin{equation}
1 = \det(\mathbf{t},\mathbf{v},\mathbf{b}) =  \det(\mathbf{t},\alpha\mathbf{n},\mathbf{b}) =\mu\,.
\end{equation}
\begin{remark}
Using the definition of the Galilean trigonometric functions, i.e., $\cosg\,\phi=1$ and $\sing\,\phi=\phi$ \cite{Yaglom1979}, we can write an RM vector field $\mathbf{v}$ in terms of the Frenet frame as
\begin{equation}
\left\{
\begin{array}{l}
\mathbf{v}  =  \cosg(\theta)\,\mathbf{n}-\sing(\theta)\,\mathbf{b}\\
\theta'  =  \tau
\end{array}
\right.\,.
\end{equation}
This is analogous to RM frames in both Euclidean and Lorentz-Minkowski spaces \cite{BishopMonthly,OzdemirMJMS2008}.\label{remark::GalTrigFunc}
\end{remark}
\begin{prop}
Let $\mathbf{v}$ be a unit normal vector field along $\alpha$. If $\mathbf{v}$ is RM and $\{\mathbf{t},\mathbf{v},\mathbf{b}\}$ has the same orientation as the Frenet frame, then
\begin{equation}
 \mathbf{v}(s)=\mathbf{n}(s)-\left(\int_{s_0}^s \tau(x)\mathrm{d}x+\tau_0\right)\,\mathbf{b}(s)\,,
\end{equation}
where $\tau_0$ is a constant and we shall define $\sing\,\theta(s)=\theta(s)=\int_{s_0}^s \tau(x)\mathrm{d}x+\tau_0$.
\end{prop}

\begin{theorem}
A rotation minimizing frame $\{\mathbf{t},\mathbf{n}_1,\mathbf{n}_2=\mathbf{b}\}$ in isotropic space $\mathbb{I}^3$ satisfies
\begin{equation}
\frac{\mathrm{d}}{\mathrm{d}s}\left(
\begin{array}{c}
\mathbf{t}\\
\mathbf{n}_1\\
\mathbf{n}_2\\
\end{array}
\right)=\left(
\begin{array}{ccc}
0 & \kappa_1 & \kappa_2\\
-\kappa_1 & 0 & 0\\
0 & 0 & 0\\
\end{array}
\right)\left(
\begin{array}{c}
\mathbf{t}\\
\mathbf{n}_1\\
\mathbf{n}_2\\
\end{array}
\right),
\end{equation}
where the natural curvatures are $\kappa_1=\kappa$ and $\kappa_2=\kappa\,\theta$.
\end{theorem}
\textit{Proof.} The equation for $\mathbf{b}'$ is obvious. For the derivative of $\mathbf{n}_1$ we have
\begin{equation}
\mathbf{n}_1'=-\kappa\,\mathbf{t}+\tau\,\mathbf{b}-\tau\,\mathbf{b}=-\kappa\,\mathbf{t}\,.
\end{equation}
Finally, taking into account that $\mathbf{n}=\mathbf{n}_1+\theta\,\mathbf{b}$, we find
\begin{equation}
\mathbf{t}' = \kappa\,\mathbf{n}=\kappa\,\mathbf{n}_1+\kappa\theta\,\mathbf{b}\,.
\end{equation}
From the equalities above we find the desired equations of motion for the trihedron $\{\mathbf{t},\mathbf{n}_1,\mathbf{b}\}$.
\qed

Using the definition for the Galilean trigonometric functions again, we can relate the RM frame curvatures $(\kappa_1,\kappa_2)$ with the Frenet ones $(\kappa,\tau)$ according to
\begin{equation}
\left\{
\begin{array}{c}
\kappa_1(s) = \kappa(s)\,\cosg\,\theta(s)\\[4pt]
\kappa_2(s) = \kappa(s)\,\sing\,\theta(s)\\[4pt]
\theta'(s)= \tau(s)\\
\end{array}
\right.\,.\label{eq::IsoRelBetweenRMandFrenet}
\end{equation}

\subsubsection{\textit{Moving bivectors}}

In $\mathbb{I}^3$ it is not possible to define a vector product with the same invariance significance as in Euclidean space. However, one can still do some interesting investigations by employing the usual vector product $\vecpe$ from Euclidean space $\mathbb{E}^3$ in isotropic space $\mathbb{I}^3$. Associated with the isotropic Frenet frame, we introduce a (moving) bivector frame  
\begin{equation}
\left\{
\begin{array}{ccl}
\mathcal{T} &=& \mathbf{n}\vecpe\mathbf{b}\\
\mathcal{H} &=& \mathbf{b}\vecpe\mathbf{t}\\
\mathcal{B} &=& \mathbf{t}\vecpe\mathbf{n}\\
\end{array}
\right.\,,
\end{equation}
which satisfies the equation
\begin{equation}
\frac{\mathrm{d}}{\mathrm{d}s}\left(
\begin{array}{c}
\mathcal{T}\\
\mathcal{H}\\
\mathcal{B}\\
\end{array}
\right)=\left(
\begin{array}{ccr}
0 & \kappa & \,\,0\\
-\kappa & 0 & \,\,0\\
0 & -\tau & \,\,0\\
\end{array}
\right)\left(
\begin{array}{c}
\mathcal{T}\\
\mathcal{H}\\
\mathcal{B}\\\end{array}
\right)\,
\end{equation}
and \cite{Sachs1990}, Eqs. (7.43a-c), p. 130,
\begin{equation}
\det(\mathcal{T},\mathcal{H},\mathcal{B})=\det(\mathbf{t},\mathbf{n},\mathbf{b})=1\,\,\mbox{ and }\,\,\mathcal{T}=\tilde{\mathbf{t}}.
\end{equation}

Analogously, we shall introduce  the following (moving) RM bivector frame associated with an RM frame $\{\mathbf{t},\mathbf{n}_1,\mathbf{n}_2=\mathbf{b}\}$
\begin{equation}
\left\{
\begin{array}{lcl}
\mathcal{T} &=& \mathbf{n}_1\vecpe\mathbf{n}_2 = \mathbf{n}\vecpe\mathbf{b}\\
\mathcal{N}_1 &=& \mathbf{n}_2\vecpe\mathbf{t}\\
\mathcal{N}_2 &=& \mathbf{t}\vecpe\mathbf{n}_1
\end{array}
\right.\,.
\end{equation}
\begin{lemma}
The moving frame $\{\mathcal{T},\mathcal{N}_1,\mathcal{N}_2\}$ forms a basis for $\mathbb{R}^3$.
\end{lemma}
\textit{Proof.} We have 
\begin{equation}
\mathbf{n}_1\vecpe\mathbf{n}_2=(\mathbf{n}-\theta\,\mathbf{b})\vecpe\mathbf{b}=\mathbf{n}\vecpe\mathbf{b}=\mathcal{T}.
\end{equation}
In addition
\begin{equation}
\mathcal{N}_1=\mathbf{b}\vecpe\mathbf{t}=\mathcal{H},
\end{equation}
and
\begin{equation}
\mathcal{N}_2=\mathbf{t}\vecpe(\mathbf{n}-\theta\,\mathbf{b})=\mathcal{B}+\theta\,\mathcal{H}\,.
\end{equation}

Then, we find $\det(\mathcal{T},\mathcal{N}_1,\mathcal{N}_2)=\det(\mathcal{T},\mathcal{H},\mathcal{B})=1$.
\qed

\begin{prop}
A moving RM bivector frame satisfies the equation
\begin{equation}
\frac{\mathrm{d}}{\mathrm{d}s}\left(
\begin{array}{c}
\mathcal{T}\\
\mathcal{N}_1\\
\mathcal{N}_2\\
\end{array}
\right)=\left(
\begin{array}{ccr}
0 & \kappa_1 & \,\,0\\
-\kappa_1 & 0 & \,\,0\\
-\kappa_2 & 0 & \,\,0\\
\end{array}
\right)\left(
\begin{array}{c}
\mathcal{T}\\
\mathcal{N}_1\\
\mathcal{N}_2\\
\end{array}
\right)\,,
\end{equation}
where $\kappa_1=\kappa$ and $\kappa_2=\kappa\,\theta$.
\end{prop}
\textit{Proof.} Using the definitions of the bivectors, we have
\begin{equation}
\mathcal{T}'=\mathbf{n}_1'\vecpe\mathbf{n}_2+\mathbf{n}_1\vecpe\mathbf{n}_2'=\kappa_1\mathcal{N}_1,
\end{equation}
\begin{equation}
\mathcal{N}_1'=\mathbf{n}_2'\vecpe\mathbf{t}+\mathbf{n}_2\vecpe\mathbf{t}'=\mathbf{n}_2\vecpe(\kappa_1\mathbf{n}_1+\kappa_2\mathbf{n}_2)=-\kappa_1\mathcal{T},
\end{equation}
and
\begin{equation}
\mathcal{N}_2'=\mathbf{t}'\vecpe\mathbf{n}_1+\mathbf{t}\vecpe\mathbf{n}_1'=(\kappa_1\mathbf{n}_1+\kappa_2\mathbf{n}_2)\vecpe\mathbf{n}_1=-\kappa_2\mathcal{T}.
\end{equation}
\qed

\subsection{Spherical curves in the simply isotropic space}

Our approach to spherical curves is based on order of contact. More precisely, we first investigate spheres in $\mathbb{I}^3$ with an order 3 contact with $\alpha$ by using RM frames and their associated bivector frames, i.e., we describe the osculating spheres. Then, we use that a curve is spherical when its osculating spheres are all equal to the sphere that contains the curve. We refer to the proof of theorem \ref{thr::charSphCurves}, p. \ref{thr::charSphCurves}, for a similar approach in the simpler setting of Euclidean spherical curves.

Defining a function $F(\mathbf{x})=\lambda\langle\mathbf{x}-\alpha_0,\mathbf{x}-\alpha_0\rangle_z+\langle\mathbf{u},\mathbf{x}-\alpha_0\rangle$, where $\alpha_0=\alpha(s_0)$, $\langle\cdot,\cdot\rangle$ is the usual inner product in $\mathbb{E}^3$, and $\lambda,\mathbf{u}$ are constants to be determined. We have for the derivatives of $(F\circ\alpha)(s)$
\begin{eqnarray}
\left\{
\begin{array}{ccc}
F' & = & 2\lambda\langle\alpha(s)-\alpha_0,\mathbf{t}\rangle_z+\langle\mathbf{u},\mathbf{t}\rangle,\\[6pt]
F'' & = & 2\lambda\langle\mathbf{t},\mathbf{t}\rangle_z+2\lambda\langle\alpha(s)-\alpha_0,\sum_i\kappa_i\mathbf{n}_i\rangle_z+\langle\mathbf{u},\sum_i\kappa_i\mathbf{n}_i\rangle,\\[6pt]
F''' & = & 2\lambda\langle\alpha(s)-\alpha_0,-\kappa_1^2\mathbf{t}+\sum_i\kappa_i'\mathbf{n}_i\rangle_z+\langle\mathbf{u},-\kappa_1^2\mathbf{t}+\sum_i\kappa_i'\mathbf{n}_i\rangle\,\\
\end{array}
\right..
\end{eqnarray}
Imposing the condition $(F\circ\alpha)'(s_0)=(F\circ\alpha)''(s_0)=(F\circ\alpha)'''(s_0)=0$ (contact of order 3) gives
\begin{equation}
\left\{
\begin{array}{c}
\langle\mathbf{u},\mathbf{t}(s_0)\rangle = 0\\[5pt]
2\lambda = - \langle\mathbf{u},\sum_i\kappa_i(s_0)\mathbf{n}_i(s_0)\rangle\\[5pt]
\langle\mathbf{u},\sum_i\kappa_i'(s_0)\mathbf{n}_i(s_0)\rangle = 0\\
\end{array}
\right..
\end{equation}
From the first and third equations above, we find that
\begin{equation}
\mathbf{u}=\rho\,[\mathbf{t}\vecpe(\kappa_1'\mathbf{n}_1+\kappa_2'\mathbf{n}_2)](s_0)=\rho\,[\kappa_1'\mathcal{N}_2-\kappa_2'\mathcal{N}_1](s_0),
\end{equation}
for some constant $\rho\not=0$. On the other hand, from the second equation, we find that
\begin{equation}
2\lambda+\rho\,[\kappa_1'\kappa_2\langle\mathbf{n}_2,\mathcal{N}_2\rangle-\kappa_1\kappa_2'\langle\mathbf{n}_1,\mathcal{N}_1\rangle](s_0)=0.
\end{equation}
It can be easily verified that $\langle\mathbf{n}_i,\mathcal{N}_i\rangle = \det(\mathbf{t},\mathbf{n}_1,\mathbf{n}_2)=1$, and then we can rewrite the expression above as
\begin{equation}
2\lambda = \rho\,[\kappa_1\kappa_2'-\kappa_1'\kappa_2](s_0)=\rho\,\tau(s_0)\kappa^2(s_0),
\end{equation}
where in the last equality one should use the expressions form $(\kappa_1,\kappa_2)$ in terms of $(\kappa,\tau)$, see remark \ref{remark::GalTrigFunc} and Eq.  (\ref{eq::IsoRelBetweenRMandFrenet}).

In short, the equation for the isotropic osculating sphere (\ref{eq::IsoOscSphere}), with respect to an RM frame and its associated bivector frame, can be written as
\begin{equation}
\tilde{\alpha}^2-2\left\langle\alpha,\tilde{\alpha}_0+\frac{\kappa_2'\mathcal{N}_1}{\tau\kappa^2}\vert_{s_0}-\frac{\kappa_1'\mathcal{N}_2}{\tau\kappa^2}\vert_{s_0}\right\rangle+2\left[\frac{\tilde{\alpha}_0^2}{2}-\left\langle\alpha_0,\frac{\kappa_1'\mathcal{N}_2-\kappa_2'\mathcal{N}_1}{\tau\kappa^2}\vert_{s_0}\right\rangle\right]=0\,.
\end{equation}
\begin{theorem}
An admissible regular curve $\alpha:I\to \mathbb{I}^3$ lies on the surface of a sphere if and only if its normal development, i.e., the curve $(\kappa_1(s),\kappa_2(s))$, lies on a line not passing through the origin. In addition, $\alpha$ is a spherical curve of cylindrical type with radius $r$ if and only if $\kappa$ is constant and equal to $r^{-1}$. 
\end{theorem}
\textit{Proof.} The condition of being spherical implies that the isotropic osculating spheres are constant (and equal to the sphere that contains the curve). This condition demands
\begin{equation}
\frac{\mathrm{d}}{\mathrm{d}s}\left[\tilde{\alpha}+\frac{\kappa_2'\mathcal{N}_1}{\tau\kappa^2}-\frac{\kappa_1'\mathcal{N}_2}{\tau\kappa^2}\right]=0\label{eq::CondToBeSphericalVector}
\end{equation}
and
\begin{equation}
\frac{\mathrm{d}}{\mathrm{d}s}\left[\frac{\tilde{\alpha}^2}{2}-\left\langle\alpha,\frac{\kappa_1'\mathcal{N}_2-\kappa_2'\mathcal{N}_1}{\tau\kappa^2}\right\rangle\right]=\frac{\mathrm{d}}{\mathrm{d}s}\left[\left\langle\alpha,\frac{\tilde{\alpha}}{2}-\frac{\kappa_1'\mathcal{N}_2-\kappa_2'\mathcal{N}_1}{\tau\kappa^2}\right\rangle\right]=0\,,\label{eq::CondToBeSphericalScalar}
\end{equation}
where we used that $\tilde{\alpha}^2=\langle\alpha,\alpha\rangle_z=\langle\alpha,\tilde{\alpha}\rangle$.

The first condition gives
\begin{eqnarray}
0 & = & \tilde{\mathbf{t}}+\left(\frac{\kappa_2'}{\tau\kappa^2}\right)'\mathcal{N}_1-\left(\frac{\kappa_1'}{\tau\kappa^2}\right)'\mathcal{N}_2+\frac{\kappa_2'}{\tau\kappa^2}(-\kappa_1\mathcal{T})-\frac{\kappa_1'}{\tau\kappa^2}(-\kappa_2\mathcal{T})\nonumber\\
& = & \left(\frac{\kappa_2'}{\tau\kappa^2}\right)'\mathcal{N}_1-\left(\frac{\kappa_1'}{\tau\kappa^2}\right)'\mathcal{N}_2,
\end{eqnarray}
which, by taking into account the linear independence of $\{\mathcal{N}_1,\mathcal{N}_2\}$, implies that 
\begin{equation}
a_1:=-\frac{\kappa_2'}{\tau\kappa^2}=\mbox{constant};\,a_2:=\frac{\kappa_1'}{\tau\kappa^2}=\mbox{constant}.
\end{equation}
On the other hand, condition (\ref{eq::CondToBeSphericalScalar}) implies
\begin{eqnarray}
0 & = & \left\langle\alpha,\frac{\mathrm{d}}{\mathrm{d}s}\left(\frac{\tilde{\alpha}}{2}-\sum_ia_i\mathcal{N}_i\right)\right\rangle+\left\langle\mathbf{t},\frac{\tilde{\alpha}}{2}-\sum_ia_i\mathcal{N}_i\right\rangle\nonumber\\
& = & (1+a_1\kappa_1+a_2\kappa_2)\langle\alpha,\mathbf{t}\rangle_z,
\end{eqnarray}
where we used that $\langle\tilde{\alpha},\mathbf{t}\rangle=\langle\alpha,\tilde{\mathbf{t}}\rangle=\langle\alpha,\mathbf{t}\rangle_z$ to obtain the second equality. If the curve is not of cylindrical type, we can not have $\langle\alpha,\alpha\rangle_z=\mbox{constant}$, and then we conclude that for a parabolic spherical curve the normal development $(\kappa_1,\kappa_2)$ lies on a line not passing through the origin.

On the other hand, if the curve is of cylindrical type $\langle\alpha(s)-P,\alpha(s)-P\rangle_z=r^2$, taking the derivative gives
\begin{equation}
\langle\mathbf{t},\alpha-P\rangle_z=0\,.\label{eq::InnerProdTangwithalphaP}
\end{equation}
Then $\alpha-P=a_1\mathbf{n}_1+a_2\mathbf{n}_2$. We have that $a_1=\langle\alpha-P,\mathbf{n}_1\rangle_z$ and, therefore, $a_1'=\langle\mathbf{t},\mathbf{n}_1\rangle_z+\langle\alpha-P,-\kappa_1\mathbf{t}\rangle_z=0$ and $a_1$ is a constant.

Taking the derivative of Eq. (\ref{eq::InnerProdTangwithalphaP}) gives
\begin{eqnarray}
0 & = & \langle\mathbf{t},\mathbf{t}\rangle_z+\langle\kappa_1\mathbf{n}_1+\kappa_2\mathbf{n}_2,\alpha-P\rangle_z\nonumber\\
& = & 1 + \langle\kappa_1\mathbf{n}_1+\kappa_2\mathbf{n}_2,a_1\mathbf{n}_1+a_2\mathbf{n}_2\rangle_z\nonumber\\
& = & 1+a_1\kappa.
\end{eqnarray}
Then, the curvature $\kappa=\kappa_1$ is a constant and, in addition, $r^2=\langle\alpha-P,\alpha-P\rangle_z=\langle a_1\mathbf{n}_1+a_2\mathbf{n}_2,a_1\mathbf{n}_1+a_2\mathbf{n}_2\rangle_z=a_1^2$.

Reciprocally, if $\kappa$ is a (non-zero) constant, define $P=\alpha+\kappa^{-1}\mathbf{n}_1$. Taking the derivative gives $P'=\mathbf{t}+\kappa^{-1}(-\kappa\mathbf{t})=0$ and then $P$ is a constant. Clearly we have $\langle\alpha-P,\alpha-P\rangle_z=1/\kappa^2$. 
\qed

To complete our analysis, let us describe the case where the normal development curve is a straight line passing through the origin.

\begin{prop}
An admissible regular curve $\alpha:I\to \mathbb{I}^3$ lies on a plane if and only if its normal development $(\kappa_1(s),\kappa_2(s))$ lies on a line passing through the origin. 
\end{prop}
\textit{Proof.} It is known that $\alpha$ is a plane curve if and only if all its osculating planes are equal to the plane that contains the curve. Define a function $F(\mathbf{x})=\langle\mathbf{x}-\alpha_0,\mathbf{u}\rangle$, where $\alpha_0=\alpha(s_0)$ and $\langle\mathbf{u},\mathbf{u}\rangle=1$. Taking the derivatives of $F\circ\alpha$ twice and demand a contact of order 2, we have
\begin{equation}
\left\{
\begin{array}{ccccc}
(F\circ\alpha)'(s_0) & = & \langle\mathbf{t}(s_0),\mathbf{u}\rangle & = & 0\\[6pt]
(F\circ\alpha)''(s_0) & = & \langle[\kappa_1\mathbf{n}_1+\kappa_2\mathbf{n}_2]\vert_{s_0},\mathbf{u}\rangle & = & 0\\
\end{array}
\right..
\end{equation}
From these equations we deduce that
\begin{equation}
\mathbf{u}=\mathbf{u}(s_0) = \rho(s_0)\,[\mathbf{t}\times_e(\kappa_1\mathbf{n}_1+\kappa_2\mathbf{n}_2)]\vert_{s_0}=\rho(s_0)\,[\kappa_1\mathcal{N}_2-\kappa_2\mathcal{N}_1]\vert_{s_0}\,,
\end{equation}
where, by applying the definition of the Frenet and RM bivectors, we can write $\rho=(\kappa_1\Vert\mathcal{B}\Vert)^{-1}$. 

The condition of being a plane curve is equivalent to $\mathrm{d}\mathbf{u}/\mathrm{d}s=0$. This leads to
\begin{eqnarray}
\frac{\mathrm{d}\mathbf{u}}{\mathrm{d}s} & = & -\frac{1}{\kappa_1\Vert\mathcal{B}\Vert}\left[\left(\frac{\kappa_1\kappa_2'-\kappa_1'\kappa_2}{\kappa_1}+\frac{\kappa_2\tau\langle\mathcal{B},\mathcal{H}\rangle}{\langle\mathcal{B},\mathcal{B}\rangle}\right)\mathcal{N}_1+\frac{\kappa_1\tau\langle\mathcal{B},\mathcal{H}\rangle}{\langle\mathcal{B},\mathcal{B}\rangle}\mathcal{N}_2\right]\nonumber\\
& = &-\frac{\tau}{\Vert\mathcal{B}\Vert}\left[\left(1+\frac{\theta\langle\mathcal{B},\mathcal{H}\rangle}{\langle\mathcal{B},\mathcal{B}\rangle}\right)\mathcal{N}_1+\frac{\langle\mathcal{B},\mathcal{H}\rangle}{\langle\mathcal{B},\mathcal{B}\rangle}\mathcal{N}_2\right],
\end{eqnarray}
where we used that $\tau\kappa^2=\kappa_1\kappa_2'-\kappa_1'\kappa_2$, $\mathcal{N}_1=\mathcal{H}$, $\langle\mathcal{H},\mathcal{H}\rangle=1$, and $\langle\mathcal{T},\mathcal{T}\rangle=1$.

Finally, it is easy to see that the planarity condition, i.e., $\tau=0\Leftrightarrow (\kappa_2/\kappa_1)'=\kappa_1\kappa_2'-\kappa_1'\kappa_2=0$, is equivalent to $\mathbf{u}'=0$, from which we deduce that it is equivalent to $\kappa_2/\kappa_1=\mbox{constant}$ and then $(\kappa_1,\kappa_2)$ lies on a line passing through the origin.
\qed
\chapter{ROTATION MINIMIZING FRAMES AND (GEODESIC) SPHERICAL CURVES IN RIEMANNIAN GEOMETRY}
\label{chap::RMFinRiemGeom}

It is possible to introduce Frenet frames along curves in Riemannian manifolds, see e.g. \cite{BolcskeiBAG2007,GutkinJGP2011,LucasJMAA2015,LucasMJM2016,SzilagyiSUZ2003}, and also rotation minimizing frames \cite{Etayo2016,EtayoTJM2017}(\footnote{It is worth mentioning that a Frenet-like theorem is valid only for manifolds of constant curvature ($\mathbb{E}^{m+1}$, $\mathbb{S}^{m+1}(r)$, and $\mathbb{H}^{m+1}(r)$ are the prototypes of such spaces \cite{doCarmo1992}), i.e., two curves are congruent if and only if they have the same curvatures \cite{CastrillonLopezAM2015,CastrillonLopezDGA2014}.}). In order to introduce such concepts, one should take covariant derivatives in the direction of the unit tangent instead of the ordinary derivation. More precisely, let $M^{m+1}$ be a Riemannian manifold with Levi-Civita connection $\nabla$ and $\alpha:I\to M$ a regular curve parametrized by an arc-length $s$. The unit tangent and the curvature are defined as usual, i.e.,
\begin{equation}
\mathbf{t}=\alpha'\,\,\mbox{ and }\,\,\kappa =\Vert\nabla_{\mathbf{t}}\,\mathbf{t}\Vert\,,
\end{equation}
respectively. If $\kappa\not=0$, we define the principal normal as
\begin{equation}
\mathbf{n} = \frac{1}{\kappa}\nabla_{\mathbf{t}}\,\mathbf{t}\,.
\end{equation}
The binormal vector $\mathbf{b}$ is such that $\{\mathbf{t},\mathbf{n},\mathbf{b}\}$ is a positively oriented orthonormal frame along $T_{\alpha(s)}M$. The torsion is given by
\begin{equation}
\tau = -\langle\nabla_{\mathbf{t}}\,\mathbf{b},\mathbf{n}\rangle\,.
\end{equation}
The Frenet equations in a Riemannian manifold can be written as
\begin{equation}
\nabla_{\mathbf{t}}\left(
\begin{array}{c}
\mathbf{t}\\
\mathbf{n}\\
\mathbf{b}\\
\end{array}
\right)=\left(
\begin{array}{ccc}
0 & \kappa & 0\\
-\kappa & 0 & \tau\\
 0 & -\tau & 0\\ 
\end{array}
\right)\left(
\begin{array}{c}
\mathbf{t}\\
\mathbf{n}\\
\mathbf{b}\\
\end{array}
\right).\label{eq::RiemFrenetEqs}
\end{equation}

Finally, we say that $\mathbf{x}\in\mathfrak{X}(M)$ is an RM vector field along a regular curve $\alpha:I\to M^{m+1}$ if for some real function $\lambda$ one has $\nabla_{\mathbf{t}}\,\mathbf{x}=\lambda\,\mathbf{t}$, where $\mathfrak{X}(M)$ is the module of vector fields in $M$ \cite{Etayo2016}.

\begin{remark}
For an $m$-dimensional hypersurface $M^m\subseteq\mathbb{E}^{m+k}$, the covariant derivative may be constructed as follows. Given a tangent vector $\mathbf{v}$ along a curve $\alpha$ on $M^{m}$, we take the (usual) derivative of $\mathbf{v}$ along the curve and then projected it on the tangent plane $T_{c(s)}M$. This furnishes the so called \emph{(induced) covariant derivative}:
\begin{equation}
\frac{\mathrm{D}\mathbf{v}}{\mathrm{d}s}\stackrel{\mbox{def}}{=} \mathrm{Proj}_{\,T_cM}\,\frac{\mathrm{d}\mathbf{v}}{\mathrm{d}s}\,.
\end{equation}
The \emph{induced Levi-Civita connection} $\nabla$ on $M$ is such that $\mathrm{D}\mathbf{v}/\mathrm{d}s=\nabla_{\alpha'}\mathbf{v}$. Since we will work with the sphere and hyperbolic space modeled as hypersurfaces in flat space and we will be mostly interested in computing $\nabla_{\alpha'}\mathbf{v}$, this definition does not represent a severe restriction and the reader may keep it in mind as a guide to intuition.
\end{remark}

\section{Preliminaries}

In this work we will be primarily interested in the $(m+1)$-dimensional sphere $\mathbb{S}^{m+1}(r)$ and in the hyperbolic space $\mathbb{H}^{m+1}(r)$. We will use them modeled as submanifolds of $\mathbb{E}^{m+2}$ and $\mathbb{E}_1^{m+2}$, respectively\footnote{Here, $\mathbb{E}_1^{m+2}$ denotes the Lorentz space equipped with the index 1 metric $\langle \mathbf{x},\mathbf{y}\rangle_1=-x_1y_1+\sum_{i=2}^{m+2}x_iy_i$.}:
\begin{equation}
\mathbb{S}^{m+1}(r) = \{q\in\mathbb{R}^{m+2}\,:\,\langle q,q\rangle_e=r^2\}
\end{equation}
and
\begin{equation}
\mathbb{H}^{m+1}(r) = \{q\in\mathbb{R}^{m+2}\,:\,\langle q,q\rangle_1=-r^2,\,x_1>0\},
\end{equation}
equipped with the induced metric, denoted by $\langle\cdot,\cdot\rangle$ (the context will make clear if we are using the Euclidean $\langle\cdot,\cdot\rangle_e$ or Lorentzian $\langle\cdot,\cdot\rangle_1$ metric).

Denoting by $\nabla$ and $\nabla^0$ the Levi-Civita connections on $\mathbb{S}^{m+1}(r)$ (or $\mathbb{H}^{m+1}(r)$) and $\mathbb{E}^{m+2}$ (or $\mathbb{E}^{m+2}_1$), they are related by the Gauss formula according to
\begin{equation}
\nabla^0_{\mathbf{x}}\,\mathbf{y}=\nabla_{\mathbf{x}}\,\mathbf{y}\mp\frac{1}{r^2}\langle\mathbf{x},\mathbf{y}\rangle \,q\,,\label{eq::CovDerAndUsualDer}
\end{equation}
where $q$ denotes the position vector, i.e., the canonical immersion $q:\mathbb{S}^{m+1}(r)\to\mathbb{E}^{m+2}$ for the minus sign and $q:\mathbb{H}^{m+1}(r)\to\mathbb{E}^{m+2}_1$ for the plus sign.
\begin{remark}
The models above do not represent the only choices. Another common way of looking at the spherical geometry is the intrinsic model based on stereographic projection \cite{doCarmo1992}. On the other hand, besides the hyperboloid model above, other common models for the hyperbolic space are the Poincar\'e ball and half-plane models \cite{BenedettiPetronio,ReynoldsMonthly1993}. Anyway, the important fact is that these models are all isometric. Then, intrinsically speaking they are all the same, the choice between them being a matter of convenience.
\end{remark}

The concept of normal curves will play an important role in our work. In Euclidean space we say that $\alpha$ is a \emph{normal curve} if
\begin{equation}
\alpha(s)-p \in \mbox{span}\{\mathbf{t}(s)\}^{\perp},
\end{equation}
where $p$ is a fixed point (the center of the normal curve). We can straightforwardly prove that normal curves in $\mathbb{E}^{m+1}$ are precisely the spherical ones (in this case, $p$ is the center of the respective sphere): $\langle \alpha-p,\mathbf{t}\rangle=0\Leftrightarrow \langle \alpha-p,\alpha-p\rangle=$ constant. This definition makes sense due to the double nature of $\mathbb{E}^{m+1}$ as both a manifold and a tangent space. In order to extend it to a Riemannian manifold $M^{m+1}$, we should replace $\alpha-p$ by a geodesic connecting $p$ to a point $\alpha(s)$ on the curve, as done in \cite{LucasJMAA2015,LucasMJM2016} for the study of rectifying curves:
\begin{definition}
A regular curve $\alpha:I\to M^{m+1}$ is a \emph{normal curve} with center $p$ if the geodesics $\beta_s$ connecting $p$ to $\alpha(s)$ are orthogonal to $\alpha$, i.e., $\langle\mathbf{t}_{\alpha},\mathbf{t}_{\beta}\rangle=0$ along $\alpha$. 
\end{definition}

The equivalence between spherical and normal curves can be extended to a Riemannian manifold by applying the Gauss lemma for the exponential map \cite{doCarmo1992}:

\begin{prop}
On a sufficiently small neighborhood of $p\in M^{m+1}$, a curve $\alpha:I\to M^{m+1}$ is normal (with center $p$) if and only if it lies on a geodesic sphere (with center $p$). \label{prop::EquivNormalAndSphCurves}
\end{prop}

Straight lines in $\mathbb{S}^{m+1}(r)$ or $\mathbb{H}^{m+1}(r)$ can be constructed by intersecting planes in $\mathbb{E}_i^{m+1}$ passing through the origin with $\mathbb{S}^{m+1}(r)$ or $\mathbb{H}^{m+1}(r)$ \cite{ReynoldsMonthly1993}. Then, given $p\in\mathbb{S}^{m+1}(r)$, $\mathbf{v}\in \mathbb{S}^m(1)\subset T_p\mathbb{S}^{m+1}(r)$ or $p\in\mathbb{H}^{m+1}(r)$, $\mathbf{v}\in\mathbb{S}^m(1)\subset T_p\mathbb{H}^{m+1}(r)$, the exponential map is
\begin{equation}
\mbox{exp}_p(u\mathbf{v})=\cos\left(\frac{u}{r}\right)p+r\sin\left(\frac{u}{r}\right)\mathbf{v}\,
\end{equation}
or
\begin{equation}
\mbox{exp}_p(u\mathbf{v})=\cosh\left(\frac{u}{r}\right)p+r\sinh\left(\frac{u}{r}\right)\mathbf{v}\,,
\end{equation}
respectively. Observe that the geodesics $\beta(u)$ above are defined for any value $u\in\mathbb{R}$. So, the equivalence in Proposition \ref{prop::EquivNormalAndSphCurves} is valid globally.

\section{Spherical curves in $\mathbb{S}^{m+1}(r)$ and $\mathbb{H}^{m+1}(r)$}

\begin{figure}[tbp]
\centering
    {\includegraphics[width=0.33\linewidth]{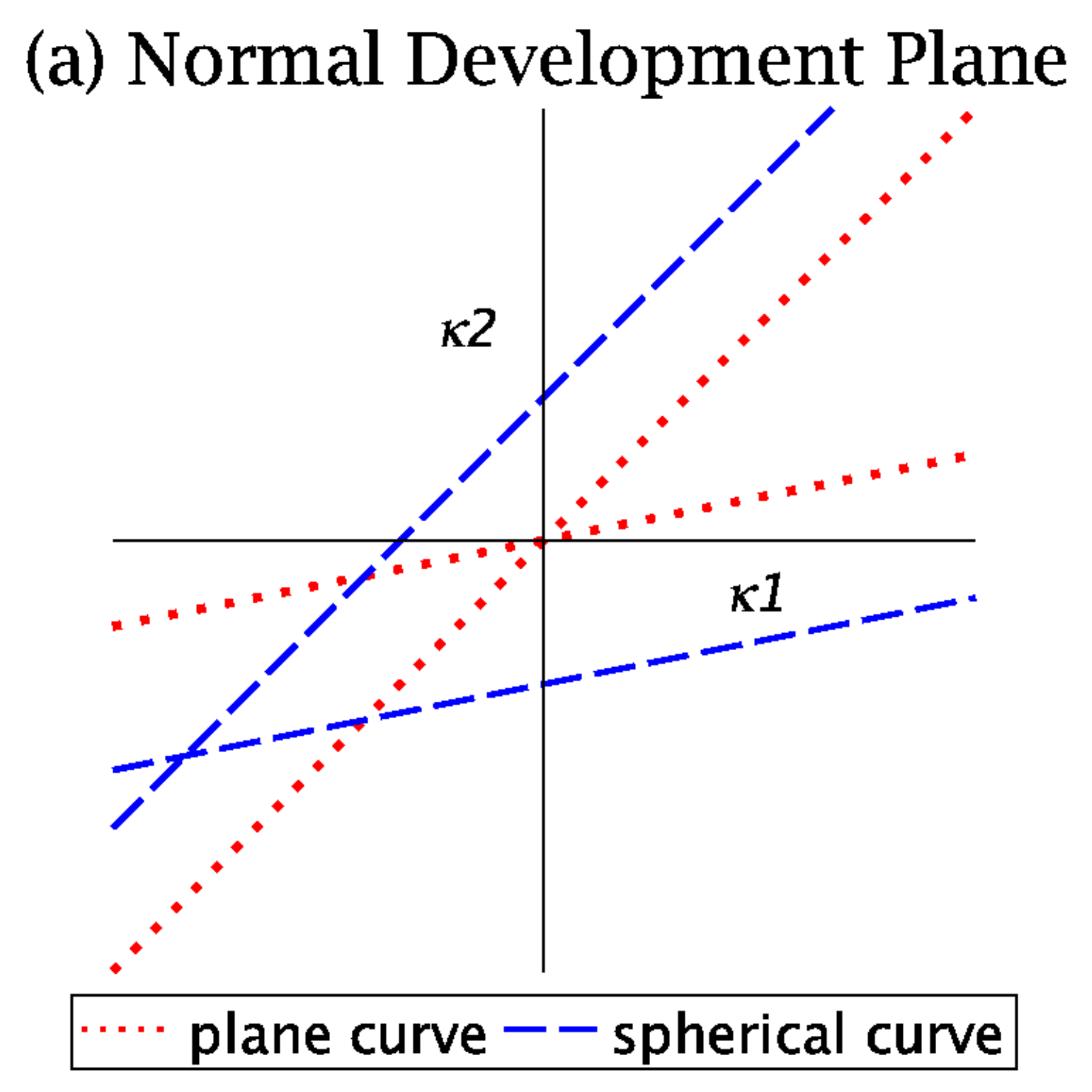}}
    {\includegraphics[width=0.32\linewidth]{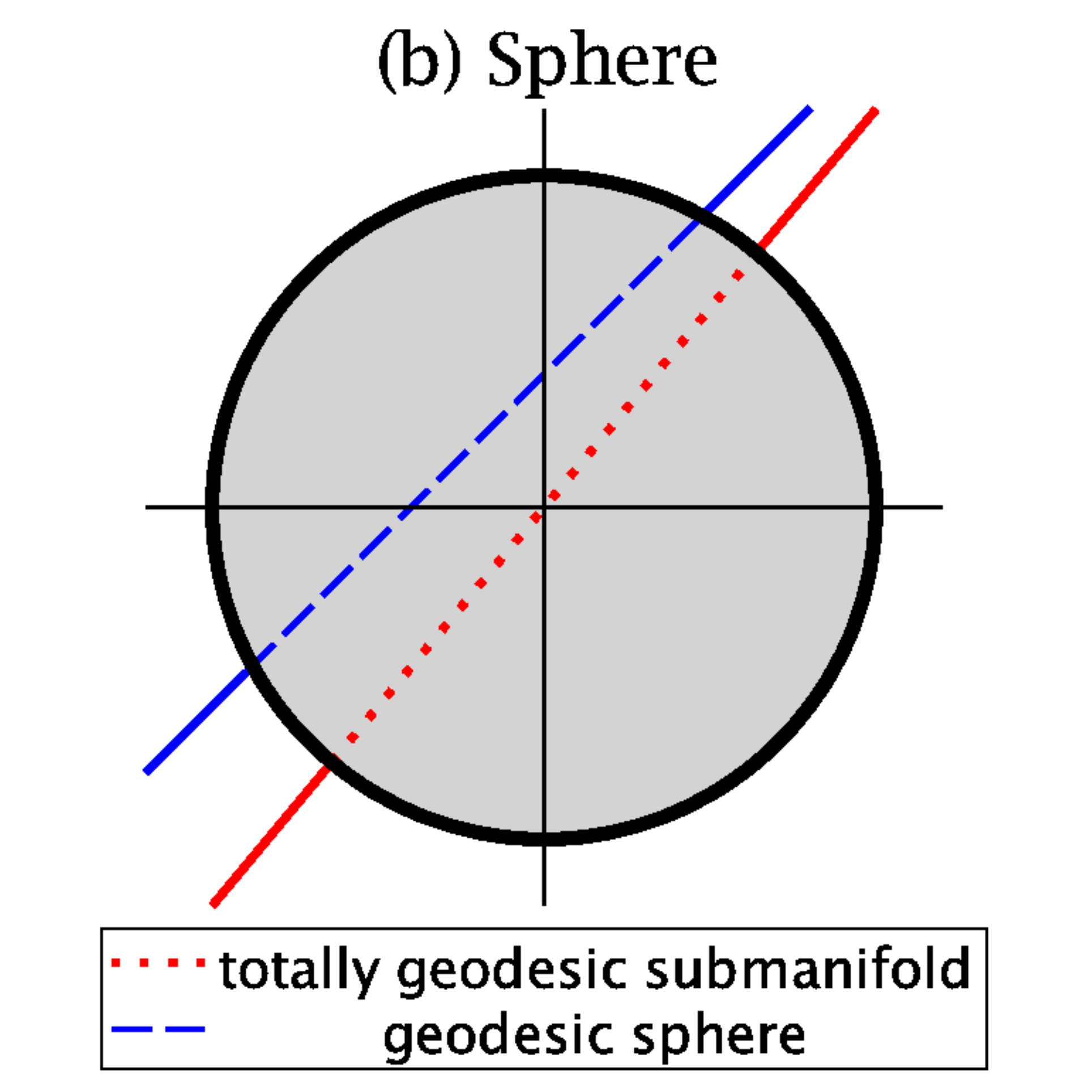}}
    {\includegraphics[width=0.32\linewidth]{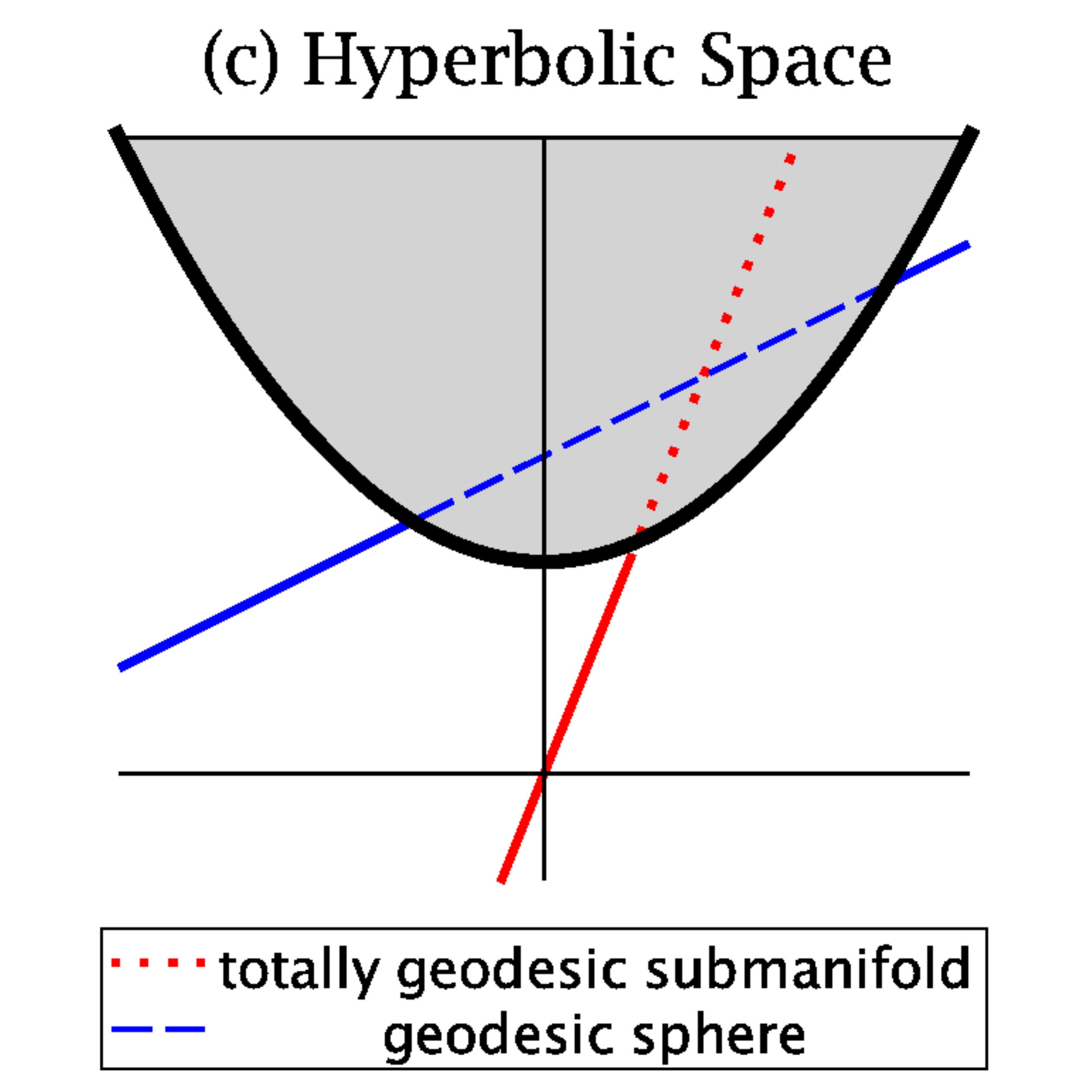}}
          \caption{The geometry of the normal development, geodesic spheres, and totally geodesics submanifolds in $\mathbb{S}^{m+1}(r)$ and $\mathbb{H}^{m+1}(r)$: \textbf{(a)} ($m=2$ in the figure) Lines not passing through the origin (dashed blue line) represent geodesic spherical curves (Theorem \ref{thr::RMCharSphGeodCurves}) and lines through the origin (dotted red line) represent plane curves, i.e., curves on totally geodesic submanifolds (Theorem \ref{thr::NormDevelopS3andH3PlaneCurv}); \textbf{(b)} and \textbf{(c)} Lines passing through the origin (dotted red line) represent hyperplanes passing through the origin and, when intersected with $\mathbb{S}^{m+1}(r)$  or $\mathbb{H}^{m+1}(r)$, give rise to totally geodesic submanifolds, while lines not passing through the origin (dashed blue line) represent hyperplanes not passing through the origin and, when intersected with $\mathbb{S}^{m+1}(r)$ or $\mathbb{H}^{m+1}(r)$, give rise to geodesic spheres (an intersection with hyperplanes forming smaller angles with the hyperboloid axis of the hyperboloid give rise to equidistant surfaces and horospheres) \cite{Spivak1979v4}.}
          \label{fig::DiagramSphPlaneCurv}
\end{figure}

As investigated in previous chapters, by equipping a curve in $\mathbb{E}^{m+1}$ with an RM frame it is possible to characterize spherical curves by means of a linear relation involving the coefficient which dictate the frame motion. We now extend these characterizations for curves on geodesic spheres of $\mathbb{S}^{m+1}(r)$ and $\mathbb{H}^{m+1}(r)$ (see Figure \ref{fig::DiagramSphPlaneCurv}).

\begin{theorem}
Let $\alpha$ be a regular curve in $\mathbb{S}^{m+1}(r)$ or $\mathbb{H}^{m+1}(r)$. Then, $\alpha$ lies on a geodesic sphere if and only if
\begin{equation}
\left\{\begin{array}{lcc}
\displaystyle\sum_{i=1}^m a_i\,\mathbf{\kappa}_i+\displaystyle\frac{1}{r}\,\cot\left(\frac{z_0}{r}\right) &=0\,, &\mbox{ if }\,\alpha\subseteq\mathbb{S}^{m+1}(r)\\[8pt]
\displaystyle\sum_{i=1}^m a_i\,\mathbf{\kappa}_i+\displaystyle\frac{1}{r}\,\coth\left(\frac{z_0}{r}\right) &=0\,, &\mbox{ if }\,\alpha\subseteq\mathbb{H}^{m+1}(r)\\
\end{array}
\right.,
\end{equation}
for some constants $z_0< \pi r/2$ (the radius of the geodesic sphere\footnote{The restriction $z_0<\pi r/2$ guarantees that the center of the geodesic sphere is well defined: if $z_0=\pi r/2$, both $p$ and its antipodal $-p$ are equidistant from the geodesic sphere.}) and $a_i$, $1\leq i\leq m$. \label{thr::RMCharSphGeodCurves}
\end{theorem}
\textit{Proof. } We will do the proof for $\alpha\subseteq\mathbb{S}^{m+1}(r)$ only, the case for $\mathbb{H}^{m+1}(r)$ being analogous (one just needs to use the hyperbolic versions of the trigonometric functions). 

If $\alpha:I\to\mathbb{S}^{m+1}(r)$ is a normal curve parametrized by arc-length $s$, let us write
\begin{equation}
\alpha(s) =  \exp_p(z_0\,\mathbf{v}(c_0\,s))\,,
\end{equation}
where $z_0$ and $c_0=[r\sin(z_0/r)]^{-1}$ are constants\footnote{The constant $c_0$ is here in order to have $s$ as an arc-length parameter, as can be easily checked.} and $\mathbf{v}:I\to\mathbb{S}^{m}(1)\subseteq T_p\mathbb{S}^{m+1}(r)$ is a unit speed curve such that $\langle p,\mathbf{v}\rangle_e=0$. If $\{\mathbf{t}_{\alpha},\mathbf{n}_1,...,\mathbf{n}_m\}$ is an RM frame along $\alpha$, the unit tangent of $\alpha$ can be written as
\begin{equation}
\mathbf{t}_{\alpha}(s) = \mathbf{v}'(c_0s).\label{eqTalphaInTermsOfV}
\end{equation}
On the other hand, the unit speed geodesic $\beta_s$ connecting $p$ to a point $\alpha(s)$ is
\begin{equation}
\beta_s(u) = \cos\left(\frac{u}{r}\right)\,p+r\sin\left(\frac{u}{r}\right)\,\mathbf{v}(c_0\,s)\,\Rightarrow \beta_s(z_0)=\alpha(s)\,.
\end{equation}
The normality condition $\langle \mathbf{t}_{\alpha},\mathbf{t}_{\beta}\rangle=0$ implies
\begin{equation}
``\mathbf{t}_{\beta_s}\mbox{ at }\alpha(s)" = \mathbf{t}_{\beta_s}(z_0)=\displaystyle\sum_{i=1}^ma_i(s)\mathbf{n}_i(s)\,.
\end{equation}

The derivative of the coefficients $a_i=\langle \mathbf{t}_{\beta_s},\mathbf{n}_i\rangle$ gives
\begin{equation}
a'_i = \langle\nabla_{\mathbf{t}_{\alpha}}\mathbf{t}_{\beta_s},\mathbf{n}_i\rangle+\langle \mathbf{t}_{\beta_s},\nabla_{\mathbf{t}_{\alpha}}\mathbf{n}_i\rangle= \langle \nabla^0_{\mathbf{t}_{\alpha}}\mathbf{t}_{\beta_s},\mathbf{n}_i\rangle\,\label{eq::DerivativeOfai}
\end{equation}
where the last equality is a consequence of the fact that $\mathbf{n}_i$ is RM and also that $\nabla_{\mathbf{x}}\,\mathbf{y}=\nabla_{\mathbf{x}}^0\,\mathbf{y}$ for two orthogonal vectors $\mathbf{x},\,\mathbf{y}$ in Eq. (\ref{eq::CovDerAndUsualDer}). Now, using that $\mathbf{t}_{\beta}$ along $\alpha$ can be also written as
\begin{equation}
\mathbf{t}_{\beta_s}(z_0) = -\frac{1}{r}\sin\left(\frac{z_0}{r}\right)\,p+\cos\left(\frac{z_0}{r}\right)\,\mathbf{v}(c_0\,s),\label{eq::TbetaAlongAlpha}
\end{equation}
we have
\begin{equation}
\nabla_{\mathbf{t}_{\alpha}}^0\mathbf{t}_{\beta} = \frac{1}{r}\frac{\cos(z_0/r)}{\sin(z_0/r)}\,\mathbf{v}'(c_0\,s)=\frac{\cot(z_0/r)}{r}\,\mathbf{t}_{\alpha}\,.\label{eq::NablaTbAlongTaIsTa}
\end{equation}
Inserting the expression above in Eq. (\ref{eq::DerivativeOfai}) shows that $a_i'=0$ and, therefore , the coefficients $a_i$, $1\leq i\leq m$, are all constants. Finally, taking the derivative of $\langle \mathbf{t}_{\beta},\mathbf{t}_{\alpha}\rangle=0$ along $\alpha$ gives
\begin{eqnarray}
0 & = & \langle\nabla_{\mathbf{t}_{\alpha}}\mathbf{t}_{\beta},\mathbf{t}_{\alpha}\rangle+\langle \mathbf{t}_{\beta},\nabla_{\mathbf{t}_{\alpha}}\mathbf{t}_{\alpha}\rangle\nonumber\\
& = & \left\langle\frac{\cot(z_0/r)}{r}\,\mathbf{t}_{\alpha},\mathbf{t}_{\alpha}\right\rangle+\left\langle \sum_{i=1}^ma_i\mathbf{n}_i,\sum_{j=1}^m\kappa_j\mathbf{n}_i\right\rangle\nonumber\\
& = & \frac{1}{r}\cot\left(\frac{z_0}{r}\right)+\sum_{i=1}^ma_i\kappa_i.
\end{eqnarray}

Conversely, suppose that $\alpha$ is a regular curve satisfying $\sum_{i}a_i\kappa_i+\cot(z_0r^{-1})/r=0$. The proof is based on the following observation: for a spherical curve, if we invert the direction of the motion of $\beta_s$ we have a geodesic connecting $\alpha(s)$ to $p$, whose initial velocity vector according to Eq. (\ref{eq::TbetaAlongAlpha}) should be $-\mathbf{t}_{\beta}$. Now, let us define
\begin{equation}
\mathbf{w}(s)=-\sum_{i=1}^ma_i\,\mathbf{n}_i(s)\,\mbox{ and }\, 
P(s) = \cos\left(\frac{z_0}{r}\right)\alpha(s) - r\sin\left(\frac{z_0}{r}\right)\mathbf{w}(s).
\end{equation}
Taking the derivative of the last equation, we find $P'(s)=0$ and then $P$ is a constant. Consequently, it means that the geodesics with initial conditions $\alpha(s)$ and $\mathbf{w}(s)$ travel always the same distance in order to arrive at $P$, i.e., $\alpha$ is a spherical curve. \qed

Finding RM frames along a curve may be a difficult problem and in general one must resort to some kind of numerical method, see e.g. \cite{WangACMTOG}. However, for a curve $\alpha$ in $\mathbb{S}^2(r,p)\subseteq\mathbb{R}^3$, computing RM frames is not difficult: $\mathbf{u}=(\alpha(s)-p)/r$ is RM \cite{WangACMTOG}. This result can be extended for other ambient spaces by taking into account Eq. (\ref{eq::NablaTbAlongTaIsTa}) in the proof above. So, we have

\begin{corollary}
For a curve $\alpha(s)$ on a geodesic sphere of $\mathbb{S}^{m+1}(r)$, or $\mathbb{H}^{m+1}(r)$, the tangents of the geodesics connecting the center of the geodesic sphere to a point on the curve is a rotation minimizing vector field. 
\end{corollary}

The previous theorem was obtained by expressing $\mathbf{t}_{\beta}$ in terms of an RM basis for the normal plane $\mbox{span}\{\mathbf{t}_{\alpha}\}^{\perp}$. If we use the Frenet frame instead, then we can extend a classical characterization result for spherical curves in $\mathbb{R}^3$.

\begin{theorem}
Let $\alpha$ be a regular curve with non-zero torsion in $\mathbb{S}^3(r)$ or $\mathbb{H}^3(r)$. If $\alpha$ lies on a geodesic sphere, then
\begin{equation}
\frac{\mathrm{d}}{\mathrm{d}s}\left[\frac{1}{\tau}\frac{\mathrm{d}}{\mathrm{d}s}\left(\frac{1}{\kappa}\right)\right]+\frac{\tau}{\kappa}=0\,.\label{eq::FrenetCharSphCurves}
\end{equation}
\label{thr::FrenetChar3DSphCurv}
\end{theorem}
\textit{Proof. } We will do the proof for $\alpha\subseteq\mathbb{S}^3(r)$ only, the case for $\mathbb{H}^3$ being analogous. 

Let $\alpha$ be a spherical curve and $\{\mathbf{t}_{\alpha},\mathbf{n},\mathbf{b}\}$ its Frenet frame, then there exists a point $p$ such that the geodesic $\beta_s$ connecting $p$ to $\alpha(s)$ satisfies $\langle\mathbf{t}_{\beta},\mathbf{t}_{\alpha}\rangle=0$. Let us write
\begin{equation}
\mathbf{t}_{\beta} = c_1\mathbf{n}+c_2\mathbf{b}\,.
\end{equation}

Taking the derivative gives
\begin{eqnarray}
\nabla_{\mathbf{t}_{\alpha}}\mathbf{t}_{\beta} & = & c'_1\mathbf{n}+c'_2\mathbf{b}+c_1\nabla_{\mathbf{t}_{\alpha}}\mathbf{n}+c_2\nabla_{\mathbf{t}_{\alpha}}\mathbf{b}\nonumber\\
\frac{1}{r}\cot\left(\frac{z_0}{r}\right)\mathbf{t}_{\alpha} & = & -c_1\kappa\mathbf{t}_{\alpha}+(c'_1-\tau c_2)\mathbf{n}+(c'_2+\tau c_1)\mathbf{b}\,,
\end{eqnarray}
where we used Eq. (\ref{eq::NablaTbAlongTaIsTa}) to arrive at the second equality above. Now, comparing the coefficients of the last equation leads to 
\begin{equation}
\left\{
\begin{array}{ccc}
-\kappa c_1 & = & \frac{1}{r}\cot\left(\frac{z_0}{r}\right)\\
c'_1 - \tau c_2 & = & 0\\
c'_2 + \tau c_1 & = & 0\\
\end{array}
\right..\label{eq::auxFrenetCharSphCurves}
\end{equation}
 
 From the 1st and 2nd equations we find
\begin{equation}
c_1=-\frac{1}{r\kappa}\cot\left(\frac{z_0}{r}\right)\Rightarrow \tau c_2 = c'_1 = -\frac{\mathrm{d}}{\mathrm{d}s}\left[\frac{1}{r\kappa}\cot\left(\frac{z_0}{r}\right)\right].
\end{equation}
Now, using the expression above in combination with the 3rd equation of (\ref{eq::auxFrenetCharSphCurves}), furnishes
\begin{equation}
-\tau c_1 = c_2'=-\frac{\mathrm{d}}{\mathrm{d}s}\left\{\frac{1}{\tau}\frac{\mathrm{d}}{\mathrm{d}s}\left[\frac{1}{r\kappa}\cot\left(\frac{z_0}{r}\right)\right]\right\}\,.
\end{equation}
Finally, the desired result follows from the finding above and the 1st equation of (\ref{eq::auxFrenetCharSphCurves}).

Conversely, let $\alpha$ be a regular curve satisfying Eq. (\ref{eq::FrenetCharSphCurves}). As in the proof for the characterization of spherical curve via RM frames, the idea is to find a (fixed) point $P$ and a vector field $\mathbf{w}$ such that all the geodesics emanating from $\alpha$ with initial velocity $W$ reach $P$ after traveling the same distance. Define the following vector field along $\alpha(s)$
\begin{equation}
\mathbf{w}(s)= -\frac{1}{r\kappa(s)}\cot\left(\frac{z_0}{r}\right)\mathbf{n}(s)-\frac{1}{\tau(s)}\frac{\textrm{d}}{\textrm{d}s}\left[\frac{1}{r\kappa(s)}\cot\left(\frac{z_0}{r}\right)\right]\,\mathbf{b}(s),
\end{equation}
which satisfies $\nabla_{\mathbf{t}_{\alpha}}\mathbf{w}=r^{-1}\cot(z_0r^{-1})\,\mathbf{t}_{\alpha}\,$. Now define 
\begin{equation}
P(s) = \cos\left(\frac{z_0}{r}\right)\alpha(s)-r\sin\left(\frac{z_0}{r}\right)\mathbf{w}(s)\,.
\end{equation}
Taking the derivative of $P$ shows that $P'(s)=0$ and, therefore, $P$ is constant and will be the center of the geodesic sphere that contains $\alpha$.
\qed

\begin{remark}
One can also equip a curve with a Frenet frame in higher dimensional Riemannian manifolds \cite{Spivak1979v4}, p. 29, and use them to characterize (geodesic) spherical curves. One can follow the same steps as in the previous theorem, i.e., use that a spherical curve must be normal and then investigate the coefficients $c_i$ of $\mathbf{t}_{\beta}$ in terms of the Frenet frame. The expressions however are quite cumbersome and will not attempt to write it here. We just remark that, as in happens in 3d, the values of $r$ and of the geodesic sphere radius do not appear in the expression characterizing spherical curves. Note in addition that the curve must be of class $C^{m+2}$, in contrast with the $C^2$ requirement in theorem \ref{thr::RMCharSphGeodCurves} via RM frames.
\end{remark}

\section{Plane curves: totally geodesic submanifolds}

The so called totally geodesic submanifolds in a Riemannian ambient space have the simplest shape and play the role of (hyper)planes\footnote{It is worth mentioning that despite their simplicity, in general, Riemannian manifolds do not have non trivial totally geodesic submanifolds \cite{MurphyArXiv2017,Tsukada1996}. The existence of such submanifolds imposes severe restrictions on the geometry of the ambient manifold, see e.g., \cite{NikolauevskyIJM2015}.}. 

\begin{definition}
A submanifold $N$ of a Riemannian manifold $M$ is a \emph{totally geodesic submanifold} if any geodesic on the submanifold $N$ with the induced Riemannian metric is also a geodesic on $M$. One dimensional totally geodesic submanifolds are geodesics.
\end{definition}

There are many equivalent ways of characterizing a totally geodesic submanifold. Indeed, all the conditions below are equivalent \cite{Cartan1946}, p. 114,
\begin{enumerate}
\item $N\subset M$ is totally geodesic;

\item the principal curvatures vanish in every point of $N$;
\item the normal field to $N$ remains normal if parallel transported along any curve on $N$;
\item any tangent field to $N$ remains tangent if parallel transported along any curve on $N$.
\end{enumerate}
Note that property 3 essentially says that the normal field of a totally geodesic submanifold is constant, which is a crucial feature of Euclidean planes: $\pi$ is a plane if and only if there exist $\mathbf{u}_0$ and $x_0$ constants such that $\pi=\{x:\langle x-x_0,\mathbf{u}_0\rangle=0\}$. Thus, we say that \emph{plane curves} in a Riemannian manifold are those curves on totally geodesic submanifolds.

In Euclidean space it is known that normal development curves $(\kappa_1,...\,,\kappa_m)$ which are lines passing through
the origin characterize plane curves. Here we (partially) extend this result to totally geodesic curves on any Riemannian manifold.
\begin{theorem}
Let $\alpha:I\to N^n\subset M^{m+1}$ be a regular curve and $\{\mathbf{t},\mathbf{n}_1,...\,,\mathbf{n}_m\}$ a rotation minimizing frame along it. If $\alpha$ lies on a totally geodesic submanifold $N$, then its normal development curve 
$(\kappa_1,...\,,\kappa_m)$ lies on a line passing through the origin.
\label{thr::NormDevelopRiemPlaneCurv}
\end{theorem}
\textit{Proof. } Let $\mathbf{u}$ be a normal vector field on $N$. Since $N$ is totally geodesic, we can use that $\nabla_{\mathbf{t}}\,\mathbf{u}=0$. In addition, we can also write $\mathbf{u}=\sum_{i=1}^ma_i\mathbf{n}_i$ for the normal $\mathbf{u}$ along $\alpha$. The coefficient $a_i=\langle\mathbf{u},\mathbf{n}_i\rangle$ satisfies
\begin{equation}
a'_i = \langle\nabla_{\mathbf{t}}\,\mathbf{u},\mathbf{n}_i\rangle+\langle\mathbf{u},\nabla_{\mathbf{t}}\,\mathbf{n}_i\rangle=0\,.
\end{equation}
Then, for all $i\in\{1,...\,,m\}$, $a_i$ is a constant. Finally,
\begin{equation}
0=\nabla_{\mathbf{t}}\,\mathbf{u}=\sum_{i=1}^ma_i\nabla_{\mathbf{t}}\,\mathbf{n}_i=\sum_{i=1}^m(-a_i\kappa_i\,\mathbf{t})
\end{equation}
and, therefore, $\sum a_i\kappa_i=0$ represents the equation of a line passing through the origin.
\qed

Let us now discuss the reciprocal of the theorem above. Given a curve $\alpha:I\to M$ satisfying $\sum_{i=1}^ma_i\,\kappa_i=0$ for some constants $a_1,...\,,a_m$, we may define $\mathbf{u}(s)=\sum a_i\,\mathbf{n}_i$. Then, it follows that
\begin{equation}
\nabla_{\mathbf{t}}\,\mathbf{u}=\sum -a_i\kappa_i\,\mathbf{t}=0\,.
\end{equation}
Thus, $\mathbf{u}$ is parallel transported along $\alpha$. The problem now is to find a totally geodesic submanifold containing $\alpha$ and whose normal field when restricted to $\alpha$ is equal to $\mathbf{u}$. A candidate to solution is the submanifold given by the following parametrization
\begin{equation}
X(s_1,...\,,s_m) = \exp_{\alpha(s_1)}\left(\sum_{i=2}^{m}s_i\,\mathbf{u}_i\right),
\end{equation}
where $\{\mathbf{u}_i(s)\}_{i=2}^{m}$ is an orthonormal basis for  $\mbox{span}\{\mathbf{t}(s),\mathbf{u}(s)\}^{\perp}$ for all $s=s_1$. 

Observe however, the fact that $X$ is geodesic along $\alpha$ does not implies that it will also be geodesic in all its points. In fact, the existence of totally geodesic submanifolds is an exceptional fact. On the other hand, in both $\mathbb{S}^{m+1}(r)$ and $\mathbb{H}^{m+1}(r)$ the situation is easier, since that totally geodesic submanifolds are precisely the intersection of linear subspaces of $\mathbb{R}^{m+2}$ with $\mathbb{S}^{m+1}(r)$ and $\mathbb{H}^{m+1}(r)$ \cite{Spivak1979v4} (see Figure \ref{fig::DiagramSphPlaneCurv}). Then, we have
\begin{theorem}
Let $\alpha$ be a $C^2$ regular curve in $\mathbb{S}^{m+1}(r)$, or $\mathbb{H}^{m+1}(r)$, and $\{\mathbf{t},\mathbf{n}_1,...\,,\mathbf{n}_m\}$ an RM frame along it. Then, $\alpha$ is a plane curve, i.e., it lies on a totally geodesic submanifold, if and only if the normal development $(\kappa_1,...\,,\kappa_m)$ is a line passing through the origin.
\label{thr::NormDevelopS3andH3PlaneCurv}
\end{theorem}
\textit{Proof.} The direction ``plane curve $\Rightarrow$ $\sum_{i=1}^ma_i\kappa_i=0$ ($a_i$ constant)'' is a consequence of the previous theorem. For the reciprocal, define a vector field along $\alpha$ as $\mathbf{u}(s)=\sum_{i=1}^ma_i\mathbf{n}_i(s)$. Using that the normal development is a line passing through the origin, we have 
\begin{equation}
\frac{\mathrm{d}\mathbf{u}}{\mathrm{d}s}=\nabla_{\mathbf{t}}\,\mathbf{u} = \sum_{i=1}^m-a_i\kappa_i\,\mathbf{t} = 0,
\end{equation}
where for the first equality we used that $\langle\mathbf{t},\mathbf{u}\rangle=0$ in Eq. (\ref{eq::CovDerAndUsualDer}). Therefore, $\mathbf{u}$ is a constant vector in $\mathbb{R}^{m+2}$ and it follows that $\alpha$ is contained in the plane, in $\mathbb{R}^{m+2}$, given by $\{x\in\mathbb{R}^{m+2}:\langle x,\mathbf{u}\rangle=0\}$. In fact, 
\begin{equation}
\langle\alpha,\mathbf{u}\rangle'=\langle\mathbf{t},\mathbf{u}\rangle=0\Rightarrow \langle\alpha,\mathbf{u}\rangle=c\,\,\mbox{ constant}\,.
\end{equation}
The constant $c$ must be zero. Otherwise, $\alpha$ would be contained on an intersection of $\mathbb{S}^{m+1}(r)$, or $\mathbb{H}^{m+1}(r)$, with a plane not passing through the origin, which is a geodesic sphere \cite{Spivak1979v4}. Since the normal development of a spherical curve does not pass through the origin, we conclude that $c=0$.
\qed
\part{APPLICATIONS IN PHYSICS: ON THE QUANTUM MECHANICS OF A CONSTRAINED PARTICLE}

\chapter{QUANTUM CONSTRAINED SYSTEMS}
\label{chap::Vgip}

The study of new material properties due to its shape has gained importance since the birth of nanoscience. The experimental techniques have evolved to a stage where various examples of nanostructures whose shape resembles planes, spheres, cylinders, and other non-trivial geometries, have been synthesized \cite{TerronesNewJPhys,CastroNetoRepProgPhys}. However, it is still difficult to establish a relation between the geometry and the quantum behavior of such systems. In face of these developments, writing the appropriate equations for a particle confined somewhere is essential to a proper understanding and modeling of these phenomena. In the 1950s De Witt addressed the problem of describing a confinement in a curved space through a quantization procedure, which resulted however in an ordering ambiguity \cite{DeWittRMP1957}. Later on, an approach which does not suffer from this ambiguity was devised by Jensen and Koppe \cite{JensenKoppeAnnPhys} in the 1970s and by Da Costa \cite{DaCostaPRA1981,DaCostaPRA1982} in the 1980s, showing that a geometry-induced potential (GIP) acts upon the dynamics\footnote{For more rigorous studies, see e.g. \cite{FroeseCMP2001,MitchellPRA2001,SchusterJaffeAnnPhys2003}; for studies allowing the confinement to vary along the constraint region, see e.g. \cite{WachsmuthPRA2010,SchmelcherPRA2014}, in which case analogies with the Born-Oppenheimer approximation is an important tool \cite{JeckoJMP2014}; and for studies taking into account that, from a realistic viewpoint, a particle can not move on a truly 2D region, which gives rise to corrections to the GIP, see e.g. \cite{IkegamiPTP1991,IkegamiSurfScience,GravesenJMathPhys}.}. Since then, some research on the subject has been reported, such as a path integral formulation \cite{MatsutaniJPhysSocJapan,MatsutaniPRa}, effects on the eigenstates of nanostructures \cite{EncinosaPRA,GravesenPRA}, interaction with an electromagnetic potential \cite{IkegamiSurfScience,FerrariPRL,deOliveiraJMathPhys,BrasileirosAnnPhys}, modeling of bound states on conical surfaces \cite{MoraesAnnPhys,BrasileirosJMP,ChinesesPhysicaE}, spin-orbit interaction \cite{EntinPRB,OrtixSpin,OrtixPRB2015}, electronic transport on surfaces \cite{MarchiEtAlPRB2005,FernandosEfranceses,WangJPD2016}, and bent waveguides \cite{KrejcirikJGP2003,delCampoSciRep2014,SchmelcherPRA2014,HaagAHP2015}, just to name a few.

For surfaces in $\mathbb{R}^3$, Encinosa and Etemadi found that the shift in the ground-state energy may be of sufficient order to be observable in quantum nanostructures \cite{EncinosaPRA}. In addition, taking into account the effects of the GIP may lead to qualitative changes in the transmission characteristics of a semiconductor two dimensional electron gas forming a Y junction when compared to the case with no GIP \cite{CuoghiPRB2009}. More recently, on the experimental side, Onoe \textit{et al.} reported on the observation of Riemannian geometric effects through the GIP on the Tomonaga-Luttinger liquid exponent in a 1D metallic $C_{60}$ polymer with an uneven periodic peanut-shaped structure \cite{OnoePRB,OnoeEPL}. In addition, Szameit \textit{et al.} described the experimental realization of an optical analogue of the GIP \cite{SzameitEtAlPRL}: for more on the interplay between geometry and optics, see e.g. \cite{SchultheissPRL2010,PedersenJMO2015}.

\section{Geometry-induced potential for constrained systems}

Let a mass $m^*$ in a space $M$ be confined to some $d$-dimensional region $N^d\subseteq M^{d+k}$ (the usual case being $M^{d+k}=\mathbb{R}^{d+k}$). Then, what are the ``correct'' equations that govern the (constrained) dynamics on $N^d$? A first approach would be to use the intrinsic coordinates of $N^d$ and write the equations according to them\footnote{For example, the dynamics governed by a differential operator $L_M$ in $M$, such as the Laplacian $-\Delta_M$, is then described by the respective operator $L_N$ written on the $N^d$-coordinates.}. According to such an intrinsic scheme, the ambient space $M^{d+k}$ plays no relevant role at all. On the other hand, a different and more realistic approach would be to appeal to an explicit confining mechanism. In other words, one imposes that some kind of confining potential is responsible for the constraining, e.g., a strong force acting in the normal direction to $N$. Here the ambient space $M^{d+k}$ may play some role, since the confining potential ``sees'' the directions normal to  $N^d$, and then the constrained equations may depend on the way $N^d$ is embedded on $M^{d+k}$. Finally, one can also imagine a third different approach. Namely, one writes the equations in $M^{d+k}$ according to some coordinate system adapted to $N^d$, i.e., coordinates $(u^1,...,u^{d+k})$ such that $N^d=\{u\in M\,:\,u^{d+1}=u^{d+1}_0,...,u^{d+k}=u^{d+k}_0\}$ for some constants $u^{d+i}_0$, $i=1,...,k$, and then one takes the constrained dynamics on $N^d$ as the dynamics in $M$ after the last $k$ coordinates being fixed\footnote{We mention that, by the definition of a submanifold, it is always possible to find an adapted coordinate system in a certain
  neighborhood of a point of $N^d$; naturally, $\protect \mathaccentV
  {bar}016{u}=(u^1,...,u^d)\DOTSB \mapstochar \rightarrow
  (u^1,...,u^d,u^{d+1}_0,...,u_0^{d+k})\in M^{d+k}$,  for some constants $u^{d+i}_0$, $i=1,...,k$, is a (local)
  parametrization of $N^d$ into $M^{d+k}$.}: e.g., spheres in spherical coordinates. Generally, this approach is not equivalent to a confining potential one \cite{JensenKoppeAnnPhys,BernardAndLewYanVoonEurJPhys}. Indeed, since the equation $L_M(u)=0$, which describes the dynamics of the particle in $M^{d+k}$ according to a differential operator $L_M$, may involve derivatives with respect to $u^{d+1},...,u^{d+k}$, in general it does not follow that the solutions of $L_M(u\,;\{u^{d+i}=u^{d+i}_0\})$ are equivalent to the solutions of the respective operator $L_N(\bar{u})$ on $N$ written according to the adapted coordinate system.

In the classical mechanics picture, the approaches described above are shown to be equivalent, the choice between them being a matter of convenience. However, on the quantum mechanical counterpart, the dynamics must obey the uncertainty relations and, since any kind of confinement involves the fully knowledge of some degrees of freedom, it is not clear if  different approaches would lead to equivalent results for the constrained dynamics. In the 1970s Jensen and Koppe \cite{JensenKoppeAnnPhys} showed how the many available approaches discussed above would lead to non-equivalent results through the illustrative example of a circle of radius $R$. More recently, Bernard and Lew Yan Voon also discussed the non-equivalence for the case of spheroidal surfaces in $\mathbb{R}^3$ \cite{BernardAndLewYanVoonEurJPhys}, while Filgueiras \textit{et al.} discussed the difference between intrinsic and confining potential approaches for conical surfaces \cite{BrasileirosJMP}. Finally, we also mention that, by approaching the problem through a quantization procedure in the intrinsic coordinates of $N^d$, the resulting equations suffer from an ordering ambiguity \cite{DeWittRMP1957}. On the other hand, a confining potential approach does not suffer from such a problem: the confining potential approach gives a unique effective Hamiltonian to the confined dynamics \cite{DaCostaPRA1981}.

In order to find the equations for the constrained dynamics in a surface $\Sigma\subset\mathbb{R}^3$, Jensen and Koppe \cite{JensenKoppeAnnPhys} devised an approach which consists in describing the confinement by starting from the dynamics in the region between two neighboring parallel surfaces and imposing homogeneous boundary conditions along them. So, taking the limit as the distance between the neighboring surfaces goes to zero, one obtains the equations that govern the constrained dynamics. Some years later, Da Costa devised an approach which consists in applying an explicit strong confining potential to restrict the motion of the particle to the desired surface (or curve) \cite{DaCostaPRA1981}. As expected, both formalisms coincide \cite{JensenKoppeAnnPhys,DaCostaPRA1981}; for surfaces one finds \cite{JensenKoppeAnnPhys,DaCostaPRA1981}    
\begin{equation}
\textrm{i}\hbar\frac{\partial \psi}{\partial t} = -\frac{\hbar^2}{2m^*}\left[\Delta_{\Sigma}+(H^2-K)\right]\psi;\label{eq:SchrodingerEqExtrinsic}
\end{equation}
while for curves one has \cite{DaCostaPRA1981}
\begin{equation}
\textrm{i}\hbar\frac{\partial \psi}{\partial t} = -\frac{\hbar^2}{2m^*}\left(\Delta_{\alpha}+\frac{\kappa^2}{4}\right)\psi,\label{eq:SchrodingerEqExtrinsicCurves}
\end{equation}
where $\Delta_{\alpha}=\textrm{ d}^2/\textrm{d}s^2$ is the Laplacian on the curve in terms of its arc-length parameter. The above equations show that in general the intrinsic and confining potential approaches do not lead to the same constrained dynamics. In the former, the dynamics is governed by the Laplacian operator only, while in the latter the Laplacian is coupled to a scalar geometry-induced potential. So, in order to do a more realistic study, where the global geometry should be taken into account, an extrinsic scheme would be more appropriate. Additionally, the equations will be exactly the same only for (regions) of the plane or spheres, since these are the only surfaces where $H^2-K\equiv0$, while the equality for curves occurs uniquely for line segments, since it is demanded $\kappa^2\equiv0$.

Finally, it is worth mentioning that these results for the constrained dynamics are based on the assumption that the confining potential $V_{\mathrm{conf}}$ is uniform, i.e., its equipotential level sets only depend on the distance from the constraint region $N\subseteq M$: $V_{\mathrm{conf}}(q)=V_{\mathrm{conf}}(\mbox{dist}(q,N))$. The confinement is then put forward through a limiting procedure, i.e., one considers a sequence of potentials $\{V_{\lambda}\}_{\lambda}$ that approximates the confining one $V_{\mathrm{conf}}$ for $\lambda\to\infty$ \cite{DaCostaPRA1981}: 
\begin{equation}
V_{\mathrm{conf}}(q)=\lim_{\lambda\to\infty} V_{\lambda}(q) = \left\{
\begin{array}{ccc}
0 & , & q\in N\\
\infty & , & q\not\in N\\
\end{array}
\right.,
\end{equation}
which allows for the decoupling between the tangential and normal degrees of freedom in the limit $\lambda\to\infty$. So, one separates the Hamiltonian into a term that governs the low energy motion in the tangent direction, which is the effective Hamiltonian along the constraint region, and a high energy motion in the normal direction. However, in some context this hypothesis is no longer realistic and one can not suppose that the equipotentials are equidistant. As a consequence, the tangential and normal degrees of freedom are coupled \cite{WachsmuthPRA2010,SchmelcherPRA2014}. In addition, let us mention that from a realistic viewpoint, a particle can not move on a truly two dimensional region and then the constrained dynamics is meant to be a mathematical construct. In this case, one must consider that a particle may move on a tubular neighborhood around the constraint region, which gives rise to corrections to the geometry-induced potential \cite{IkegamiPTP1991,IkegamiSurfScience,GravesenJMathPhys}. In what follows we will not consider such possibilities, i.e., we shall admit that the confinement is uniform and realizable.

Let us finish this section by making some remarks concerning the role played by the torsion for curves. Interestingly, the torsion of a curve does not appear in the GIP \cite{DaCostaPRA1981}. Nonetheless, Takagi and Tanzawa put forward an investigation for a particle confined to a thin tube, which is twisted and curved to form a closed loop \cite{TwistedRingProgTheorPhys}, and described the effect of both curvature and torsion of the loop up to second order. They then observed that the torsion may give rise to a geometry-induced Aharonov-Bohm effect. On the other hand, in the study of the Schr\"odinger-Pauli equation for a spin-orbit coupled electron constrained to a space curve \cite{OrtixPRB2015}, it was found that the torsion of the curve generates an additional quantum geometry-induced potential, adding to the known curvature-induced one. In short, besides making the integration of the Frenet equations more difficult, these studies suggest that by considering other effects, in addition to the constraining for the Schr\"odinger equation, the torsion naturally appears in the discussion. Moreover, by noticing that the torsion has to do with the derivative of the binormal vector $\mathbf{b}$, which can be expressed as $\mathbf{b}=(\alpha'\times\alpha'')/\Vert\alpha'\times\alpha''\Vert$ \cite{Manfredo,Struik}, one would say that the torsion is somehow related to an angular momentum. So, it seems natural to expect that the torsion appears in those contexts where the angular momentum plays a role.
\chapter{CONSTRAINED QUANTUM DYNAMICS \, ON \, A \, TUBULAR \, SURFACE}
\label{chap::TubularSurf}

\section{Parametrization of tubular surfaces}

From now on, we will focus on the application of the theoretical framework discussed in the previous chapter to generic tubular surfaces. 
Given a regular curve $\alpha:[0,L]\to\mathbb{R}^3$ of length $L$ parametrized by arc-length $s$, i.e., $\Vert \alpha'(s)\Vert=1$, we may consider a tubular surface of radius $r$ around it. A tube of radius $r$ with generating curve $\alpha(s)$ is parametrized via the Frenet frame $\{\mathbf{t},\mathbf{n},\mathbf{b}\}$ as:
\begin{equation}
X_F(s,\phi)=\alpha(s)+r\cos\phi\,\mathbf{n}(s)+r\sin\phi\,\mathbf{b}(s),\label{eqFrenetParTube}
\end{equation}
which means that a circle with radius $r$ in the plane normal to the tangent vector $\mathbf{t}(s)=\alpha'(s)$ moves along $\alpha(s)$ and then generates the surface of the tube. 
 The first and second fundamental forms of the tube $g_{ij}$ and $h_{ij}$, respectively, read \cite{TubularSurfaces} 
\begin{equation}
g_{ij}=\left(
\begin{array}{cc}
[f_F(s,\phi)^2+\tau(s)^2 r^2] & \tau(s) r^2\\[6pt]
\tau(s) r^2 & r^2\\
\end{array}
\right)\label{eqFirstFormFranetParTube}
\end{equation}
and 
\begin{equation}
h_{ij}=
\left(
\begin{array}{cc}
-[f_F(s,\phi)\kappa(s)\cos\phi-\tau(s)^2r] & \tau(s) r\\
\tau(s) r & r\\
\end{array}
\right),
\label{eqSecondFormFrenetParTube}
\end{equation}
where $i=1$ (or $2$) represents $s$ (or $\phi$) and we have defined
\begin{equation}
f_F(s,\phi)=1-r\kappa(s)\cos\phi;
\end{equation}
$\kappa(s) = \Vert\alpha''\Vert$ and $\tau(s)=\langle\alpha',\alpha''\times \alpha'''\,\rangle\Vert\alpha'\times \alpha''\Vert^{-2}$ are the curvature and the torsion of $\alpha(s)$, respectively \cite{Struik} (if $\alpha$ is not parametrized by an arc-length, then $\kappa=\Vert\alpha'\times\alpha''\Vert\Vert\alpha'\Vert^{-3}$). 

The condition for the parametrization be regular is $g=\det g_{ij}>0$. So, we must have $f_F(s,\phi)\not=0$.  This condition can be satisfied for every coordinate pair $(s,\phi)$ by choosing $r$ smaller than the radius of curvature $R_c(s)=1/\kappa(s)$ for every $s$.

In addition, the Gaussian and Mean curvatures are
\begin{equation}
K=-\frac{\kappa(s)\cos\phi}{rf_F(s,\phi)},\,M=\frac{1}{2r}+\frac{rK}{2},\label{eqGaussAndMeanCurvFrenetParTube}
\end{equation}
respectively. Then we see that the Gaussian and Mean curvatures do not depend on the torsion $\tau$ and thus, the same property holds for the geometry-induced potential in Eq. (\ref{eq:SchrodingerEqExtrinsic}):
\begin{equation}
V_{gip}=-\frac{\hbar^2}{2{m^*}}(M^2-K)=-\frac{\hbar^2}{2{m^*}}\left[\frac{1}{2rf_F(s,\phi)}\right]^2\,.
\label{eqGeoPotencialFrenetParTube}
\end{equation}
Moreover, unless the curve is planar, i.e., with zero torsion, the parametrization $X_F(s,\phi)$ in Eq. (\ref{eqFrenetParTube}) is not orthogonal: $g_{12}=0\Leftrightarrow\tau=0$. Therefore, in the expression for the Laplacian it will appear terms involving derivatives with respect to both $s$ and $\phi$, i.e., terms proportional to $\partial_s\partial_{\phi}$. Nonetheless, we may avoid it by considering a rotation minimizing (RM) frame instead of the Frenet one $\{\mathbf{t},\mathbf{n},\mathbf{b}\}$. An RM frame\footnote{See chapter \ref{chap::RMframes} for more details about RM frames.} is an alternative orthonormal moving frame $\{\mathbf{t},\mathbf{n}_1,\mathbf{n}_2\}$ along $\alpha$, which is written by rotating $\mathbf{n}$ and $\mathbf{b}$ in the plane normal to $\mathbf{t}$:
\begin{equation}
\left(
\begin{array}{c}
\mathbf{t}\\
\mathbf{n}_1\\
\mathbf{n}_2\\
\end{array}
\right)=
\left(
\begin{array}{ccc}
1 & 0 & 0\\
0 & \cos\theta(s) & -\sin\theta(s)\\
0 & \sin\theta(s) & \cos\theta(s)\\
\end{array}
\right)
\left(
\begin{array}{c}
\mathbf{t}\\
\mathbf{n}\\
\mathbf{b}\\
\end{array}
\right).
\end{equation}
Assuming that $\kappa_1(s)$ and $\kappa_2(s)$ are the Bishop parameters \cite{TubularSurfaces,BishopMonthly}, one finds the following set of differential equations which specify (up to an additive constant in $\theta$ \cite{BishopMonthly}) a new frame $\{\mathbf{t},\mathbf{n}_1,\mathbf{n}_2\}$
\begin{equation}
\left\{\begin{array}{ccc}
\mathbf{t}'(s) & = & \kappa_1(s)\,\mathbf{n}_1(s)+\kappa_2(s)\,\mathbf{n}_2(s)\\[3pt]
\mathbf{n}_i'(s) & = & -\kappa_i(s)\,\mathbf{t}(s)\\[3pt]
\kappa_1(s) & = & \kappa(s) \cos\theta(s)\\[3pt]
\kappa_2(s) & = & \kappa(s) \sin\theta(s)\\[3pt]
\theta'(s) & = & \tau(s)\\
\end{array}
\right..
\label{BishopFrame}
\end{equation}
Besides the advantages of the metric obtained from an RM frame (as we will made clear in the following), we mention that if $\alpha$ is singular somewhere, i.e., $\kappa=0$, the Frenet frame is not well-defined. On the other hand,  we can still define $\{\mathbf{t},\mathbf{n}_1, \mathbf{n}_2\}$ even at the singular points \cite{BishopMonthly}. In other words, the Frenet frame is globally defined only if $\alpha$ is a convex curve, while an RM one can be globally defined. On the other hand, if the curve $\alpha$ is closed, the vector $\mathbf{N}_i(0)$ will differ from $\mathbf{N}_i(L)$ ($i=1,2$) by an amount of $\Delta\theta=\int_0^L \tau(u)\mathrm{d}u$, which is a Berry phase \cite{ColinDeVerdiere2006}, as we will show in section \ref{sec::SEqAndGeomPhase} below: see also remark \ref{remark::TerminologyRM}, page \pageref{remark::TerminologyRM}.

The parametrization of a tubular surface of radius $r$ according to an RM frame takes the form
\begin{equation}
X_{RM}(s,\phi)=\alpha(s)+r\cos\phi\,\mathbf{n}_1(s)+r\sin\phi\,\mathbf{n}_2(s)\,,\label{eqBishopParTube}
\end{equation}
with the first and second fundamental forms respectively written as
\begin{equation}
g_{ij}=\left(
\begin{array}{cc}
f_{RM}(s,\phi)^2 & 0 \\
0 & r^2\\
\end{array}
\right)
\label{eqFirstFormBishopParTube}
\end{equation}
and 
\begin{equation}
h_{ij}=\left(
\begin{array}{cc}
-f_{RM}(s,\phi)\,\kappa(s)\cos(\phi-\theta(s)) & 0 \\
0 & r\\
\end{array}
\right)\,,
\label{eqSecondFormBishopParTube}
\end{equation}
where we have defined 
\begin{equation}
f_{RM}(s,\phi)=1-r\kappa(s)\cos(\phi-\theta(s)).
\end{equation}

As before, the condition for the parametrization be regular reads $f_{RM}(s,\phi)\not=0$, which could be satisfied for all values of parameters $(s, \phi)$ by imposing $r<R_c=1/\kappa(s)$. Moreover, the parametrizations $X_F$ and $X_{RM}$ in Eqs. (\ref{eqFrenetParTube}) and (\ref{eqBishopParTube}) are related by the coordinate change $(s,\phi_{RM})\mapsto (s,\phi_F=\phi_{RM}-\theta(s)).$ 

A direct consequence of choosing this new frame is that the parametrization is orthogonal. Then, we can take advantage that the Laplacian operator does not have any term involving derivatives with respect to both $s$ and $\phi$. The only drawback is that for curves with nontrivial torsion the phase $\theta$ includes an integral term $\theta=\int\,\tau$. In this case we must use a numerical method in order to compute $\theta$, see e.g. \cite{WangACMTOG}.

In this new coordinate system the Gaussian and Mean curvatures read
\begin{equation}
K=-\frac{\kappa(s)\cos(\phi-\theta)}{rf_{RM}(s,\phi)},\,M=\frac{1}{2r}+\frac{rK}{2},\label{eqGaussAndMeanCurvBishopParTube}
\end{equation}
respectively. From these expressions one would conclude that the Gaussian and Mean curvatures depend on the torsion. But this dependence is only apparent, since the coordinates according to the Frenet, $(s,\phi_F)$, and to the RM frames, $(s,\phi_{RM})$, are related according to $\phi_F=\phi_{RM}-\theta$. Nonetheless, the dynamics still depends on the torsion (see Section \ref{sec::SEqAndGeomPhase}). 

The geometry-induced potential $V_{gip}$ in the coordinate system $X_{RM}$, Eq. (\ref{eqBishopParTube}), reads
\begin{equation}
-\frac{\hbar^2}{2{m^*}}(M^2-K)=-\frac{\hbar^2}{2{m^*}}\left[\frac{1}{2r f_{RM}(s,\phi)}\right]^2\,,
\label{eqGeoPotencialBishopParTube}
\end{equation}
which just differs from the previous expression, Eq. (\ref{eqGeoPotencialFrenetParTube}), by a torsion-dependent phase $\theta=\int \tau$.

\section{Energy landscape of the geometry-induced potential}

To finish the geometric considerations about tubular surfaces, let us investigate how the landscape of the geometry-induced potential, i.e., the surface $(s,\phi,V_{gip}(s,\phi))$, relate to that of $\kappa$, the curvature function of the reference curve. Such an understanding may play a major role in the study of the quantum dynamics on a tubular surface\footnote{We will study in more detail the problem of finding a surface with a prescribed $V_{gip}$ in chapter \ref{chap::PrescVgip}, where we take into account surfaces with symmetries.}.

To find the critical points of $V_{gip}$, i.e., maxima, minima, and saddle points, one must solve the following system
\begin{equation}
\left\{\begin{array}{ccccc}
\displaystyle\frac{\partial V_{gip}}{\partial s} & = & -\displaystyle\frac{2r}{f_F^3}\,\kappa'\cos\phi & = & 0\\[9pt]
\displaystyle\frac{\partial V_{gip}}{\partial \phi}& =& \displaystyle\frac{2r}{f_F^3}\,\kappa\sin\phi & = & 0\\[7pt]
\end{array}
\right.\,.\label{eq:PartialDerivativesVgeo}
\end{equation}
Since $\cos\phi$ and $\sin\phi$ can not be simultaneously zero, $\kappa$ or $\kappa'$ must vanish. At $\kappa=0$ the generating curve is singular and the analysis of the critical points is subtle. On the other hand, if $\kappa\not=0$ for all points of the generating curve, a critical point $(s_c,\phi_c)$ of $V_{gip}$ satisfies $\kappa'(s_c)=0$. Now, using the second expression in (\ref{eq:PartialDerivativesVgeo}), we find $\phi_c=0$ or $\pi$. The Hessian at a critical point is 
\begin{equation}
\mbox{Hess }V_{gip}(s_c,\phi_c=m\pi) = \frac{2(-1)^mr}{f(s_c,m\pi)^3}\left(
\begin{array}{cc}
-\kappa''(s_c) & 0\\
0 & \kappa(s_c)\\
\end{array}
\right)\,.
\end{equation}
So, we have the following classification scheme for the critical points $(s_c,\phi_c)$ of the geometry-induced potential $V_{gip}$ ($\kappa>0$):
\begin{enumerate}
\item $(s_{\min},0) \Rightarrow$ saddle point;
\item $(s_{\min},\pi) \Rightarrow$ maximum point;
\item $(s_{\max},0) \Rightarrow$ mininum point;
\item $(s_{\max},\pi) \Rightarrow$ saddle point,
\end{enumerate}
where $s_{\min}$ ($s_{\max}$) denotes a minimum (maximum) of $\kappa(s)$. Note that if $s_c$ is an inflexion point of $\kappa$, i.e., $\kappa'(s_c)=\kappa''(s_c)=0$, then $\det\mbox{Hess}\,V_{gip}(s_c,m\pi)=0$ and then $(s_c,m\pi)$ is a degenerate critical point of $V_{gip}$ (i.e., a two-dimensional inflexion point).
\newline
\newline
\begin{remark}
When $\kappa(s^*)=0$, the curve is singular at $s=s^*$ and the Frenet frame is not well-defined at this point. Additionally, if $\kappa$ vanishes along an interval $(s_1,s_2)$, then $\alpha:(s_1,s_2)\to\mathbb{R}^3$ is a line segment and, therefore, the tubular surface is a cylinder and $V_{gip}$ is constant: $-\hbar^2/8{m^*}r^2$.
\end{remark}

\section{Schr\"odinger equation and geometric phase}
\label{sec::SEqAndGeomPhase}

In the parametrization via Frenet frame, $X_{F}$ in Eq. (\ref{eqFrenetParTube}), the eigenvalue problem $\Delta_{LB}\psi+(M^2-K)\psi=-2{m^*}E\psi/\hbar^2$ reads
$$\frac{1}{f_F^2}\left[\frac{\partial^2\psi}{\partial s^2}+\frac{\tau \partial_{\phi}f_{F}-\partial_sf_{F}}{f_F}\frac{\partial\psi}{\partial s}\right]- \displaystyle\frac{2\tau}{f_F^2}\frac{\partial^2\psi}{\partial s\partial\phi}+\frac{f_F^2+r^2\tau^2}{r^2f_F^2}\frac{\partial^2\psi}{\partial\phi^2}+$$
\begin{equation}
+\displaystyle\frac{r^2[\tau \partial_sf_F-\tau'f_F]+[f_F^2-\tau^2r^2]\partial_{\phi}f_F}{r^2f_F^3}\frac{\partial\psi}{\partial\phi}+\frac{\psi}{4r^2f_F^2}= -\displaystyle\frac{2{m^*}E}{\hbar^2}\psi\,.
\label{eqEigenvalueFrenetParTub}
\end{equation}
Observe that the geometry-induced potential does not depend on the curve torsion $\tau$, Eq. (\ref{eqGeoPotencialFrenetParTube}), which implies that a planar and non-planar tube with the same curvature profile will have the same $V_{gip}$. On the other hand, it is easy to see from the above expression for the Schr\"odinger equation that the torsion does play a role on the quantum dynamics of a particle constrained to a tubular surface.

Now, using the parametrization via an RM frame, $X_{RM}$ in Eq. (\ref{eqBishopParTube}), we obtain a simpler eigenvalue problem than the previous one using Frenet frame. However, we mention that it involves an integration $\theta(s)=\int\,\tau$, which may turn the approach analytically more difficult. The corresponding eigenvalue problem via an RM then reads
\begin{eqnarray}
\frac{1}{f_{RM}^2}\left(\frac{\partial^2\psi}{\partial s^2}-\frac{\partial_s f_{RM}}{f_{RM}}\frac{\partial \psi}{\partial s}\right) & + & \frac{1}{r^2}\left(\frac{\partial^2\psi}{\partial \phi^2} + \frac{\partial_{\phi}f_{RM}}{f_{RM}}\frac{\partial\psi}{\partial \phi}\right)+\frac{\psi}{4r^2f_{RM}^2}= -\frac{2{m^*}E\psi}{\hbar^2}
\label{eqEigenvalueBishopParTub}
\end{eqnarray}

Now, we may pass from the ``curvilinear'' equation above to an ``rectangular'' one by   rescaling the wave function as $\Psi=\sqrt{f_{RM}}\psi$. This conserves the probability density $\int\vert\Psi\vert^2\,\mathrm{d}s\,\mathrm{d}\phi=\int\vert\psi\vert^2\,\mathrm{d}S=1$ and the following eigenvalue problem follows
\begin{equation}
-\frac{1}{f_{RM}^2}\frac{\partial^2\Psi}{\partial s^2}-\frac{1}{r^2}\frac{\partial^2\Psi}{\partial \phi^2}+2\frac{\partial_s f_{RM}}{f_{RM}^3}\frac{\partial\Psi}{\partial s}+V_{eff}\,\Psi=\frac{2{m^*}E}{\hbar^2}\Psi,
\label{eqEigenvalueBishopRescaled}
\end{equation}
where
\begin{equation}
V_{eff} = \frac{\partial^2_{\phi} f_{RM}}{2r^2f_{RM}}-\frac{(\partial_{\phi}f_{RM})^2}{4r^2f_{RM}^2}+\frac{\partial^2_sf_{RM}}{2f_{RM}^2}-\frac{5(\partial_sf_{RM})^2}{4f_{RM}^4}-\frac{1}{4r^2f_{RM}^2}.
\end{equation}

It is known that the torsion of a curve may give rise to a geometric phase (Berry phase) for the propagation along a waveguide \cite{KuglerPRD1988,TwistedRingProgTheorPhys}: $\theta=\int\tau$. This is precisely the phase introduced in previous sections in order to relate an RM to the Frenet frame, Eq. (\ref{BishopFrame}). In fact, this can be formalized by noting that an RM is composed by vectors parallel-transported along the reference curve with respect to the normal bundle connection \cite{ColinDeVerdiere2006,Etayo2016}. 

We can directly observe the appearance of the Berry phase for the propagation on a tubular surface through the study of a thin tube. In fact, let us assume the limiting behavior $r\kappa,r\kappa',r\kappa''\ll1$ and expand the coefficients in the Schr\"odinger equation up to first order. Then
\begin{equation}
\frac{\partial^2\xi}{\partial s^2}+\frac{1}{r^2}\frac{\partial^2\xi}{\partial \phi^2}+2\,V_{eff}^{(2)}\,\xi=-\left(\frac{2{m^*}E}{\hbar^2}+\frac{1}{4r^2}\right)\xi,
\end{equation}
where we have denoted $r\kappa=\epsilon$, $r\kappa'=\eta$, $r\kappa''=\nu\ll1$ and the (expanded) effective potential reads
$$
V_{eff}^{(2)} 
  =  \frac{(\nu-\epsilon\tau^2)}{4}\cos\left(\phi-\theta\right)+\frac{(2\eta\tau+\epsilon\tau')}{4}\sin\left(\phi-\theta\right)+$$
  \begin{equation}
  \left[\eta\cos\left(\phi-\theta\right)+
  \epsilon\tau\sin\left(\phi-\theta\right)\right]\frac{\partial}{\partial s}+ \epsilon\cos\left(\phi-\theta\right)\frac{\partial^2}{\partial^2 s}\,.
\end{equation}
Now, collecting the terms of zero order gives
\begin{equation}
-\frac{\partial^2\Psi}{\partial s^2}-\frac{1}{r^2}\frac{\partial^2\Psi}{\partial \phi_{RM}^2}=\left(\frac{2{m^*}E}{\hbar^2}+\frac{1}{4r^2}\right)\Psi,
\end{equation}
where $(s,\phi)\in[0,L]\times[0,2\pi]$. 
We may assume an angular periodicity for $\phi$ and,
depending if the reference curve is open or closed, homogeneous or periodic boundary conditions in $s$, respectively. Once we have fixed the type of boundary conditions in $s$, the eigenvalues for a tubular surface with and without the $V_{gip}$ (in thin cylinder approximation) are the same. In particular, the eigenvalues do not depend on the curve torsion on such an approximation. On the other hand, the eigenfunctions do depend on the curve torsion through its geometric phase $\theta=\int \tau$:
\begin{equation}
\Psi_{n,\ell} = \psi_n(s)\,\mathrm{e}^{\mathrm{i}\,r\,\ell\,\phi_{RM}}= \psi_n(s)\,\mathrm{e}^{\mathrm{i}\,r\,\ell\,\phi_F}\mathrm{e}^{\mathrm{i}\,r\,\ell\,\theta},
\end{equation}
where $\psi_n(s)$ is an eigenfunction relative to $s$ and we used the relation $\phi_F=\phi_{RM}-\theta$. From this expression we see that the solution of a non-planar case ($\tau\not=0$) differs from the planar one $\Psi_{\tau=0}=\psi_n(s)\,\mathrm{e}^{\mathrm{i}\,r\,\ell\,\phi_F}$ by an amount corresponding to a geometric phase $r\,\ell\,\theta=r\,\ell\int\tau$.

The appearance of the geometric phase in the eigenfunctions shows that the curve torsion is an important ingredient in the study of transport and interference problems on tubular surfaces.
\chapter{CURVES AND SURFACES WITH A \,\,PRESCRIBED GEOMETRY-INDUCED POTENTIAL}
\label{chap::PrescVgip}

Exploiting the effects of an extra contribution to the Hamiltonian due to a confining potential approach is essential and in this respect an important problem is that of a prescribed geometry-induced potential, i.e., the inverse problem of finding a curved region with a potential given \textit{a priori}. Inverse problems constitute an important subject from both  experimental and theoretical viewpoints, a classical problem being that of hearing the shape of a drum \cite{KacMonthly1966}, i.e., the determination of information about the geometry of a region that gives rise to a prescribed spectrum. More recently, we can mention the success in the detection of gravitational waves \cite{AbbottPRL2016}, which allows one to infer information about the spacetime geometry via measurements of an interferometric gravitational-wave detector. In the context of the constrained quantum dynamics this kind of problem offers the possibility of engineering surfaces and curves with a quantum behavior prescribed \textit{a priori} through their geometry-induced potential and has been already investigated for curves \cite{delCampoSciRep2014} and a class of revolution surfaces \cite{AtanasovPLA2007}. Nonetheless, a comprehensive understanding of such an inverse problem for the constrained dynamics is still absent.

In this work we address the problem of prescribed GIP for curves and for surfaces in Euclidean space $\mathbb{R}^3$. The former can be easily solved by integrating the Frenet equations, while the latter involves the solution of a non-linear 2nd order PDE. We restrict ourselves to the study of surfaces invariant by a 1-parameter group of isometries of $\mathbb{R}^3$, which turns the PDE for the prescribed GIP into an ODE and leads us to the study of cylindrical, revolution, and helicoidal surfaces. The latter class is particularly important due to the fact that, by screw-rotating a curve clock- and counterclockwisely, one can generate pairs of enantiomorphic surfaces, which turn these objects the natural candidates to test and exploit a link between chirality and the effects of a GIP. We show how to find helicoidal surfaces associated with a given non-negative function and further specialize to the study of helicoidal minimal surfaces. For this class of minimal surfaces we prove the existence of localized and geometry-induced bound states, then generalizing known results for the dynamics on a helicoid \cite{AtanasovPRB2009}, and also the possibility of controlling the change in the distribution of the probability density when the surface is subjected to an extra charge.

\section{Curves with prescribed geometry-induced potential}

For the confinement on a curve $\alpha:I\to\mathbb{R}^3$ the problem reduces to that of finding a curve with a prescribed curvature function. The curve is then obtained after integration of the Frenet equations (\ref{eq::FrenetFrameEqs}) as $\alpha(s)=\alpha_0+\int_{s_0}^s\,\mathbf{t}(u)\,du.$ It is worth to mention that in the case of plan curves, i.e., $\tau\equiv0$, the parametrization of the solution curve for the Frenet equations are \cite{Struik}
\begin{equation}
\left\{
\begin{array}{c}
x(s)=z_{1}\,C(s)-z_{2}\,S(s)+x_0\\[6pt]
y(s)=z_{1}\,S(s)+z_{2}\,C(s)+y_0\\
\end{array}
\right.\,,
\end{equation}
where $x_0,\,y_0,$ and $z_{i}$ are constants to be specified by the initial conditions and
\begin{equation}
\left\{
\begin{array}{c}
S(s)=+\displaystyle\int_{s_0}^s\cos\Big(\int_{s_0}^v\kappa(u)\,\mathrm{d}u\Big)\mathrm{d}v\\[8pt]
C(s)=-\displaystyle\int_{s_0}^s\sin\Big(\int_{s_0}^v\kappa(u)\,\mathrm{d}u\Big)\mathrm{d}v\\
\end{array}
\right.\,.\label{eq::SenoCosFrenetFrame}
\end{equation}
Recently, del Campo \textit{et al.} exploited such an explicit solution in order to find pair of curves with the same scattering properties \cite{delCampoSciRep2014}: see also \cite{MatsutaniJPSJ1991}. Finally, for $\tau\not=0$, it is possible to find the general solution for the Frenet equations by writing them in term of a complex Riccati equation \cite{Struik}.

\subsection{The Hydrogen atom in a curve: power-law curvature functions}

The 1D Hydrogen atom is characterized by the following Hamiltonian
\begin{equation}
\hat{H} = -\frac{\hbar^2}{2m^*}\frac{\mathrm{d}^2}{\mathrm{d}x^2}-\frac{\mathrm{e}^2}{4\pi\varepsilon_0\vert x\vert}\,.
\end{equation}
Then, we can \textit{model} the Hydrogen atom through a confinement on a curve by considering the curvature function  
\begin{equation}
\kappa(s)=\sqrt{\frac{8m^*}{\hbar^2}\frac{\mathrm{e}^2}{4\pi\varepsilon_0}}\,\frac{1}{\sqrt{s}}\,,
\end{equation}
where $s$ denotes the arc-length parameter of $\alpha$.

In order to find a plane curves whose geometry-induced potential is Hydrogen-like we may consider curves with a power law curvature function, i.e., $\tau\equiv0$ and $\kappa(s)=c_0/s^p$, where $c_0>0$ and $p\in\mathbb{R}$ are constants. So, a planar Hydrogen curve is that curve with $p=1/2$ and $c_0=\sqrt{8m^*\mathrm{e}^2/4\pi\varepsilon_0\hbar^2}$. The parametrization for a Hydrogen curve may be written as (see subsection \ref{subsec::PowerLawCurvature}, p. \pageref{subsec::PowerLawCurvature}, for details)
\begin{equation}
\alpha(s)=\left(
\begin{array}{cc}
\frac{1}{2c_0^2} & \frac{\sqrt{s}}{c_0}\\
-\frac{\sqrt{s}}{c_0} & \frac{1}{2c_0^2}\\
\end{array}
\right)\left(
\begin{array}{c}
\cos\left(2c_0\sqrt{s}\right)\\
\sin\left(2c_0\sqrt{s}\right)\\
\end{array}
\right)\,,
\end{equation}

In the 1950s Loudon solved the 1D-Hydrogen atom on the line \cite{LoudonAmJPhys}. Then, since the eigenfunction along the Hydrogen curve will be a function of the arc-length parameter $s>0$, we have the following wave function along the curve
\begin{equation}
\psi = B\mathrm{e}^{-\frac{z}{2}}zL_{N}^1(z),\,\,z=\frac{2s}{Na_0},
\end{equation}
where $B$ is a normalizing constant, $a_0=\hbar^2/m^*\mathrm{e}^2$, and $L_a^b(z)$ denotes an associated Laguerre polynomial. This solution is not equal to the radial solution of the 3D Hydrogen:
\begin{equation}
R_{N\ell}(r) = B_{N\ell}\,\mathrm{e}^{-\frac{z}{2}}z^{\ell}L_{N}^{2\ell+1}(z),\,\,z=\frac{2r}{Na_0},
\end{equation}
where $B_{N\ell}$ is a normalizing constant. However, taking into account the use of spherical coordinates to describe the radial part, one obtains the same probability density in both cases: $\mathrm{d}P_{1D}=\vert\psi_{1D}\vert^2\mathrm{d}s=\mathrm{d}P_{3D}=r^2\vert\psi_{1D}\vert^2\mathrm{d}r$, where one must take $\ell=0$ in the 3D solution in order to properly compare the solutions in both dimensions. As expected, this means that in the 1D solution only $s$ orbitals make sense. Then, a 1D periodic table will have 2 columns only \cite{Chemistry1D,Chem1DPackage}.

\subsection{Intrinsic approach to the 1D constrained dynamics}

Before leaving the study of curves, let us investigate what happens in an intrinsic approach for 1D dynamics. Recently, it was observed that the energy spectrum of a free non-relativistic particle (intrinsically) confined to a curve $\alpha:I\to\mathbb{R}^2$ only depends on the length $L$ of the curve and on the imposed boundary conditions \cite{BastosEtAlPhysEd2012}. In other words, there exist essentially two types of a 1D Particle in a Box Model (PIB model), the open and closed boxes, i.e., homogeneous (HBC) or periodic (PBC) boundary conditions\footnote{One may impose general boundary conditions on the PIB model \cite{GeneralPIBPRD1990}, but the problem may be no longer exactly soluble \cite{GeneralPIBPRA1995}.}. Here we make the important observation that this result does not depend neither on the ambient space nor on the dimensions considered, i.e., given any Riemannian manifold $M$, the spectrum of the Schr\"odinger equation on a curve $\alpha:[0,L]\to M$ only depends on the length and imposed boundary conditions. Indeed, this is based on the fact that the Laplace operator on $\alpha$ is simply the second derivative with respect to its arc-length $s$:
\begin{equation}
\Delta_{\alpha} = \frac{\mathrm{d}^2}{\mathrm{d}\,s^2}\,.\label{eq::LaplaceOnCurves}
\end{equation}
So, we have the following general result for the spectrum on a curve in any ambient space:
\begin{enumerate}
\item If $\alpha$ is an open curve, i.e., assuming HBC,
 \begin{equation}
 E_n(HBC) = \frac{h^2n^2}{8m^*L^2}\,\,,n=1,2,...;
\end{equation}
\item If $\alpha$ is a closed curve, i.e., assuming PBC,
 \begin{equation}
 E_n(PBC) = \frac{h^2n^2}{2m^*L^2} = 4E_n(HBC)\,\,,n=1,2,...\,.
\end{equation}
\end{enumerate}

Since the effect of an external scalar potential $\mathcal{V}$ in an intrinsic scheme is simply given by its restriction to $\alpha$, $\mathcal{V}\equiv \mathcal{V}\vert_{\alpha}$, we see that the specific shape of a 1D box is immaterial, the important feature being the fact that the particle is confined somewhere. Indeed, given a curve $\alpha(u)$ in $M$ subjected to an external potential function $\mathcal{V}$ (this potential is not a confining one), we can apply a re-parametrization by arc-length and then the Hamiltonian operator changes as
\begin{equation}
\hat{H}_{u}=-\frac{\hbar^2}{2m^*}\Delta_{\alpha}+\mathcal{V}(\alpha(u))\,\mapsto\, \hat{H}_{s}=-\frac{\hbar^2}{2m^*}\frac{\mathrm{d}^2}{\mathrm{d}s^2}+\mathcal{V}(\alpha(s))\,.\label{eq::changeOfHbyarclength} 
\end{equation} 

Although simple, in many contexts the PIB model applies nicely \cite{NanotubeAs1DBoxPRL1999,BilayerGrapheneNanoRes}. This is the case because in such models one is primarily interested on the existence of a confinement. Naturally, an improved version of the PIB model is necessary in a more realistic context, e.g., if one wants to take into account the role played by the surfaces of distinct nanostructures on the 1D constrained dynamics. 

Finally, let us comment that the situation in dimension greater than 1 is much more complex and, to the best of our knowledge, no simple characterization of the Laplacian operator, which would allow for a fully description of the eigenvalue problem in higher dimensions is available \cite{GeometryLaplacian,GilkeyMatContemp}. We should necessarily restrict ourselves to some particular class of manifolds, e.g., spherical space forms \cite{GilkeyMatContemp} (constant positive curvature) or generalized cylinders \cite{BastosEtAlPreprint}, just to name a few.

\section{Surfaces with prescribed geometry-induced potential}

For surfaces, the situation is more complex. Indeed, the prescribed GIP problem generally demands the solution of a 2nd order non-linear PDE. For example, assuming the surface to be the graph of a smooth function $Z(x,y)$, i.e., the parametrization is given by $\mathbf{r}(x,y) = (x,y,Z(x,y))$, the Gaussian and Mean curvatures are written as
\begin{equation}
K(x,y) = \displaystyle\frac{Z_{xx}Z_{yy}-Z_{xy}^2}{(1+Z_{x}^2+Z_{y}^2)^2}=\frac{\det(\mbox{Hess}\,Z)}{(1+\Vert\nabla Z\Vert^2)^2}
\end{equation}
and
\begin{equation}
H(x,y) = \displaystyle\frac{Z_{xx}(1+Z_{y}^2)-2Z_{xy}Z_{x}Z_{y}+Z_{yy}(1+Z_{x}^2)}{2(1+Z_{x}^2+Z_{y}^2)^{3/2}}=\frac{1}{2}\nabla\cdot\left(\frac{\nabla Z}{\sqrt{1+\Vert\nabla Z\Vert^2}\,}\right),
\end{equation}
respectively. The equation for $K$ is a nonlinear elliptic PDE of Hessian type (also referred as Monge-Amp\`ere equation) \cite{Gutierrez2001}, while the equation for $H$ is a nonlinear elliptic PDE of divergent type \cite{GilbargTrudinger1977}. 

A general study of the PDE associated with the prescribed GIP $H^2-K$ is not a trivial task. In addition it can encode in its generality useless examples. In this respect, the study of particular classes can turn to be more useful and insightful than a general analysis. So, instead of studying the prescribed potential problem in general, which would lead us to the realm of non-linear analysis \cite{Gutierrez2001,GilbargTrudinger1977}, here we restrict ourselves to the simpler, but still important and difficult, context of invariant surfaces (continuous symmetries). To be more precise, we assume the surfaces to be invariant by a 1-parameter group of isometries of $\mathbb{R}^3$ \cite{DoCarmoTohoku1982,MedeirosRMU1991}. This allows us to avoid the study of a non-linear PDE, since the symmetry turns the equation into an (non-linear) ODE along the so called generating curve.

\subsection{Surfaces invariant by a 1-parameter subgroup of isometries}

Basically there exist three types of surfaces invariant by a 1-parameter subgroup of isometries of $\mathbb{R}^3$, namely (i) cylindrical surfaces (translation symmetry), (ii) surfaces of revolution (rotational symmetry), and (iii) helicoidal surfaces (screw rotation symmetry, i.e., a combination of a translation and a rotation). Due to their appealing symmetry, these surfaces are commonly encountered in applications and theoretical studies in the context of a constrained dynamics: e.g., cylindrical surfaces to model rolled-up nanotubes \cite{OrtixPRB2010} and $\pi$ electron energies of aromatic molecules \cite{BastosEtAlPreprint,MiliordosPRA2010}; surfaces of revolution as tractable examples to test the validity and potentialities of an extrinsic confinement approach \cite{EncinosaPRA,GravesenJMathPhys}; and helicoidal surfaces to study geometry-induced charge separation \cite{AtanasovPRB2009} and the relation to the concept of chirality \cite{AtanasovPRB2015}, just to name a few.
\newline

A function $T:\mathbb{R}^3\to\mathbb{R}^3$ is an isometry of $\mathbb{R}^3$ if it satisfies for all $q,\tilde{q}\in\mathbb{R}^3$ the relation $\langle T(q),T(\tilde{q})\rangle=\langle q,\tilde{q}\rangle$. These functions form the so called group of rigid motions of $\mathbb{R}^3$, which are composed by translations $T_a(q)=q+a$ and rotations $R\in O(3)$ (or $SO(3)$ if one imposes that $T$ preserves orientation). By a one-parameter subgroup of isometries we mean an action of the additive group $(\mathbb{R},+)$ on the symmetry group (rigid motions) of $\mathbb{R}^3$. In other words, a 1-parameter subgroup of isometries is a smooth map $\gamma:\mathbb{R}\times\mathbb{R}^3\to\mathbb{R}^3$ such that
\begin{enumerate}
\item For all $t\in\mathbb{R}$ the map $q\mapsto \gamma(t,q)$, denoted by $\gamma_t$, is a rigid motion;
\item For all $t,s\in\mathbb{R}$, $\gamma_t\circ\gamma_s=\gamma_{t+s}$ and $\gamma_0=Id$ is the identity map.
\end{enumerate}
Up to a change of variables, every 1-parameter subgroup can be written as \cite{MedeirosRMU1991}
\begin{equation}
\gamma(t,q)=(q^1\cos t+q^2\sin t,-q^1\sin t+q^2\cos t,q^3+ht),
\end{equation}
or as
\begin{equation}
\gamma(t,q)=(q^1,q^2,q^3+ht),
\end{equation}
where $h\in\mathbb{R}$ is a constant, equal to zero for rotational symmetry in the former or equal to zero for the identity map in the latter.
\newline
\newline
\textit{Remark:} When discussing the constrained dynamics on a helicoidal surface it will prove useful to adopt a different notation. More precisely, we will assume the 1-parameter subgroup of isometries to be $$\gamma(t,q)=(q^1\cos (\omega t)+q^2\sin(\omega t),-q^1\sin(\omega t)+q^2\cos(\omega t),q^3+t),$$ where $\omega$ is a constant.
\newline

A surface $\Sigma\subseteq \mathbb{R}^3$ invariant by a 1-parameter subgroup of isometries of $\mathbb{R}^3$ is characterized by 
\begin{equation}
\forall\,t\in\mathbb{R},\,\Sigma = \gamma_t(\Sigma)\,.
\end{equation}
Intuitively, we can approximate an invariant surface by successive applications, to a given curve $\alpha(s)$, of a certain kind of rigid motion: $$\Sigma\cong\{\gamma_{t_0}(\alpha(s)),\gamma_{t_0+\Delta t}(\alpha(s)),\cdots,\gamma_{t_0+n\Delta t}(\alpha(s))\}.$$ So, in the limit $\Delta t\to 0$, we generate the surface by continuously moving the curve $\alpha$ by the action of a 1-parameter subgroup $\gamma_t$. We call such a curve the \textit{generating curve}, which can be assumed to be planar.

It follows that the values of the Gaussian and Mean curvature only depend on the values assumed along the generating curve. As a corollary of the invariance of $K$ and $H$, the prescribed GIP problem demands the solution of an ODE instead of a PDE. In the next section we present a study of this problem for each type of invariant surface. 

\section{Cylindrical surfaces with prescribed geometry-induced potential}

Now we focus on the simplest instance of surfaces invariant by a 1-parameter subgroup of isometries, namely, surfaces with translation symmetry. A cylinder is the standard example, it is just the surface obtained by translating a circle. More generally, a cylindrical surface is obtained by taking a generating curve (the cross section) which can be any planar curve $\alpha(s):I\to\mathbb{R}^2\subset\mathbb{R}^3$ (for a study of cylindrical surfaces with a varying cross section see \cite{BastosEtAlPreprint}). We then translate this curve in the direction of a unit vector $\textbf{a}=(a_1,a_2,a_3)$, where we assume $a_3\not=0$ in order to have a regular surface, i.e., $\textbf{a}$ is out of the $xy$ plane\footnote{We could have assumed $\mathbf{a}$ to be $(0,0,1)$, but we decided to work with an arbitrary vector in order to include inclined cylinders in our discussion.}. By denoting $\alpha(s)=(x(s),y(s),0)$, where $s$ is an arc-length parameter, we have the following parametrization for a cylindrical surface
\begin{equation}
\mathbf{x}(s,t) = \alpha(s)+t\,\textbf{a}.\label{eq::parCylindricalSurface}
\end{equation}
Observe that the generating curve does not need to be closed.

The coefficients of the first and second fundamental forms are given by
\begin{equation}
g_{11}(s,t)=1,\,g_{12}(s,t)=\cos\theta,\,g_{22}(s,t)=1,
\end{equation}
and
\begin{equation}
h_{11}(s,t)=h_{12}(s,t)=0,h_{22}(s,t)=\langle\alpha'\times\alpha'',\mathbf{a}\rangle,
\end{equation}
respectively; where we have adopted the unit normal $\mathbf{n}=\textbf{a}\times\alpha'$ and $\theta=\cos^{-1}\langle \alpha',\textbf{a}\rangle$  is the (constant) angle between $\textbf{a}$ and $\alpha'$. Now we can compute the Gaussian and Mean curvatures of a cylindrical surface as
\begin{equation}
K\equiv0\,\mbox{ and }\,H=\frac{a_3[x'(s)\,y''(s)-x''(s)\,y'(s)]}{2\,\sin^2\theta}\,,
\end{equation}
respectively. Notice, as expected, that due to the translation symmetry the Gaussian and Mean curvatures are functions of $s$ only. On the other hand, since $K\equiv0$, the problem of a prescribed GIP $H^2-K$ is equivalent to the problem of finding cylindrical surfaces with prescribed Mean curvature. Then, given a function $H(s)$, one must solve the following system of 2nd order nonlinear ODEs
\begin{equation}
\left\{
\begin{array}{c}
x'\,y''-x''\,y'=\frac{2\,{\sin^2\theta}}{a_3}\,H(s)\\[5pt]
(x')^2+(y')^2=1\\
\end{array}
\right.\,,\label{eq::ODEPresc_H_ForCylindricalSur}
\end{equation}
where the second equation comes from the parametrization by arc-length. 

For a planar curve $\alpha(s)=(x(s),y(s))$, we can write the curvature function as \cite{Manfredo,Struik}
\begin{equation}
\kappa = \frac{x'\,y''-x''\,y'}{[(x')^2+(y')^2]^{3/2}}.
\end{equation}
Then, we have the following result 
\begin{prop} 
The Mean curvature $H(s)$ of a cylindrical surface and the curvature function $\kappa(s)$ of its generating curve (cross section) are related according to
\begin{equation}
\kappa(s) = \frac{2\,{\sin^2\theta}}{a_3}H(s)\,,\label{eq::RelationCurvCurveAndHcylindSurf}
\end{equation}
where $\theta$ is the (constant) angle between the direction of translation $\mathbf{a}=(a_1,a_2,a_3)$ and the plane which contains the generating curve. Moreover, it follows that Eq. (\ref{eq::FranetFrameEqs}) solves the problem of prescribed Mean curvature, i.e., there is an equivalence between finding curves with prescribed curvature and finding cylindrical surfaces with prescribed Mean curvature.
\end{prop}
\textbf{Example} (Right cylindrical surfaces with constant Mean curvature): A right cylindrical surface is given by the condition $\textbf{a}=(0,0,\pm1)$, which implies $\cos\theta\equiv\pm1$. Now, assume that $H(s)\equiv H_0\not=0$ is a constant. Also assume for simplicity all signs equal to $+$ and $s_0=0$ (the other cases are analogous). Then, one finds
\begin{equation}
\left(
\begin{array}{c}
x-x_0\\
y-y_0\\
\end{array}
\right)
=\frac{1}{2H_0}\left(
\begin{array}{cc}
\cos(2H_0\,s) & - \sin(2H_0\,s)\\
\sin(2H_0\,s) & \cos(2H_0\,s)\\
\end{array}
\right)\left(
\begin{array}{c}
z_{1}\\
z_{2}\\
\end{array}
\right)- \frac{1}{2H_0}\left(
\begin{array}{c}
z_{1}\\
z_{2}\\
\end{array}
\right),
\end{equation}
which represents a right cylinder with radius $R=1/2H_0$. On the other hand, if $H_0\equiv0$, then $(x(s),y(s))=(-z_{02},z_{01})\,s+(x_0,y_0)$, which represents a line segment that generates a  cylindrical surface which is a (piece of a) plane. 


\section{Surfaces of revolution with prescribed geometry-induced potential}


A first attempt to solve the prescribed GIP problem for surfaces of revolution was devised by Atanasov and Dandoloff \cite{AtanasovPLA2007}. They considered surfaces of revolution whose generating curve, to be rotated around the $z$ axis, is a graph on the $xz$ plane. They also investigated the existence of bound states and surfaces in the form of circular strips around the symmetry axis.

In the following, we consider surfaces of revolution without imposing any restriction on the generating curve. We show that the equation for the prescribed GIP can be rewritten as a first order complex equation. Further, we specialize to surfaces whose generating curve is a graph on the $xz$ plane that can be rotated around either the $x$ or the $z$ axis. 

Given a curve $\alpha(s)=(x(s),0,z(s))$ on the $xz$ plane ($s$ being its arc-length parameter), the surface of revolution obtained by rotating $\alpha$ around the $z$ axis is parametrized by
\begin{equation}
\mathbf{x}(s,\phi)=(x(s)\cos\phi,x(s)\sin\phi,z(s)),
\end{equation}
where we must assume $x(s)>0$ for all $s$.

The coefficients of the first and second fundamental forms are given by
\begin{equation}
g_{11} = 1,\,\,g_{12}=0,\,\,g_{22}=x^2(s),
\end{equation}
and
\begin{equation}
h_{11} = x'(s)z''(s)-x''(s)z'(s),\,\,h_{12}=0,\,\,h_{22}=x(s)z'(s),
\end{equation}
respectively. From these expressions we find
\begin{equation}
U=\sqrt{H^2-K} = \frac{x(x'z''-x''z')-z'}{2x}.\label{eqPrescribedGIPFullRevolSurf}
\end{equation}
Observe the similarity of this expression with that of the Mean curvature:
\begin{equation}
H = \frac{x(x'z''-x''z')+z'}{2x}.
\end{equation}
Indeed, they are the same except for the exchange of the sign in front of $z'$. This similarity will be exploited in the following.

The equation of prescribed $U$, or $H$, is a 2nd order non-linear ODE\footnote{In fact, since we are assuming $x'^2+z'^2=1$, the prescribed curvature problem is given by a system of 2nd order non-linear ODE's.}. In the 80's, Kenmotsu solved the prescribed Mean curvature equation by transforming it in a 1st order complex linear ODE \cite{KenmotsuTohokuMathJ}: $Z'-2\,\textrm{i}\,H\,Z+1=0$. This technique can be applied to our problem, i.e., we can write the equation for $U$ as a 1st order complex ODE, which in our case is non-linear: $Z'-2\,\textrm{i}\,U\,Z+\vert Z\vert^2=0$.

Multiplying Eq. (\ref{eqPrescribedGIPFullRevolSurf}) by $x'$, and using $x'^2+z'^2=1$ and its derivative, we find
\begin{equation}
0=2xx'U+x'z'-xz''=2\frac{x'}{x}U-\left(\frac{z'}{x}\right)'\,.
\end{equation}
On the other hand, multiplying Eq. (\ref{eqPrescribedGIPFullRevolSurf}) by $z'$, and using $x'\,^2+z'\,^2=1$ and its derivative, we have
\begin{equation}
0=2xz'U+xx''-x'\,^2+1=2\frac{z'}{x}U+\left(\frac{x'}{x}\right)'+\frac{1}{x^2}\,.
\end{equation}
Finally, defining $Z(s)=x^{-1}(s)[x'(s)+\textrm{i}\,z'(s)]$, we can glue the above equations together and write
\begin{equation}
Z'(s)-2\,\textrm{i}\,U(s)\,Z(s)+\vert Z(s)\vert^2=0\,.
\end{equation}

In the next subsections we will study some particular classes of revolution surfaces where the equation for the prescribed potential can be effectively solved. 

\subsection{Surfaces whose generating curve is a graph rotated around a vertical axis}
 
In the end of the 1990s, Baikoussis and Koufogiorgos \cite{BaikoussisJGeom} studied the problem of finding helicoidal surfaces with prescribed Mean or Gaussian curvatures. They assumed a parametrization
given by
\begin{equation}
\mathbf{x}(\rho,\phi) = (\rho\cos \phi,\rho\sin \phi,\lambda(\rho)+h\,\phi),\,\rho>0,\label{eq::ParamHelicoidalSur}
\end{equation}
where $h$ is a constant and $\lambda(\rho)$ a smooth function, which represents the generating curve $(\rho,0,\lambda(\rho))$. As natural, $\phi$ stands for the rotation angle around the $Oz$ axis, the screw axis, and $\rho$ for the distance from it. 

If $h=0$, the helicoidal surface is just a surface of revolution, while if $\lambda\equiv0$ and $h\not=0$ one has the usual helicoid surface. In addition, since the generating curve $\lambda$ is supposed to be a graph, cylinders are not covered by (\ref{eq::ParamHelicoidalSur}) (such an example will be covered in the following subsection by allowing a rotation around the $x$ axis).  

The problem of prescribed Mean or Gaussian curvatures is then solved by writing the curvatures of the given surface in terms of the parameters $h$ and $\lambda(\rho)$. This leads to an ODE that, if properly manipulated, can be written as 
\begin{equation}
\frac{\rho}{2}A'(\rho)+A(\rho)=H_0(\rho) \mbox{ and }\frac{1}{2\rho}(B^2(\rho))'=K_0(\rho),\label{eq::MeanAndGaussODE}
\end{equation}
where
\begin{equation}
A = \frac{\lambda'}{\sqrt{\rho^2(1+\lambda'\,^2)+h^2}}\,;\,B^2=\frac{\rho^2\lambda'\,^2+h^2}{\rho^2(1+\lambda'\,^2)+h^2}\,.\label{eq::AandBfromMeanGaussEDO}
\end{equation}

We now apply these ideas to surfaces of revolution by imposing  $h=0$. It follows that $B^2=\rho^2A^2$, which gives us the following ODE in terms of $U$ ($=\sqrt{H^2-K}$)
\begin{equation}
\frac{\rho^2}{4}(A')^2=U^2\Rightarrow A(\rho)=\pm\left(2\int U(\rho)\,\frac{\textrm{d}\rho}{\rho}+a_1\right)\,,\label{eq::willmoreEDO}
\end{equation}
where $a_1$ is a constant of integration. Using this in Eq. (\ref{eq::AandBfromMeanGaussEDO}) under the condition $h=0$, one obtains an ODE for the generating curve $\lambda(\rho)$:
\begin{equation}
\lambda'\,^2 = A^2\,\rho^2\,(1+\lambda'\,^2)\Rightarrow [1-\rho^2A^2]\,\lambda'\,^2=\rho^2\,A^2\geq0.
\end{equation}
By continuity, if $1-\rho_0^2A(\rho_0)>0$ at some $\rho_0\in\mathbb{R}-\{0\}$, then $1-\rho^2A^2(\rho)>0$ on a neighborhood of $\rho_0$. So, one gets the general solution in the neighborhood of $\rho_0$
\begin{equation}
\lambda(\rho) = \pm\int\frac{\rho A(\rho)}{\sqrt{1-\rho^2A^2(\rho)}}\,\textrm{d}\rho+a_2,\label{eq::CurvaSolucaoPotPrescrito}
\end{equation}
where $A(\rho)$ is given by Eq (\ref{eq::willmoreEDO}) and $a_2$ is another constant of integration.

In short, given a smooth function $U(\rho)$, we can define a 2-parameter family of curves
\begin{equation}
\gamma(\rho;U(\rho),a_1,a_2) = \pm\int\displaystyle\frac{\rho\,\Big(2\int\,U\,\frac{\textrm{d}\rho}{\rho}+a_1\Big)}{[1-\rho^2(2\int\,U\,\frac{\textrm{d}\rho}{\rho}+a_1)^2]^{1/2}}\,\textrm{d}\rho+a_2\,.\label{eq::familyCurvesLambda}
\end{equation}
which furnishes a 2-parameter family of surfaces of revolution with a GIP $\sqrt{H^2(\rho)-K(\rho)}=U(\rho)$ by applying a rotation around the $z$-axis.
\newline
\newline
\begin{example} (vanishing geometry-induced potential)
For $U\equiv0$, Eq. (\ref{eq::willmoreEDO}) gives $A(\rho)=a_1$ constant and, from Eq. (\ref{eq::familyCurvesLambda}), one has 
\begin{equation}
\lambda(\rho) =\left\{
\begin{array}{ccc}
 \pm\sqrt{a_1^{-2}-\rho^2}+a_2 & , & a_1\not=0\\
 a_2 & , & a_1 = 0\\
\end{array}
\right..
\end{equation}
Then, for $a_1\not=0$, one has a sphere of radius $R=1/a_1$, and if $a_1=0$ one has a region of a plane. By a well known result, the only surfaces satisfying $H^2-K\equiv0$ are (pieces of a) sphere or plane (see \cite{Manfredo}, p. 147). In this way we recovered the two cases of surfaces where $H^2-K\equiv0$.
\end{example}
\begin{example} (constant geometry-induced potential)
Remember that for a cylinder of radius $R$, the geometry-induced potential is $U\equiv(2R)^{-1}$. However, a cylinder can not be obtained from the parametrization in Eq. (\ref{eq::ParamHelicoidalSur}); for a cylinder $\mathbf{x}(\rho,\phi)=(R\,\cos \phi,R\,\sin \phi,\rho)$. Now we show that there are other examples of surfaces of revolution, which are not a cylinder, with $U\equiv U_0\not=0$ constant. The importance of such examples lies in the fact that surfaces with a constant GIP have the same set of eigenfunctions of the problem without the GIP \footnote{Indeed, two Hamiltonians $\hat{H}_i=-\hbar^2/2m^*\,\Delta_g+\mathbb{G}_i$ differ by a constant, i.e., $\mathbb{G}_1-\mathbb{G}_2\equiv$ constant, if and only if they have the same set of eigenfunctions when subjected to the same boundary conditions. In this case, if $E^{(1)}_n$ and $E^{(2)}_n$ denote the respective eigenvalues for the same eigenfunction $\psi_n$, we have $E^{(1)}_n-E^{(2)}_n = \mathbb{G}_1-\mathbb{G}_2$ (notice that the gap between the eigenvalues satisfies $E_{n+k}^{(2)}-E_{n}^{(2)}=E_{n+k}^{(1)}-E_{n}^{(1)}$).}. 

Indeed, assuming $U(\rho)=U_0$ constant, Eq. (\ref{eq::familyCurvesLambda}) gives 
\begin{equation}
\lambda(\rho)=\pm\int_{\rho_0}^{\rho}\frac{x\Big(2\,U_0\ln\left(\frac{x}{\rho_0}\right)+a_1\Big)}{[1-x^2(2\,U_0\ln\left(\frac{x}{\rho_0}\right)+a_1)^2\,]^{1/2}}\,\textrm{d}x+a_2\,. 
\end{equation} 
The rotation of this curve around the $z$ axis generates a non-cylindrical surface with constant GIP $U_0$.
\end{example}

\subsection{Surfaces whose generating curve is a graph rotated around a horizontal axis}
 
Now we focus on another class of surfaces of revolution. In the previous analysis, the curve on the $xz$ plane to be rotated around the $z$ axis was supposed to be a graph, i.e., of the form $z=z(x)$. In this way, the surfaces obtained do not include cylinders and, more generally, do not include the surface of deformed nanotubes \cite{FernandosEfranceses} also. To include such examples, we can enlarge our class of surfaces by allowing a rotation of a curve $z=z(x)$ around the $x$ axis. We can parametrize these surfaces according to
\begin{equation}
\mathbf{x}(q,\phi) = (q,\rho(q)\sin\phi,\rho(q)\cos\phi),
\end{equation}
where $\rho(q)>0$ is a function which represents the distance to the rotation axis and defines the generating curve $(q,0,\rho(q))$ in the $xz$ plane to be rotated around the $x$ axis. As usual, $\phi$ is the angle of rotation. 

The geometry-induced potential of such surfaces can be written as \cite{FernandosEfranceses}
\begin{equation}
V_{gip}=-\frac{\hbar^2}{2m}\frac{[1+\rho'(q)^2+\rho(q)\rho''(q)]^2}{4\rho(q)^2[1+\rho'(q)^2]^{3}}\,,\label{eq::PotGeomTuboEnrugado}
\end{equation}
which furnishes for $U=\sqrt{H^2-K}$ the expression
\begin{equation}
\pm U=\frac{1+\rho'(q)^2+\rho(q)\rho''(q)}{2\rho(q)[1+\rho'(q)^2]^{3/2}}=-\frac{\rho}{2\rho'}\frac{\textrm{d}A}{\textrm{d}q},
\end{equation}
where
\begin{equation}
A=\frac{1}{\rho(q)[1+\rho'(q)^2]^{1/2}}\,.
\end{equation}
Then, we have the following differential equation for $A$
\begin{equation}
\rho\frac{\textrm{d}A}{\textrm{d}q}+2(\pm U)\frac{\textrm{d}\rho}{\textrm{d}q}=\left[\rho\frac{\textrm{d}A}{\textrm{d}\rho}+2(\pm U)\right]\frac{\textrm{d}\rho}{\textrm{d}q}=0.
\end{equation}
If $\rho'\equiv0$, then $\rho=$ constant and we have a cylinder. Otherwise, we find the following ODE in terms of $\rho$
\begin{equation}
\rho\frac{\textrm{d}A}{\textrm{d}\rho}+2(\pm U)=0\Rightarrow A(\rho)=\pm\left(2\int\frac{\textrm{d}\rho}{\rho}\,U(\rho)+a_1\right)\,,
\end{equation}
where $a_1$ is a constant of integration. Notice that this last equation is identical to Eq. (\ref{eq::willmoreEDO}), with the difference that here $\rho=\rho(q)$ is the function that we are trying to find.

Now, by using the definition of $A$, we find
\begin{equation}
\frac{\textrm{d}\rho}{\textrm{d}q}=\pm\sqrt{\frac{1-\rho^2A^2}{\rho^2A^2}}\Rightarrow q(\rho)=\pm\int \frac{\rho A}{\sqrt{1-\rho^2A^2}}\,\textrm{d}\rho+q_0\,.
\end{equation}
This equation is identical to Eq. (\ref{eq::CurvaSolucaoPotPrescrito}), but instead of obtaining the function which gives the generating curve, we obtained its inverse. This result reveals a certain duality between the surface of revolution obtained by rotating a curve $z=z(x)$ around the $x$ or the $z$ axes. In other words, 

\begin{prop}
Let $U$ be a smooth function of one variable, then each curve of the 2-parameter family given in (\ref{eq::CurvaSolucaoPotPrescrito}) generates a surface of revolution whose geometry-induced potential is $\sqrt{H^2-K}=U$ when rotated around the $x$ or the $z$ axis.
\end{prop}

\chapter{CONSTRAINED \,\, DYNAMICS \,\, ON \,\, HELICOIDAL \,\, SURFACES}
\label{chap_ConstDynHelSurf}

The definition of chirality comes from the fact that some objects can not be transformed into their mirror image under applications of rigid motions. This idea is present in many scientific areas and is of fundamental importance \cite{CahnAngewChem1966}. It appears in nature, such as tendrils and gastropod shells, and more fundamentally in the structure of DNA molecules. The study of chiral molecules is an important branch of stereochemistry with many applications in inorganic, organic, and physical chemistry, and also with several implications for the pharmaceutical industry. The concept of chirality is also present in particle physics and condensed matter \cite{KondepudiSciAmer1990}. In particular, this concept has proved to be useful in understanding some recent experimental results related to electronic, mechanical, and optical properties of nanotubes \cite{Yakobson2014}.

Recently, a link between chirality and the constrained particle dynamics was observed in the study of a particle on a helicoid \cite{AtanasovPRB2009, AtanasovPRB2015}. A helicoid is a particular instance of a helicoidal surface. These surfaces form the natural candidates to investigate a link with the concept of chirality. Indeed, given a curve $\alpha(\rho)=(\rho,0,\lambda(\rho))$ on the $xz$ plane, we can obtain enantiomorphic surfaces by screw-rotating $\alpha$ around the $z$ axis clock and counterclockwisely:
\begin{equation}
(\rho\cos(\omega\phi),\rho\sin(\omega\phi),\lambda(\rho)+\phi)\leftrightarrow(\rho\cos(\omega\phi),-\rho\sin(\omega\phi),\lambda(\rho)+\phi)\,.
\end{equation}
Observe that the sign of the constant $\omega$ can be used in order to control the chirality of the respective surface.

In the following, we study the geometric properties of helicoidal surfaces and comment on the existence of the so-called natural parameters, which allows for a better understanding and unified approach to such surfaces. The study of the respective Schr\"odinger equation under the influence of the GIP in such a coordinate system, along with some comparisons with known results for the dynamics on a helicoid, are present in the next section.  

\section{Parametrization by natural parameters}

Helicoidal surfaces are invariant by a rotation in combination with a translation (screw-rotation), the standard example being a helicoid, whose generating curve is just a line segment $(\rho,0,0)$:
\begin{equation}
\mathbf{x}_{helic}(\rho,\phi) = (\rho\cos(\omega\phi),\rho\sin(\omega\phi),\phi),\label{eq:ParHelicoid}
\end{equation}
where $\omega$ is a constant. If $L$ is the height of the helicoid, then we can write $\omega=2\pi n/L$, where $n$ is the number of twists around the screw-rotation axis. Moreover, the sign of $\omega$ governs the distinct chiralities states exhibited by helicoidal surfaces and has some consequences for the dynamics \cite{AtanasovPRB2009,AtanasovPRB2015}.

\begin{remark} In the previous Section we have already encountered helicoidal surfaces, Eq. (\ref{eq::ParamHelicoidalSur}), but here we adopt a different notation in order to ease comparisons with known results for the dynamics on a helicoid. As a consequence, surfaces of revolution are not allowed, since a translation in the direction of the screw axis is always present. However, surfaces of revolution can be formally obtained by changing $\omega\phi\mapsto \phi$ and then taking $\omega\to\infty$.
\end{remark}

For the helicoid, the coordinate system $(\rho,\phi)$ allows for a simple interpretation: $\phi$ represents the rotation angle (observe that the translation in the direction of the screw axis is proportional to the angular rotation), while the $\rho$-constant curves are helices; $\rho$ is the distance from the screw axis. On the other hand, for a general helicoidal surface, the translation along the screw axis has an extra contribution, which depends on the height of the generating curve $\alpha(\rho)=(\rho,0,\lambda(\rho)),$ $\rho>0$. Then, we have
\begin{equation}
\mathbf{x}(\rho,\phi) = (\rho\cos(\omega\phi),\rho\sin(\omega\phi),\lambda(\rho)+\phi).\label{eq:ParHelicSurf}
\end{equation}

In the above parametrization of a helicoidal surface the coordinate system $(\rho,\phi)$ does not have the same interpretation as happens for a helicoid. Indeed, in order to achieve that one could use a coordinate system composed of \textit{natural parameters}. More precisely, we say that a helicoidal surface $\Sigma$ is parametrized by \textit{natural parameters} $(\xi,\chi)$ if:
\begin{enumerate}
\item $\xi$-curves ($\chi$ constant) are parametrized by the arc-length parameter; and
\item $\chi$-curves ($\xi$ constant) are helices orthogonal to the $\xi$-curves. 
\end{enumerate}
In other words, since $\xi$ is the arc-parameter of a $\chi$-curve, the parameter $\xi$ represents a distance from the screw axis, while $\chi$ denotes the parameter along the orbits of the screw rotation symmetry, i.e., helices. This is precisely what happens for a helicoid, where $\xi_{helic}=\rho$ and $\chi_{helic}=\phi$ (here $\lambda_{helic}\equiv0$).

A useful consequence of using natural parameters is that the metric can be written in a simpler form:
\begin{equation}
\textrm{d}s^2 = \textrm{d}\xi^2+\mathcal{U}^2(\xi)\,\textrm{d}\chi^2,
\end{equation}
for some function $\mathcal{U}$.

 It is possible to show that every helicoidal surface admits a reparametrization by natural parameters \cite{DoCarmoTohoku1982}. Indeed, from the line element
\begin{eqnarray}
\textrm{d}s^2 & = & (1+\lambda'\,^2)\textrm{d}\rho^2+2\lambda'\textrm{d}\rho\,\textrm{d}\phi+(1+\omega^2\rho^2)\textrm{d}\phi^2\\[4pt]
& = &\left(1+\frac{\omega^2\rho^2\lambda'\,^2}{1+\omega^2\rho^2}\right)\textrm{d}\rho^2+(1+\omega^2\rho^2)\left(\textrm{d}\phi+\frac{\lambda'}{1+\omega^2\rho^2}\textrm{d}\rho\right)^2\,,
\end{eqnarray}
one finds the desired coordinate system $(\xi,\chi)=(\xi(\rho,\phi),\chi(\rho,\phi))$ by solving
\begin{equation}
\left\{
\begin{array}{ccc}
\textrm{d}\xi & = & \displaystyle\left(1+\frac{\omega^2\rho^2\lambda'\,^2}{1+\omega^2\rho^2}\right)^{1/2}\textrm{d}\rho\\ [10pt]
\textrm{d}\chi & = & \displaystyle\frac{\lambda'}{1+\omega^2\rho^2}\,\textrm{d}\rho+\textrm{d}\phi
\end{array}
\right.\,.
\end{equation}
Observe that our notations are slightly distinct from that of Do Carmo and Dajczer \cite{DoCarmoTohoku1982}: $(\xi,\chi,\omega\phi,\omega,a)_{ours}\mapsto(s,t,\phi,1/h,m)_{theirs}$. 

Using natural parameters $(\xi,\chi)$ to write the line element gives
\begin{equation}
\textrm{d}s^2 = \textrm{d}\xi^2+(1+\omega^2\rho^2)\textrm{d}\chi^2,
\end{equation}
which, by taking into account that $\rho$ does not depends on $\chi$, i.e., $\rho=\rho(\xi)$, and consequently also $\lambda=\lambda(\xi)$, can be rewritten as
\begin{equation}
\textrm{d}s^2 = \textrm{d}\xi^2+\mathcal{U}^2(\xi)\textrm{d}\chi^2,
\end{equation}
where $\mathcal{U}^2(\xi)=1+\omega^2\rho^2(\xi)$. For a helicoid, the map  $(\rho,\phi)\mapsto(\xi,\chi)$ is just the identity and, therefore, one has $\mathcal{U}^2_{helicoid}(\rho=\xi)=1+\omega^2\rho^2$. 

The function $\mathcal{U}$ encodes all the geometric information of its associated helicoidal surface and, consequently, both the Gaussian and the Mean curvatures are written in terms of $\mathcal{U}$. Further, we mention that $\mathcal{U}$ also determines the geometry-induced potential which governs the behavior of a quantum particle confined on the associated helicoidal surface.
\newline

A natural question now is if we can associate a helicoidal surface with a given non-negative function $\tilde{\mathcal{U}}(\xi)$, i.e., 
\newline
\textbf{Problem:} Given a function $\tilde{\mathcal{U}}(\xi)>0$, is it possible to find a constant $\tilde{\omega}$ and some functions $\rho,\phi$, and $\lambda$, such that the helicoidal surface $$\mathbf{x}(\rho,\phi)=(\rho\cos(\tilde{\omega}\phi),\rho\sin(\tilde{\omega}\phi),\lambda(\rho)+\phi)$$ has its line element written in natural coordinates as $\textrm{d}s^2=\textrm{d}\xi^2+\tilde{\mathcal{U}}^2\,\textrm{d}\chi^2$?

This problem do admit a solution to any given function $\mathcal{U}>0$. In fact, it is always possible to find a 2-parameter family of helicoidal surfaces associated with it. This is precisely the content of the \textit{Bour Lemma} \cite{DoCarmoTohoku1982}. It states that for every non-zero function $\mathcal{U}$ there exists a 2-parameter family of isometric helicoidal surfaces associate with it. The functions $(\rho,\phi)$ and $\lambda(\rho)$ which characterize the helicoidal surface can be written as \cite{DoCarmoTohoku1982}
\begin{equation}
\left\{
\begin{array}{l}
\rho=\rho(\xi)=\displaystyle\frac{1}{\omega}\sqrt{a^2\,\mathcal{U}^2-1}\\[10pt]
\lambda=\lambda(\xi)=\displaystyle\frac{1}{\omega}\int \textrm{d}\xi\,\frac{a\,\mathcal{U}}{a^2\,\mathcal{U}^2-1}\sqrt{a^2\mathcal{U}^2\,[\omega^2-a^2\,\dot{\mathcal{U}}^2]-\omega^2}\\[8pt]
\phi=\phi(\xi,\chi)=\displaystyle\frac{\chi}{a}-\frac{1}{\omega}\displaystyle\int \textrm{d}\xi\,\frac{\sqrt{a^2\,\mathcal{U}^2\,[\omega^2-a^2\,\dot{\mathcal{U}}^2]-\omega^2}}{a\,\mathcal{U}[a^2\,\mathcal{U}^2-1]}
\end{array}
\right.,\label{Eq::solucaoLemaBour}
\end{equation}
where a dot represent the derivative with respect to $\xi$: $\dot{\mathcal{U}}=\textrm{d}\,\mathcal{U}/\textrm{d}\xi$. By varying the constants $a$ and $\omega$ above, we generate a 2-parameter family of isometric helicoidal surfaces associated with the $\mathcal{U}$ given \textit{a priori}\footnote{If we choose $\mathcal{U}=\mathcal{U}_0$ to be a constant function, then we obtain a 2-parameter family of helicoidal surfaces which are contained on a cylinder of radius $\rho=\sqrt{a^2\mathcal{U}^2_0-1}$.}.

Finally, the Gaussian and Mean curvatures are written as \cite{DoCarmoTohoku1982}
\begin{equation}
K=K(\xi) = - \frac{\ddot{\mathcal{U}}}{\mathcal{U}},\label{eqGaussCurvHelicoidalSurfaces} 
\end{equation}
and
\begin{equation}
H=H(\xi)=\frac{a^2\, \mathcal{U}\,\ddot{\mathcal{U}}+a^2\,\dot{\mathcal{U}}^2-\omega^2}{2\sqrt{a^2\,\mathcal{U}^2\,[\omega^2-a^2\,\dot{\mathcal{U}}^2]-\omega^2}}\label{eqMeanCurvHelicoidalSurfaces}
\end{equation}
respectively, where he have adopted the surface normal $\mathbf{N}=\mathcal{U}^{-1}(\partial_{\chi}X\times\partial_{\xi}X)$.

According to the Bour lemma, we have for each function $\mathcal{U}$ a 2-parameter family of isometric helicoidal surfaces $[\mathcal{U},\omega,a]$. This means that the metric, and also the Gaussian curvature, is the same for all the helicoidal surfaces in the family. However, since the Mean curvature is not a bending invariant, the parameters $\omega$ and $a$ can give rise to different values of $H$. It follows that these parameters can be of physical relevance, since the geometry-induced potential also depends on the Mean curvature $H$.
\newline
\newline
\begin{example} (helicoidal minimal surfaces) Imposing the condition $H=0$ to Eq. (\ref{eqMeanCurvHelicoidalSurfaces})
gives
\begin{eqnarray}
a^2\mathcal{U}\,\ddot{\mathcal{U}}+a^2\dot{\mathcal{U}}^2=\omega^2\,\Rightarrow\,\mathcal{U}^2(\xi)=\frac{1}{a^2}(\omega^2\xi^2+2\,\omega_1\omega\xi+\omega_0),
\end{eqnarray}
where $\omega_0$, $\omega_1$ are constants  satisfying $b=\omega_0-\omega_1^2\geq1$, since $a^2\mathcal{U}^2-1>0$. In short, helicoidal minimal surfaces are characterized by a quadratic polynomial (for the particular case of a helicoid, we have $a=\omega_0=1$ and $\omega_1=0$). The Gaussian curvature of a helicoidal minimal surface is given by
\begin{equation}
K(\xi)=-\frac{\omega^2(\omega_0-\omega_1^2)}{a^4\mathcal{U}^4}=-\frac{b\,\omega^2}{[(\omega\xi+\omega_1)^2+b]^2}<0\,.\label{eqGaussianCurvMinHelSurf}
\end{equation}

The solution of Eq. (\ref{Eq::solucaoLemaBour}) for a helicoidal minimal surface is
\begin{equation}
\left\{
\begin{array}{l}
\rho=\rho(\xi)=\displaystyle\frac{1}{\omega}\sqrt{\omega^2\xi^2+2\,\omega_1\omega\xi+\omega_0-1}\\[8pt]
\lambda(\xi)=\sqrt{b-1}\displaystyle\int \textrm{d}\xi\,\frac{1}{\sqrt{\omega^2\xi^2+2\,\omega_1\omega\xi+\omega_0}\,(\omega^2\xi^2+2\,\omega_1\omega\xi+\omega_0-1)}\\[10pt]
\phi=\phi(\xi,\chi)=\displaystyle\frac{\chi}{a}-\sqrt{b-1}\displaystyle\int \textrm{d}\xi\,\frac{\sqrt{\omega^2\xi^2+2\,\omega_1\omega\xi+\omega_0}}{\omega^2\xi^2+2\,\omega_1\omega\xi+\omega_0-1}
\end{array}
\right..
\end{equation}
The parameter $a$ plays no relevant role. Indeed, by doing $\chi\mapsto a\chi$, we see that all the surfaces with distinct $a$ have the same image but different parametrizations. 
\end{example}

\section{Schr\"odinger equation on invariant surfaces}

In the previous sections, we have introduced and studied the geometry of surfaces invariant by a 1-parameter subgroup of isometries of $\mathbb{R}^3$. In this section, we devote our attention to the Schr\"odinger equation for a constrained particle on such surfaces.

All the three types of invariant surfaces have in common the following property: they admit the existence of a coordinate system $(u,v)$ such that the respective line element can be written as
\begin{equation}
\textrm{d}s^2 = \textrm{d}u^2 + f^2(u)\,\textrm{d}v^2,\,(u,v)\in[u_0,u_1]\times[v_0,v_1],
\end{equation}
where $f$ is a positive smooth function. Since such a metric has $g_{12}=0$, the Gaussian curvature (which only depends on the coefficients $g_{ij}$) can be expressed as
\begin{equation}
K = -\frac{1}{2\sqrt{g_{11}g_{22}}}\left[\frac{\partial}{\partial v}\left(\frac{g_{11,2}}{\sqrt{g_{11}g_{22}}}\right)+\frac{\partial}{\partial u}\left(\frac{g_{22,1}}{\sqrt{g_{11}g_{22}}}\right)\right] = -\frac{\ddot{f}}{f},
\end{equation}
where $g_{ij,k}=\partial g_{ij}/\partial q^k$, with $q^1=u$, $q^2=v$, and a dot denotes the derivative with respect to $u$. Naturally, the Mean curvature does not admit such a unified description, since it is not a bending invariant.

Let $\Sigma$ be an invariant surface with coordinate system $(u,v)$ as above. The Hamiltonian reads
\begin{equation}
\hat{H} = -\frac{\hbar^2}{2m^*}\Delta_g+V_{gip}=-\frac{\hbar^2}{2m^*f}\left[\frac{\partial}{\partial u}\left(f\frac{\partial}{\partial u}\right)+\frac{1}{f}\frac{\partial^2}{\partial v^2}\right]+V_{gip}\,.
\end{equation}
Now, rescaling the wave function as $\psi \mapsto \psi/g^{1/4}=\psi/\sqrt{f}$ (the Hamiltonian $\hat{H}$ should be rescaled as $f^{\frac{1}{2}}\,\hat{H}\,f^{-\frac{1}{2}}$), we have
\begin{eqnarray}
\hat{H} & = & -\frac{\hbar^2}{2m^*}\left[\frac{\partial^2}{\partial u^2}+\frac{1}{f^2}\frac{\partial^2}{\partial v^2}\right]+V_{eff}\,,
\end{eqnarray}
where
\begin{eqnarray}
V_{eff} &=&-\frac{\hbar^2}{2m^*}\left(-\frac{\ddot{f}}{2f}+\frac{\dot{f}\,^2}{4f^2}\right)-\frac{\hbar^2}{2m^*}(H^2-K),\\
& = & -\frac{\hbar^2}{2m^*}\left(\frac{\dot{f}\,^2}{4f^2}+\frac{\ddot{f}}{2f}\right)-\frac{\hbar^2}{2m^*}\,H^2\,.
\end{eqnarray}

As a corollary, it follows that the stationary Schr\"odinger equation can be solved by separation of variables. Indeed, writing $\psi(u,v)=A(u)B(v)$, we have
\begin{equation}
\frac{1}{B(v)}\frac{\textrm{d}^2B(v)}{\textrm{d}v^2}=-\Big(U(u)+k^2\Big)f^2(u)-\frac{f^2(u)}{A(u)}\frac{\textrm{d}^2A(u)}{\textrm{d}u^2},
\end{equation}
where  $k^2=2m^*E/\hbar^2$ and $U(u)=-2m^*\,V_{eff}(u)/\hbar^2$. This procedure   furnishes the following equations
\begin{equation}
\left\{
\begin{array}{c}
B''(v) =-\lambda B(v) \\[5pt]
A''(u) + \left(U(u)+k^2-\displaystyle\frac{\lambda}{f^2(u)}\right)A(u) = 0\\
\end{array}
\right.\,,\label{eq::EqSchrDesacopladaSuperInvar}
\end{equation}
whose solutions depends on the imposed boundary conditions. 

The above equations clearly show that for an invariant surface the stationary Schr\"odinger equation decouples into an equation along the orbits of the 1-parameter subgroup ($v$-curves) and an effective equation along the direction orthogonal to the orbits ($u$-curves), i.e., an effective equation along the generating curve.

\subsection{Schr\"odinger equation for cylindrical surfaces}

For a cylindrical surface one has 
\begin{equation}
\left\{\begin{array}{c}
\textrm{d}s_{cyl}^2 = \textrm{d}u^2+\textrm{d}v^2\\[5pt]
K_{cyl}\equiv0,\,H_{cyl} = \displaystyle\frac{\dot{x}\ddot{y}-\ddot{x}\dot{y}}{2} = \frac{\kappa}{2}\\
\end{array}
\right.,
\end{equation}
where $\kappa(u)$ is the curvature function of the cross section $\alpha(u)=(x(u),y(u),0)$ (generating curve), $u$ being its arc-length, which is translated in the direction of $(0,0,1)$. Thus, the decoupled equations (\ref{eq::EqSchrDesacopladaSuperInvar}) read
\begin{equation}
\left\{
\begin{array}{c}
B''(v) =-\lambda B(v)\\[5pt]
A''(u) +\frac{\kappa^2}{8}A(u)+(k^2-\lambda)A(u) = 0\\
\end{array}
\right.\,.\label{eq::EqSchrDesacopladaSuperCylind}
\end{equation}
For a cylindrical surface we may assume homogeneous boundary conditions for the $v$-directions. Then, the energy spectrum is given by
\begin{equation}
E_{cyl}(n_u,n_v) = \frac{h^2n_v^2}{8m^*L_v^2} + E_{\kappa,n_u}\,,
\end{equation}
where $L_v$ is the height of the cylindrical surface, with $n_v\in\{1,2,...\}$, and $E_{\kappa,n_u}$ is the $n_u$-th eigenenergy of a constrained particle in a 1D box of length $L_u$ under a potential $V_{gip}=-\hbar^2\kappa^2/8m^*$: a box with homogeneous or periodic boundary conditions if $\alpha$ is open or closed, respectively.\footnote{In an intrinsic approach, i.e., in the absence of $V_{gip}=-\hbar^2\kappa^2/8m^*$, one would find $E_{cyl}(n_u,n_v) = h^2n_v^2/8m^*L_v^2 + h^2n_u^2/8m^*L_u^2$, with $n_u,n_v\in\{1,2,...\}$, for an open cross section or $E_{cyl}(n_u,n_v) = h^2n_v^2/8m^*L_v^2 + h^2n_u^2/2m^*L_u^2$, with $n_u\in\{1,2,...\}$ and $n_v\in\{0,1,2,...\}$, for a closed cross section \cite{BastosEtAlPreprint}.}.

\subsection{Schr\"odinger equation for surfaces of Revolution}

For a surface of revolution one has 
\begin{equation}
\left\{\begin{array}{c}
\textrm{d}s_{rev}^2 = \textrm{d}u^2+x^2(u)\textrm{d}v^2\\[5pt]
K_{rev}=-\frac{\ddot{x}}{x},\,H_{rev} = \displaystyle\frac{x(\dot{x}\ddot{z}-\ddot{x}\dot{z})+\dot{z}}{2x}\\
\end{array}
\right.,
\end{equation}
where $\alpha(u)=(x(u),0,z(u))$, with $x>0$, is the generating curve which is rotated around the $z$ axis, with arc-length parameter $u$. Then, the decoupled equations (\ref{eq::EqSchrDesacopladaSuperInvar}) read
\begin{equation}
\left\{
\begin{array}{c}
B''(v) =-\lambda B(v)\\[5pt]
A''(u) +\left(\displaystyle\frac{\dot{x}^2}{4x^2}+\frac{\ddot{x}}{2x}+H_{rev}^2\right)A(u)+ \left(k^2-\displaystyle\frac{\lambda}{x^2}\right)A(u) = 0\\
\end{array}
\right.\,.\label{eq::EqSchrDesacopladaSuperRevol}
\end{equation}
For a revolution surface we may assume an angular periodicity for the $v$-curves, which gives $\lambda=m_{\chi}^2$, $m_{\chi}\in\mathbb{Z}$. Then, the effective dynamics in the $u$-direction is
\begin{equation}
-\frac{\hbar^2}{2m^*}A''(u) -\frac{\hbar^2}{2m^*}\left(\displaystyle\frac{\dot{x}^2+2x\ddot{x}-4m_{\chi}^2}{4x^2}+H_{rev}^2\right)A(u) = E\,A(u)\,.
\end{equation}
Depending on the concavity of $x(u)$ and on the values of the angular momentum quantum number $m_{\chi}$, the contribution of $\dot{x}^2+2x\ddot{x}-4m_{\chi}^2$ in the effective potential for the $u$-direction can be attractive or repulsive, then changing the way it favors the existence of geometry-induced bound states \cite{MoraesAnnPhys,AtanasovPLA2007}.

\subsection{Schr\"odinger equation for helicoidal surfaces}

For a helicoidal surface one has 
\begin{equation}
\left\{\begin{array}{c}
\textrm{d}s_{hel}^2 = \textrm{d}u^2+\mathcal{U}^2(u)\textrm{d}v^2\\[5pt]
K_{hel}=-\displaystyle\frac{\ddot{\mathcal{U}}}{\mathcal{U}},\,H_{hel}=\displaystyle\frac{a^2\, \mathcal{U}\,\ddot{\mathcal{U}}+a^2\,\dot{\mathcal{U}}^2-\omega^2}{2\sqrt{a^2\,\mathcal{U}^2\,[\omega^2-a^2\,\dot{\mathcal{U}}^2]-\omega^2}}\\
\end{array}
\right.,
\end{equation}
where $u=\xi$ and $v=\chi$ are natural parameters of the helicoidal surface introduced in section 6. Then, the decoupled equations (\ref{eq::EqSchrDesacopladaSuperInvar}) read
\begin{equation}
\left\{
\begin{array}{c}
B''(v) =-\lambda B(v)\\[5pt]
A''(u) +\left(\displaystyle\frac{\dot{\mathcal{U}}^2}{4\,\mathcal{U}^2}+\frac{\ddot{\mathcal{U}}}{2\,\mathcal{U}}+H_{hel}^2\right)A(u)+ \left(k^2-\displaystyle\frac{\lambda}{\mathcal{U}^2}\right)A(u) = 0\\
\end{array}
\right.\,.\label{eq::EqSchrDesacopladaSuperRevol}
\end{equation}

The standard example of a helicoidal surface is that of a helicoid. For such a surface it is known that particles with distinct angular quantum numbers tend to localize in distinct parts of the helicoid and also that there exist geometry-induced bound states \cite{AtanasovPRB2009}. In the following we extend these findings to all helicoidal minimal surfaces (the helicoid being the simplest example) and, due to the existence of other parameters associated to a helicoidal minimal surface, we show in addition the possibility of controlling the change in the distribution of the probability density when the surface is subjected to an extra charge, i.e., where the particles are find with greatest probability.

\section{Constrained dynamics on helicoidal minimal surfaces} 

\begin{figure*}[tbp]
\centering
  \includegraphics[width=0.45\linewidth]{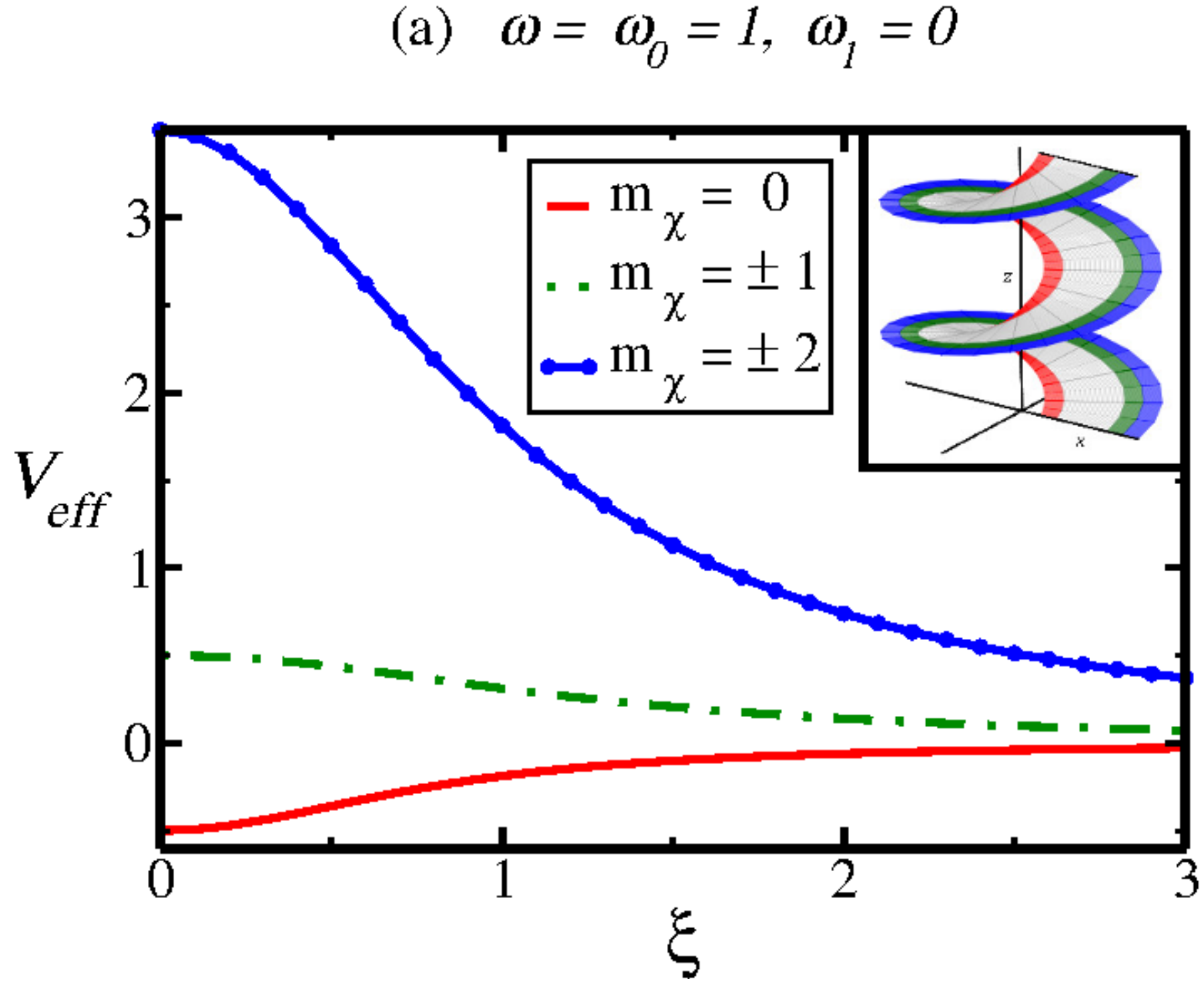}  {\includegraphics[width=0.45\linewidth]{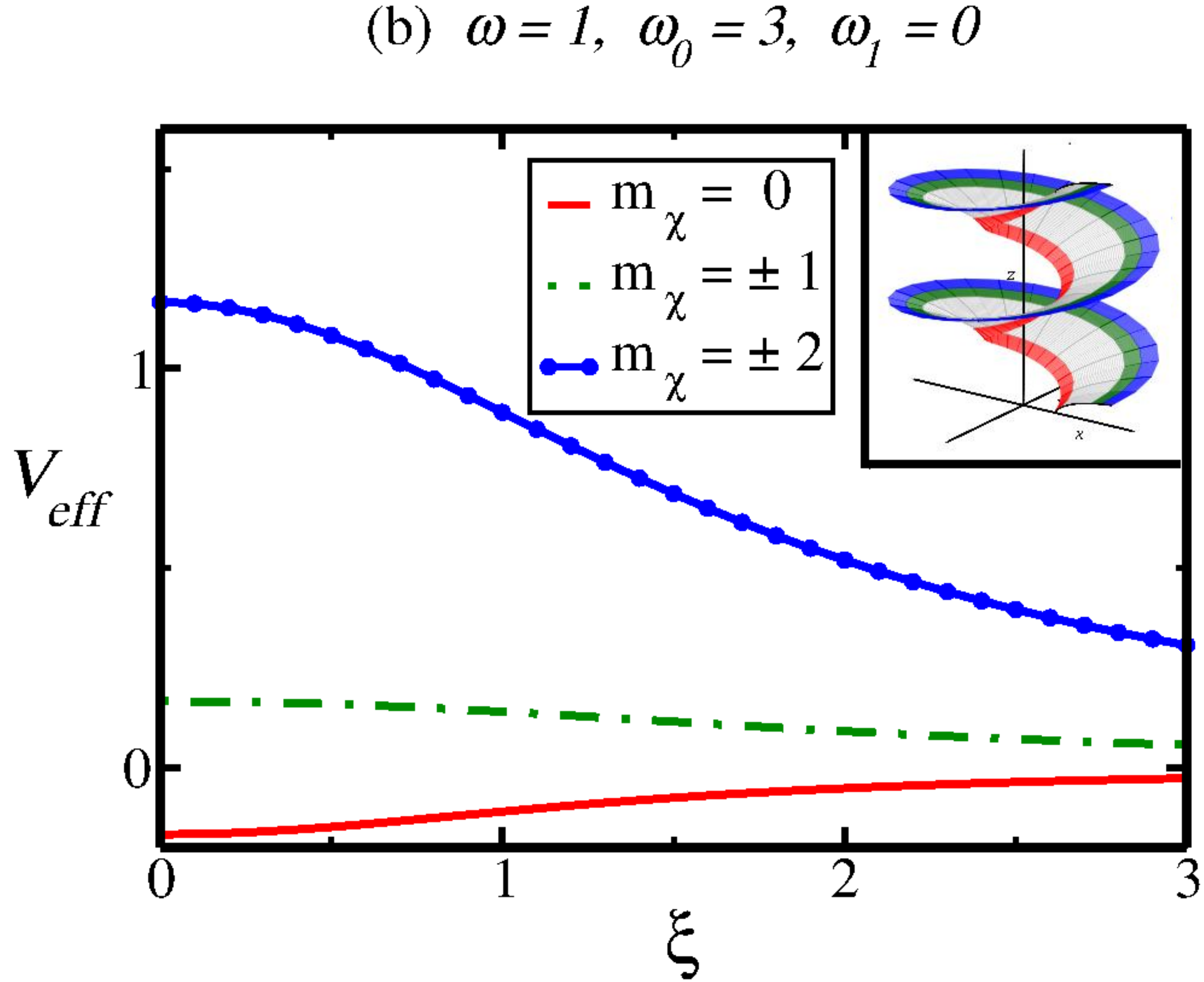}}
  
          \caption{The behavior of the effective potential $V_{eff}$ as a function of $\xi$ with $\hbar^2/2m^*=1$, $\omega=1$, $\omega_1=0$, and different values of $m_{\chi}$. Insets:  helicoidal surfaces for \textbf{(a)} $\omega_0=1$ and \textbf{(b)} $\omega_0=3$, respectively. The insets illustrate the fact that particles with distinct values of $m_{\chi}$ tend to localize in different parts of the surface.
          }
\label{fig:Fig1}
\end{figure*}

For a helicoidal minimal surface one has $H\equiv0$ and $\mathcal{U}^2=(\omega\xi+\omega_1)^2+b$ as seen in Example 6.1 (without loss of generality, we set $a=1$). Then, the decoupled Schr\"odinger equation reads ($u=\xi$, $v=\chi$)
\begin{equation}
B''(\chi) =-\lambda B(\chi)
\end{equation}
and
\begin{equation}
A''(\xi) + \frac{\omega^2}{4}\left\{\displaystyle\frac{1-\lambda}{b+(\omega\xi+\omega_1)^2}+\frac{b}{[b+(\omega\xi+\omega_1)^2]^2}\right\}A(\xi)
+k^2A(\xi) = 0\,.\label{eqfuncAforMHS}
\end{equation}

Writing the solution for $B$ as $B(\chi)=\textrm{e}^{\textrm{i}k_{\chi}\,\chi}$ furnishes
\begin{equation}
\lambda = k_{\chi}^2,
\end{equation}
with $k_{\chi}$ being the partial moment in the $\chi$ direction. The canonical momentum associated to the coordinate $\chi$, $L_{\chi}=-\textrm{i}\,\hbar/\omega\,\partial_{\chi}$, has the same eigenfunctions as the equation for $B$. The momentum $k_{\chi}$ is quantized according to 
\begin{equation}
k_{\chi} = m_{\chi}\,\omega\,,m_{\chi}\in\mathbb{Z}\,.
\end{equation}

Using the expression for $\lambda$ in Eq. (\ref{eqfuncAforMHS}) shows that the equation for the $\xi$ direction is subjected to the following effective potential
\begin{equation}\label{Veff}
V_{eff}(\xi)=-\frac{\hbar^2}{2m^*}\frac{\omega^2}{4}\left\{\frac{b}{[b+(\omega\xi+\omega_1)^2]^2}+\displaystyle\frac{1-4\,m_{\chi}^2}{b+(\omega\xi+\omega_1)^2}\right\}\,.
\end{equation}

This effective potential displays two terms with distinct contributions. The first term contributes attractively, while the second one depends on the sign of $m_{\chi}^2-1/4$, acting attractively for $m_{\chi}^2<1/4$, i.e., $m_{\chi}=0$, or repulsively for $m_{\chi}^2>1/4$, i.e., $m_{\chi}\not=0$. The effect of this variable part is of a centrifugal potential character for $m_{\chi}\not=0$ (repulsive), it pushes a particle to the outer border of the surface. On the other hand, when $m_{\chi}=0$ (attractive), the contribution of this variable part is of a anticentrifugal character and it concentrates the particles in the inner border of the minimal helicoidal surface, i.e., around the screw axis. This analysis is in agreement with what happens for the particular case of a helicoid \cite{AtanasovPRB2009}, where $\omega_0=1$ and $\omega_1=0$ ($b=1$). 

Now, we benchmark our analytical expression for the effective potential [see Eq. (\ref{Veff})] with the one derived by Atanasov \textit{et al.} \cite{AtanasovPRB2009}. For that purpose, we can assume the following set of parameters: $\omega_1 =0$ and $\omega_0 =1$ (then $b=1$). After substitution of values, we find that
\begin{equation}\label{VeffAtanosov}
V_{eff}(\xi)=-\frac{\hbar^2}{2m^*}\frac{\omega^2}{4}\left[\frac{1}{(1+\omega^2\xi^2)^2}+\displaystyle\frac{1-4\,m_{\chi}^2}{1+\omega^2\xi^2}\right]\,.
\end{equation}
The evolution of $V_{eff}$, in Eq. (\ref{VeffAtanosov}), as a function of $\xi$ is depicted in Fig. \ref{fig:Fig1} (a). As one readily sees, the behavior of $V_{eff}$ is strongly affected by the angular momentum quantum number $m_{\chi}$. When
$m_{\chi} =0$, we have $V_{eff}(\xi) < 0$, which leads to the existence of bound states\footnote{A globally attractive potential $V$ satisfying the criterion $\int V(x)\,\textrm{d}x^n<0$ do admit the existence of bound states for $n=1$ or $n=2$ \cite{Chadan2003}.}. On the other hand, for nonvanishing angular momentum quantum numbers, we observe that $V_{eff}(\xi) > 0$ (the energy spectrum is positively valued) and no bound state is allowed \cite{AtanasovPRB2009}.

The above analysis is still valid for other values of the parameters $\omega$, $\omega_1$, and $\omega_0$. In other words, the existence of geometry-induced bound and localized states previously verified for a helicoid \cite{AtanasovPRB2009} can be extended to any helicoidal minimal surface.
\newline

Finally, let us comment that other results established for a helicoid can be extended to all helicoidal minimal surfaces with some additional advantages. Indeed, applying a change of variables
\begin{equation}
\xi(\tilde{\xi})=\sqrt{b}\,\,\tilde{\xi}-\frac{\omega_1}{\omega}=\sqrt{(\omega_0-\omega_1^2)}\,\,\tilde{\xi}-\frac{\omega_1}{\omega},\label{eq::MappingFromHelMinSurfToAhelicoid}
\end{equation}
which implies $\textrm{d}A/\textrm{d}\tilde{\xi}=\sqrt{b}\,\textrm{d}A/\textrm{d}\xi$, we can map our effective equation for the $\xi$-direction into that of a helicoid \cite{AtanasovPRB2009}
\begin{equation}
-\frac{\hbar^2}{2m^*}\frac{\textrm{d}^2A}{\textrm{d}\,\tilde{\xi}^2}-\frac{\hbar^2}{2m^*}\frac{\omega^2}{4}\left\{\frac{1}{(1+\omega\tilde{\xi} ^2)^2}+\displaystyle\frac{1-4\,m_{\chi}^2}{1+\omega\tilde{\xi}^2}\right\}A(\tilde{\xi}) = b\,E\,A(\tilde{\xi})\,.
\end{equation}
For example, when analyzing the distribution of the probability density for a constrained particle on a helicoidal minimal surface subjected to some charge distribution, as analyzed by Atanasov \textit{et al.} \cite{AtanasovPRB2009}, we can use the correspondence above to map the problem for a helicoidal minimal surface into an equivalent problem for a helicoid and then, by inverting Eq. (\ref{eq::MappingFromHelMinSurfToAhelicoid}), solve the original problem. So, by tuning the parameters $\omega_0$ and $\omega_1$, we can control the changing in the distribution of the probability density, e.g., we can govern the location where the particle will be found with greatest probability when the outer border of the helicoidal minimal surface is uniformly charged (i.e., where the extra charge will concentrate).

\chapter{CONCLUSIONS}
This thesis was devoted to the differential geometry of curves and surfaces along with some applications in quantum mechanics. In its first part we dedicated our attention to moving frames along curves. First, we introduced the well known Frenet frame and discussed on plane curves whose curvature is a power-law function, showing that they are related to spiral curves. We also described the curvature and torsion of space curves in terms of osculating spheres. This showed that, as happens for plane curves, spherical analogs may be used to describe space curves. Later, we introduced a general framework to describe adapted frames along curves and proved that the curvature function can be seen as a lower bound for the scalar angular velocity of any moving frame. This allowed us to define Rotation Minimizing (RM) frames as those frames that achieve this minimum, which can be done by somehow eliminating the contribution for the frame rotation related to the curve torsion. Interestingly, RM frames apply very well in the study of spherical curves and allow one to characterize curves on spheres through a linear equation involving the coefficients that dictate the frame motion. This happens even in higher dimensions, in contrast with a Frenet-like approach. Indeed, we discussed such characterization via a Frenet approach and the difficulties one may encounter in generalizing the characterization to higher dimensions.

We also applied these ideas to characterize curves that lie on a level surface, $\Sigma=F^{-1}(c)$. This was done by reinterpreting the problem in the context of a metric induced by the Hessian of $F$, which may fail to be positive or non-degenerate and naturally led us to the study of curves in Lorentz-Minkowski and isotropic spaces. We developed a systematic approach to the construction of RM frames and characterization of spherical curves in a Lorentz-Minkowski and isotropic spaces, and furnished a general criterion for a curve to lie on a level set surface. As a particular instance of this problem, we were able to completely characterize curves that lie on an Euclidean quadric.

We also extended the previous investigations in order to characterize curves that lie on the (hyper)surface of geodesics spheres in a Riemannian manifold. We discussed on the concept of normal curves, which are precisely the curves whose geodesics connecting a fixed point to points on the curve induce a normal vector field (along the curve), and mentioned that, as a consequence of the Gauss lemma for the exponential map in a Riemannian manifold $M$, on a sufficiently small neighborhood of a given point $p\in M$ the condition of being a normal curve (with center $p$) is equivalent to be a curve on a geodesic sphere. We then used this equivalence to characterize geodesic spherical curves in hyperbolic and spherical Riemannian geometries through a linear relation involving the coefficients that dictate the frame motion. For completeness, we also discussed the characterization of geodesic spherical curves in terms of a Frenet frame and show that the characterization of (geodesic) spherical curves is the same as in Euclidean space. Finally, we showed that if a Riemannian manifold contains totally geodesic submanifolds, which play the role of planes, then their curves are associated with a normal development that lies on a line passing through the origin. For the reciprocal, we proved in $\mathbb{S}^{m+1}(r)$ and $\mathbb{H}^{m+1}(r)$ that a curve lies on a totally geodesic submanifold if and only if its normal development is a line passing through the origin.

In the second half of this thesis we applied some of the theoretical framework developed in the first part in the quantum dynamics of a constrained particle. After describing the confining potential formalism to the constrained dynamics, from which emerges a geometry-induced potential (GIP) acting upon the dynamics, we devoted our attention to tubular surfaces as a mean to model the particle quantum dynamics on curved nanotubes. The use of RM frames offered a simpler description for the (tubular) Schrödinger equation, when compared with an approach based on the Frenet frame. In addition, it allowed us to show that the torsion of the centerline of a curved tube gives rise to a geometric phase, which is of fundamental importance in future applications in connection with phenomena like the Aharanov-Bohm effect.

Later, we studied the problem of prescribed GIP for curves and surfaces in Euclidean space. We showed that the problem for curves is easily solved by integrating Frenet equations, while the problem for surfaces involves a non-linear 2nd order partial differential equation (PDE). We exemplified the prescribed GIP problem through the study of curves whose geometry-induced potential is that of a Hydrogen atom. We further restricted ourselves to the prescribed GIP problem for surfaces invariant by a 1-parameter group of isometries. Due to their appealing symmetry, i.e., translation, rotation, and screw (helicoidal) symmetry, these surfaces are commonly encountered in applications and theoretical studies of quantum mechanics and do not constitute any severe restriction to the investigation of a constrained dynamics on surfaces. Besides, this simplifying hypothesis turns the study of the PDE for the prescribed potential into that of an ODE and discloses many potentialities of invariant surfaces in applications. In addition, the invariance property also allows for a unified description of the Schrödinger equation under the effect of a geometry-induced potential. We completely solved the problem for cylindrical and revolution surfaces. For the class of helicoidal surfaces we presented the concept of natural parameters, which allows for a unified description of such surfaces and also the association of a 2-parameter family of isometric helicoidal surfaces with a given positive function. These surfaces are particularly important due to the fact that, by screw-rotating a curve clockwisely and counterclockwisely, one can easily generate pairs of enantiomorphic surfaces, which naturally turns helicoidal surfaces an adequate setting to test and exploit a link between chirality and the effects of a geometry-induced potential. Finally, for the family of helicoidal minimal surfaces we proved the existence of geometry-induced bound and localized states, then generalizing known results for the particular case of a helicoid, and in addition we also showed the possibility of controlling the change in the distribution of the probability density when the surface is subjected to an extra charge. This control is a fundamental step toward future applications of this formalism.

Naturally, some questions are still open and the seek for solutions may be seen as future perspectives. On the geometric side, a problem that remains open is that of computing RM frames for a generic curve. To the best of our knowledge, no exact solution exists and in general one must resort to a numerical approach (e.g., the double reflection method \cite{WangACMTOG}). In this respect we showed that the use of osculating spheres allowed us to write an alternative expression for the torsion, whose integration constitute a way to address the problem of computing RM frames. We hope our investigations may give some hints in this direction. In addition, the approach to RM frames in Riemannian geometry was restricted to constant curvature ambient spaces. These extensions are presently under investigation for some homogeneous spaces and will be the subject of a follow-up work.

On the physics side, it is worth mentioning that for a more realistic description of the constrained dynamics, the approach presented here must be extended to others contexts that take into account spin and relativistic effects, such as the Pauli-Schrödinger and Dirac equations. Due to the importance of nanotubes, it would be also interesting to investigate theses extensions for tubular surfaces. Finally, there is the obvious need of considering the prescribed geometry-induced problem for non-invariant surfaces.

\postextual

\addcontentsline{toc}{chapter}{\,\,\,\,\,\,\,\,\,\,\,\,\,\,\,\,\,\,\,\,\,\,\,\,\,\,\,REFERENCES}
\nocite{}
\clearpage
\markboth{REFERENCES}{REFERENCES}


\providecommand{\abntreprintinfo}[1]{%
 \citeonline{#1}}
\setlength{\labelsep}{0pt}

\printindex

\end{document}